%% file: main.tex
\documentclass[11pt]{article}

\usepackage{setspace}
\usepackage[utf8]{inputenc}

\usepackage{bbold}
\usepackage{dsfont}
\usepackage{mathtools}
\usepackage{stmaryrd}
\usepackage{verbatim}
\usepackage{amssymb}
\usepackage{titlesec}
\usepackage{mathrsfs}

\usepackage{amsmath}
\DeclareFontFamily{U}{mathx}{}
\DeclareFontShape{U}{mathx}{m}{n}{<-> mathx10}{}
\DeclareSymbolFont{mathx}{U}{mathx}{m}{n}
\DeclareMathAccent{\widehat}{0}{mathx}{"70}
\DeclareMathAccent{\widecheck}{0}{mathx}{"71}

\usepackage{float}
\usepackage{amsfonts}
\usepackage{comment}

\usepackage{mhchem}
\graphicspath{ {images/} }
\usepackage{array}
\newcolumntype{P}[1]{>{\centering\arraybackslash}p{#1}}
\usepackage{longtable}
\usepackage{mathrsfs}
\usepackage{caption}
\DeclareCaptionLabelFormat{adja-page}{\hrulefill\\#1 #2 \emph{(previous page)}}
\usepackage{subcaption}
\usepackage{xcolor}
\allowdisplaybreaks

\usepackage{cite}
\usepackage{tikz}

\usepackage{tabulary}
\usepackage{booktabs}
\usepackage{array} 
\usepackage[flushleft]{threeparttable}
\usepackage{blkarray}

\usepackage{hyperref}
\hypersetup{
    colorlinks=true,
    linkcolor=blue,
    filecolor=blue,      
    urlcolor=blue,
    citecolor=blue,
    pdftitle={Overleaf Example},
    pdfpagemode=FullScreen,
    }

\usepackage[toc,page,titletoc]{appendix}
\usepackage[capitalise,nameinlink]{cleveref}

\usepackage[T3,T1]{fontenc}
\DeclareSymbolFont{tipa}{T3}{cmr}{m}{n}
\DeclareMathAccent{\invbreve}{\mathalpha}{tipa}{16}

\usepackage{soul,xcolor}
\setstcolor{red}

\crefname{supp}{Supplement}{Supplements}

\usepackage{amsthm}
\newtheorem{theorem}{Theorem}[section]
\newtheorem{corollary}{Corollary}[section]
\newtheorem{lemma}{Lemma}[section]

\newtheorem{proposition}{Proposition}[section]

\newtheorem{assumption}{Assumption}[section]
\numberwithin{equation}{section}

\theoremstyle{definition}
\newtheorem{definition}{Definition}[section]

\theoremstyle{definition}
\newtheorem{remark}{Remark}[section]

\theoremstyle{definition}

\usepackage{authblk}

\newcommand{\initialx}{x^\circ}
\newcommand{\initialxbreve}{\Breve{x}^\circ}

\newcommand{\initialz}{z^\circ}
\newcommand{\initialzbreve}{\Breve{z}^\circ}

\newcommand{\labelDTMC}{\Breve}

\renewcommand{\l}{\ell}

\newcommand{\eps}{\varepsilon}

\newcommand{\Ss}{\mathscr{S}}
\newcommand{\Rs}{\mathscr{R}}

\newcommand{\PP}{\mathds{P}}
\newcommand{\R}{\mathds{R}}
\newcommand{\N}{\mathds{N}}
\newcommand{\Z}{\mathds{Z}}
\newcommand{\E}{\mathds{E}} 

\newcommand{\C}{\mathds{C}}

\newcommand{\A}{\mathcal{A}}
\newcommand{\D}{\mathcal{D}}
\newcommand{\F}{\mathcal{F}}
\newcommand{\B}{\mathcal{B}}

\renewcommand{\L}{\mathcal{L}}

\newcommand{\CC}{\mathcal{C}}

\newcommand{\Bc}{\mathcal{B}^c} 
\newcommand{\nBc}{\mathcal{B}}

\newcommand{\X}{\mathcal{X}}
\newcommand{\Y}{\mathcal{Y}}
\newcommand{\T}{\mathcal{T}}
\newcommand{\V}{\mathcal{V}}
\newcommand{\NN}{\mathcal{N}}
\newcommand{\one}{\mathds{1}}

\newcommand{\DR}{\text{D}^\text{R}}
\newcommand{\DA}{\text{D}^\text{A}}

\newcommand{\Dtot}{\text{D}_{\text{tot}}}

\newcommand{\prp}{\Lambda}  
\newcommand{\rate}{\Upsilon} 
\newcommand{\upb}{K} 

\DeclareMathOperator{\spr}{spr}
\DeclareMathOperator{\spec}{sp}
\DeclareMathOperator{\diag}{diag}
\DeclareMathOperator{\nullity}{nullity}

\newcommand{\interior}[1]{%
  {\kern0pt#1}^{\mathrm{o}}%
}

\DeclarePairedDelimiterX{\abs}[1]{\lvert}{\rvert}{#1}
\DeclarePairedDelimiterX{\norm}[1]{\lVert}{\rVert}{#1}
\DeclarePairedDelimiterX{\floor}[1]{\lfloor}{\rfloor}{#1}
\DeclarePairedDelimiterX{\qvar}[1]{\langle}{\rangle}{#1}
\DeclarePairedDelimiterX{\inn}[1]{\langle}{\rangle}{#1}

\DeclareCaptionLabelFormat{cont}{#1~#2\alph{ContinuedFloat}}
\title{Analysis of singularly perturbed stochastic chemical reaction networks motivated by applications to epigenetic cell memory}
\author[1,*]{Simone Bruno}
\author[2,*]{Felipe A. Campos}
\author[2]{Yi Fu}
\author[1]{Domitilla Del Vecchio}
\author[2]{Ruth J. Williams}
\affil[1]{Department of Mechanical Engineering, Massachusetts Institute of Technology, 77 Massachusetts Avenue, Cambridge, MA 02139. Emails: {\tt\small sbruno@mit.edu, ddv@mit.edu}}
\affil[2]{Department of Mathematics, University of California, San Diego, 9500 Gilman Drive, La Jolla CA 92093-0112. Email: {\tt\small (yif064,fcamposv,rjwilliams)@ucsd.edu}}
\affil[*]{These authors contributed equally: S. Bruno and F. Campos}
\date{}                     
\begin{document}
	\maketitle

\begin{abstract}

Epigenetic cell memory, the inheritance of gene expression patterns across subsequent cell divisions, is a critical property of multi-cellular organisms. In recent work \cite{BrunoDelVecchio}, a subset of the authors observed in a simulation study how the stochastic dynamics and time-scale differences between establishment and erasure processes in chromatin modifications (such as histone modifications and DNA methylation) can have a critical effect on epigenetic cell memory. 
In this paper, we provide a mathematical framework to rigorously validate and extend beyond these computational findings. Viewing our stochastic model of a chromatin modification circuit as a singularly perturbed, finite state, continuous time Markov chain, we 
extend beyond existing theory in order to characterize the leading coefficients in the series expansions of stationary distributions and mean first passage times. In particular, we characterize the limiting stationary distribution in terms of a reduced Markov chain, provide an algorithm to determine the orders of the poles of mean first passage times, and determine how changing erasure rates affects system behavior.
%
The theoretical tools developed in this paper not only allow us to set a rigorous mathematical basis for the computational findings of our prior work, highlighting the effect of chromatin modification dynamics on epigenetic cell memory, but they can also be applied to other singularly perturbed Markov chains beyond the applications in this paper, especially those associated with chemical reaction networks.

\end{abstract}

\input{Introduction.tex}

\input{Motivatingexamples.tex}

\input{AnalyticalpertCTMC.tex}		
	
\input{StationaryDistributions.tex}

\input{MFPT.tex}		

\input{Monotonicity.tex}

\input{Applications.tex}

\input{Conclusion.tex}

\noindent
{\bf Supplementary information (SI) file:} file containing the proofs of the theoretical tools developed in this paper, and detailed mathematical derivations for some of the chromatin modification circuit models analyzed.\\

\noindent
{\bf Acknowledgements:} S.B. was supported by NSF Collaborative Research grant MCB-2027949 (PI: D.D.V.). R.J.W., F.C and Y.F. were supported in part by NSF Collaborative Research grant MCB-2027947 (PI: R.J.W.) and by the Charles Lee Powell Foundation (PI: R.J.W.).\\ 

\noindent
{\bf Ethics declarations:} The authors declare that they have no conflicts of interest.\\

\noindent
{\bf Data availability:} Data sharing not applicable to this article as no datasets were generated or analysed during the current study.\\

	\newpage

\appendixpageoff
\appendixtitleoff
\renewcommand{\appendixtocname}{Supplementary Information}
\begin{appendices}
  \setcounter{section}{19}
  \crefalias{section}{supp}
  \setcounter{figure}{0}
\renewcommand{\thefigure}{S.\arabic{figure}}

  \setcounter{equation}{0}
\renewcommand{\theequation}{S.\arabic{equation}}

\pagenumbering{arabic}
\renewcommand*{\thepage}{\arabic{page}}	
	
\input{Appendix_title}

\input{Appendix_Probability}

\input{Appendix_StatDist}

\input{Appendix_orderCTMCalgorithm}
\input{Appendix_DeviationMatrix}

\input{Appendix_1DModel}

\input{Appendix_2DModel}
\input{Appendix_3DModel}

\input{Appendix_4DModel}

\vfill\break

\end{appendices}

\end{document}

%% file: Introduction.tex
\section{Introduction}
\label{sec:introduction}



	        
	        
	    
	
\subsection{Background}

	   Epigenetic cell memory, the inheritance of gene expression patterns across subsequent cell divisions \cite{Blakey}, is a critical property of multi-cellular organisms of intense interest in the field of systems biology \cite{Ptashne2013, Ferrell}. It has previously been discovered that chromatin modifications, such as DNA methylation and histone modifications, are key mediators of epigenetic cell memory \cite{Carey2013, EpiReview2015, Allis,Kim2017} (see references in \cite{BrunoDelVecchio} for more biological background). More precisely, it was found via simulations of stochastic models that the time scale separation between establishment (fast) and erasure (slow) of these modifications extends the duration of cell memory, and that different asymmetries between erasure rates of chromatin modifications can lead to different gene expression patterns \cite{ECC2022,CDC2022,BrunoDelVecchio}. Here, we provide a mathematical framework to rigorously validate these computational findings and to further explore models of chromatin modification circuits. We do this in a way that the results obtained and the tools developed can be applied to other mathematical models beyond the applications in this paper, especially stochastic models of chemical reaction networks.
%


\subsection{Focus of our work}	   
        In this paper, we consider different versions of the chromatin modification circuit proposed in \cite{BrunoDelVecchio}. In particular, we start with simpler circuits that include histone modifications only and then we consider more elaborate circuits that include also DNA methylation. All of these circuits can be viewed as examples of Stochastic Chemical Reaction Networks (SCRNs). A SCRN is a continuous time Markov chain living in the non-negative integer lattice in $d$-dimensions, where the components of the Markov chain track the number of molecules of each of $d$ species in the network over time, and each jump of the Markov chain corresponds to the firing of a reaction in the network \cite{Kurtz:11}. A more precise description is given in Section \ref{sec:StochasticChemicalReactionNetworks}.

        In order to analyze these stochastic models, we first determine how the stationary distributions and mean first passage times between states vary when a small parameter $\eps$ (non-dimensional parameter that scales the speed of the basal erasure of all the chromatin modifications) tends to zero. To this end, we show that the stationary distributions and the mean first passage times of these singularly perturbed Markov chains admit series expansions in $\eps$ and we develop theoretical tools to determine the coefficients in these expansions.  Then, we focus on determining how the different erasure rates of chromatin modifications affect the behavior of the chromatin modification circuit models. This latter study is conducted by exploiting comparison theorems for Markov chains recently developed in \cite{Monotonicitypaper}.
        
One of the key features of our work is that these tools and the associated mathematical results are not only applicable to the chromatin modification models, but they can also be used to analyze other models that respect the same set of assumptions.

\subsection{Related work}

As mentioned in the previous paragraph, the stochastic behavior of the chromatin modification circuit models can be described by singularly perturbed continuous time Markov chains. There is some literature on discrete and continuous time, singularly perturbed Markov chains, especially by Avrachenkov et al. \cite{AvrachenkovFilarHowlett}, Hassin \& Haviv \cite{HH}, Beltrán and Landim \cite{BeltranLandim,BeltranLandim2}, and Yin \& Zhang  \cite{yinzhang2013}. Avrachenkov et al. \cite{AvrachenkovFilarHowlett} gave general characterizations of series expansions for the stationary distribution and mean first passage times of a singularly perturbed \textit{discrete time} Markov chain with finite state space. While their theory can be in principle translated to continuous time Markov chains, our work mostly deals directly with the singularly perturbed continuous time Markov chains and provides more concrete theoretical results for the leading coefficients of the stationary distribution series expansion and the orders of the poles of the mean first passage times. For the leading coefficients in the series expansion for the mean first passage times, we use in part the results of Avrachenkov \& Haviv \cite{AvrachenkovHaviv} and Avrachenkov et al.
\cite{AvrachenkovFilarHowlett} and adapt their work to the continuous time Markov chain setting. We treat in detail the case where the chain for $\eps = 0$ has more than one absorbing state and at least one transient state.
Furthermore, we also provide an interpretation of leading coefficients in the series expansion of the stationary distribution in terms of a certain restricted Markov chain. An algorithm we give to determine the order of the pole of the mean first passage time extends the work of Hassin \& Haviv \cite{HH} from discrete time to continuous time. We also extend the original algorithm’s scope to treat mean first passage times to a subset of states, instead of just a single state. 
Beltrán and Landim \cite{BeltranLandim,BeltranLandim2} study metastable and tunneling behavior for a sequence $\{\eta^N\}_{N=1}^\infty$ of  time-homogeneous continuous time Markov chains with  countable state spaces. Under an acceleration of time by a factor $\theta_N$, 
they give conditions under which the trace of the accelerated process on the metastates is asymptotically Markovian as $N \rightarrow \infty$. For our case, this would correspond to accelerating time for $\eta^N=X^\eps$ by $\theta_N \approx \frac{1}{\eps}$. Beltrán and Landim identified the transition rates for the limiting Markov chain and proved that its stationary distribution can be obtained as a limit from the stationary distribution for $\eta^N$. While this work is potentially related to what we did, it requires knowing the stationary distribution for $\eta^N$ {\it a priori}. Our approach does not need to know that stationary distribution explicitly and we also study mean first passage times, giving explicit asymptotics for both.
Finally, Yin \& Zhang \cite{yinzhang2013} conducted an extensive study focused on determining matched asymptotic expansions for the marginal distributions at time $t$ of singularly perturbed continuous time Markov chains. Their infinitesimal generators, generalizing those of Phillips \& Kokotovic \cite{Phillips1981} and Pan \& Basar \cite{Pan1999}, are of the form $Q(\eps)=\frac{1}{\eps}Q^{(0)}+Q^{(1)}$, and can be time dependent. 
For the time independent case, this would correspond to studying the marginal distributions of our Markov chain $X^{\eps}$ in the "linear" case and at time $\frac{t}{\eps}$ as $\eps \rightarrow 0$, i.e., $\lim_{\eps \rightarrow 0}X^{\eps}(\frac{t}{\eps})$.
Thus, while their work potentially might provide information about stationary distributions as $\eps \rightarrow 0$, we directly study the power series expansion (in $\eps$) of the stationary distribution of $X^{\eps}$, and we also study series expansions of mean first passage times for $X^{\eps}$, and we develop more concrete analyses for both.

\subsection{Outline of the paper}
    
    In Section \ref{motivexamples} we introduce two simplified models for the chromatin modification circuit that do not include DNA methylation. Through these examples, we introduce the mathematical setting and questions we address in this paper. We describe the basic setup and definitions needed for this paper in Section \ref{sec:BasicDefinitions}. We present our main results in Section \ref{sec:MainResults}. Some proofs are given there, whilst others are in the Supplementary Information (SI). Further applications of the theoretical tools developed in Section \ref{sec:MainResults} for chromatin modification circuits that include DNA methylation are presented in Section \ref{sec:Applications}. Concluding remarks are given in Section \ref{sec:conclusion}.

\subsection{Preliminaries and notation}
\label{sec:PreliminariesAndNotation}

    Denote the set of integers by $\Z$. For an integer $d \geq 2$ we denote by $\Z^d$ the set of $d$-dimensional vectors with entries in $\Z$. Denote by $\Z_+ = \{0,1,2, \ldots \}$, the set of non-negative integers. For an integer $d \geq 2$ we denote by $\Z_+^d$ the set of $d$-dimensional vectors with entries in $\Z_+$. We denote by $\one$ a vector of any dimension where all entries are $1$'s. The size of $\one$ will be understood from the context. The set of real numbers will be denoted by $\R$, $\R_+=[0,\infty)$, $\R_{>0}=(0,\infty)$, and $d$-dimensional Euclidean space will be denoted by $\R^d$ for $d \geq 2$. For integers $n,m \geq 1$, the set of $n \times m$ matrices with real-valued entries will be denoted by $\R^{n \times m}$. The set of complex numbers will be denoted by $\C$.
    
    Let $\X$ be a finite set. If needed, we will enumerate the entries of $\X$ by $\{1,\ldots,|\X|\}$. For a matrix $A = (A_{x,y})_{x,y \in \X}$ with real-valued entries, we denote the kernel of $A$ by $\ker(A) := \{ x \in \R^{|\X|}:\: Ax=0 \}$ and the nullity of $A$ by $\nullity(A) := \dim(\ker(A))$. We denote the spectrum of $A$ by $\spec(A)$ and the spectral radius by $\spr(A) = \max\{|\lambda|:\: \lambda \in \spec(A) \}$. A matrix $Q = (Q_{x,y})_{x,y \in \X}$ will be called a $Q$-matrix if $Q_{x,y} \geq 0$ for every $x \neq y \in \X$ and $Q\one = 0$. We denote the identity matrix, which has $1$'s on the diagonal and zeros elsewhere, by $I = (I_{x,y})_{x,y \in \X}$. For a vector $v=(v_x)_{x \in \X}$ we denote by $\diag((v_x)_{x \in \X})$ the diagonal matrix in $\X$ with entries given by $v$. 
Vectors are column vectors unless indicated otherwise and a superscript of $T$ will denote the transpose of a vector or matrix.
    For integers $n,m \geq 1$ and a matrix $A \in \R^{n \times m}$, we denote by $\norm{A} = (\sum_{i=1}^n\sum_{j=1}^m|A_{i,j}|^2)^{1/2}$ the Frobenius norm of $A$\footnote{Here, we chose to fix a particular norm on $\R^{n \times m}$, although other choices of norm will often work.}. For a vector $v \in \R^n$, we denote the Euclidean norm of $v$ by $\norm{v}=(\sum_{i=1}^n|v_{i}|^2)^{1/2}$.

    \begin{definition}
    \label{def:RealAnalyticPerturbation}
    Given a matrix $A^{(0)}$ in $\R^{n \times m}$, a \textbf{real-analytic perturbation} of $A^{(0)}$ is a matrix-valued function $A:[0,\eps_0) \longrightarrow \R^{n \times m}$, where $\eps_0 > 0$, and
        \begin{equation}
        \label{eq:AepsExpansion}
        A(\eps) = \sum_{k=0}^{\infty} \eps^kA^{(k)}, \qquad 0 \leq \eps < \eps_0,
        \end{equation}
    in which $\{A^{(k)} :\: k \geq 0\}$ is a sequence of matrices in  $\R^{n \times m}$ such that
        \begin{equation}
        \label{eq:AbsoluteConvergenceMatrixSeries}
        \sum_{k = 0}^{\infty} \eps^k\norm{A^{(k)}} < \infty, \qquad \text{for every } 0 \leq \eps < \eps_0.
        \end{equation}
    Such a perturbation is called \textbf{linear} if $A(\eps) = A^{(0)}+\eps A^{(1)}$ for $0 \leq \eps < \eps_0$. 
    \end{definition}
    
    
    By \eqref{eq:AbsoluteConvergenceMatrixSeries}, a real-analytic perturbation of $A^{(0)}$ can be extended to a function $F(z) := \sum_{k=0}^{\infty} z^kA^{(k)}$ defined on $B(0,\eps_0) = \{z\in \C \::\: |z| < \eps_0 \}$. The function $F$ will be called an \textbf{analytic perturbation} or \textbf{complex-analytic perturbation} of $A^{(0)}$. This extension will allow us to invoke results in complex analysis in order to study real-analytic perturbations. An example of this is the following result.
  
    \begin{proposition}
    \label{prop:InverseOfAnalytic}
     Let $A:[0,\eps_0) \longrightarrow \R^{n \times n}$ be a real-analytic perturbation of $A^{(0)}$ such that $A^{-1}(\eps)$ exists for every $0 < \eps < \eps_0$. Then, there is $\eps_1 \in (0,\eps_0)$ and $p \in \Z_+$ such that
        \begin{equation}
        \label{eq:MatrixValuedLaurentSeries}
        A^{-1}(\eps) = \sum_{k=-p}^{\infty} \eps^kB^{(k)}, \quad  0 < \eps < \eps_1,
        \end{equation}
        where $\sum_{k = -p}^{\infty} \eps^k\norm{B^{(k)}} < \infty$ for every $0 < \eps < \eps_1$, $\{B^{(k)}:\: k \geq -p\}$ is a sequence of matrices in $\R^{n \times n}$, $B^{(-p)}$ is not the identically zero matrix and $p$ is called the \textbf{order of the pole} at $\eps=0$.
    \end{proposition}
	
	This result is given in the analytic setting as Theorem 2.4 in \cite{AvrachenkovFilarHowlett}. Proposition \ref{prop:InverseOfAnalytic} follows by extending $A(\cdot)$ to a complex disk, then using Cramer's rule as in the proof of Theorem 2.4 in \cite{AvrachenkovFilarHowlett} and checking that the matrices $\{B^{(k)}: \: k \geq -p\}$ obtained are real-valued.


%% file: Motivatingexamples.tex
\section{Motivating Example: Chromatin Modification Circuit}
\label{motivexamples}
	
	In order to understand how the interactions among known chromatin modifications influence epigenetic cell memory, we consider the	chemical reaction model of the gene’s inner chromatin modification circuit introduced in \cite{BrunoDelVecchio}. This model has the nucleosome with DNA wrapped around it, D, as a basic unit that can be modified either with activating marks, such as H3K4 methylation (H3K4me3) or H3K4 acetylation (H3K4ac), or repressive marks, such as H3K9 methylation (H3K9me3) or DNA methylation. H3K4me3 and H3K4ac are two histone modifications that promote a less compact DNA around the nucleosomes and they are then associated with gene activation (see Chapter 3 of \cite{Allis} and \cite{H3K4meMaintenance2015}). In the model, it is assumed that H3K4me3 and H3K4ac co-exist and the nucleosome with either of these modifications is represented by $\mathrm{D^A}$. On the contrary, both the histone modification H3K9me3 and DNA methylation cause the DNA to be tightly wrapped around the nucleosome and therefore, they are associated with gene repression \cite{Blakey}. A nucleosome with DNA methylation only, H3K9 methylation (H3K9me3) only or both is represented by $\mathrm{D^R_1}$, $\mathrm{D^R_2}$ and $\mathrm{D^R_{12}}$, respectively.
	
	One of the key parameters of the system is $\eps>0$, a non-dimensional parameter that scales the speed of basal erasure of all chromatin modifications. We are interested in studying the behavior of the system in the limiting regime $\eps \rightarrow 0$, in which the chromatin modification system has a bimodal limiting stationary distribution \cite{BrunoDelVecchio}. One peak corresponds to the active chromatin state (most of the nucleosomes are modified with activating marks) and the other one is in the repressed chromatin state (most of the nucleosomes are modified with repressive marks). We aim to derive formulas that characterize, as $\eps$ goes to $0$, the behavior of the stationary distribution and the ``time to memory loss'' of the active (repressed) state, defined as the mean first passage time to reach the repressed (active) state, starting from the active (repressed) state. 
	
	Two other critical parameters of the system are $\mu$ and $\mu'$: they capture the asymmetry between the erasure rates of repressive and activating chromatin modifications. More precisely, $\mu$ ($\mu'$) quantifies the asymmetry between erasure rates of repressive histone modifications (DNA methylation) and activating histone modifications. Part of our study is to analytically determine how $\mu$ and $\mu'$ affect the stationary distribution and the time to memory loss of the active and repressed states.
	
	In this section, we introduce two simplified models of the chromatin modification circuit in which, compared to the full model described above, DNA methylation is not included and the only chromatin marks are histone modifications. We will use these simpler models in Section \ref{sec:MainResults} to directly apply and then better understand the theory developed in this paper. Then, in Section \ref{sec:Applications} we deal with more elaborate models that also include DNA methylation.
Note that, for consistency, we use the same notation for the species and the reaction rate constants as the one used in the paper where these models were introduced \cite{BrunoDelVecchio}.
	
	
	\subsection{1D model}
	\label{SS1}


	We first consider a simplified model in which a gene has a total of $\Dtot \ge 2$ nucleosomes, where each nucleosome either has an activating histone modification, $\mathrm{D^A}$, or a repressive histone modification, $\mathrm{D^R}$, and there are no unmodified nucleosomes in this simplified model. If the amounts of nucleosomes having repressive ($\mathrm{D^R}$) and activating ($\mathrm{D^A}$) modifications are denoted as $n_{\mathrm{D^R}}$ and $n_{\mathrm{D^A}}$, respectively, then we have the conservation law $n_{\mathrm{D^R}}+n_{\mathrm{D^A}}=\Dtot$. We call this the 1D model because it suffices to keep track of the amount of $\mathrm{D^R}$ (for example), since the amount of $\mathrm{D^A}$ can be deduced by the conservation law. Furthermore, the basal and recruited erasure of $\mathrm{D^A}$ ($\mathrm{D^R}$) coincide with the basal {\it de-novo} establishment and maintenance of $\mathrm{D^R}$ ($\mathrm{D^A}$). The chemical reaction system for this 1D model is the following:
\begin{equation}\label{eq:CRN_1D_Model}
	\begin{aligned}
 &{\large \textcircled{\small 1}}\;\ce{D^A + D^{R} ->[$k^A_E$] D^R + D^{R}},\;\;\;{\large \textcircled{\small 2}}\;\ce{D^A ->[$\delta +\bar k^A_E$] D^{R}},\\
 &{\large \textcircled{\small 3}}\;\ce{D^R + D^{A} ->[$k^R_E$] D^{A} + D^{A}},\;\;\;{\large \textcircled{\small 4}}\;\ce{D^R ->[$\delta+\bar k^R_E$] D^{A}},
	\end{aligned}
 \end{equation}
    where $\delta, k^A_E,\bar k^A_E, k^R_E,\bar k^R_E > 0$. Here, the form of the reaction rate constants is due to the fact that reactions with the same reactants and products have been combined. We denote the reaction volume by $V$, and let $\eps := \frac{\delta + \bar k_E^A}{k_E^A(\Dtot/V)} = \frac{\delta_A}{k_E^A(\Dtot/V)}$, where $\delta_A:=\delta + \bar k_E^A$. We also consider the constant $\mu := \frac{k^R_E}{k^A_E}$, which captures the asymmetry between the erasure rates of repressive and activating histone modifications.We introduce the constant $b$ such that $\mu b = \frac{\delta_{R}}{\delta_{A}}$, with $\delta_R:=\delta + \bar k_E^R$. Then, $\delta_{A}=\eps\frac{k^A_E \Dtot}{V}$ and $\delta_R := \delta_A\mu b = \eps \frac{k_E^A\Dtot}{V}\mu b$. So, as $\eps \to 0$, both $\delta_A$ and $\delta_R$ go to $0$, with $\Dtot, \frac{k_E^A}{V}, \mu, $ and $b$ fixed.     
        	        \begin{figure}[t]
            \centering
            \includegraphics[scale=0.42]{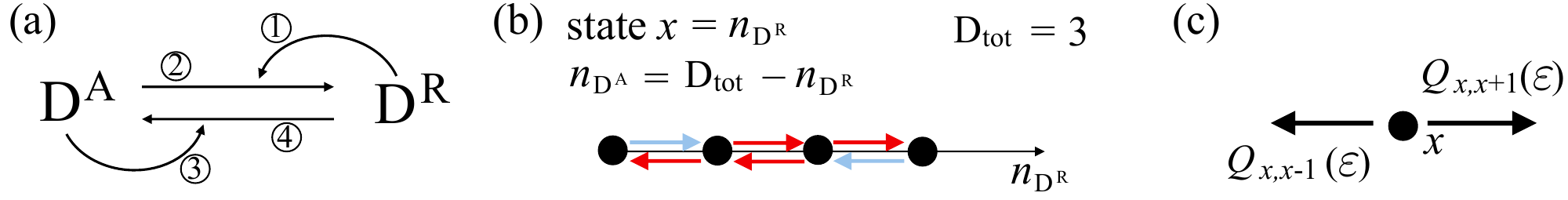}
            \caption{\small { \bf 1D model.} (a) Chemical reaction system. The numbers on the arrows correspond to the chemical reactions associated with the arrows as described in (\ref{eq:CRN_1D_Model}) in the main text. (b) Markov chain graph. Here, we consider $\mathrm{D_{tot}}=3$ and we use black dots to represent the states, red arrows to represent transition rates that are $O(1)$, and blue arrows to represent transition rates that are $O(\eps)$. (c) Directions of the potential transitions of $X^{\eps}$ starting from a state $x$, whose rates are given in equation (\ref{eq:Qmatrix_eps}).
            }
   \label{fig:1Dmodel}
        \end{figure}
Now, consider a continuous time Markov chain $X^{\eps}$, with state space $\X := \{0,\ldots,\Dtot\}$, where $\Dtot \geq 2$ is an integer, which keeps track of $n_{\mathrm{D^R}}$ through time. Given that we have the conservation law $n_{\mathrm{D^R}}+n_{\mathrm{D^A}}=\Dtot$, $n_{\mathrm{D^A}}$ can be obtained as a function of $n_{\mathrm{D^R}}$, that is $n_{\mathrm{D^A}}=\Dtot-n_{\mathrm{D^R}}$. Assuming stochastic mass-action kinetics (including the usual volume scaling of rate constants \cite{BFS}), the infinitesimal generator $Q(\eps)$ \footnote{Note that $Q(\eps)$ is sometimes called an infinitesimal transition matrix. The entries $Q_{x,y}(\eps)$ for $x \neq y$ are the infinitesimal transition rates of going from $x$ to $y$: $\PP[X^{\eps}(t+h) = y | X^{\eps}(t)=x] = Q_{x,y}(\eps)h +o(h)$ as $h \to 0$.} for $X^{\eps}$ is given by:
    \begin{equation}
    \label{eq:Qmatrix_eps}
    Q_{x,x+\l}(\eps) = \begin{cases} \left(\frac{k^A_E}{V}x+\eps \frac{k^A_E}{V} \mathrm{D_{tot}}\right)(\Dtot - x) & \text{ if } \l=1 \\
    \mu \left(\frac{k^A_E}{V}(\Dtot - x)+b \eps \frac{k^A_E}{V} \mathrm{D_{tot}}\right)x  & \text{ if } \l=-1 \\
               0 & \text{otherwise,} 
         \end{cases}
    \end{equation}
     for $x \in \X$, $\l \in \Z \setminus \{0\}$ and $x+\l \in \X$, and $Q_{x,x}(\eps)=-\sum_{y\in \X \setminus \{x\}}Q_{x,y}(\eps)$ for $x \in \X$. We extend this definition to $\eps=0$ by defining $Q_{x,y}(0) := \lim_{\eps \to 0} Q_{x,y}(\eps)$ for $x,y \in \X$. We will follow a similar convention for other examples. We consider $X^{0}$ to be the continuous time Markov chain with infinitesimal generator given by $Q(0)$. The process $X^{0}$ corresponds to a SCRN model associated with the autocatalytic reactions $\large \textcircled{\small 1}$ and $\large \textcircled{\small 3}$ in \eqref{eq:CRN_1D_Model}, alone. Note that
        \begin{align}
        \label{eq:LinearPerturbationQMatrices_1D}
            Q(\eps) = Q^{(0)} + \eps Q^{(1)}, \qquad  \eps \geq 0,
        \end{align}
            for appropriate matrices $Q^{(0)}$ and $Q^{(1)}$ in $\R^{|\X|\times|\X|}$. By \eqref{eq:LinearPerturbationQMatrices_1D}, we can see that $Q(\cdot)$ is a \textit{real-analytic} (and moreover \textit{linear}) perturbation of $Q^{(0)}$ (see Section \ref{sec:PreliminariesAndNotation} for definitions). Note that for every $\eps > 0$, $X^{\eps}$ is irreducible, while $X^{0}$ has a transient communicating class $\{1, \ldots, \Dtot-1\}$ and two absorbing states ($0$ and $\Dtot$) (see SI - Section \ref{sec:Appendix1Dmodel}). Because of this discontinuity at $\eps =0$, we say that $Q(\cdot)$ is a \textit{singular perturbation} of $Q^{(0)}$ (see Section \ref{sec:pertCTMC} for a precise definition).
    
    We first want to determine the probability for the gene to be in the active state $a$ ($x=0$), repressed state $r$ ($x=\mathrm{D_{tot}}$) or one of the intermediate states ($x \in \{1, \ldots, \Dtot-1\}$) after a long time (life-time of the organism), as a function of $\eps$.
    We are especially interested in the limit of the stationary distribution for the system, $\pi(\eps)$, as $\eps \rightarrow 0$ (i.e., the basal erasure rate of the chromatin modifications is much lower than their maintenance rate). Since $X^{\eps}$ is irreducible for $\eps >0$ (and it has a finite state space), it has a unique stationary distribution $\pi(\eps)$. In Section \ref{sec:pertCTMC} we show that $\pi(0) := \lim_{\eps \to 0} \pi(\eps)$ exists and the function $\pi(\cdot)$ admits a convergent power series expansion:
        \begin{equation}
        \label{eq:pi_expansion_1D}
        \pi(\eps) = \sum_{k = 0}^{\infty} \eps^k\pi^{(k)}\;\; \text{for } 0 \leq \eps < \eps_1,
        \end{equation}
    \noindent for some $\eps_1 > 0$. In order to determine $\pi(0)$, we can take limits and observe that $\pi(0)Q(0)=0$ and so $\pi(0)$ is a stationary distribution for $Q(0)$. Indeed, $\pi(0)$ is a specific mixture of atoms on the two absorbing states ($0$ and $\Dtot$) for $X^{0}$.
    
\begin{figure}
        \centering
        \includegraphics[scale=0.36]{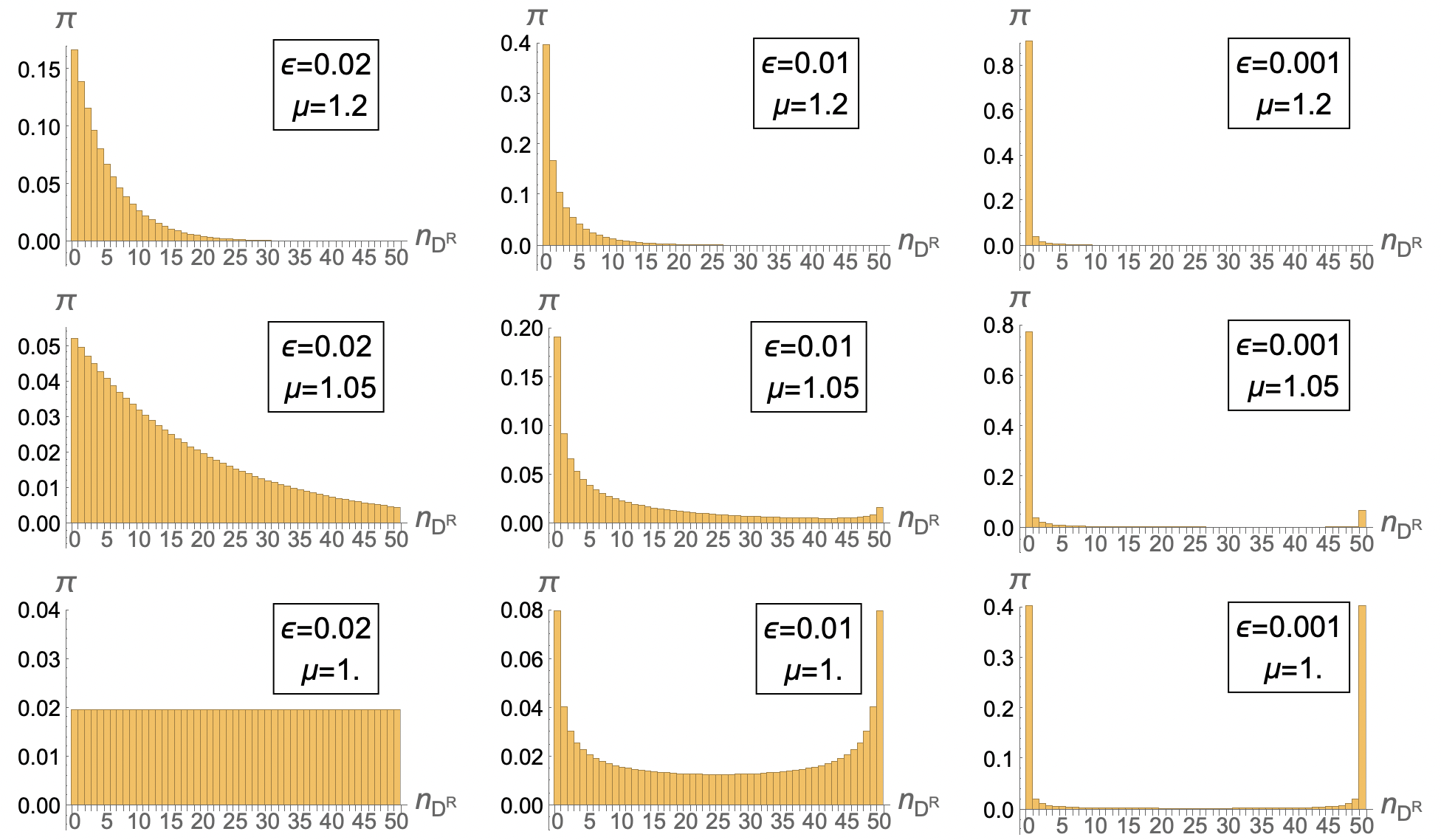}\\
        \caption{Histograms for the stationary distribution $\pi(\eps)$ of the Markov chain $X^{\eps}$ defined by \eqref{eq:Qmatrix_eps}, for different values of $\eps$ and $\mu$. The plot was generated by numerically solving $\pi(\eps)Q(\eps) = 0$ using the \textit{Eigenvector} function in Mathematica. The parameters used were $\Dtot = 50, k^A_E/V = 1,$ and $b = 1$. }
        \label{fig:StationaryDistribution1D}
    \end{figure}

    In Figure \ref{fig:StationaryDistribution1D} we see how the function $\pi(\eps)$ changes as $\eps \to 0$ for several values of $\mu$ with $\Dtot$, $\frac{k_E^A}{V}$ and $b$ fixed. Furthermore, for this simpler chromatin modification circuit, because of the birth-death structure of $X^{\eps}$, we can obtain explicit formulas for $\pi(\eps)$ when $\eps > 0$ (see SI - Section \ref{sec:Appendix1Dmodel}). On letting $\eps \to 0$, we obtain:
        \begin{equation}
        \label{pi01Dmodel}
         \pi_x(0) = \begin{cases}
                    \frac{b\mu^{\mathrm{D_{tot}}}}{1 + b\mu^{\mathrm{D_{tot}}}} & \text{ if } x=0 \\
                    0 & \text{ if } x \in \{1,\ldots, \Dtot-1\} \\
                    \frac{1}{1 + b\mu^{\mathrm{D_{tot}}}} & \text{ if } x = \Dtot.
        \end{cases}
        \end{equation}
    %
%
   \noindent Thus, $\pi_x(0)\ne 0$ only for $x=0$ and $x=\Dtot$ and $\pi_0(0)$ increases as $\mu$ increases, while $\pi_{\Dtot}(0)$ decreases as $\mu$ increases. 
    
    For continuous time Markov chains beyond the one-dimensional birth-death processes seen here, determining $\pi(0)$ will be a considerable task. In Section \ref{sec:StationaryDistributions}, we address the problem of determining $\pi(0)$, together with the whole expansion \eqref{eq:pi_expansion_1D}, in a systematic way, for a class of singularly perturbed Markov chains that includes our models of chromatin modification circuits. For the 1D model considered here, the derivation of the first two terms in the expansion is given in Section \ref{sec:IllustrativeStationaryDistribution}.

    Now, in order to evaluate the time to memory loss of the active and repressed states, let us define the first passage time as $\tau_y^{\eps} = \inf \{t \geq 0:\: X^{\eps}(t)=y\}$ for a state $y \in \X$. We will see in \eqref{eq:ExpansionTauEps} that the mean first passage time (MFPT) for $X^\eps$ starting from $x \in \X$, $h_{x,y}(\eps)=\E_x [\tau_y^{\eps}]$, has a Laurent series expansion of the form:
    \begin{equation}\label{TTMLexpansion}
         h_{x,y}(\eps) =   \frac{c_{-p}}{\eps^p} + \ldots + \frac{c_{-1}}{\eps} + c_0 + \eps c_1 + \ldots \qquad \text{for } 0 < \eps < \eps_{\{y\}},
    \end{equation}
    for some $\eps_{\{y\}} > 0$, for some natural number $p \geq 0$ and where $c_{-p} \neq 0$. Then, considering the repressed state $r=\mathrm{D_{tot}}$ and the active state $a=0$, we define the time to memory loss of the repressed state as $h_{r,a}(\eps)$ and the time to memory loss of the active state as $h_{a,r}(\eps)$. Now, we are interested in the derivation of analytical formulas for $h_{r,a}(\eps)$ and $h_{a,r}(\eps)$. This will allow us to understand how the time to memory loss changes as $\eps \rightarrow 0$, and how the asymmetry of the system, represented by $\mu$, affects this limit. For this case study, exploiting its birth-death structure, we can directly derive relevant formulas (see SI - Section \ref{sec:Appendix1Dmodel}, SI - Equations (\ref{formulaDto0})-(\ref{formula3})). In particular, defining $\lambda^{\eps}_x=Q_{x,x+1}(\eps)$, $\gamma^{\eps}_x=Q_{x,x-1}(\eps)$, with $Q_{x,x+1}(\eps)$ and $Q_{x,x-1}(\eps)$ defined in (\ref{eq:Qmatrix_eps}), and $r^{\eps}_j=\frac{\lambda^{\eps}_1 \lambda^{\eps}_2...\lambda^{\eps}_j}{\gamma^{\eps}_1 \gamma^{\eps}_2 ... \gamma^{\eps}_j}$, for $j=1,2,...,\Dtot -1$, the time to memory loss of the repressed state is given by
	\begin{equation}\label{formulaDto0INTRO}
	h_{r,a}(\eps) =\frac{r^{\eps}_{\mathrm{D_{tot}}-1}}{\gamma^{\eps}_{\mathrm{D_{tot}}}}\left(1+\sum_{i=1}^{\mathrm{D_{tot}}-1}\frac{1}{r^{\eps}_i}\right)+\sum_{i=2}^{\mathrm{D_{tot}}-1}\left[\frac{r^{\eps}_{i-1}}{\gamma^{\eps}_i}\left(1+\sum_{j=1}^{i-1}\frac{1}{r^{\eps}_j}\right)\right]+\frac{1}{\gamma^{\eps}_1}. 
	\end{equation}
	Similarly, defining $\tilde r^{\eps}_j=\frac{\gamma^{\eps}_{\mathrm{D_{tot}}-1}\gamma^{\eps}_{\mathrm{D_{tot}}-2}...\gamma^{\eps}_{\mathrm{D_{tot}}-j}}{\lambda^{\eps}_{\mathrm{D_{tot}}-1}\lambda^{\eps}_{\mathrm{D_{tot}}-2}...\lambda^{\eps}_{\mathrm{D_{tot}}-j}}$, for $j=1,2,...,\Dtot -1$, the time to memory loss of the active state is given by
	\begin{equation}\label{formula3INTRO}
	h_{a,r}(\eps)=\frac{\tilde r^{\eps}_{\mathrm{D_{tot}}-1}}{\lambda^{\eps}_{0}}\left(1+\sum_{j=1}^{\mathrm{D_{tot}}-1}\frac{1}{\tilde r^{\eps}_i}\right)+\sum_{i=2}^{\mathrm{D_{tot}}-1}\left[\frac{\tilde r^{\eps}_{i-1}}{\lambda^{\eps}_{\mathrm{D_{tot}}-i}}\left(1+\sum_{j=1}^{i-1}\frac{1}{\tilde r^{\eps}_j}\right)\right]+\frac{1}{\lambda^{\eps}_{\mathrm{D_{tot}}-1}}. 
	\end{equation}
	Since $\lambda^{\eps}_0$ and $\gamma^{\eps}_\mathrm{D_{tot}}$ are the only transition rates that are $O(\eps)$ with the rest being $O(1)$, the time to memory loss of both the active and repressed states are $O(\eps^{-1})$, that is, $p=1$, and as $\eps\rightarrow 0$, these mean times tend to infinity.
    
    Furthermore, $\gamma^{\eps}_{x}$, with $x\in \{1,2,..., \Dtot\}$, are the only rates that depend on $\mu$ (they are linear in $\mu$). Examining \eqref{formulaDto0INTRO} and \eqref{formula3INTRO} with this observation in mind, we see that, if $\mu$ is increased (that is, the erasure rate of the repressive histone modification is increased compared to that of the active histone modification), $h_{a,r}(\eps)$ increases, while $h_{r,a}(\eps)$ decreases. The opposite happens when $\mu$ is decreased.
    
    More complicated situations arise when we do not have a birth-death structure to work with, as in the model of the next example. To evaluate how critical system parameters affect the time to memory loss for such more elaborate systems, in Section \ref{sec:MainResults}, we develop an algorithm to determine $p$ (see Section \ref{orderMFPT}), we give an expression for the leading term in the series expansion of the mean first passage time, and we exploit theoretical results developed in our paper \cite{Monotonicitypaper} for comparing continuous time Markov chains, to determine how the asymmetry of the system affects the time to memory loss (see Section \ref{sec:Mon}). 
    

\begin{figure}[t]
            \centering
            \includegraphics[scale=0.42]{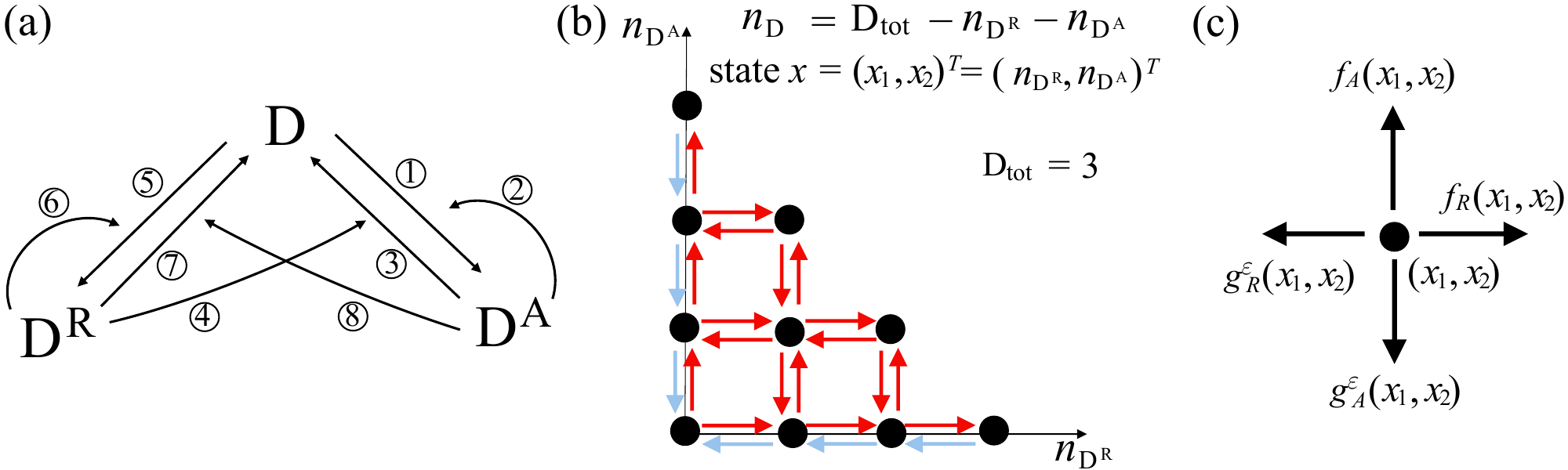}
            \caption{\small { \bf 2D model.} (a) Chemical reaction system.  The numbers on the arrows correspond to the reactions associated with the arrows as described in (\ref{reacs2D}) in the main text. (b) Markov chain graph. Here, we consider $\mathrm{D_{tot}}=3$ and we use black dots to represent the states, red arrows to represent transition rates that are $O(1)$, and blue arrows to represent transition rates that are $O(\eps)$. (c) Directions of the possible one step transitions for $X^{\eps}$ starting from a state $x=(x_1,x_2)^T$, whose rates are given in equation (\ref{rates2D}).
            }
   \label{fig:2Dmodel}
        \end{figure}
        
	\subsection{2D model}\label{exp2}
		
	Let us consider a model in which, compared to the previous one, we assume that a nucleosome can also be unmodified. More precisely, in this case we denote the number of nucleosomes unmodified (D), modified with repressive modifications ($\mathrm{D^R}$), and modified with activating modifications ($\mathrm{D^A}$) by $n_{\mathrm{D}}$, $n_{\mathrm{D^R}}$ and $n_{\mathrm{D^A}}$, respectively, and we have that $n_{\mathrm{D}}+n_{\mathrm{D^R}}+n_{\mathrm{D^A}}=\Dtot$, with $\Dtot$ representing the total number of nucleosomes within the gene. Furthermore, each histone modification autocatalyzes its own production and promotes the erasure of the other one \cite{BrunoDelVecchio,2007DoddCellPaper}. The chemical reaction system is the following:
    
    \begingroup
    \small
	\begin{equation}\label{reacs2D}
	\begin{aligned}
	&{\large \textcircled{\small 1}}\;\ce{D ->[$k^A_{W0}+k^A_W$] D^A },\;\;{\large \textcircled{\small 2}}\;\ce{D +D^A ->[$k^A_M$] D^A + D^A},\;\;{\large \textcircled{\small 3}}\;\ce{D^A ->[$\delta +\bar k^A_E$] D},\;\;{\large \textcircled{\small 4}}\;\ce{D^A + D^{R} ->[$k^A_E$] D + D^{R}},\\
	&{\large \textcircled{\small 5}}\;\ce{D ->[$k^R_{W0}+k^R_{W}$] D^R },\;\;{\large \textcircled{\small 6}}\;\ce{D +D^R ->[$k^R_M$] D^R + D^R},\;\;{\large \textcircled{\small 7}}\;\ce{D^R ->[$\delta+\bar k^R_E$] D},\;\;{\large \textcircled{\small 8}}\;\ce{D^R + D^{A} ->[$k^R_E$] D + D^{A}},\\
	\end{aligned}
	\end{equation}
	\endgroup
    where $k^A_{W0},k^A_W,k^A_M,\delta, \bar k^A_E,  k^A_E, k^R_{W0}, k^R_{W}, k^R_M, \bar k^R_E, k^R_E > 0$. Here, the form of the reaction rate constants is due to the fact
    that reactions with the same reactants and products have been combined. Now, similarly to what we did for the previous model, let us denote the reaction volume by $V$, and let $\eps := \frac{\delta + \bar k_E^A}{k_M^A(\Dtot/V)} = \frac{\delta_A}{k_M^A(\Dtot/V)}$, with $\delta_A:=\delta + \bar k_E^A$, and $\mu := \frac{k^R_E}{k^A_E}$. Additionally, consider the constant $b$ such that $\mu b = \frac{\delta_{R}}{\delta_{A}}$, with $\delta_R:=\delta + \bar k_E^R$. Then $\delta_R = \delta_A \mu b = \eps \frac{k_M^A\Dtot}{V} \mu b$. So, as $\eps \to 0$, both $\delta_A$ and $\delta_R$ go to $0$ with $\Dtot, \frac{k_M^A\Dtot}{V}, \mu$ and $b$ fixed. 
    
    We consider the continuous time Markov chain $X^{\eps} =\{(X_1^{\eps}(t),X_2^{\eps}(t))^T, t\ge0\}$, which keeps track of $(n_{\mathrm{D^R}},n_{\mathrm{D^A}})$ through time. Since the total number of nucleosomes $\Dtot$ is constant, the state space is $\mathcal{X}= \{x=(x_1,x_2)^T \in \Z_+^2 :\: x_1 + x_2 \leq \Dtot \}$. The potential one step transitions for $X^{\eps}$ from $x \in \X$ are shown in Figure \ref{fig:2Dmodel}(c), where the associated transition vectors are given by $v_1=-v_2=(0,1)^T$ and $v_3=-v_4=(1,0)^T$ and the infinitesimal transition rates (assuming mass-action kinetics with the usual volume scaling of rate constants) are given by

    \begingroup
    \small
    \begin{equation}\label{rates2D} 
    \begin{aligned}
&Q_{x,x+v_1}(\eps)=f_A(x) = (\Dtot -(x_1+x_2))\left(k_{W0}^A+k_{W}^A + \frac{k_{M}^A}{V}x_2\right),\\
&Q_{x,x+v_3}(\eps)=f_R(x) = (\Dtot -(x_1+x_2))\left(k_{W0}^R+k_{W}^R + \frac{k_{M}^R}{V}x_1\right),\\
&Q_{x,x+v_2}(\eps)=g^{\eps}_A(x) = x_2\left(\eps \frac{k_{M}^A}{V}\Dtot + x_1\frac{k^A_E}{V}\right),\;\;
Q_{x,x+v_4}(\eps)=g^{\eps}_R(x) = x_1\mu\left(\eps \frac{k_{M}^A}{V}\Dtot b + x_2\frac{k^A_E}{V}\right).
       \end{aligned}
    \end{equation}
\endgroup
%

\noindent   This is a more complicated model compared to the previous example and, in order to study its stationary distribution and mean first passage times, we will exploit the theory developed in this paper, as shown in Section \ref{sec:MainResults}.

%% file: AnalyticalpertCTMC.tex

\section{Basic Setup and Definitions}
\label{sec:BasicDefinitions}

    In Section \ref{sec:pertCTMC} we provide basic definitions for singularly perturbed continuous time Markov chains and describe some key properties for them. In particular, we describe the form of series expansions for their stationary distributions and mean first passage times. We will study these quantities and apply our results to a class of continuous time Markov chains called Stochastic Chemical Reaction Networks (SCRNs) which are defined in Section \ref{sec:StochasticChemicalReactionNetworks}. Our models of chromatin modification circuits will be SCRNs. All of the models considered will have a finite state space.

\subsection{Singularly perturbed, finite state, continuous time Markov chains}
\label{sec:pertCTMC}

    Suppose $\mathcal{X}$ is a finite set and $|\mathcal{X}| > 1$. For a value $\eps_0 > 0$, consider a family $\{X^{\eps}:\: 0 \leq \eps < \eps_0\}$ of continuous time Markov chains with state space $\X$ and infinitesimal generators $\{Q(\eps):\: 0 \leq \eps < \eps_0\}$ where $\eps \mapsto Q(\eps)$ is a real-analytic perturbation of $Q(0)$. Thus,
	\begin{equation}
	\label{eq:QmatrixPerturbedAPMC}
	Q(\eps) = Q^{(0)} + \eps Q^{(1)} + \eps^2 Q^{(2)} + \cdots,    
	\end{equation}
	where $\{Q^{(k)}: \: k \geq 0\}$ is a family of $|\X|\times|\X|$ real-valued matrices such that $\sum_{k=0}^\infty \eps^k\norm{Q^{(k)}} < \infty$ for every $0 \leq \eps < \eps_0$. Assume that the continuous time Markov chains $X^{\eps}$ are irreducible for $0 < \eps < \eps_0$. In this context, the perturbation is \textbf{singular} when $X^{0}$ has more than one recurrent class. This notion of singular will be the focus of our attention although some of our work applies for the regular (non-singular) case too. All of our chromatin modification circuit models have associated singular continuous time Markov chains, where the perturbation is \textbf{linear}, i.e., $Q^{(k)} = 0$ for every $k \geq 2$.
	
	When $0 < \eps < \eps_0$, there is an equivalent characterization of $X^{\eps}$ using holding times with exponential parameters $\{q_x (\eps)\}_{x \in \mathcal{X}}$ and a transition probability matrix $P(\eps)$ for the \textbf{embedded discrete time Markov chain}. Specifically, for each $x \in \mathcal{X}$, $q_x(\eps) = - Q_{x,x}(\eps) \neq 0$, since $X^{\eps}$ is irreducible, and for all $x,y \in \X$, $P_{x,x}(\eps) = 0$, $P_{x,y}(\eps) = \frac{Q_{x,y}(\eps)}{q_x(\eps)}$, for $y \neq x$ in $\X$. Note that $Q(\eps)=\diag(q(\eps))(P(\eps)-I)$. The matrix $P(\eps)$ has a power series expansion in $\eps$ for sufficiently small $0 \leq \eps < \eps_{P}$ for some $\eps_P > 0$ (the justification is similar to that for \eqref{eq:ExpansionPiEps} below). 
 
	The first quantities we are interested in studying are mean first passage times. Consider a nonempty set $\B \subseteq \X$ such that $\B \neq \X$ and let 
    \begin{equation*}
    \tau_\B^{\eps} := \inf \{t \geq 0: X^{\eps}(t) \in \B\}.
    \end{equation*}
    We define the \textbf{mean first passage time (MFPT)} (for $X^{\eps}$) from $x \in \X$ to $\B$ as 
    \begin{equation*}
    h_{x,\B}(\eps)=\E[\tau_\B^{\eps} \;|\; X^{\eps}(0) = x].
    \end{equation*}
    If $\B = \{y\}$ for some $y \in \X$, we adopt the notation: $h_{x,y}(\eps) := h_{x,\{y\}}(\eps)$. Using first step analysis (see (3.1) in \cite{Norris}), for $0 < \eps < \eps_0$,
	\begin{equation}
	h_{x,\B}(\eps)=\begin{cases}\label{MFPTsystem}
	0 & \mbox{if } x \in \B \\ 
	\frac{1}{q_x (\eps)}+\sum_{y\in \mathcal{X}}P_{x,y}(\eps)h_{y,\B}(\eps) & \mbox{if } x \in \B^c.
	\end{cases}
	\end{equation}
    Now, define $P^{\Bc} (\eps)$ and $Q^{\Bc}(\eps)$ as the matrices obtained by removing the columns and rows of $P(\eps)$ and $Q(\eps)$, respectively, corresponding to states in $\B$. Then, by noting that $I-P^{\Bc}(\eps)$ is invertible (see SI - Lemma \ref{lem:IMinusPInvertible}) and that $Q^{\Bc}(\eps)=-\diag((q_x(\eps))_{x \in \B^c})(I - P^{\Bc}(\eps))$ is invertible, from \eqref{MFPTsystem}, we obtain 
        \begin{equation}
        \label{eq:h_A_with_Q_inverse}
         h_{\B}(\eps) = -(Q^{\Bc}(\eps))^{-1}\one,
        \end{equation}
    where $h_{\B}(\eps) := (h_{x,\B}(\eps))_{x \in \B^c}$, $I$ is the identity matrix of dimension $|\B^c|$, and $\one$ is the vector of all $1$'s, of size $|\B^c|$. Proposition \ref{prop:InverseOfAnalytic}, yields that there is $0 < \eps_{\B} < \eps_0$ such that $-(Q^{\Bc}(\eps))^{-1}$ can be expanded as a matrix-valued Laurent series as in \eqref{eq:MatrixValuedLaurentSeries} for $0 < \eps < \eps_{\B}$, and then for each $x \in \B^c$,
        \begin{equation}
        \label{eq:ExpansionTauEps}
        \E_x[\tau_\B^{\eps}] = h_{x,\B}(\eps)= \sum_{k=-p(x)}^{\infty} \rho_x^{(k)}\eps^k, \quad 0 < \eps < \eps_{\B}, 
        \end{equation}
    where $p(x) \geq 0$ is an integer, $\rho_x^{(-p(x))} > 0$, $\rho_x^{(k)}\in \R$ for $k>-p(x)$, and the convergence is absolute convergence for $0 < \eps < \eps_{\B}$. The quantity $p(x)$ will be called the \textbf{order of the pole} of \eqref{eq:ExpansionTauEps}. In Section \ref{orderMFPT} we will show how to find $p(x)$ by using an algorithm that uses the order, with respect to $\eps$, of the transitions of the Markov chain $X^{\eps}$.

 


	
	A second quantity of interest is the stationary distribution for $X^\eps$. For $0 < \eps < \eps_0$, since $X^{\eps}$ is assumed to be irreducible and has finite state space, there is a unique \textbf{stationary distribution} $\pi(\eps) = (\pi_x (\eps))_{x \in \mathcal{X}}$, which is the unique probability row vector satisfying $\pi(\eps)Q(\eps) = 0$. We are interested in studying $\pi(\eps)$ as $\eps \to 0$. For this, first consider $\eta_{x}^{\eps} = \inf \{t \geq 0: X^{\eps} (t) \neq x\}$ and $\zeta_{x}^{\eps} = \inf \{t > \eta_{x}^{\eps}: X^{\eps} (t) = x\}$, $x \in \X$. Note that $\E_y[\zeta_x^{\eps}] = h_{y,x}(\eps)$ for $y \neq x$. For each $x \in \mathcal{X}$, $\E_x[\zeta_x^{\eps}]$ is called the \textbf{mean return time} to the state $x$, and for $0 < \eps < \eps_0$ satisfies
	\begin{equation}
	\label{eq:ReturnTimes_expansion}
	\E_x[\zeta_x^{\eps}] = \frac{1}{q_x (\eps)} + \sum_{y \neq x} P_{xy} (\eps) \E_y [\zeta_x^{\eps}] = \frac{1}{q_x (\eps)} + \sum_{y \neq x} P_{xy} (\eps)h_{y,x}(\eps).    
	\end{equation}
	It is well known (see Theorem 3.8.1 in \cite{Norris}) that for $0 < \eps < \eps_0$,
	\begin{equation}
	\label{eq:Pi_related_to_returntimes}
	\pi_x (\eps) = \frac{1}{\E_x[\zeta_x^{\eps}]} \cdot \frac{1}{q_x (\eps)}, \quad x \in \mathcal{X}.    
	\end{equation}

   From \eqref{eq:ExpansionTauEps} and \eqref{eq:ReturnTimes_expansion}, we can see that $\eps \mapsto q_x(\eps)\E_x[\zeta_x^{\eps}]$ can be extended to an analytic function on a punctured disk about $0$ in $\C$, with a Laurent series expansion having at most a pole of finite order at $0$. The radius of the punctured disk may be smaller than $\eps_0$. This, together with \eqref{eq:Pi_related_to_returntimes}, implies that $\eps \mapsto \pi_x (\eps)$ can be extended to an analytic function on a punctured disk about $0$ in $\C$, also with a Laurent series expansion. Since this function is bounded by one when restricted to sufficiently small positive values of $\eps$, we can remove the singularity at $0$ and obtain that $\pi(0) := \lim_{\eps \to 0} \pi(\eps)$ exists and furthermore $\eps \mapsto \pi(\eps)$ is a real-analytic perturbation of $\pi(0)$. In other words,
        \begin{equation}
        \label{eq:ExpansionPiEps}
        \pi(\eps) = \sum_{k = 0}^{\infty} \eps^k\pi^{(k)}, \qquad 0 \leq \eps < \eps_1,
        \end{equation}
    for sufficiently small $\eps_1 > 0$ and where $\{\pi^{(k)}: \: k \geq 0 \}$ is a sequence of real-valued $|\X|$-dimensional vectors such that $\sum_{k = 0}^{\infty} \eps^k\norm{\pi^{(k)}} < \infty$ for every  $0 \leq \eps < \eps_1$.

\subsection{Stochastic Chemical Reaction Networks (SCRNs)}
\label{sec:StochasticChemicalReactionNetworks}

    In this section, we provide some background on Stochastic Chemical Reaction Networks. The reader is referred to Anderson \& Kurtz \cite{AndersonKurtzBook} for a more in depth introduction to this subject.

    We assume there is a finite non-empty set $\Ss = \{\mathrm{S}_1,\ldots,\mathrm{S}_d\}$ of $d$ \textbf{species}, and a finite non-empty set $\Rs \subseteq \Z_+^d \times \Z_+^d$ that represents chemical \textbf{reactions}. We assume that $(w,w) \notin \Rs$ for every $w \in \Z^d_+$. The set $\Ss$ represents $d$ different molecular species in a system subject to reactions $\Rs$ which change the number of molecules of some species. For each $(v^{-},v^+) \in \Rs$, the $d$-dimensional vector $v^{-}$ (the \textbf{reactant vector}) counts how many molecules of each species are consumed in the reaction, while $v^{+}$ (the \textbf{product vector}) counts how many molecules of each species are produced. The reaction is usually written as
    \begin{equation}
	   \label{eq:ReactionNotation}
	   \sum_{i=1}^d (v^{-})_{i}\mathrm{S}_i \longrightarrow  \sum_{i=1}^d (v^{+})_{i}\mathrm{S}_i.
    \end{equation}
    To avoid the use of unnecessary species, we will assume that for each $1 \leq i \leq d$, there exists a vector $w=(w_1, \ldots,w_d)^T \in \Z_+^d$ with $w_i >0$ such that $(w,v)$ or $(v,w)$ is in $\Rs$ for some $v \in \Z^d_+$, i.e., each species is either a reactant or a product in some reaction.

    The net change in the quantity of molecules of each species due to a reaction $(v^{-},v^{+}) \in \Rs$ is described by $v^{+}-v^{-}$ and it is called the associated \textbf{reaction vector}. We denote the set of reaction vectors by $\V := \{ v \in \Z^d \:|\: v = v^{+}- v^{-} \text{ for some } (v^{-},v^{+}) \in \Rs\}$, we let $n := |\V|$ the size of $\V$ and we enumerate the members of $\V$ as $\{v_1,\ldots,v_n\}$. Note that $\V$ does not contain the zero vector because $\Rs$ has no elements of the form $(w,w)$. Different reactions might have the same reaction vector. For each $v_j \in \V$ we consider the set $\Rs_{v_j} := \{(v^{-},v^{+}) \in \Rs \:|\: v_j =v^{+}-v^{-} \}$. The matrix $S \in \R^{d \times |\Rs|}$ whose columns are the elements $v^+-v^-$ for $(v^-,v^+) \in \Rs$ will be called the \textbf{stoichiometric matrix}. 

    Given $(\Ss,\Rs)$ we will consider an associated continuous time Markov chain $X=(X_1,$ $\ldots,X_d)^T$, with a state space $\X$ contained in $\Z^d_+$, which tracks the number of molecules of each species over time. Roughly speaking, the dynamics of $X$ will be given by the following: given a current state  $x=(x_1,\ldots,x_d)^T \in \X \subseteq \Z_+^{d}$, for each reaction $(v^{-},v^{+}) \in \Rs$, there is a clock which will ring at an exponentially distributed time (with rate $\Lambda_{(v^{-},v^{+})}(x)$). The clocks for distinct reactions are independent of one another. If the clock corresponding to $(v^{-},v^{+})\in \Rs$ rings first, the system moves from $x$ to $x+v^{+}- v^{-}$ at that time, and then the process repeats. We now define the continuous time Markov chain in more detail.

    Consider sets of species $\Ss$ and reactions $\Rs$, a non-empty set $\X \subseteq \Z^d_+$ and a collection of functions $\prp= \{\prp_{(v^{-},v^{+})}:\X \longrightarrow \R_+\}_{(v^{-},v^{+}) \in \Rs}$ such that for each $x \in \X$ and $(v^{-},v^{+}) \in \Rs$, if $x+v^{+}-v^{-} \notin \X$, then $\Lambda_{(v^{-},v^{+})}(x)=0$. Now, for $1 \leq j \leq n$, define
    \begin{equation}\label{defrate}
	   \rate_j(x) := \sum_{(v^{-},v^{+}) \in \Rs_{v_j}} \prp_{(v^{-},v^{+})}(x).
    \end{equation}
    Note that for each $x \in \X$ and $1 \leq j \leq n$, if $x +v_j \notin \X$, then $\rate_j(x) = 0$. The functions $\{\prp_{(v^{-},v^{+})}:\X \longrightarrow \R_+\}_{(v^{-},v^{+}) \in \Rs}$ are called \textbf{propensity} or \textbf{intensity} functions. A common form for the propensity functions is the following, which is associated with \textbf{mass action kinetics}:
    \begin{equation}
	   \label{eq:PropensityFunctions}
	   \prp_{(v^{-},v^{+})}(x) = \kappa_{(v^{-},v^{+})}\prod_{i=1}^{d}(x_i)_{(v^{-})_i},
        \end{equation}
    where $\{\kappa_{(v^{-},v^{+})}\}_{(v^{-},v^{+}) \in \Rs}$ are non-negative constants and for $m,\l \in \Z_+$, the quantity $(m)_\l$ is the falling factorial, i.e., $(m)_0 := 1$ and $(m)_\l := m(m-1)\ldots(m-\l+1)$.
    
    A \textbf{stochastic chemical reaction network (SCRN)} (associated with $(\Ss,\Rs,\X,\prp)$) is a continuous time Markov chain $X$ with state space $\X$ and infinitesimal generator $Q$ given for $x,y \in \X$ by
    \begin{equation}
	   \label{eq:TransitionMatrixQ}
	   Q_{x,y} = \begin{cases}
		\rate_j(x) & \text{ if } y-x = v_j \text{ for some } 1 \leq j \leq n, \\
		- \sum_{j=1}^n\rate_j(x) & \text{ if } y = x, \\
		0 & \text{ otherwise.}
	\end{cases}
    \end{equation}
    %
%
  A SCRN associated with $(\Ss,\Rs,\X,\prp)$ is said to satisfy a \textbf{conservation law} if there is a $d$-dimensional non-zero vector $m$ such that $m^TS=0$, and hence $m^TX(t)=x_{\mathrm{tot}}$ for every $t \geq 0$, for some constant $x_{\mathrm{tot}}$.
  Consequently, we can reduce the dimension of the continuous time Markov chain describing the system by one. For example, if $m=(1,\ldots,1)^T$, then the projected process $(X_1,\ldots,X_{d-1})^T$ is again a continuous time Markov chain with state space $\{(x_1,\ldots,x_{d-1})^T \in \Z^{d-1}_+ \;|\; (x_1,\ldots,x_{d-1},x_{\mathrm{tot}} - \sum_{i=1}^{d-1} x_i)^T \in \X \}$. In our examples, we will often use this type of reduction.


    
    


%% file: StationaryDistributions.tex
\section{Main Results}
\label{sec:MainResults}

    In this section we describe the main theoretical results of this paper, under assumptions that go beyond those of our models of chromatin modification circuits. More precisely, we present results on stationary distributions and mean first passage times in Sections \ref{sec:StationaryDistributions} and \ref{sec:MFPT}. Then, in Section \ref{sec:Mon} we exploit theoretical results developed in our companion work \cite{Monotonicitypaper} to study monotonic dependence on parameters for a class of continuous time Markov chains related to chromatin modification circuits and other SCRNs.

\subsection{Stationary distributions}
\label{sec:StationaryDistributions}

    This section focuses on characterizing the terms in the series expansion \eqref{eq:ExpansionPiEps}. In Section \ref{sec:ZerothTerm} we focus on determining the term $\pi^{(0)}=\pi(0)$, while in SI - Section \ref{sec:StatDistHigherOrder} we provide a result which enables computation of all of the higher order terms $\pi^{(k)}$, for $k>0$, under additional assumptions. In Section \ref{sec:IllustrativeStationaryDistribution} we apply these results to the examples introduced in Section \ref{motivexamples}. Additional characterizations of $\pi(0)$ and $\pi^{(1)}$ are given in the SI - Sections \ref{sec:AdditionalCharactStatDist} and \ref{sec:PartialBalance}. Further examples for higher dimensional models of the chromatin modification circuits will be given in Section \ref{sec:Applications}. We remind the reader that to ease notation, we have adopted the convention that stationary distribution vectors will be row vectors, even though we do not use the transpose notation $T$ to indicate this.
    
\subsubsection{The zeroth order term}
\label{sec:ZerothTerm}
    
    As in Section \ref{sec:pertCTMC}, consider a family $\{X^\eps:\: 0 \leq \eps < \eps_0\}$ of continuous time Markov chains on a  finite state space $\mathcal{X}$, with infinitesimal generators $\{Q(\eps):\: 0 \leq \eps < \eps_0\}$ where $\eps \mapsto Q(\eps)$ is a real-analytic perturbation of $Q(0)$ with coefficients $\{Q^{(k)}:\: k \geq 0\}$ and additionally $Q(\eps)$ is irreducible for every $0 < \eps < \eps_0$. The matrix $Q(0) = Q^{(0)}$ is a $Q$-matrix for which $\X$ decomposes into recurrent (or ergodic) states $\A$ and transient states $\T$. From now on, we assume the following.
    
    \begin{assumption}
    \label{assumption:transient_absorbing}
     The set $\A$ consists of $|\A| \geq 1$ absorbing states for $Q(0)$, while $\T$ consists of $|\T| \geq 1$ transient states for $Q(0)$.
    \end{assumption}
    
    In other words, in the dynamics of $Q(0)$ there is at least one transient state, at least one recurrent state and all the recurrent states are absorbing. Now, we label the state space starting with the states in $\A$ and followed by the ones in $\T$. For every $k \geq 0$, we can write $Q^{(k)}$ as
     \begin{equation}
     \label{eq:Qdecomposition}
     Q^{(k)} =\left( \begin{array}{c | c }
                A_{k} & S_{k} \\
                \hline
                R_{k} & T_{k}
    \end{array}    \right),
     \end{equation}
    where $A_{k} \in \R^{|\A| \times |\A|}$, $S_{k} \in \R^{|\A| \times |\T|}$, $R_{k} \in \R^{|\T| \times |\A|}$ and $T_{k} \in \R^{|\T| \times |\T|}$. In a similar fashion, we can write
        \begin{equation}
        \label{eq:Q_eps_decomposition}
          Q(\eps) = \left( \begin{array}{c | c }
                A(\eps) & S(\eps) \\
                \hline
                R(\eps) & T(\eps)
    \end{array}    \right),
        \end{equation}
    for $0 \leq \eps < \eps_0$, where $A(\eps) \in \R^{|\A| \times |\A|}$, $S(\eps) \in \R^{|\A| \times |\T|}$, $R(\eps) \in \R^{|\T| \times |\A|}$ and $T(\eps) \in \R^{|\T| \times |\T|}$. From Assumption \ref{assumption:transient_absorbing}, we obtain that
    \begin{equation}
    \label{matrixQ00}
     Q^{(0)}= Q(0) =\left( \begin{array}{c | c }
                0 & 0 \\
                \hline
                R_0 & T_0
    \end{array}    \right),
     \end{equation}
     where $T_0$ is an invertible matrix (see SI - Lemma \ref{lem:InvertibilityQ^B}).
     
    
    For each $0 < \eps < \eps_0$, we denote by $\pi(\eps) = (\pi_x (\eps))_{x \in \mathcal{X}}$ the stationary distribution for $Q(\eps)$. In Section \ref{sec:pertCTMC}, we showed that the limit $\pi(0) := \lim_{\eps \to 0} \pi(\eps)$ exists and that $\eps \mapsto \pi(\eps)$ is a real-analytic perturbation of $\pi(0)$ with expansion given by \eqref{eq:ExpansionPiEps} for $0 \leq \eps < \eps_1$. For convenience, decompose the row vector $\pi(\eps)$ as $\pi(\eps) =[\alpha(\eps),\beta(\eps)]$ for $0 \leq \eps < \eps_1$ where $\alpha(\eps) \in \R^{|\A|}$ and  $\beta(\eps) \in \R^{|\T|}$. From \eqref{eq:ExpansionPiEps}, letting $\pi^{(k)}=[\alpha^{(k)},\beta^{(k)}]$, we have
        \begin{equation*}
        \alpha(\eps) = \sum_{k=0}^{\infty} \eps^k\alpha^{(k)} \quad \text{ and } \quad \beta(\eps)= \sum_{k=0}^{\infty} \eps^k\beta^{(k)}   
        \end{equation*}
    for $0 \leq \eps < \eps_1$. Since $\pi(\eps)$ is a probability distribution for every $0 \leq \eps < \eps_1$, we have that $\sum_{k=0}^{\infty} \eps^k(\pi^{(k)}\one) = 1$, which yields that  $\pi(0)\one=1$  and $\pi^{(k)}\one = 0$ for every $k \geq 1$. Since $\pi(0)Q(0) = 0$, 
    $\pi(0)$ is a stationary distribution for $X^{0}$ and so, by Assumption \ref{assumption:transient_absorbing}, it must be supported on $\A$ and so $\beta^{(0)}=0$. In the next result we establish an equation that is satisfied by $\alpha^{(0)}=\alpha(0)$ and introduce a key matrix for our analysis. For convenience, let $\alpha:= \alpha(0)$. 
    
    \begin{lemma}
    \label{lem:alphaL_A0}
    Under Assumption \ref{assumption:transient_absorbing}, $\pi(0)=[\alpha,0]$, where $0$ is the zero row vector of size $|\T|$ and $\alpha$ is an $|\A|$-dimensional probability vector satisfying the equation:
        \begin{equation}
        \label{eqsatisfied}
             \alpha(A_1 +  S_1(-T_0)^{-1}R_0) = 0.
        \end{equation}
    In addition,
        \begin{equation}
        \label{eq:Beta1InTermsOfAlpha}
        \beta^{(1)} = \alpha S_1(-T_0)^{-1}.
        \end{equation}
    \end{lemma}
%
%
\noindent See SI - Section \ref{proofs} for the proof of Lemma \ref{lem:alphaL_A0}. For convenience, we adopt the notation:
        \begin{equation}
        \label{eq:QA}
        Q_{\A} := A_1 +  S_1(-T_0)^{-1}R_0.   
        \end{equation}
    In SI - Lemma \ref{lem:Q_A_Qtilde_are_Qmatrix}, we show that $Q_{\A}$ is a $Q$-matrix of size $|\A| \times |\A|$. As a consequence, there exists a continuous time Markov chain with state space $\A$ and infinitesimal generator $Q_{\A}$. In general, a probability vector satisfying \eqref{eqsatisfied} needs not be unique. The following condition will imply uniqueness.
       
    \begin{assumption}
    \label{assumption:SingleRecurrentClass}
    The Markov chain associated with $Q_{\A}$ has a single recurrent class.
    \end{assumption}
    
    By SI - Lemma \ref{lem:nullityQmatrix}, Assumption \ref{assumption:SingleRecurrentClass} is equivalent to the condition $\dim(\ker(Q^T_{\A}))=1$. The next result then follows from Lemma \ref{lem:alphaL_A0}.
    
    \begin{theorem}
    \label{thm:alpha_unique}
    Suppose Assumptions \ref{assumption:transient_absorbing} and \ref{assumption:SingleRecurrentClass} hold. Then, $\pi(0)=[\alpha,0]$, where $\alpha$ is the unique probability vector on $\A$ such that $\alpha Q_{\A} = 0$.
    \end{theorem}
    
    As we will see, all of the chromatin modification circuit models presented in this work satisfy both Assumptions \ref{assumption:transient_absorbing} and \ref{assumption:SingleRecurrentClass}. Also note that Lemma \ref{lem:alphaL_A0} yields a characterization of $\beta^{(1)}$ by means of \eqref{eq:Beta1InTermsOfAlpha}.
    
    Theorem \ref{thm:alpha_unique} is simple to state, yet less easy to use since simple formulas for $Q_{\A}$ can be seldom obtained, making Assumption \ref{assumption:SingleRecurrentClass} hard to verify directly using \eqref{eq:QA}. In this regard, we now introduce an auxiliary continuous time Markov chain $\tilde{X}$ and use it to construct (via time-change) a realization $\hat{X}_{\A}$ of the continuous time Markov chain with infinitesimal generator $Q_{\A}$. This will enable us to give assumptions on $\tilde{X}$ that will imply Assumption \ref{assumption:SingleRecurrentClass} and which can sometimes be easier to verify. Also, this explicit realization for $\hat{X}_{\A}$ can lead to alternative ways to verify Assumption \ref{assumption:SingleRecurrentClass}. Under Assumption \ref{assumption:transient_absorbing}, consider the matrix     
    \begin{equation}
    \label{eq:DefTildeQ}
     \tilde{Q} :=\left( \begin{array}{c | c }
                A_1 & S_1 \\
                \hline
                R_0  & T_0 
    \end{array}    \right).
     \end{equation}
     
\begin{figure}[h]
        \centering
        \includegraphics[scale=0.9]{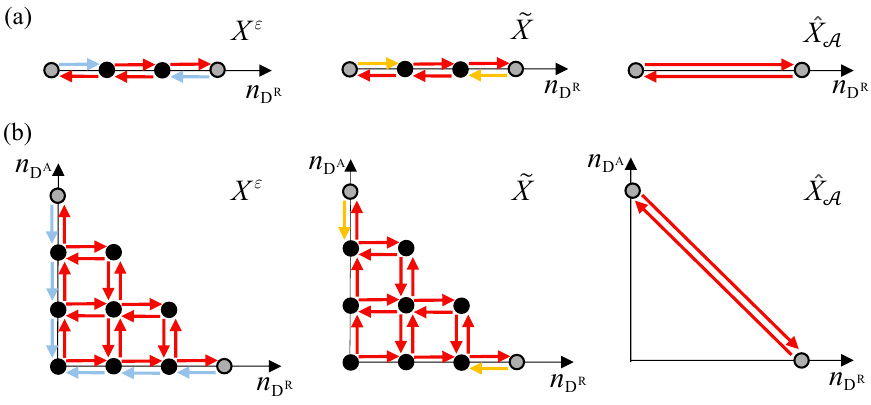}
            \caption{\small { \bf Graphs for the one-step transitions of $X^{\eps}$, $\tilde{X}$ and $\hat X_{\A}$ for the (a) 1D model and (b) 2D model.} Here, we consider $\mathrm{D_{tot}}=3$ and we use gray dots to represent the states belonging to $\A$ and black dots to represent all the other states, red arrows to represent transitions that are $O(1)$ for $X^{\eps}$, $\tilde{X}$ and $\hat X_{\A}$, blue arrows to represent transitions that are $O(\eps)$ for $X^{\eps}$, and golden arrows to represent the transitions for $\tilde{X}$ that were $O(\eps)$ for $X^{\eps}$ and became $O(1)$ for $\tilde{X}$.
            }%
            \label{fig:1D2DtildeX}
        \end{figure}

    In SI - Lemma \ref{lem:Q_A_Qtilde_are_Qmatrix} we prove that $\tilde{Q}$ is a $Q$-matrix. Let $\tilde{X}$ be a continuous time Markov chain with infinitesimal generator $\tilde{Q}$. For the purpose of illustration, if we assume that the perturbation is linear (as in \eqref{eq:LinearPerturbationQMatrices_1D}) and $\eps_0 > 1$, then the transitions of $\tilde{X}$ consist of the transitions of $X^0$ augmented by the transitions of $X^1$ that emanate from $\A$. See Figure \ref{fig:1D2DtildeX}(a)-(b) for an illustration related to the 1D and 2D models, respectively, introduced in Section \ref{motivexamples}.
    
    
    Now, consider the \textbf{occupation time} of $\A$ by the Markov chain $\tilde{X}$ up to time $t \geq 0$, given by $\chi_{\A}(t) := \int_{0}^{t} \one_{\A}(\tilde{X}(s))ds$ for $t \geq 0$. Denote by $\chi_{\A}(\infty) = \lim_{t \to \infty} \chi_{\A}(t) = \int_0^{\infty} \one_{\A}(\tilde{X}(s))ds$. Since $T_0$ is invertible, SI - Lemmas \ref{lem:InvertibilityQ^B} and \ref{lem:OccupationTimesAreInfinity} yield that $\PP_x[\chi_{\A}(\infty)=\infty]=1$ for all $x \in \mathcal{X}$. Additionally, consider the right-continuous inverse of $\chi_{\A}$, $\tau(s) := \inf\{ t \geq 0:\: \chi_\A(t) > s \}$, defined for $s \geq 0$. We define the \textbf{restriction} process $\hat{X}_{\A}$ as
        \begin{equation}
        \label{eq:RestrictiontoAbsorbing}
         \hat{X}_{\A}(s) := \tilde{X}(\tau(s)),\;\;\;s \ge 0.
        \end{equation}
   By properties of the right-continuous inverse (see Problem 4.5 in \cite{KaratzasShreve}, for example), the reader may verify that $\hat{X}_{\A}$ corresponds to observing $\tilde{X}$ only on the time intervals where $\tilde{X}$ is in $\A$. Roughly speaking, we are \textit{erasing} the times where $\tilde{X}$ is outside of $\A$. In the language of Blumenthal \& Getoor \cite{BlumenthalGetoor}, $\chi_{\A}$ is a continuous additive functional for $\tilde{X}$, and by Exercise V.2.11 in \cite{BlumenthalGetoor}, we obtain that $\hat{X}_{\A}$ is a continuous time Markov chain  with state space $\A$. In the next result, we prove that $\hat{X}_{\A}$ is a realization of the continuous time Markov chain associated with $Q_{\A}$. See Figure \ref{fig:1D2DtildeX}(a)-(b) for a representation of $\hat{X}_{\A}$ associated with the 1D and 2D models, respectively.

    \begin{lemma}
    \label{lem:welldefined_QA}
    Suppose Assumption \ref{assumption:transient_absorbing} holds. Then, $\hat{X}_{\A}$ has infinitesimal generator $Q_{\A}$.
    \end{lemma}

\noindent The proof of Lemma \ref{lem:welldefined_QA} is given in SI - Section \ref{proofs}.
    We now introduce some assumptions that imply that Assumption \ref{assumption:SingleRecurrentClass} holds. In addition, these assumptions will allow for some refinements (see SI - Section \ref{sec:AppendixHigherOrder}).
    
    \begin{assumption}
    \label{assumption:A_is_recurrent}
    For $\tilde{X}$, there exists a communicating class $\CC$ such that $\A \subseteq \CC$.
    \end{assumption}
    
    We note that, if such a class $\CC$ exists, then it has to be recurrent. In fact, if it was transient then $\chi_{\A}(\infty) < \infty$ with positive probability under $\PP_x$, $x \in \A$, which is a contradiction.

    \begin{assumption}
    \label{assumption:tildeX_irreducible}
    The Markov chain $\tilde{X}$ is irreducible.
    \end{assumption}
    
    We note that Assumption \ref{assumption:tildeX_irreducible} implies Assumption \ref{assumption:A_is_recurrent}. Moreover, they are both related to Assumption \ref{assumption:SingleRecurrentClass} in the following way.
    
    \begin{lemma}\label{lemma3}
    Suppose Assumptions \ref{assumption:transient_absorbing} and \ref{assumption:A_is_recurrent} hold. Then, the process $\hat{X}_{\A}$ is irreducible. As a consequence, either of Assumptions \ref{assumption:tildeX_irreducible} or \ref{assumption:A_is_recurrent} implies that Assumption \ref{assumption:SingleRecurrentClass} holds.
    \end{lemma}
    

\noindent The proof of Lemma \ref{lemma3} is given in SI - Section \ref{proofs}. The next result follows from Lemmas \ref{lem:welldefined_QA}, \ref{lemma3} and Theorem \ref{thm:alpha_unique}.
    
    \begin{theorem}
    \label{thm:alpha_rho}
    Suppose Assumptions \ref{assumption:transient_absorbing} and \ref{assumption:A_is_recurrent} hold. Then, $\pi(0)=[\alpha,0]$ where $\alpha$ is the unique stationary distribution for the process $\hat{X}_{\A}$ and all entries of $\alpha$ are strictly positive.
    \end{theorem}
    
    Assumptions \ref{assumption:A_is_recurrent} and \ref{assumption:tildeX_irreducible} can be understood graphically in some cases. For example, Figure \ref{fig:1D2DtildeX} illustrates that for the 1D-model, Assumption \ref{assumption:tildeX_irreducible} is satisfied. For the 2D-model, we can see that while Assumption \ref{assumption:tildeX_irreducible} is not satisfied (since the state $(0,0)$ forms its own (transient) class for $\tilde{X}$), Assumption \ref{assumption:A_is_recurrent} does indeed hold. 
    In Section \ref{sec:Applications} we will see that neither Assumption \ref{assumption:tildeX_irreducible} nor \ref{assumption:A_is_recurrent} is satisfied by the 3D or 4D model. However, the weaker Assumption \ref{assumption:SingleRecurrentClass} does hold.



    
      

In the SI, we give recursive formulae for the higher order terms $\pi^{(k)}$,  $k=1,2,\dots$, under the following additional assumption (see SI - Theorem \ref{thm:HigherOrderTerms_for_stationary distribution}).

\begin{assumption}
    \label{assumption:linear_case}
    The perturbation is linear, i.e., $Q(\eps) = Q^{(0)} + \eps Q^{(1)}$ for $0 \leq \eps < \eps_0$.
    \end{assumption}


\subsubsection{Illustrative examples: 1D and 2D model}
\label{sec:IllustrativeStationaryDistribution}
    
     \textbf{1D model.} We use the tools developed in the preceding section to derive the terms $\pi^{(0)}$ and $\pi^{(1)}$ in the expansion \eqref{eq:pi_expansion_1D} for the 1D model introduced in Section \ref{SS1}. Fix $\Dtot \ge 2$ and let $X^\eps$ with infinitesimal generator $Q(\eps)$ be as in Section \ref{SS1}, with the expression for $Q^{(k)}$ given in \eqref{eq:Qdecomposition}. By \eqref{eq:LinearPerturbationQMatrices_1D}, Assumption \ref{assumption:linear_case} holds.
%
  %
    Moreover, for each $0 < \eps < \eps_0$, with $\eps_0$ being a fixed, positive constant, $Q(\eps)$ is irreducible, while $Q(0)$ has a non-empty set of transient states $\T= \{1, \ldots, \Dtot-1\}$ and a set of two absorbing states $\A = \{a,r\}$, with $a=0$ representing the fully active state ($n_{\DA}=\Dtot$) and $r=\Dtot$ representing the fully repressed state ($n_{\DR}=\Dtot$). Then, Assumption \ref{assumption:transient_absorbing} holds (see SI - Section \ref{sec:Appendix1Dmodel}). Furthermore, by defining $f(x):=x(\Dtot-x)$ for $x \in \X$, we can write the matrices $R_0$ and $T_0$ in the matrix $Q^{(0)}$ as follows:
    
 %
\begin{equation*}
\scalebox{0.82}{$
\begin{aligned}
    R_0 &= \begin{pmatrix}
            \mu\frac{k^A_E}{V}f(1) & 0 \\
            0 & 0 \\
            \vdots & \vdots\\
            0 & 0 \\
            0 & \frac{k^A_E}{V}f(1) \\
        \end{pmatrix},\;\;
    T_0 &= \begin{pmatrix}
            -(\lambda^{0}_1+\gamma^{0}_1) & \lambda^{0}_1 & 0 & \ldots & \ldots & \ldots & 0\\
            0 & \ddots & \ddots & \vdots & \vdots & \vdots &  \vdots\\
            \vdots & \ddots & \ddots & \ddots & 0 &  \vdots &  \vdots\\
            \ldots & 0 & \gamma^{0}_x & -(\lambda^{0}_x+\gamma^{0}_x) & \lambda^{0}_x & 0 & \ldots\\
            \vdots & 0 & \vdots & \ddots & \ddots & \ddots & \vdots\\
            \vdots & \vdots & \vdots & \vdots & \ddots & \ddots & 0\\
            0 & \ldots & \ldots & \ldots & 0 & \gamma^{0}_{\Dtot-1} & -(\lambda^{0}_{\Dtot-1}+\gamma^{0}_{\Dtot-1})\\
        \end{pmatrix}
\end{aligned}$}
\end{equation*}
 
    \noindent where $R_0$ is a $(\Dtot -1)\times 2$ matrix and $T_0$ is a $(\Dtot -1)\times (\Dtot -1)$ tridiagonal matrix, and
    $\gamma^{0}_x=\mu\frac{k^A_E}{V}f(x)$, $\lambda^{0}_x=\frac{k^A_E}{V}f(x)$, and $f(\Dtot-1) = f(1)=(\Dtot-1)$. In addition, we can write $A_1$ and $S_1$ of $Q^{(1)}$ as follows:

    \begingroup
    \small 
    \begin{equation}\nonumber
    A_1 =  \left( \begin{array}{c c }
                -\frac{k^A_E}{V} \mathrm{D^2_{tot}} & 0 \\
                0 & -b\mu\frac{k^A_E}{V} \mathrm{D^2_{tot}} \\
    \end{array}    \right),\;\;\;    
    S_1 =  \left( \begin{array}{c c c c c}
                \frac{k^A_E}{V} \mathrm{D^2_{tot}} & 0 & \ldots & \ldots & 0 \\
                0 & \ldots & \ldots & 0 & b\mu\frac{k^A_E}{V} \mathrm{D^2_{tot}}\\
    \end{array}    \right).\;\;\;  
    \end{equation}
    \endgroup
    \normalsize
    The process $\tilde{X}$, whose infinitesimal generator is defined in \eqref{eq:DefTildeQ}, is irreducible (see SI - Section \ref{sec:Appendix1Dmodel}). This is illustrated in Figure \ref{fig:1D2DtildeX}(a). Thus, Assumption \ref{assumption:tildeX_irreducible} holds. Then, Assumption \ref{assumption:A_is_recurrent} is also satisfied and Theorem \ref{thm:alpha_rho} can be applied. This yields that $\pi(0)=\pi^{(0)}=[\alpha,0] = \left[\alpha_a,\alpha_r, 0 \ldots ,0 \right]$ where $\alpha$ is the unique stationary distribution for the restriction process $\hat{X}_{\A}$ (defined by \eqref{eq:RestrictiontoAbsorbing}), whose infinitesimal generator is $Q_{\A} = A_1 +  S_1(-T_0)^{-1}R_0$ by Lemma \ref{lem:welldefined_QA} and \eqref{eq:QA}. 
    %
%
Now,
        \begingroup
    \small 
    	\begin{equation}\label{QA1D}
    	    Q_\A = 
	\left( \begin{array}{c c}
	    -\frac{1-\mu}{1-\mu^{\Dtot}} \frac{k_E^A}{V} \Dtot^2 & \frac{1-\mu}{1-\mu^{\Dtot}} \frac{k_E^A}{V} \Dtot^2 \\
	    b \mu^{\Dtot} \frac{1-\mu}{1-\mu^{\Dtot}} \frac{k_E^A}{V} \Dtot^2 & - b \mu^{\Dtot} \frac{1-\mu}{1-\mu^{\Dtot}} \frac{k_E^A}{V} \Dtot^2 \\
    \end{array} \right) = \frac{1-\mu}{1-\mu^{\Dtot}} \frac{k_E^A}{V} \Dtot^2 \left( \begin{array}{c c}
	    -1 & 1 \\
	    b \mu^{\Dtot} & - b \mu^{\Dtot} \\
    \end{array} \right),
    \end{equation}
    \endgroup
    \normalsize
    \noindent and since $\alpha$  is the unique probability vector satisfying $\alpha Q_{\A} = 0$, we have
\begingroup
\small
    \begin{equation}\nonumber
    \begin{aligned}
        \alpha_a=\frac{b\mu^{\mathrm{D_{tot}}}}{1 + b\mu^{\mathrm{D_{tot}}}},\;\;\alpha_r=\frac{1}{1 + b\mu^{\mathrm{D_{tot}}}}.
        \end{aligned}
    \end{equation}
    \endgroup
 \normalsize
These results are in agreement with \eqref{pi01Dmodel} in Section \ref{SS1}, where we explicitly computed the stationary distribution $\pi(\eps)$ and let $\eps \rightarrow 0$ (see SI - Section \ref{sec:Appendix1Dmodel}).

Now, since Assumptions \ref{assumption:transient_absorbing}, \ref{assumption:A_is_recurrent}, and \ref{assumption:linear_case} hold, we can apply SI - Theorem \ref{thm:HigherOrderTerms_for_stationary distribution} to derive an expression for $\beta^{(1)}$. For the transient states $\T=\{1,...,\Dtot-1\}$, we have $\beta^{(1)}=[\pi^{(1)}_{1},...,\pi^{(1)}_{\Dtot-1}] = \alpha S_1(-T_0)^{-1}$, and so for $x \in \T$
         \begin{equation}        \beta^{(1)}_x=\frac{b\mu^{\mathrm{D_{tot}}}}{1 + b\mu^{\mathrm{D_{tot}}}}\frac{k^A_E}{V}\mathrm{D^2_{tot}}(-T_0)^{-1}_{1,x}+\frac{1}{1 + b\mu^{\mathrm{D_{tot}}}}b \mu \frac{k^A_E}{V}\mathrm{D^2_{tot}}(-T_0)^{-1}_{\Dtot-1,x},
        \end{equation}   
in which $(-T_0)^{-1}_{1,x}$ and $(-T_0)^{-1}_{\Dtot-1,x}$, for $x \in \T$, are the elements indexed by $(1,x)$ and $(\Dtot-1,x)$ of the matrix $(-T_0)^{-1}$, respectively. After some calculations, we obtain 

\begingroup
\small
         \begin{equation}\nonumber
         \begin{aligned}
         (-T_0)^{-1}_{1,x}&=\frac{\left(\frac{k^A_E}{V}\right)^{\Dtot-2}\frac{\prod_{i=1}^{\Dtot-1}f(i)}{f(x)}\left(1+\prod_{i=1}^{\Dtot-1-x}\mu^i\right)}{\left(\frac{k^A_E}{V}\right)^{\Dtot-1}\prod_{i=1}^{\Dtot-1}f(i)\left(1+\prod_{i=1}^{\Dtot-1}\mu^i\right)}= \frac{\left(1+\prod_{i=1}^{\Dtot-1-x}\mu^i\right)}{\frac{k^A_E}{V}f(x)\left(1+\prod_{i=1}^{\Dtot-1}\mu^i\right)},\\         
         (-T_0)^{-1}_{\Dtot-1,x}&=\frac{\left(\frac{k^A_E}{V}\right)^{\Dtot-2}\frac{\prod_{i=1}^{\Dtot-1}f(i)}{f(x)}\left(1+\prod_{i=1}^{x-1}\mu^i\right)\mu^{\Dtot-1-x}}{\left(\frac{k^A_E}{V}\right)^{\Dtot-1}\prod_{i=1}^{\Dtot-1}f(i)\left(1+\prod_{i=1}^{\Dtot-1}\mu^i\right)}=\frac{\left(1+\prod_{i=1}^{x-1}\mu^i\right)\mu^{\Dtot-1-x}}{\frac{k^A_E}{V}f(x)\left(1+\prod_{i=1}^{\Dtot-1}\mu^i\right)} ,
         \end{aligned}
        \end{equation}  
\endgroup
 \normalsize
and then $\beta^{(1)}_x$, $x\in \T$, can be written as follows:
        \begin{equation}\label{eq:beta1}
        \beta^{(1)}_x=\frac{\mathrm{D^2_{tot}}}{f(x)}\frac{b\mu^{\Dtot-x}}{1+b\mu^{\Dtot}} =\frac{\mathrm{D^2_{tot}}}{x(\Dtot-x)}\frac{b\mu^{\Dtot-x}}{1+b\mu^{\Dtot}}.
        \end{equation}

    \textbf{2D model.} In this section we analyze the stationary distribution for the 2D model introduced in Section \ref{exp2}.
    %
    %
    Fix $\Dtot \ge 2$ and let $X^\eps$ with infinitesimal generator $Q(\eps)$ be as in Section \ref{SS1}, with the expression for the $Q^{(k)}$ given by \eqref{eq:Qdecomposition}. By \eqref{rates2D}, for this model Assumption \ref{assumption:linear_case} holds. Furthermore, $Q(0)$ has a non-empty set of transient states $\T= \{i_1, \ldots, i_m\}$ where $m = \frac{(\Dtot+2)(\Dtot+1)}{2}-2$, $i_1=(0,\Dtot-1)^T$, $i_m=(\Dtot-1,0)^T$, and absorbing states $\A = \{a,r\}$, with $a=(0,\Dtot)^T$ corresponding to the fully active state ($n_{\DA}=\Dtot$) and with $r=(\Dtot,0)^T$ corresponding to the fully repressed state ($n_{\DR}=\Dtot$), respectively. Then, Assumption \ref{assumption:transient_absorbing} holds (see SI - Section \ref{sec:Appendix2Dmodel}).
From \eqref{rates2D}, we see that $A_0 =0$, $S_0=0$ and 
\begin{equation}\nonumber
     A_1 = \left( \begin{array}{cc}
        - \frac{k^A_M}{V}\mathrm{D^2_{tot}} & 0  \\
            0 & -\frac{k^A_M}{V}\mathrm{D^2_{tot}}\mu b
    \end{array}    \right),\;
         S_1 = \left( \begin{array}{ccccc}
                \frac{k^A_M}{V}\mathrm{D^2_{tot}} & 0 & \ldots & \ldots & 0 \\
                0 & \ldots & \ldots & 0 & \frac{k^A_M}{V}\mathrm{D^2_{tot}}\mu b
    \end{array}    \right).
    \end{equation}
Furthermore, $R_0\in \R^{m \times 2}$ is given by
\begin{equation*}
\scalebox{0.89}{$
         \left( \begin{array}{cc}
                f_A(0,\Dtot-1) & 0 \\
                0 & 0                \\
                \vdots & \vdots       \\
                0 & 0                \\
                0 & f_R(\Dtot-1,0)
    \end{array}    \right) = 
    \left( \begin{array}{cc}
                k_{W0}^A+k_{W}^A + \frac{k_{M}^A}{V}(\Dtot-1) & 0 \\
                0 & 0                \\
                \vdots & \vdots        \\
                0 & 0                \\
                0 &  k_{W0}^R +k_{W}^R + \frac{k_{M}^R}{V}(\Dtot-1) 
    \end{array}    \right),$}
\end{equation*}  
    \noindent and $R_1=0$. The matrices $T_0$ and $T_1$ are more complex and examples of them, for $\Dtot=2$, are provided in SI - Section \ref{sec:Appendix2Dmodel}. For $\tilde{X}$, $\CC = \X \setminus \{(0,0)\}$ is a communicating class such that $\A \subseteq \CC$. This implies that Assumption \ref{assumption:A_is_recurrent} is satisfied. Given that Assumptions \ref{assumption:transient_absorbing} and \ref{assumption:A_is_recurrent} are satisfied, Theorem \ref{thm:alpha_rho} can be applied and we obtain that $\pi(0)=\pi^{(0)}=[\alpha,0] = \left[\alpha_a,\alpha_r, 0 \ldots ,0 \right]$ where $\alpha$ is the unique stationary distribution for the process $\hat{X}_{\A}$, whose infinitesimal generator is $Q_{\A} = A_1 +  S_1(-T_0)^{-1}R_0$. This means that $\alpha$ is the unique probability vector such that $\alpha(A_1 +  S_1(-T_0)^{-1}R_0) = 0$. Furthermore, given that Assumption \ref{assumption:linear_case} is satisfied, we can apply SI - Theorem \ref{thm:HigherOrderTerms_for_stationary distribution} to derive an expression for $\beta^{(1)}=[\pi^{(1)}_{i_1},...,\pi^{(1)}_{i_m}] = \alpha S_1(-T_0)^{-1}$.
    For example, if $\Dtot=2$, the matrix $Q_\A$ is given by
    
    \begingroup
    \small 
    	\begin{equation}\label{QA2D}
    	    Q_\A = \frac{4}{K} \frac{k_M^A}{V} \left( \begin{array}{c c}
	    -(k_{W0}^R+k_{W}^R)(k_{W0}^R+k_{W}^R + \frac{k_{M}^R}{V}) & (k_{W0}^R+k_{W}^R)(k_{W0}^R+k_{W}^R + \frac{k_{M}^R}{V}) \\
	    b\mu^2(k_{W0}^A+k_{W}^A)(k_{W0}^A+k_{W}^A+\frac{k_{M}^A}{V}) & - b\mu^2(k_{W0}^A+k_{W}^A)(k_{W0}^A+k_{W}^A+\frac{k_{M}^A}{V}) \\
    \end{array} \right),
    \end{equation}
    \endgroup
    \normalsize    
    with 
    \begingroup
    \small 
    \begin{equation}\label{QA2DcoefK}
    K=(k_{W0}^A+k_{W}^A+\frac{k_{M}^A}{V}+k_{W0}^R+k_{W}^R)(k_{W0}^R+k_{W}^R + \frac{k_{M}^R}{V})+\mu(k_{W0}^R+k_{W}^R + \frac{k_{M}^R}{V}+k_{W0}^A+k_{W}^A)(k_{W0}^A+k_{W}^A+\frac{k_{M}^A}{V}),
    \end{equation}
    \endgroup
    \normalsize
    \noindent and then $\pi^{(0)}$ is given by
    \begin{equation}
        \label{pi2DmodelRESULTS}
         \pi^{(0)}_x = \begin{cases}
         \frac{b\mu^2(k_{W0}^A+k_{W}^A)(k_{W0}^A+k_{W}^A+\frac{k_{M}^A}{V})}{b\mu^2(k_{W0}^A+k_{W}^A)(k_{W0}^A+k_{W}^A+\frac{k_{M}^A}{V})+(k_{W0}^R+k_{W}^R)(k_{W0}^R+k_{W}^R + \frac{k_{M}^R}{V})} & \text{ if } x=(0,\Dtot)^T \\
         0 & \text{ if } x \in \T \\     \frac{(k_{W0}^R+k_{W}^R)(k_{W0}^R+k_{W}^R + \frac{k_{M}^R}{V})}{b\mu^2(k_{W0}^A+k_{W}^A)(k_{W0}^A+k_{W}^A+\frac{k_{M}^A}{V})+(k_{W0}^R+k_{W}^R)(k_{W0}^R+k_{W}^R + \frac{k_{M}^R}{V})} & \text{ if } x = (\Dtot,0)^T.
        \end{cases}
        \end{equation}
    See SI - Section \ref{sec:Appendix2Dmodel} for the evaluation of $\pi^{(1)}_{x}$ for the transient states $x \in \T$ when $\Dtot=2$. For this value of $\mathrm{D_{tot}}$, we see from \eqref{pi2DmodelRESULTS} that $\pi^{(0)}_x$ depends monotonically on $\mu$ for each fixed $x$. As $\mathrm{D_{tot}}$ increases, the algebraic complexity of a full parameter representation of $\pi^{(0)}_x$ increases very rapidly. Thus, to investigate monotonic dependence on parameters for biologically relevant values of $\mathrm{D_{tot}}$ (of the order of 50, considering an average gene length of 10,000 bp \cite{Cooper} and one nucleosome per 200 bp \cite{Groudine03}), we shall use comparison theorems developed in \cite{Monotonicitypaper}, without calculating any explicit formula (Section \ref{sec:Mon}).

%% file: MFPT.tex
\subsection{Mean first passage times (MFPTs)}
\label{sec:MFPT}
In this section we develop a theoretical framework to study mean first passage times for continuous time Markov chains. We first develop an algorithm to determine the order of the pole of MFPTs for singularly perturbed Markov chains (Section \ref{orderMFPT}). In Section \ref{sec:leadingMFPT}, we focus on determining the leading coefficient for MFPTs, under some assumptions introduced in Section \ref{sec:StationaryDistributions}.
%
%
In Section \ref{sec:MFPTexamples}, we apply these results to the examples introduced in Section \ref{motivexamples}.
 


\subsubsection{Algorithm to find the order of the poles for MFPTs}
\label{orderMFPT}

	
	
	Our algorithm is adapted from an algorithm developed by Hassin and Haviv \cite{HH} for discrete time Markov chains. The idea used in \cite{HH} was to consider transitions between subsets of states and to keep track of the sojourn times in the sets of states. This is used to define a coarser version of the process, which may not be a Markov process and which moves between groups of states of the original Markov chain. This idea can be adapted to the continuous time setting as well. For this, we introduce stopping times to more explicitly track the sojourn times than was done in \cite{HH}. In addition, we extend the original algorithm's scope to consider the mean first passage time to a subset of states, instead of just a single state. The paper \cite{HH} uses r-cycles and notes that these could be replaced by more general r-components. Here, we focus on using the latter and call the set of vertices in such an r-component an r-connected set.
	
	In this section, we consider a singularly perturbed, finite-state, continuous time Markov chain $X^\eps$ on $\X$ with infinitesimal generator $Q(\eps)$ as described in Section \ref{sec:pertCTMC}. We provide an algorithm for finding the orders $\{p(v):\: v \in \B^c \}$ of the poles for the mean first passage times to $\B \subset \mathcal{X}$ for $X^\eps$ starting from states in $\B^c$, where $\B \neq \emptyset$ is a strict subset of $\X$. We begin with a few definitions and some notation and then present the algorithm.
	
	\begin{definition}
		\label{def:BigTheta}
		Given $\eps_0 > 0$ and a function $f: (0,\eps_0) \rightarrow \R_{>0}$, we say $f = \Theta (\eps ^k)$ if there exist $k \in \Z$ and strictly positive $m, M \in \R_{> 0}$ such that, for all $0 < \eps < \eps_0$, 
	\begin{equation*}
        m \eps^k \leq f(\eps) \leq M \eps^k.
        \end{equation*}
		If $f = \Theta (\eps ^{k})$ for some $k \in \Z$, we say the \textbf{order} (at the origin) of $f$ is $k$. If $f = \Theta (\eps ^{-k})$ where $k \in \Z_+$, we say that the \textbf{order of the pole} of $f$ is $k$.
	\end{definition}

	Because the perturbation of $X^\eps$ is real analytic, $|\X|>1$ and $X^\eps$ is irreducible for $\eps>0$, there exists $\eps_{\max} > 0$ such that for each $x \neq y \in \X$, either $Q_{x,y}(\eps) \equiv 0$ for all $\eps \in (0,\eps_{\max})$ or $Q_{x,y}(\eps) > 0$ for all $\eps \in (0,\eps_{\max})$. In the latter case, the order of $Q_{x,y}(\eps)$ is a non-negative integer, which we denote by $k_{xy}$. We let $E_0 = \{(x,y):\: Q_{x,y}(\eps) > 0 \text{ for all } \eps \in (0,\eps_{\max}) \}$. As the algorithm progresses, states of $\X$ are gathered together to form composite nodes and the graph of the states of $X^\eps$ progresses through a series of reduced graphs. If $u$ is a node in one of the graphs, then $S(u) \subset \X$ consists of the states in $\X$ that are collapsed to form the (reduced) node $u$. In Steps 2 and 3 of the algorithm, the function $\mathcal{K}$ and the initial values of $p$ are inductively determined for all of these graphs. The final values of $p$ for nodes in $\B^c$ are then determined in Step 4. With $\mathcal{K}_{uv}$ being defined, a directed edge $(u,v)$ in one of the graphs is called an \textbf{r-edge}, where r is for regular, if $\mathcal{K}_{uv} = 0$, and an \textbf{r-path} is a directed path in the graph consisting of r-edges only. A set $C$ in one of the graphs is called an \textbf{r-connected set} if $|C|>1$ and there exists an r-path from $u$ to $v$ for any $u \neq v \in C$. The \textbf{order of the pole of the expected sojourn time} spent in an r-connected set $C$ depends only on the set $C$ and is denoted by $p(c)$ where $c$ is a node representing the set $C$. For any node $w$ outside of $C$, $\mathcal{K}_{cw}$ and $\mathcal{K}_{wc}$ are the \textbf{the order of the probabilities of a one-step transition} from $c$ to $w$ and from $w$ to $c$, respectively. In Step 4 of the algorithm, $p(\cdot)$ keeps being updated but will stay finite and eventually fixate. The algorithm statement and related proof can be found in the SI - Sections \ref{sec:Algstat} - \ref{sec:justification}.
	
	\subsubsection{Leading coefficient in MFPT series expansion}
	\label{sec:leadingMFPT}
	
        In Section \ref{sec:pertCTMC}, we have shown that for each $0 < \eps < \eps_0$, the unique stationary distribution $\pi(\eps)$ for $X^\eps$ admits a real-analytic expansion in powers of $\eps$. By \eqref{eq:ReturnTimes_expansion} and \eqref{eq:Pi_related_to_returntimes}, for $x \in \X$,
	\begin{equation}
	\frac{1}{\pi_x(\eps)} = q_x(\eps)\E_x[\zeta^{\eps}_x] = 1 + \sum_{y \neq x} Q_{x,y}(\eps) h_{y,x}(\eps).
	\end{equation}
	Recall that $E_0 = \{(x,y):\: Q_{x,y}(\eps) > 0 \text{ for all } \eps \in (0,\eps_{\max}) \}$ and $k_{xy}$ is the order of $Q_{x,y} (\eps)$ for each $(x,y) \in E_0$. Using the algorithm in Section \ref{orderMFPT}, we can obtain the order of the pole, $p_x(y)$, of the mean first passage time $h_{y,x}(\eps)$ from $y$ to $x$ for all $y \neq x \in \X$. Therefore, for each $x \in \X$, the order of $\pi_x(\eps)$ is 
	\begin{equation}
	    \label{eqn:orderSD}
	    k_x = \max\{p_x(y)-k_{xy}: (x,y) \in E_0; 0 \} \geq 0,
	\end{equation} 
	and then
	\begin{equation*}
        \pi_x(\eps) = \sum_{k = k_x}^{\infty} \eps^{k} \pi_x^{(k)}.
        \end{equation*}

        The following theorem is for continuous time and builds on discrete time results of 
        
        \noindent Avrachenkov et al. \cite{AvrachenkovHaviv,AvrachenkovFilarHowlett}.
        
	\begin{theorem}
	\label{thm:MFPTcoefficient}
		Suppose Assumptions \ref{assumption:transient_absorbing}, \ref{assumption:SingleRecurrentClass} and \ref{assumption:linear_case} hold. Let $Q_\A$ be given by \eqref{eq:QA}, $\hat{X}_\A$ be as defined in \eqref{eq:RestrictiontoAbsorbing}, and $\alpha$ be the unique stationary distribution for $\hat{X}_\A$ defined in Theorem \ref{thm:alpha_unique}. Let $D=(-Q_\A + \one \alpha)^{-1} - \one \alpha$. For $y \in \X$, let $k_y$ be the order of the stationary distribution $\pi_y(\eps)$ of $X^\eps$, defined by \eqref{eqn:orderSD}. Then, for $x,y \in \A$, the mean first passage time from $x$ to $y$ for $X^\eps$ is
		\begin{equation}
		\label{eq:leadingMFPT}
		    h_{x,y}(\eps)= \frac{D_{y,y} - D_{x,y}}{\pi_y^{(k_y)}} \frac{1}{\eps^{k_y+1}} + O\left(\frac{1}{\eps^{k_y}}\right).
		\end{equation}
		Moreover, if $\hat{X}_\A$ is an irreducible Markov chain, then the order of the pole of $h_{x,y}(\eps)$ is one, i.e., $k_y=0$, and the coefficient of $\eps^{-1}$ in \eqref{eq:leadingMFPT} is equal to the mean first passage time from $x$ to $y$ for the process $\hat{X}_\A$.
	\end{theorem}
 
	\noindent The proof of Theorem \ref{thm:MFPTcoefficient} is given in SI - Section \ref{proofthmleadcoeff}.
	\begin{remark}
	    It may be possible that $D_{y,y}-D_{x,y} = 0$. In this case, 
        \begin{equation*}   
         h_{x,y}(\eps)= 0 \cdot \frac{1}{\eps^{k_y+1}} + O \left(\frac{1}{\eps^{k_y}} \right) = O \left(\frac{1}{\eps^{k_y}} \right).
         \end{equation*}
	    However, if we find that the order of the pole of $h_{x,y}(\eps)$ is $k_y+1$, using the algorithm in Section \ref{orderMFPT}, then we can rule out the possibility of $D_{y,y}-D_{x,y}$ being $zero$.
	\end{remark}

	\begin{figure}[h!]
		\centering \includegraphics[width=\textwidth]{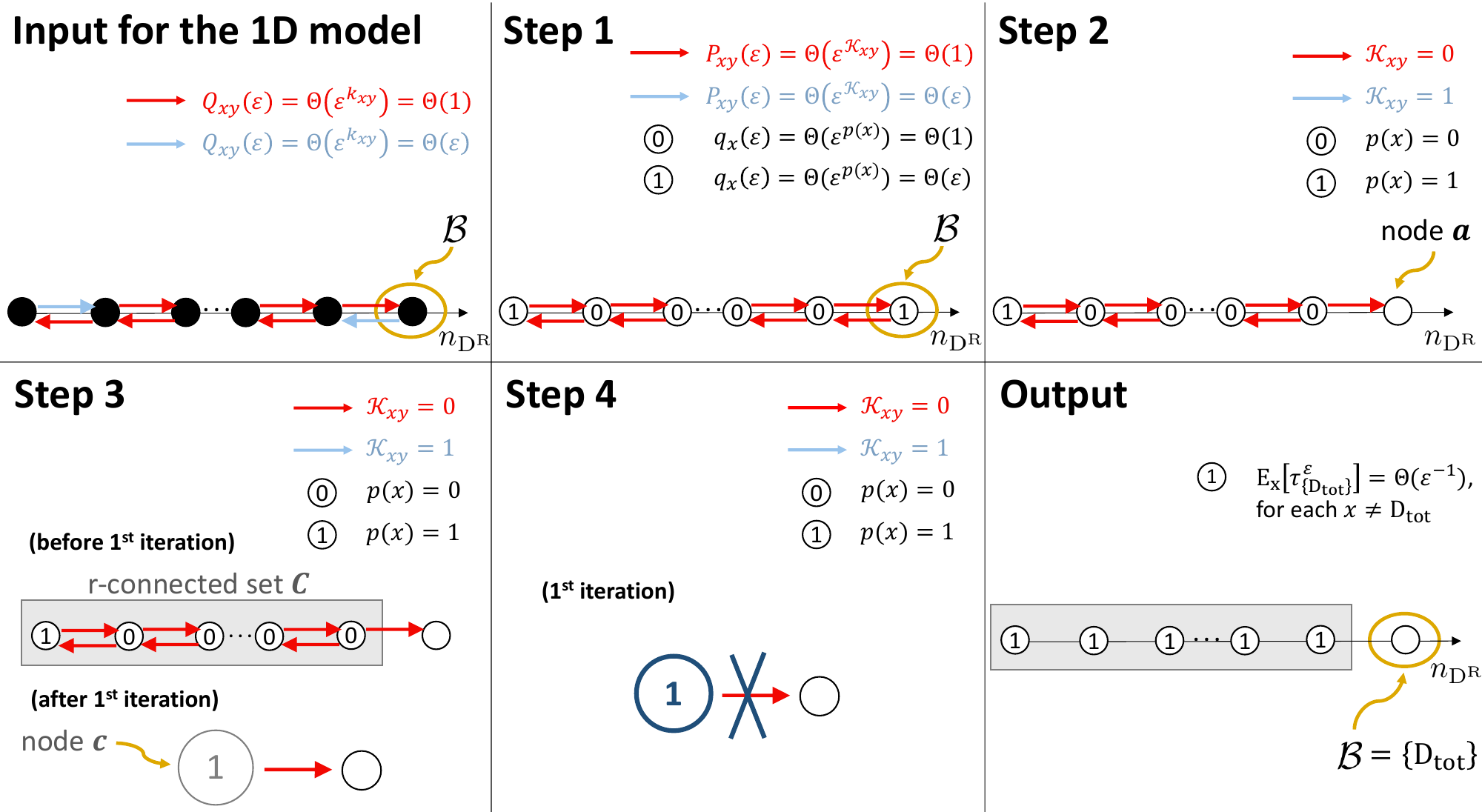}
		\caption{\small { \bf Key steps of the algorithm for the 1D model.} The algorithm is described in Section \ref{orderMFPT}, and it finds the order of the pole of the mean first passage time to $\B \subset \X$ from each state in $\B^c$. In our 1D model, the input for the algorithm is the order of each of the non-zero off-diagonal entries in $Q(\eps)$ and the set $\B=\{\mathrm{D_{tot}}\}$. The order of the non-zero entries in $Q(\eps)$ is represented by colored arrows in the graph in the ``Input'' panel. Step 1 transforms the orders in the $Q(\eps)$-matrix into the orders in the $P(\eps)$-matrix and the exponential parameters $q(\eps)$ to give an equivalent construction for the continuous time Markov chain. The order of the non-zero entries in $P(\eps)$ is represented by colored arrows in the graph, and the number in the circle at a state $x \in \B^c$ is the order of the pole $p(x)$ of $\frac{1}{q_x(\eps)}$ (the mean sojourn time at the state $x$). In Step 2, the set $\B$ is relabeled as the node $a$, and then all transitions from $a$ to $\B^c$ are removed. Step 3 for the 1D model involves only one iteration, where the collection of all nodes except the node $a$ (called an r-connected set $C$) is condensed to a single node $c$, and the order of the pole at $c$ is $p(c) = \max_{u \in C} p (u) + \min \{\mathcal{K}_{uv}: u \in C, v\notin C \text{ and } (u,v) \in E \} = 1 + 0 = 1$, where $E$ denotes the edge set of the graph in Step 3 before the $1^{st}$ iteration. Moreover, $\mathcal{K}_{ca} = \min \{\mathcal{K}_{ua}:\: u \in C \text{ and } (u,a) \in E \} - \min \{\mathcal{K}_{uv}:\: u \in C, v\notin C \text{ and } (u,v) \in E \} = 0-0 = 0$. Step 4 involves one iteration. In this iteration, the node $c$ is the only node other than $a$, so its value of $p$ is fixed, and then any edges leading to or from $c$ are removed. When all of the nodes other than $a$ have been fixed, the order of the pole of the mean first passage time from each state in $\B^c$ to $\B$ is given by the fixed value of the node to which the state belongs.
		}
		\label{fig:algorithm_1D}
	\end{figure}
 
	\subsubsection{Illustrative examples: 1D and 2D models}
    \label{sec:MFPTexamples}
    
	We first apply the algorithm given in Section \ref{orderMFPT} to find the order of the pole of the time to memory loss in the 1D and 2D models introduced in Section \ref{motivexamples}. For the 1D model, we could also directly derive the analytical expression for the time to memory loss by exploiting first step analysis \cite{Norris} and solve the system (\ref{MFPTsystem}) introduced in Section \ref{sec:pertCTMC} (see SI - Section \ref{sec:Appendix1Dmodel}). Figure \ref{fig:algorithm_1D} illustrates the key steps of the algorithm for the 1D model, which lead to the conclusion that the time to memory loss for the active state is $\Theta(\eps^{-1})$. Because of the symmetry in the input graph in Figure \ref{fig:algorithm_1D}, the time to memory loss for the repressed state is also $\Theta(\eps^{-1})$. These orders found by the algorithm are consistent with what can be directly derived by first step analysis. Similarly, SI - Figure \ref{fig:algorithm_2D} illustrates the key steps of the algorithm for the 2D model, which leads to the conclusion that the time to memory loss of both the active and the repressed states is $\Theta(\eps^{-1})$.
	
Next, we find the leading coefficient for the time to memory loss in the 1D and 2D models, which is the coefficient of the $\eps^{-1}$ term in all cases. Recall from Section \ref{sec:IllustrativeStationaryDistribution} that Assumptions \ref{assumption:transient_absorbing}, \ref{assumption:A_is_recurrent} and \ref{assumption:linear_case} hold for both 1D and 2D models and hence by Lemma \ref{lemma3}, so does Assumption \ref{assumption:SingleRecurrentClass} and $\hat{X}_\A$ is irreducible. For the 1D model, $Q_\A$ is given by \eqref{QA1D}. Thus, by Theorem \ref{thm:MFPTcoefficient}, the leading coefficient of the time to memory loss for the active state is the mean first passage time from the fully active state $a$ to the fully repressed state $r$ in $\hat{X}_\A$, which has an exponential distribution with parameter $(Q_\A)_{a,r} = \frac{1-\mu}{1-\mu^{\Dtot}} \frac{k_E^A}{V} \Dtot^2$ since $\hat{X}_\A$ has only two states. Thus, 
\begin{equation*}
 h_{a,r}(\eps)= \frac{1-\mu^{\Dtot}}{1-\mu} \frac{V}{k_E^A} \frac{1}{\Dtot^2} \eps^{-1} + O(1),
 \end{equation*}

and similarly, the time to memory loss for the repressed state is
\begin{equation*}
 h_{r,a}(\eps)= \frac{1-\mu^{\Dtot}}{1-\mu} \frac{V}{k_E^A} \frac{1}{b \mu^{\Dtot} \Dtot^2} \eps^{-1} + O(1).
 \end{equation*}

Similarly, in the 2D model, by Theorem \ref{thm:MFPTcoefficient}, 
 \begin{equation*}
 h_{a,r}(\eps) = \frac{1}{(Q_\A)_{a,r}} \eps^{-1} + O(1) \quad \text{ and } \quad h_{r,a}(\eps) = \frac{1}{(Q_\A)_{r,a}} \eps^{-1} + O(1).
 \end{equation*}
 
\noindent	As an example, when $\Dtot=2$, $Q_\A$ is shown in \eqref{QA2D} and we obtain that 
\begin{equation*}
\begin{aligned}
h_{a,r}(\eps) &= \frac{V}{k_M^A} \frac{K}{4 (k_{W0}^R+k_{W}^R)(k_{W0}^R+k_{W}^R + \frac{k_{M}^R}{V})} \eps^{-1} + O(1)\;\;\;\;\mathrm{and}\\
h_{r,a}(\eps) &= \frac{V}{k_M^A} \frac{K}{4 b\mu^2(k_{W0}^A+k_{W}^A)(k_{W0}^A+k_{W}^A+\frac{k_{M}^A}{V})} \eps^{-1} + O(1),
\end{aligned}
\end{equation*}
with $K$ defined in \eqref{QA2DcoefK}.

%% file: Monotonicity.tex
\subsection{Monotonic dependence on parameters}
\label{sec:Mon}
An important aspect to consider in the study of the stochastic behavior of the chromatin modification circuit is that the erasure rate is different for each type of chromatin modification. These differences can introduce asymmetries in the system that can affect the stationary distribution and the time to memory loss of the active state and repressed state. These asymmetries are captured by the two parameters $\mu$ and $\mu'$. In particular, $\mu$ quantifies the asymmetry between erasure rates of repressive and activating histone modifications and $\mu'$ quantifies the asymmetry between erasure rates of DNA methylation and activating histone modifications. In order to determine how the different chromatin modification erasure rates affect the stochastic behavior of the system, we study how $\mu$ and $\mu'$ affect the stationary distribution and the time to memory loss of the active and repressed gene states.

For the 1D model of the chromatin modification circuit, that does not include DNA methylation, we have an analytical expression for the stationary distribution and the time to memory loss ((\ref{pi01Dmodel}), (\ref{formulaDto0INTRO}), and  (\ref{formula3INTRO})) and we can understand the effect of $\mu$ by directly studying the formulas. However, for the higher-than-1D models we do not have an explicit expression for the stationary distribution or time to memory loss. This is the reason why for these models we exploit the comparison theory developed in \cite{Monotonicitypaper} that allows to determine how $\mu$ and $\mu'$ affect the stochastic behavior of the system through the construction of a coupling between processes with different values for these parameters. In the next subsection, we briefly summarize the relevant theory from \cite{Monotonicitypaper}. 


\subsubsection{Comparison theorems for continuous time Markov chains}
\label{sec:ComparisonTheoremsForStochasticProcesses}

    Denote by $\leq$ the usual componentwise partial order on $\R^d$, i.e., for $x,y \in \R^d$, $x \leq y$ whenever $x_i \leq y_i$ for every $1 \leq i \leq d$. Let $m,d \geq 1$ be integers, consider a matrix $A \in \R^{m \times d}$, where no row of $A$ is identically zero, and consider the following definition. 
    
    \begin{definition}[Definition 3.1 from \cite{Monotonicitypaper}]
    For $x,y \in \R^d$, we say that $x \preccurlyeq_A y$ whenever $A(y-x) \geq 0$ and we say that $x \sim_A y$ whenever $Ax=Ay$.
    \end{definition}
    
    For the matrix $A$, consider the convex cone $K_A := \{ x \in \R^d:\: Ax \geq 0\}$, and, for any $x \in \R^d$, consider the set $K_A +x =\{ y \in \R^d:\: A(y-x) \geq 0 \} = \{ y \in \R^d:\: x \preccurlyeq_A y \}$ and the sets $\partial_i(K_A+x) := \{ y \in K_A +x:\: \inn{A_{i\bullet},y} = \inn{A_{i\bullet},x} \}$ \footnote{Here, for convenience of notation, let $A_{i\bullet}$ denote the row vector corresponding to the $i$-th row of $A$, for $1 \leq i \leq m$. In this article, we will adopt the convention of considering the inner product $\inn{\cdot,\cdot}$ as a function of a row vector in its first entry and as a function of a column vector in the second entry. In particular, $\inn{A_{i\bullet},x} = \sum_{k=1}^d A_{ik}x_k$.} for $1 \leq i \leq m$. Then, the boundary of $K_A + x$ can be expressed as
        \begin{equation*}
         \partial(K_A +x) = \bigcup_{i=1}^{m}  \partial_i(K_A+x).  
        \end{equation*}
    
    Consider a non-empty set $\X \subseteq \Z_+^d$, we will say that a set $\Gamma \subseteq \X$ is \textbf{increasing} with respect to $\preccurlyeq_A$ if for every $x \in \Gamma$ and $y \in \X$, $x \preccurlyeq_A y$ implies that $y \in \Gamma$. We observe that, for $x \in \X$, the set
      \begin{equation}
            (K_A +x) \cap \X = \{ y \in \X:\: x \preccurlyeq_A y \}
        \end{equation}
   is increasing by the transitivity property of $\preccurlyeq_A$. On the other hand, we will say that a set $\Gamma \subseteq \X$ is \textbf{decreasing} with respect to $\preccurlyeq_A$ if for every $x \in \Gamma$ and $y \in \X$, $y \preccurlyeq_A x$ implies that $y \in \Gamma$. 
    
    Now, consider a non-empty set $\X \subseteq \Z_+^d$ and a finite set of distinct nonzero possible transition vectors for a pair of continuous time Markov chains on $\X$. We denote the set of vectors by $\{v_1,\ldots,v_n\} \subseteq \Z^d \setminus \{0\}$, where $0$ is the origin in $\Z^d$. Consider two collections of functions $\rate=(\rate_1, \dots,\rate_n)$ and $\Breve{\rate}= (\Breve{\rate}_1, \dots,\Breve{\rate}_n)$ from $\X$ into $\R_+$ such that $\rate_j(x) = \Breve{\rate}_j(x)= 0$ if $x+v_j \notin \X$. Assume that $Q=(Q_{x,y})_{x,y \in \X}$, given by \eqref{eq:TransitionMatrixQ}, is the infinitesimal generator for a continuous time Markov chain $X$ and $\breve{Q}$, defined by \eqref{eq:TransitionMatrixQ} but with functions $\Breve{\rate}_1, \dots,\Breve{\rate}_n$ in place of $\rate_1, \dots,\rate_n$, is the infinitesimal generator for a continuous time Markov chain $\Breve{X}$. We call $X$ and $\Breve{X}$ the continuous time Markov chains associated with $\rate$ and $\Breve{\rate}$ respectively.

    The following stochastic comparison result was proved in Campos et al. \cite{Monotonicitypaper}. The condition $(i)$ of the theorem and $A \in \Z^{n \times d}$ ensure that to go outside of $K_A+x$, the Markov chains will necessarily hit the boundary of $K_A+x$.
    

    \begin{theorem}[Theorems 3.2, 3.4, 3.5 from \cite{Monotonicitypaper}]
    \label{thm:Compartheorem}
   With $\X$, $v_1,\ldots,v_n$, $\rate$ and $\Breve{\rate}$ as described above, assume that the continuous time Markov chains associated with $\rate$ and $\Breve{\rate}$ do not explode in finite time. Consider a matrix $A \in \Z^{m \times d}$ with nonzero rows and suppose that both of the following conditions hold:
        \begin{enumerate}
            \item[(i)]
            For each $1\leq j \leq n$, the vector $Av_j$ has entries in $\{-1,0,1\}$ only.
             \item[(ii)]
             For each $x \in \X$, $1 \leq i \leq m$ and $y \in \partial_i(K_A+x) \cap \X$ we have that
            \begin{equation}
            \label{eq:CouplingConditions_InnerProduct_I}
            \Breve{\rate}_j(y) \leq \rate_j(x), \quad \text{for each } 1 \leq j \leq n \text{ such that } \inn{A_{i\bullet},v_j} < 0,
            \end{equation}
      and  
            \begin{equation}
            \label{eq:CouplingConditions_InnerProduct_II}
            \Breve{\rate}_j(y) \geq \rate_j(x), \quad \text{for each } 1 \leq j \leq n \text{ such that } \inn{A_{i\bullet},v_j} > 0.
            \end{equation}
        \end{enumerate}
    Then, for each pair $\initialx,\initialxbreve \in \X$ such that $\initialx \preccurlyeq_A \initialxbreve$, there exists a probability space $(\Omega,\F,\PP)$ with realizations of the two continuous time Markov chains $X = \{X(t):\: t \geq 0\}$ and $\Breve{X}=\{\Breve{X}(t):\: t \geq 0\}$ defined there, each having state space $\X \subseteq \Z^d_+$, with infinitesimal generators given by $Q$ and $\Breve{Q}$, associated with $\rate$ and $\Breve{\rate}$, respectively, with initial conditions $X(0)=\initialx$ and $\Breve{X}(0)=\initialxbreve$, and such that:
        \begin{equation}
        \label{eq:DominanceProperty_InnerProduct}
        \PP\left[X(t) \preccurlyeq_A \Breve{X}(t) \text{ for every } t \geq 0 \right]=1.   
        \end{equation}
Furthermore, for a non-empty set $\Gamma \subseteq \X$, consider $\tau_{\Gamma} := \inf\{ t \geq 0:\: X(t) \in \Gamma \}$ and $\Breve{\tau}_{\Gamma} := \inf\{ t \geq 0:\: \Breve{X}(t) \in \Gamma \}$. If $\Gamma$ is increasing with respect to the relation $\preccurlyeq_A$, then $\E[\Breve{\tau}_{\Gamma}] \le \E[\tau_\Gamma]$.
   If $\Gamma$ is decreasing with respect to the relation $\preccurlyeq_A$, then $\E[\tau_\Gamma] \le \E[\Breve{\tau}_{\Gamma}]$.
    Finally, suppose that the two continuous time Markov chains are irreducible and positive recurrent on $\X$, and denote the associated stationary distributions by $\pi$ and $\Breve{\pi}$, respectively. Then, if $\Gamma \subseteq \X$ is a non-empty set that is increasing with respect to $\preccurlyeq_A$, we have $\sum_{x \in \Gamma} \pi_x \leq  \sum_{x \in \Gamma} \Breve{\pi}_x$,
     or if $\Gamma \subseteq \X$ is a non-empty set that is decreasing with respect to $\preccurlyeq_A$, we have $\sum_{x \in \Gamma} \Breve{\pi}_x \leq  \sum_{x \in \Gamma} \pi_x$.
\end{theorem}    
    
	\subsubsection{Illustrative example: 2D model}\label{2DTTMLresultsmu}
	
We are interested in determining how the asymmetry of the system, represented by the parameter $\mu= k^R_E/k^A_E$ affects the stationary distribution $\pi(\eps)$ and the times to memory loss, $h_{a,r}(\eps)$ and $h_{r,a}(\eps)$, of the active ($a=(0,\mathrm{D_{tot}})^T$) and repressed ($r=(\mathrm{D_{tot}},0)^T$) states, respectively, for the continuous time Markov chain $X^{\eps}$ described in Section \ref{exp2}. For this, we use Theorem \ref{thm:Compartheorem}.
%
%
For $\eps \in (0,\eps_0)$, let $X^{\eps}$ be the continuous time Markov chain with
\begin{equation}\label{propfun2D}
\rate_1(x)=f_A(x),\;\;\rate_2(x)=g^{\eps}_A(x),\;\;\rate_3(x)=f_R(x),\;\;\rate_4(x)=g^{\eps}_R(x),\;\;\;\;x \in \X
\end{equation}
with $\X$, $v_1,...,v_4$, and $f_A(x)$, $g^{\eps}_A(x)$, $f_R(x)$, $g^{\eps}_R(x)$ as defined in Section \ref{exp2}, and introduce the continuous time Markov chain $\Breve{X}^{\eps}$ defined on $\X$, having the same transition vectors of $X^{\eps}$, and having infinitesimal transition rates $\Breve{\rate}_{1}(x),...,\Breve{\rate}_{4}(x)$ defined as for $\rate_{1}(x),...,\rate_{4}(x)$, with all the parameters having the same values except that $\mu$ is replaced by $\Breve{\mu}$, where $\mu \ge \Breve{\mu}$. Let    
    \begin{equation}\label{matrixA2D}
       A= \begin{bmatrix}
        1 & 0\\
        0 & -1
        \end{bmatrix}
    \end{equation}
    and let us consider the partial order $x \preccurlyeq_A y$. A similar example was analyzed by Campos et al. \cite{Monotonicitypaper} -  Example 4.4, using the results of Theorem \ref{thm:Compartheorem}.
    %
    %
    The only differences are that, in \cite{Monotonicitypaper}, the matrix $A$ is the negative of the matrix given in \eqref{matrixA2D} and the inequality between $\mu$ and $\Breve{\mu}$ is the opposite compared to the one considered here. The relationship between the notation in \cite{Monotonicitypaper} and our notation is $\kappa_{1a}=k_{W0}^A+k_{W}^A$, $\kappa_{1b}=(k_{M}^A/V)$, $\kappa_{2a}=k_{W0}^R+k_{W}^R$, $\kappa_{2b}=(k_{M}^R/V)$, $\kappa_{3a}=\eps (k_{M}^A/V)$, $\kappa_{3b}=(k^A_E/V)$, $c=b$.
   
   From the analysis in \cite{Monotonicitypaper}, we can directly conclude that, if $\pi(\eps)$ is the stationary distribution for $X^{\eps}$ and 
   $\Breve{\pi}(\eps)$ the stationary distribution for $\Breve{X}^{\eps}$, then $\Breve{\pi}_a(\eps) \leq \pi_a(\eps)$ and $\Breve{\pi}_r(\eps) \geq \pi_r(\eps)$. This implies that increasing $\mu$ increases the probability of the system in steady-state being in the active state $a$ to the detriment of the repressed state $r$ (and vice versa for decreasing $\mu$). We can also conclude, using natural notation for quantities associated with $X^{\eps}$ and $\Breve{X}^{\eps}$, that, defining $\tau_y^{\eps} = \inf \{t \geq 0:\: X^{\eps}(t)=y\}$ and $\Breve{\tau}_y^{\eps} = \inf \{t \geq 0:\: \Breve{X}^{\eps}(t)=y\}$, $h_{r,a}(\eps)=\E_{r}[\tau_{a}^{\eps}]\le 
   \E_{r}[\Breve{\tau}_{a}^{\eps}]=\Breve{h}_{r,a}(\eps)$ and 
   $\Breve{h}_{a,r}(\eps)=\E_{a}[\Breve{\tau}_{r}^{\eps}]\le 
   \E_{a}[\tau_{r}^{\eps}]=h_{a,r}(\eps)$, implying that the time to memory loss of the repressed state decreases for higher values of $\mu$, while the time to memory loss of the active state increases for higher values of $\mu$.

%% file: Applications.tex
\section{Further Examples}
\label{sec:Applications}

In this section, compared to the models of the chromatin modification circuit introduced in Section \ref{motivexamples}, which do not include DNA methylation, we introduce more elaborate models that include DNA methylation and we study their stochastic behavior by exploiting the theory developed in this paper.

\subsection{3D chromatin modification circuit model}\label{3Dmodelanalysis}
We now introduce a model in which DNA methylation is also a possible chromatin mark. The species involved are D (unmodified nucleosome), $\mathrm{D^R_1}$ (nucleosome with CpGme only), $\mathrm{D^R_{12}}$ (nucleosome with both H3K9me3 and CpGme) and $\mathrm{D^A}$ (nucleosome with an activating histone modification). In particular, we assume that, in order to be modified with both repressive modifications, D is first modified with DNA methylation, obtaining $\mathrm{D^R_1}$, and then with a repressive histone modification, obtaining $\mathrm{D^R_{12}}$. The opposite order of modifications is not allowed. 
    \begin{figure}[h!]
    \centering
    \includegraphics[scale=0.4]{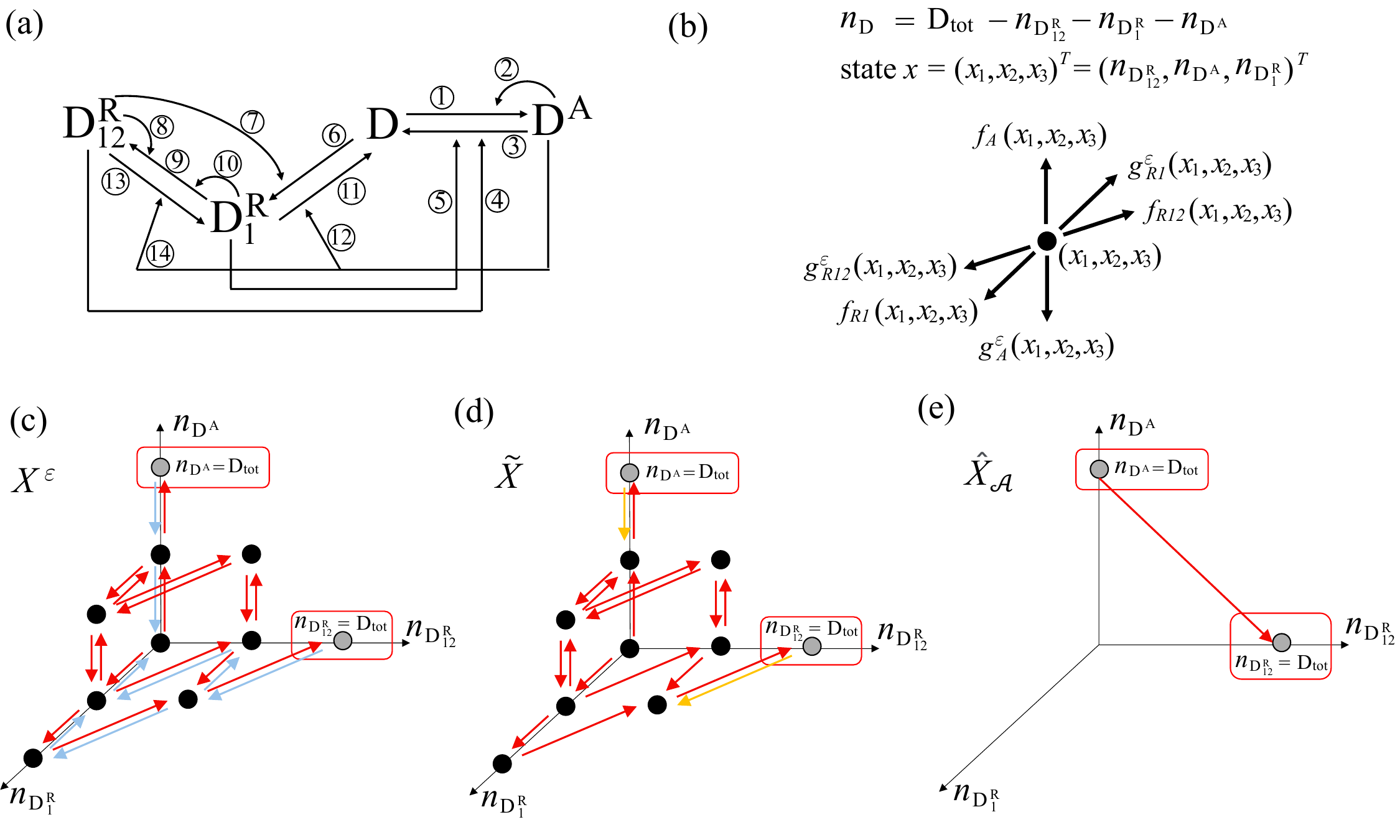}
    \caption{\small { \bf 3D Model and associated Markov chain.} (a) Original chemical reaction system.  The numbers on the arrows correspond to the reactions associated with the arrows as described in \eqref{reacs3D} in the main text. (b) Directions of the possible transitions of the Markov chain $X^{\eps}$, starting from a state $x=(x_1,x_2,x_3)^T$ and whose rates are given in equation (\ref{rates3D}). (c) Graph for $X^{\eps}$.
    Here, the red (blue) arrows correspond to $O(1)$ ($O(\eps)$) transition rates. (d) Graph for the Markov chain $\tilde{X}$. 
    Here, the gold arrows correspond to transitions that were $O(\eps)$ in $X^{\eps}$ and became $O(1)$ in $\tilde{X}$. 
    (e) Graph for the Markov chain $\hat{X}_{\mathcal{A}}$. For (c), (d), and (e) the state of the Markov chain is $x=(n_{\mathrm{D^R_{12}}},n_{\mathrm{D^A}},n_{\mathrm{D^R_1}})^T$ and we consider $\mathrm{D_{tot}}=2$. In panels (c) - (e), we use gray dots to represent the states belonging to $\A$ and black dots to represent all the other states.
            }   
   \label{fig:3Dmodel}
    \end{figure}

    \clearpage

\noindent This enables us to simplify the model and the related analysis. This assumption will be removed in the 4D model analyzed in Section {\ref{SS4}}. The chemical reaction system for the 3D model, shown in Fig. \ref{fig:3Dmodel}(a), is the following:

\begingroup
\small \label{reacs3D}
	\begin{align}
	\notag &{\large \textcircled{\small 1}}\;\ce{D ->[$k^A_{W0}+k^A_W$] D^A },\;\;\; {\large \textcircled{\small 2}}\;\ce{D + $\mathrm{D^A}$ ->[$k^A_M$] $\mathrm{D^A}$ + $\mathrm{D^A}$ },\;\;\;{\large \textcircled{\small 3}}\;\ce{$\mathrm{D^A}$ ->[$\delta + \bar k^A_{E}$] D },\\
	\notag &{\large \textcircled{\small 4}}\;\ce{$\mathrm{D^A}$ + $\mathrm{D^R_1}$ ->[$k^A_E$] D + $\mathrm{D^R_1}$ },{\large \textcircled{\small 5}}\;\ce{$\mathrm{D^A}$ + $\mathrm{D^R_{12}}$ ->[2$k^A_E$] D + $\mathrm{D^R_{12}}$ },\\
	\label{reacs3D} &{\large \textcircled{\small 6}}\; \ce{D ->[$k^1_{W0}+k^1_W$] $\mathrm{D^R_1}$ },\;\;\;{\large \textcircled{\small 7}}\;\ce{D + $\mathrm{D^R_{12}}$ ->[$k^{'}_M$] $\mathrm{D^R_1}$ + $\mathrm{D^R_{12}}$ },\;\;\;{\large \textcircled{\small 8}}\;\ce{$\mathrm{D^R_1}$ + $\mathrm{D^R_{12}}$ ->[$k_M + \bar k_M$] $\mathrm{D^R_{12}}$ + $\mathrm{D^R_{12}}$ },\\
	\notag &{\large \textcircled{\small 9}}\;\ce{$\mathrm{D^R_1}$ ->[$k^2_{W0}$] $\mathrm{D^R_{12}}$ },\;\;\;{\large \textcircled{\small 10}}\;\ce{$\mathrm{D^R_1}$ + $\mathrm{D^R_1}$ ->[$\bar k_M$] $\mathrm{D^R_{12}}$ + $\mathrm{D^R_1}$ },{\large \textcircled{\small 11}},\;\;\;\ce{$\mathrm{D^R_1}$ ->[$\delta^{'}+k^{'}_{T}$] D },\\
    \notag &{\large \textcircled{\small 12}}\;\ce{$\mathrm{D^R_1}$ + $\mathrm{D^A}$ ->[$k^{'*}_T$] D + $\mathrm{D^A}$ },{\large \textcircled{\small 13}}\;\ce{$\mathrm{D^R_{12}}$ ->[$\delta+\bar k^R_{E}$] $\mathrm{D^R_1}$ },\;\;\;{\large \textcircled{\small 14}}\;\ce{$\mathrm{D^R_{12}}$ + $\mathrm{D^A}$ ->[$k^R_E$] $\mathrm{D^R_1}$ + $\mathrm{D^A}$ },
	\end{align}
   \endgroup
   \normalsize
    where $k^A_{W0},k^A_W,k^A_M,\delta, \bar k^A_E,  k^A_E, k^1_{W0}, k^1_{W}, k^2_{W0}, k^2_{W}, k'_M, \bar k_M, k_M, \delta', k'_T, k^{'*}_T, \bar k^R_E, k^R_E > 0$ and the form of the reaction rate constants is due to the fact that reactions with the same reactants and products have been combined. As we did for the 2D model, define parameters $\eps=\frac{\delta+\bar{k}^A_E}{\frac{k^A_M}{V}\Dtot}$ and $\mu=\frac{k^R_E}{k^A_E}$, with a constant $b$ such that $\frac{\delta+\bar k^R_E}{\delta + \bar k^A_E} = b\mu$. Furthermore, since this model includes DNA methylation, we also define $\mu'=\frac{k^{'*}_T}{k^A_E}$ and a constant $\beta$ such that $\frac{\delta^{'}+k^{'}_T}{\delta + \bar k^A_E} = \beta \mu'$. The parameter $\mu'$ quantifies the asymmetry between the erasure rates of DNA methylation and activating histone modifications.
	The Markov chain $X^{\eps}$ associated with the system is a linearly perturbed finite state continuous time Markov chain with the state $x$ tracking $n_{\mathrm{D^R_{12}}}$, $n_{\mathrm{D^A}}$, $n_{\mathrm{D^R_1}}$, that is,
    the number of nucleosomes of types $\mathrm{D^R_{12}}$, $\mathrm{D^A}$, and $\mathrm{D^R_1}$, respectively. If the total number of modifiable nucleosomes is $\mathrm{D_{tot}}$, which is conserved, the state space is $\mathcal{X}= \{(x_1,x_2,x_3)^T \in \Z_+^3:\: x_1 + x_2 + x_3 \leq \Dtot \}$. The transition vectors for $X^{\eps}$ are given by $v_1=-v_2=(0,1,0)^T$, $v_3=-v_4=(0,0,1)^T$, and $v_5= - v_6=(1,0,-1)^T$. The infinitesimal transition rates are  
 
 \begingroup
     \small
   \begin{equation}
   \label{rates3D}
    \begin{aligned}
        &Q_{x,x+v_1}(\eps)=f_A(x) = (\Dtot -(x_1+x_2+x_3))\left(k_{W0}^A+k_{W}^A + \frac{k_{M}^A}{V}x_2\right),\\
        &Q_{x,x+v_2}(\eps)=g^{\eps}_A(x) = x_2\left(\eps \frac{k_{M}^A}{V}\Dtot + \frac{k^A_E}{V}(x_3+2x_1)\right),\\
        &Q_{x,x+v_3}(\eps)=f_{R1}(x) = (\Dtot -(x_1+x_2+x_3))\left(k^1_{W0}+k^1_{W} + \frac{k^{'}_{M}}{V}x_1\right),\\
        &Q_{x,x+v_4}(\eps)=g^{\eps}_{R1}(x) = x_3\mu'\left(\eps \frac{k_{M}^A}{V}\Dtot \beta + x_2\frac{k^A_E}{V}\right),\\
       &Q_{x,x+v_5}(\eps)=f_{R12}(x) = x_3\left(k^2_{W0}+ \frac{k_{M}}{V}x_1 + \frac{\bar k_{M}}{V}\left(x_1+\frac{x_3-1}{2}\right)\right),\\
       &Q_{x,x+v_6}(\eps)=g^{\eps}_{R12}(x) = x_1\mu\left(\eps \frac{k_{M}^A}{V}\Dtot b + x_2\frac{k^A_E}{V}\right).\\
    \end{aligned}
   \end{equation}
   \endgroup
   \normalsize
%
A representation of the possible transitions, with associated rates, and the Markov chain graph for $\mathrm{D_{tot}}=2$ are given in Fig. \ref{fig:3Dmodel}(b) and (c), respectively. Each rate depends on the state $x$.

\subsubsection{Stationary distribution}
\label{SD3D}

    We now focus on the expansion as a function of $\eps$ of the stationary distribution for the 3D model. In SI - Section \ref{sec:Appendix3Dmodel}, we show that, when $\eps=0$, the continuous time Markov chain associated with the 3D model has transient states $\T= \{i_1, \ldots, i_m\}$ where $m =\sum_{j=0}^{\Dtot}\left(\frac{(j+2)(j+1)}{2}\right)-2$, $i_1=(0,\Dtot-1,0)^T$, $i_m=(\Dtot-1,0,1)^T$, and absorbing states $\A = \{a,r\}$, with $a=(0,\Dtot,0)^T$ corresponding to the fully active state ($n_{\DA}=\Dtot$) and $r=(\Dtot,0,0)^T$ corresponding to the fully repressed state ($n_{\mathrm{D^R_{12}}}=\Dtot$), respectively. Then, Assumption \ref{assumption:transient_absorbing} holds (see SI - Section \ref{sec:Appendix3Dmodel}). Furthermore, $\X = \A \cup \T$ and from  \eqref{rates3D} we see that $Q(\eps)$ can be written in the form \eqref{eq:Q_eps_decomposition}, where $Q(\eps)$ is a linear perturbation of $Q(0)$.
%
%
%
Hence, Assumption \ref{assumption:linear_case} holds. Assumption \ref{assumption:SingleRecurrentClass} also holds, where the recurrent class is $\{r\}$ (see SI - Section \ref{sec:Appendix3Dmodel}). Then, we can apply SI - Theorem \ref{thm:HigherOrderTerms_for_stationary distribution}. We first obtain that $\pi(0)=\pi^{(0)}=[\alpha,0] = \left[\alpha_a,\alpha_r, 0 \ldots ,0 \right]$ where $\alpha$ is the unique stationary distribution for the process $\hat{X}_{\A}$ with infinitesimal generator $Q_{\A} = A_1 +  S_1(-T_0)^{-1}R_0$. Since the recurrent class $\{r\}$ is a singleton and $\alpha$ is supported on $\{r\}$, we must have 
%
 $\alpha_a=0$ and $\alpha_r=1$.
    
    We now derive an expression for $\pi^{(1)}$. For the transient states $\T= \{i_1, \ldots, i_m\}$, $\beta^{(1)}=[\pi^{(1)}_{i_1},...,\pi^{(1)}_{i_m}]= \alpha S_1(-T_0)^{-1}=[0,...,0,\pi^{(1)}_{i_m}]$, with  \begin{equation*}
    \pi^{(1)}_{i_m}=\frac{\mu b\frac{k_{M}^A}{V}\Dtot^2}{k^2_{W0}+ (\frac{k_{M}}{V} + \frac{\bar k_{M}}{V})(\Dtot - 1)}.
    \end{equation*}
    See SI - Section \ref{sec:Appendix3Dmodel} for the detailed mathematical derivation. Now, $\alpha^{(1)}=[\pi^{(1)}_{a},\pi^{(1)}_{r}]$ is the unique vector such that $\alpha^{(1)}Q_{\A} = -\beta^{(1)}[R_1 + T_1(-T_0)^{-1}R_0],$ $\alpha^{(1)}\one = -\beta^{(1)}\one$.
    
    As an illustration, suppose $\Dtot=2$. Then (see SI - Section \ref{sec:Appendix3Dmodel} for the detailed mathematical derivation),
    	\begin{equation}\label{QA3D}
    	    Q_\A = \frac{K_1+\mu K_2}{K_3+\mu K_4 + \mu' K_5 + \mu \mu' K_6} \left( \begin{array}{c c}
	    -1 & 1 \\
	    0 & 0 \\
    \end{array} \right),
    \end{equation}
    \begin{equation}\label{pi3D}
        \pi^{(1)}_{a}=\frac{\mu^2\mu'^2K_7}{K_8(K_9+K_{10}\mu)},\;\;\;
        \pi^{(1)}_{r}=-\pi^{(1)}_{a}-\pi^{(1)}_{i_m}=-\frac{\mu^2\mu'^2K_7+\mu K_{11}(K_{9}+K_{10}\mu)}{K_{8}(K_{9}+K_{10}\mu)},
        \end{equation}
   with $m=8$ and $K_i$, $i=1,...,11$, are non-negative constants independent of $\eps$, $\mu$ and $\mu'$ (see SI - Section \ref{sec:Appendix3Dmodel} for their precise definitions). Hence, the stationary distribution for $\Dtot=2$ satisfies
               \begin{equation}
        \label{pi3DmodelRESULTS}
         \pi_x(\eps) = \begin{cases}
         \eps\frac{(\mu\mu')^2K_7}{K_8(K_9+K_{10}\mu)} + O(\eps^2) & \text{ if } x=a=(0,2,0)^T \\
                    O(\eps^2)& \text{ if } x \in \T \backslash \{i_m\} \\    
                    \eps \frac{K_{11}}{K_8}\mu + O(\eps^2)& \text{ if } x = i_m = (1,0,1)^T \\     1-\eps\frac{(\mu\mu')^2K_7+\mu K_{11}(K_9+K_{10}\mu)}{K_8(K_9+K_{10}\mu)} + O(\eps^2)& \text{ if } x = r= (2,0,0)^T.
        \end{cases}
        \end{equation}
    For small $\eps>0$, the stationary distribution is concentrated around the active and repressed states, although more mass is concentrated around the repressed state. However, higher values of $\mu'$ increase the probability of being in the active state, while decreasing the probability of being in the repressed state.

   \subsubsection{Time to memory loss}
   \label{3DTTML}
   	
    In this section, we determine how the leakage of the system ($\eps$) and the asymmetry between activating histone modifications and DNA methylation ($\mu'$) affect the time to memory loss of the active state $h_{a,r} (\eps)$ and the time to memory loss of the repressed state $h_{r,a} (\eps)$. 
	
	Firstly, by the algorithm in Section \ref{orderMFPT}, we have that $h_{a,r} (\eps)$ is $O(\eps^{-1})$ and $h_{r,a} (\eps)$ is $O(\eps^{-2})$ (see SI - Section \ref{sec:alg_3D4D}). This means that decreasing the leakage extends the memory of both the active and repressed chromatin states, but the effect is stronger for the repressed state. This difference is influenced by the co-existence and cooperation between DNA methylation and repressive histone marks that introduce a structural bias in the 3D chromatin modification circuit towards a repressed chromatin state. 
 %
%

    These results are consistent with the ones obtained by applying Theorem \ref{thm:MFPTcoefficient}, which allows us not only to find the order of $h_{a,r} (\eps)$ and $h_{r,a} (\eps)$, but also to find an expression for their leading coefficients (see SI - Section \ref{sec:Appendix3Dmodel} for the detailed mathematical derivation).
    As an example, when $\Dtot=2$, $Q_\A$ and $\pi_a^{(1)}$ are shown in \eqref{QA3D} and \eqref{pi3D}, and we obtain from them that 
    \begin{align}
    \label{form3D1}    h_{a,r}(\eps)&=\frac{K_3+\mu K_4 + \mu' K_5 + \mu \mu' K_6}{K_1+\mu K_2} \frac{1}{\eps} + O(1),\;\;\;\mathrm{and}\\
    \label{form3D2} h_{r,a}(\eps)&=\frac{K_3+\mu K_4 + \mu' K_5 + \mu \mu' K_6}{K_1+\mu K_2} \frac{K_8(K_9+K_{10}\mu)}{\mu^2\mu'^2K_7} \frac{1}{\eps^2} + O\left(\frac{1}{\eps}\right),
    \end{align}
    where $K_i$, $i=1,...,11$, are non-negative constants independent of $\eps$, $\mu$ and $\mu'$, defined in SI - Section \ref{sec:Appendix3Dmodel}.

Now, we focus on understanding how the asymmetry between chromatin modification erasure rates affects the time to memory loss. In particular, since experimental data suggest that the asymmetry between the erasure rates of DNA methylation and activating histone modifications is more pronounced than the asymmetry between erasure rates of opposite histone modifications, in this analysis we focus only on studying the effect of $\mu'$, but a similar procedure to the one presented in the next paragraph could be applied to study the effect of $\mu$. 
%
%
To this end, we exploit the comparison Theorem \ref{thm:Compartheorem} to determine directly how $\mu'$ affects $h_{a,r}(\eps)$ and $h_{r,a}(\eps)$, without deriving an explicit expression for them. To this end, we first note that the transitions of the Markov chain $X^{\eps}(t)$ are in six possible directions, that can be written as $v_1=(0,1,0)^T$, $v_2=(0,-1,0)^T$, $v_3=(0,0,1)^T$, $v_4=(0,0,-1)^T$, $v_5=(1,0,-1)^T$, $v_6=(-1,0,1)^T$, with the associated infinitesimal transition rates that can be written as $\rate_1(x)=f_A(x)$, $\rate_2(x)=g^A_\eps(x)$, $\rate_3(x)=f_{R1}(x)$, $\rate_4(x)=g^{\eps}_{R1}(x)$, $\rate_5(x)=f_{R12}(x)$, $\rate_6(x)=g^{\eps}_{R12}(x)$.	       
%
Define the matrix
    \begin{equation}\label{matrixA3D}
       A= \begin{bmatrix}
        1 & 0 & 0\\
        0 & -1 & 0\\
        1 & 0 & 1
        \end{bmatrix}\nonumber
    \end{equation}
    and, for $x\in \X$, $(K_A+x)\cap \X=\{ w \in \X:\: x \preccurlyeq_A w \}$. Let us also introduce infinitesimal transition rates $\Breve{\rate}_{i}(x)$, $i=1,2, ..., 6$, defined as for $\rate_{i}(x)$, $i=1,2, ..., 6$, with all the parameters having the same values except that $\mu'$ is replaced by $\Breve{\mu}'$, with $\mu' \ge \Breve{\mu}'$. 
   All of the conditions of Theorem \ref{thm:Compartheorem} hold (see SI - Section \ref{sec:Appendix3Dmodel}) and so we can apply the theorem.
   %
    %
    This allows us to establish that, since $a=(0,\mathrm{D_{tot}},0)^T \preccurlyeq_A r=(\mathrm{D_{tot}},0,0)^T$, then $\Breve{h}_{a,r}(\eps)\le h_{a,r}(\eps)$ and $h_{r,a}(\eps)\le \Breve{h}_{r,a}(\eps)$, where $\Breve{}$ indicates quantities associated with $\Breve{\rate}$. Thus, we can conclude that, given that the only difference between the two systems was that $\mu' \ge\Breve{\mu}'$, the time to memory loss of the active state is monotonically increasing with $\mu'$, while the time to memory loss of the repressed state is monotonically decreasing with $\mu'$.
      \begin{figure}[h!]
    \centering
    \includegraphics[scale=0.41]{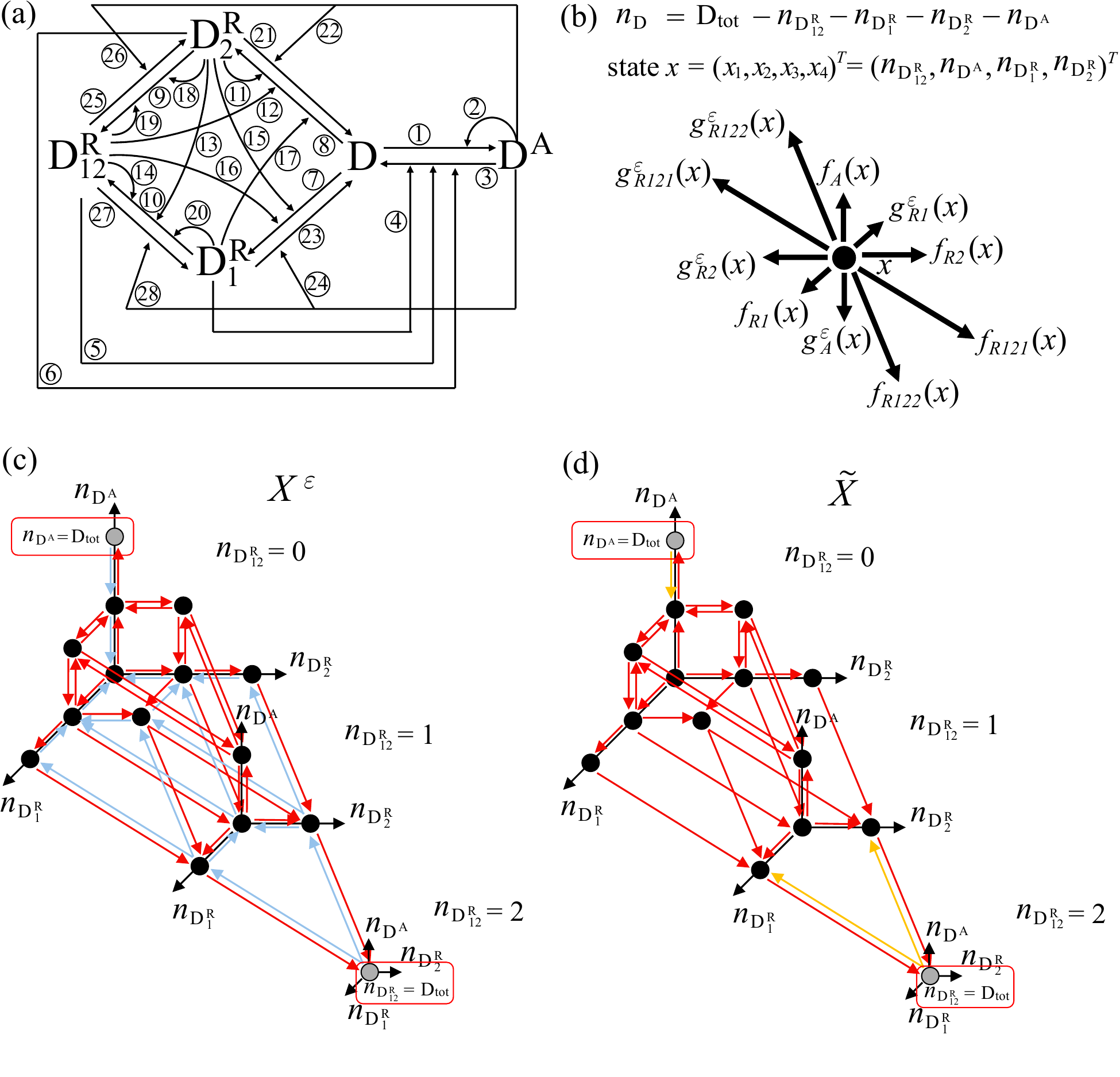}
    \caption{\small { \bf 4D Model and associated Markov chain.} (a) Chemical reaction system. The numbers on the arrows correspond to the reactions associated with the arrows as described in (\ref{reacs4D}) in the main text. (b) Directions of the possible transitions of the Markov chain $X^{\eps}$ associated with the reduced SCRN, starting from a state $x=(x_1,x_2,x_3,x_4)^T$ and whose rates are given in equation (\ref{rates4D}). (c) Graph for $X^{\eps}$.
    Here, the red (blue) arrows correspond to $O(1)$ ($O(\eps)$) transition rates. (d) Graph for the Markov chain $\tilde{X}$. 
    Here, the golden arrows correspond to the transitions that were $O(\eps)$ in $X^{\eps}$ and became $O(1)$ in $\tilde{X}$. For (c) and (d) the state of the Markov chain is $x=(n_{\mathrm{D^R_{12}}},n_{\mathrm{D^A}},n_{\mathrm{D^R_1}},n_{\mathrm{D^R_2}})^T$, we consider $\mathrm{D_{tot}}=2$, and we show three interconnected slices ($n_{\mathrm{D^R_{12}}}=0,1,2$) of the Markov chain state space. In panels (c) and (d), we use gray dots to represent the states belonging to $\A$ and black dots to represent all the other states.
            }       
   \label{fig:4Dmodel1}
    \end{figure} 

\clearpage

\subsection{4D chromatin modification circuit model}\label{SS4}
Now, we consider a complete model in which the species involved are D, $\mathrm{D^R_1}$, $\mathrm{D^R_{12}}$, $\mathrm{D^A}$ and $\mathrm{D^R_2}$ (nucleosome with H3K9me3 only). Compared to the 3D model, we assume that, in order to be modified with both repressive modifications, D can be also modified first with a repressive histone modification (H3K9me3), obtaining $\mathrm{D^R_2}$, and then with DNA methylation (CpGme), obtaining $\mathrm{D^R_{12}}$. The chemical reaction system, shown in Fig. \ref{fig:4Dmodel1}(a), is the following:

\begingroup
\small 
\begin{align}
\notag &{\large \textcircled{\small 1}}\;\ce{D ->[$k^A_{W0}+k^A_W$] D^A },\;\;\;{\large \textcircled{\small 2}}\;\ce{D + $\mathrm{D^A}$ ->[$k^A_M$] $\mathrm{D^A}$ + $\mathrm{D^A}$ },\;\;\;{\large \textcircled{\small 3}}\;\ce{$\mathrm{D^A}$ ->[$\delta + \bar k^A_{E}$] D },\\
\notag &{\large \textcircled{\small 4}}\;\ce{$\mathrm{D^A}$ + $\mathrm{D^R_1}$ ->[$k^A_E$] D + $\mathrm{D^R_1}$ },\;\;\;{\large \textcircled{\small 5}}\;\ce{$\mathrm{D^A}$ + $\mathrm{D^R_{12}}$ ->[2 $k^A_E$] D + $\mathrm{D^R_{12}}$ },\;\;\;{\large \textcircled{\small 6}}\;\ce{$\mathrm{D^A}$ + $\mathrm{D^R_2}$ ->[$k^A_E$] D + $\mathrm{D^R_2}$ },\\
\notag &{\large \textcircled{\small 7}}\;\ce{D ->[$k^1_{W0}+k^1_{W}$] $\mathrm{D^R_1}$ },\;\;\;{\large \textcircled{\small 8}}\;\ce{D ->[$k^2_{W0}+k^2_{W}$] $\mathrm{D^R_2}$ },\;\;\;{\large \textcircled{\small 9}}\;\ce{$\mathrm{D^R_2}$ ->[$k^1_{W0}$] $\mathrm{D^R_{12}}$ },\;\;\;{\large \textcircled{\small 10}}\;\ce{$\mathrm{D^R_1}$ ->[$k^2_{W0}$] $\mathrm{D^R_{12}}$ },
\\
\notag & {\large \textcircled{\small 11}}\;\ce{D + $\mathrm{D^R_2}$ ->[$k_M$] $\mathrm{D^R_2}$ + $\mathrm{D^R_2}$ },\;\;\;{\large \textcircled{\small 12}}\;\ce{D + $\mathrm{D^R_{12}}$ ->[$k_M + \bar k_M$] $\mathrm{D^R_2}$ + $\mathrm{D^R_{12}}$ },\\
\label{reacs4D} &{\large \textcircled{\small 13}}\;\ce{$\mathrm{D^R_1}$ + $\mathrm{D^R_2}$ ->[$k_M$] $\mathrm{D^R_{12}}$ + $\mathrm{D^R_2}$ },\;\;\;{\large \textcircled{\small 14}}\;\ce{$\mathrm{D^R_1}$ + $\mathrm{D^R_{12}}$ ->[$k_M + \bar k_M$] $\mathrm{D^R_{12}}$ + $\mathrm{D^R_{12}}$ },\\
\notag &{\large \textcircled{\small 15}}\;\ce{D + $\mathrm{D^R_2}$ ->[$k^{'}_M$] $\mathrm{D^R_1}$ + $\mathrm{D^R_2}$ },\;\;\;{\large \textcircled{\small 16}}\;\ce{D + $\mathrm{D^R_{12}}$ ->[$k^{'}_M$] $\mathrm{D^R_1}$ + $\mathrm{D^R_{12}}$ },\;\;\;{\large \textcircled{\small 17}}\;\ce{D + $\mathrm{D^R_1}$ ->[$\bar k_M$] $\mathrm{D^R_2}$ + $\mathrm{D^R_1}$ }\\
\notag &{\large \textcircled{\small 18}}\;\ce{$\mathrm{D^R_2}$ + $\mathrm{D^R_2}$ ->[$k^{'}_M$] $\mathrm{D^R_{12}}$ + $\mathrm{D^R_2}$ },\;\;\;{\large \textcircled{\small 19}}\;\ce{$\mathrm{D^R_2}$ + $\mathrm{D^R_{12}}$ ->[$k^{'}_M$] $\mathrm{D^R_{12}}$ + $\mathrm{D^R_{12}}$ },\\
\notag &{\large \textcircled{\small 20}}\;\ce{$\mathrm{D^R_1}$ + $\mathrm{D^R_1}$ ->[$\bar k_M$] $\mathrm{D^R_{12}}$ + $\mathrm{D^R_1}$ },\;\;\;
{\large \textcircled{\small 21}}\;\ce{$\mathrm{D^R_2}$ ->[$\delta+\bar k^R_{E}$] D },\;\;\;{\large \textcircled{\small 22}}\;\ce{$\mathrm{D^R_2}$ + $\mathrm{D^A}$ ->[$k^R_E$] D + $\mathrm{D^A}$ },\\
\notag &{\large \textcircled{\small 23}}\;\ce{$\mathrm{D^R_1}$ ->[$\delta^{'}+k^{'}_{T}$] D }\;\;\;{\large \textcircled{\small 27}}\;\ce{$\mathrm{D^R_1}$ + $\mathrm{D^A}$ ->[$k^{'*}_T$] D + $\mathrm{D^A}$ },\;\;\;
{\large \textcircled{\small 24}}\;\ce{$\mathrm{D^R_{12}}$ ->[$\delta^{'}+k^{'}_{T}$] $\mathrm{D^R_2}$ },\\
\notag &{\large \textcircled{\small 26}}\;\ce{$\mathrm{D^R_{12}}$ + $\mathrm{D^A}$ ->[$k^{'*}_T$] $\mathrm{D^R_2}$ + $\mathrm{D^A}$ },{\large \textcircled{\small 27}}\;\ce{$\mathrm{D^R_{12}}$ ->[$\delta+\bar k^R_{E}$] $\mathrm{D^R_1}$ },\;\;\;
{\large \textcircled{\small 28}}\;\ce{$\mathrm{D^R_{12}}$ + $\mathrm{D^A}$ ->[$k^R_E$] $\mathrm{D^R_1}$ + $\mathrm{D^A}$ },
\end{align}
   \endgroup
   \normalsize
in which the form of the reaction rate constants is due to the fact that reactions with the same reactants and products have been combined. As we did for the 3D model, let us define the parameter $\eps=\frac{\delta+\bar{k}^A_E}{k^A_M\frac{\mathrm{D_{tot}}}{V}}$, the parameter $\mu=\frac{k^R_E}{k^A_E}$, with a constant $b$ such that $\frac{\delta+\bar k^R_E}{\delta + \bar k^A_E} = b\mu$, and the parameter $\mu'=\frac{k^{'*}_T}{k^A_E}$, with a constant $\beta$ such that $\frac{\delta^{'}+k^{'}_T}{\delta + \bar k^A_E} = \beta \mu'$.
The Markov chain $X^{\eps}$ associated with the system is a linearly perturbed finite state continuous time Markov chain with the state $x$ representing the number of each type of modified nucleosome, i.e., $x=(n_{\mathrm{D^R_{12}}}, n_{\mathrm{D^A}}, n_{\mathrm{D^R_1}}, n_{\mathrm{D^R_2}})^T=(x_1,x_2,x_3,x_4)^T$. 
%
If the total number of nucleosomes is $\mathrm{D_{tot}}$, which is conserved, then the state space is $\mathcal{X}= \{(x_1,x_2,x_3,x_4)^T \in \Z_+^4:\: x_1 + x_2 + x_3 + x_4 \leq \Dtot \}$. The transition vectors for $X^{\eps}$ are given by $v_1=(0,1,0,0)^T$, $v_2=(0,-1,0,0)^T$, $v_3=(0,0,1,0)^T$, $v_4=(0,0,-1,0)^T$, $v_5=(0,0,0,1)^T$, $v_6=(0,0,0,-1)^T$, $v_7=(1,0,-1,0)^T$, $v_8=(-1,0,1,0)^T$, $v_9=(1,0,0,-1)^T$ and $v_{10}=(-1,0,0,1)^T$. The infinitesimal transition rates are  

  %
\begingroup
\small 
\begin{align}
\notag &Q_{x,x+v_1}(\eps)=f_A(x) = (\Dtot -(x_1+x_2+x_3+x_4))\left(k_{W0}^A+k_{W}^A + \frac{k_{M}^A}{V}x_2\right),\\
\notag &Q_{x,x+v_2}(\eps)=g^{\eps}_A(x) = x_2\left(\eps \frac{k_{M}^A}{V}\Dtot + \frac{k^A_E}{V}(x_3+x_4+2x_1)\right),\\
\notag &Q_{x,x+v_3}(\eps)=f_{R1}(x) = (\Dtot -(x_1+x_2+x_3+x_4))\left(k^1_{W0}+k^1_{W} + \frac{k^{'}_{M}}{V}(x_1+x_4)\right),\\
\notag &Q_{x,x+v_4}(\eps)=g^{\eps}_{R1}(x) = x_3\mu'\left(\eps \frac{k_{M}^A}{V}\Dtot \beta + x_2\frac{k^A_E}{V}\right),\\
\notag &Q_{x,x+v_5}(\eps)=f_{R2}(x) = (\Dtot -(x_1+x_2+x_3+x_4))\left(k^2_{W0}+k^2_{W} +  \frac{k_{M}}{V}(x_1+x_4)  +  \frac{\bar k_{M}}{V}(x_1+x_3)\right),\\
\label{rates4D} &Q_{x,x+v_6}(\eps)=g^{\eps}_{R2}(x) = x_4\mu\left(\eps \frac{k_{M}^A}{V}\Dtot b + x_2\frac{k^A_E}{V}\right),\\
\notag &Q_{x,x+v_7}(\eps)=f_{R121}(x) = x_3\left(k^2_{W0}+  \frac{k_{M}}{V}(x_1+x_4)  +  \frac{\bar k_{M}}{V}\left(x_1+\frac{x_3-1}{2}\right)\right),\\
\notag &Q_{x,x+v_8}(\eps)=g^{\eps}_{R121}(x) = x_1\mu\left(\eps \frac{k_{M}^A}{V}\Dtot b + x_2\frac{k^A_E}{V}\right), \\
\notag &Q_{x,x+v_9}(\eps)=f_{R122}(x) = x_4\left(k^1_{W0} + \frac{k^{'}_{M}}{V}\left(x_1+\frac{x_4-1}{2}\right)\right),\\
\notag &Q_{x,x+v_{10}}(\eps)=g^{\eps}_{R122}(x) = x_1\mu'\left(\eps \frac{k_{M}^A}{V}\Dtot \beta + x_2\frac{k^A_E}{V}\right).
\end{align}
   \endgroup
   \normalsize
\noindent A representation of the transition vectors and the Markov chain graph for $\Dtot=2$ are given in Fig. \ref{fig:4Dmodel1} (b) and (c), respectively. As before, each rate depends on the state $x$, but in the rest of the section we will not show this dependency to simplify the notation. Now, we determine the stochastic behavior of the full chromatin modification circuit model in terms of its stationary distribution and time to memory loss. For this study, we will consider $k^A_W=k^1_W=k^2_W=0$ (i.e., there are no external transcription factors enhancing the establishment of chromatin modifications). This assumption will not change the qualitative nature of the results focused on studying the effect of $\eps$, $\mu$, and $\mu'$ on the stochastic behavior of the chromatin modification circuit model.

\subsubsection{Stationary distribution}
\label{SD4D}

    
    


We now determine the zeroth and first order terms of the stationary distribution expansion for the 4D model. As shown in SI - Section \ref{sec:Appendix4Dmodel}, when $\eps=0$, the continuous time Markov chain associated with the 4D model has transient states $\T= \{i_1, \ldots, i_m\}$ where $m =\sum_{j=0}^{\Dtot}\sum_{k=0}^{j}\left(\frac{(k+2)(k+1)}{2}\right)-2$, $i_1=(0,\Dtot-1,0,0)^T$, $i_{m-1}=(\Dtot-1,0,0,1)^T$, $i_m=(\Dtot-1,0,1,0)^T$, and absorbing states $\A = \{a,r\}$, with $a=(0,\Dtot,0,0)^T$ corresponding to the fully active state ($n_{\DA}=\Dtot$) and $r=(\Dtot,0,0,0)^T$ corresponding to the fully repressed state ($n_{\mathrm{D^R_{12}}}=\Dtot$), respectively. Then, Assumption \ref{assumption:transient_absorbing} holds (see SI - Section \ref{sec:Appendix4Dmodel}), so that 
%
%
$\X = \A \cup \T$, and we can rewrite the infinitesimal generator $Q(\eps)$ in the form of \eqref{eq:Q_eps_decomposition} (see SI - Section \ref{sec:Appendix4Dmodel}), where the perturbation is linear and so
%
%
Assumption \ref{assumption:linear_case} holds. We can also verify Assumption \ref{assumption:SingleRecurrentClass} (see SI - Section \ref{sec:Appendix4Dmodel}). Hence, we can apply SI - Theorem \ref{thm:HigherOrderTerms_for_stationary distribution}, as was done for the 3D model.
%

In particular, we obtain $\pi(0)=\pi^{(0)}=[\alpha,0] = \left[\alpha_a,\alpha_r, 0 \ldots ,0 \right]$, with $\alpha_a=0$, $\alpha_r=1$, and
%
%
$\beta^{(1)}=[\pi^{(1)}_{i_1},...,\pi^{(1)}_{i_m}]=[0,...,0,\pi^{(1)}_{i_{m-1}},\pi^{(1)}_{i_m}]$, with 
\begin{equation}\nonumber
\pi^{(1)}_{i_{m-1}}=\frac{\mu' \beta\frac{k_{M}^A}{V}\Dtot^2}{k^1_{W0}+ \frac{k'_{M}}{V}(\Dtot - 1)},\;\;\;\;\;\; \pi^{(1)}_{i_m}=\frac{\mu b\frac{k_{M}^A}{V}\Dtot^2}{k^2_{W0}+ (\frac{k_{M}}{V} + \frac{\bar k_{M}}{V})(\Dtot - 1)}.
\end{equation}
    See SI - Section \ref{sec:Appendix4Dmodel} for the detailed mathematical derivation. Now, $\alpha^{(1)}=[\pi^{(1)}_{a},\pi^{(1)}_{r}]$ is the unique vector such that $\alpha^{(1)}Q_{\A} = -\beta^{(1)}[R_1 + T_1(-T_0)^{-1}R_0]$, $\alpha^{(1)}\one = -\beta^{(1)}\one$.
    
    As an example, suppose $\Dtot=2$ and assume $\beta=b$, $k^1_{W0}=k^2_{W0}=k^A_{W0}$ and $k'_M=\bar k_M = k_M = k^A_M$. These assumptions do not affect the final qualitative conclusions related to the effect of $\eps, \mu$ and $\mu'$ on the stationary distribution. Then (see SI - Section \ref{sec:Appendix4Dmodel} for the detailed mathematical derivation)
    	\begin{equation}\label{QA4D}
    	    Q_\A = \frac{\bar K_1(\mu,\mu')}{\bar K_2(\mu,\mu')} \left( \begin{array}{c c}
	    -1 & 1 \\
	    0 & 0 \\
    \end{array} \right),
    \end{equation}
    \begin{equation}
    \begin{aligned}\label{pi4D}
        \pi^{(1)}_{a}&=\frac{\bar K_3(\mu,\mu')}{\bar K_4(\mu,\mu')},\;\; \pi^{(1)}_{r}&=-\pi^{(1)}_{a}-\pi^{(1)}_{i_{m-1}}-\pi^{(1)}_{i_{m}}=-\frac{\bar K_3(\mu,\mu')}{\bar K_4(\mu,\mu')}-\mu' K_{18}-\mu K_{19},
     \end{aligned}   
        \end{equation}
        \noindent with
        \begin{equation}\label{constant4D}
        \begin{aligned}
\bar K_1(\mu,\mu') &=K_1((\mu')^2K_2+(\mu)^2K_3+\mu\mu'K_4+\mu'K_5+\mu K_6+K_{7}),\\
\bar K_2(\mu,\mu') &=\mu'\mu(\mu'+\mu) K_8+(\mu')^2K_{9}+ (\mu)^2K_{10}+\mu\mu'K_{11}+\mu'K_{12}+\mu K_{13}+K_{14},\\
\bar K_3(\mu,\mu')&=(\mu\mu')^2K_{15}((\mu+\mu')K_{16}+K_{17}),\\
\bar K_4(\mu,\mu') &=K_{20}((\mu')^2K_2+(\mu)^2K_3+\mu\mu'K_4+\mu'K_5+\mu K_6+K_{7}),
     \end{aligned}
     \end{equation}
   in which $m=13$ and $K_i$, $i=1,...,20$, are non-negative functions independent of $\mu$ and $\mu'$ (see SI - Section \ref{sec:Appendix4Dmodel} for their precise definitions). We then have
   
        \begin{equation}\nonumber
        \label{pi4DmodelRESULTS}
        \pi_x(\eps) = \begin{cases}
        \eps\frac{\bar K_3(\mu,\mu')}{\bar K_4(\mu,\mu')} + O(\eps^2) & \text{ if } x= a =(0,2,0,0)^T \\
        O(\eps^2)& \text{ if } x \in \T \backslash \{i_{m-1},i_m\} \\    
        \eps \mu'K_{18} + O(\eps^2)& \text{ if } x = i_{m-1} = (1,0,0,1)^T \\    
        \eps \mu K_{19} + O(\eps^2)& \text{ if } x = i_{m} = (1,0,1,0)^T \\
        1-\eps \left(\frac{\bar K_3(\mu,\mu')}{\bar K_4(\mu,\mu')}+\mu' K_{18}+\mu K_{19}\right) + O(\eps^2)& \text{ if } x = r= (2,0,0,0)^T.
        \end{cases}
        \end{equation}
        \normalsize
    For small $\eps>0$, the stationary distribution is concentrated around the active and repressed states, and higher values of $\mu'$ or $\mu$ shift the distribution towards the active state.

   \subsubsection{Time to memory loss}\label{4DTTML}
   
As was done for the 3D model, we determine for the 4D model how the parameters $\eps$ and $\mu'$ affect the time to memory loss of the active state, $h_{a,r} (\eps)$, and the time to memory loss of the repressed state, $h_{r,a} (\eps)$. Firstly, by the algorithm in Section \ref{orderMFPT}, $h_{a,r} (\eps)$ is $O(\eps^{-1})$ and $h_{r,a} (\eps)$ is $O(\eps^{-2})$ (see SI - Section \ref{sec:alg_3D4D}).
%
    %
    Then, by applying Theorem \ref{thm:MFPTcoefficient} we can obtain expressions for the leading coefficients of $h_{a,r} (\eps)$ and $h_{r,a} (\eps)$ (See SI - Section \ref{sec:Appendix4Dmodel} for the detailed mathematical derivation). As an example, when $\Dtot=2$, $Q_\A$ and $\pi_a^{(1)}$ are shown in \eqref{QA4D} and \eqref{pi4D}, and we obtain that 
    \begin{equation*}
    h_{a,r}(\eps)=\frac{\bar K_2(\mu,\mu')}{\bar K_1(\mu,\mu')} \frac{1}{\eps} + O(1),\;\;\text{and}\;\;h_{r,a}(\eps)=\frac{\bar K_2(\mu,\mu')K_4(\mu,\mu')}{\bar K_1(\mu,\mu')\bar K_3(\mu,\mu')} \frac{1}{\eps^2} + O\left(\frac{1}{\eps}\right),
    \end{equation*}
    where $\bar K_i(\mu,\mu')$, $i=1,...,4$, are defined in \eqref{constant4D}.

Now, we determine how $\mu'$, the parameter encapsulating the asymmetry between the DNA methylation erasure rate and the activating histone modification erasure rate, affects the time to memory loss.
To this end, we seek to determine directly how $\mu'$ affects $h_{a,r}(\eps)$ and $h_{r,a}(\eps)$, without deriving an explicit expression for them. To this end, we would like to exploit two theorems from \cite{Monotonicitypaper}, namely, Theorem S.2 and Theorem 3.4 there. The transitions of the Markov chain $X^{\eps}$ are in ten possible directions, $v_1=-v_2=(1,0,-1,0)^T$, $v_3=-v_4=(1,0,0,-1)^T$, $v_5=-v_6=(0,1,0,0)^T$, $v_7=-v_8=(0,0,1,0)^T$, and $v_9= -v_{10}=(0,0,0,1)^T$, with the associated infinitesimal transition rates $\rate_1(x)=f_{R121}(x)$, $\rate_2(x)=g^{\eps}_{R121}(x)$, $\rate_3(x)=f_{R122}(x)$, $\rate_4(x)=g^{\eps}_{R122}(x)$, $\rate_5(x)=f_A(x)$, $\rate_6(x)=g^{\eps}_A(x)$, $\rate_7(x)=f_{R1}(x)$, $\rate_8(x)=g^{\eps}_{R1}(x)$, $\rate_9(x)=f_{R2}(x)$, $\rate_{10}(x)=g^{\eps}_{R2}(x)$.	       
   %
Consider infinitesimal transition rates $\Breve{\rate}_{i}(x)$, $i=1, 2,..., 10$, defined as for $\rate_{i}(x)$, $i=1, 2,..., 10$,  with all the parameters having the same values except that $\mu'$ is replaced by $\Breve{\mu}'$, with $\mu' \ge \Breve{\mu}'$. 
While we have not been able to see how to exploit Theorems S.2 and 3.4 from \cite{Monotonicitypaper} for these exact rates, we have been able to do this for closely related rates. If we introduce a small approximation in the transition rates of $X^{\eps}$, namely, $\frac{x_3-1}{2}\approx x_3 $ and $ \frac{x_4-1}{2}\approx x_4$ in $f_{R121}(x)$ and $f_{R122}(x)$, respectively, then Theorems S.2 and 3.4 in \cite{Monotonicitypaper} apply with 
 \begin{equation}   \label{matrixA4D}
       A= \begin{bmatrix}
        0 & -1 & 0 & 0\\
        1 & 0 & 1 & 0\\
        1 & 0 & 0 & 1\\
        1 & 0 & 1 & 1
        \end{bmatrix},\nonumber
 \end{equation}   
     and $(K_A+x)\cap \X=\{ w \in \X:\: x \preccurlyeq_A w \}$ (see SI - Section \ref{sec:Appendix4Dmodel}). 
This approximation can be justified by introducing the reasonable assumption that each nucleosome characterized by a repressive modification ($\mathrm{D^R_1}$ and $\mathrm{D^R_2}$) has the ability to catalyze the establishment of the opposite repressive mark on itself.
With this approximation,  since $a=(0,\mathrm{D_{tot}},0,0)^T \preccurlyeq_A r=(\mathrm{D_{tot}},0,0,0)^T$, then $\Breve{h}_{a,r}(\eps)\le h_{a,r}(\eps)$ and $h_{r,a}(\eps)\le \Breve{h}_{r,a}(\eps)$. 
Thus, the time to memory loss of the active state increases with higher values of $\mu'$, while the time to memory loss of the repressed state decreases with higher values of $\mu'$.

%% file: Conclusion.tex
\section{Conclusion}
\label{sec:conclusion}
	

In this paper, we provided a mathematical formulation and rigorous proofs to validate the computational findings in \cite{BrunoDelVecchio}, showing how the time scale separation between establishment and erasure processes of chromatin modifications affects epigenetic cell memory. To this end,  we developed and adapted theory for singularly perturbed continuous time Markov chains and
%
we analyzed the behavior of stationary distributions and mean first passage times as functions of the singular perturbation parameter $\eps$. 

We first showed that $\pi(\eps)$ can be expressed as a series expansion (Section \ref{sec:pertCTMC}) for sufficiently small $\eps$. We then proved that the limit $\pi(0) = \lim_{\eps \to 0} \pi(\eps)$ is unique and we determined an expression for it (Section \ref{sec:ZerothTerm}). We also provided an iterative procedure for computing all of the higher order terms in the expansion of $\pi(\eps)$ (SI -Section \ref{sec:StatDistHigherOrder}). 
Similarly, for the mean first passage time (MFPT) between states, we first showed there is a Laurent series expansion for sufficiently small $\eps$ (Section \ref{sec:pertCTMC}, Eq. \eqref{eq:ExpansionTauEps}). We then developed a graph based algorithm to identify the order of the leading term in the series expansion (Section \ref{orderMFPT}), and we also determined the leading coefficient there (Section \ref{sec:leadingMFPT}). 

We then applied these tools to the chromatin modification circuit models proposed in \cite{BrunoDelVecchio}, to provide a rigorous basis for the computational findings given there (Sections \ref{motivexamples}, \ref{sec:MainResults}, and \ref{sec:Applications}).

Our rigorous derivations of the analytical expressions for the stationary distributions and time to memory loss, and our results on monotonic dependence on parameters, lead to a mechanistic understanding of how $\eps$, $\mu$ and $\mu'$ affect the stochastic behavior of chromatin modification systems.
As an example, our results suggest that higher values of $\mu$ and $\mu'$ shift mass of the stationary distribution more towards the active state (Sections \ref{SD3D} and \ref{SD4D}). This finding is consistent with recent experimental results demonstrating that transfection of the DNA methylation eraser enzyme TET1 (represented in our model by higher $\mu'$ \cite{BrunoDelVecchio}) into Chinese hamster ovary (CHO-K1) repressed cells causes them to shift towards the active state \cite{PalaciosBruno2024}.
More generally, the mechanistic understanding of how $\eps$, $\mu$, and $\mu'$ affect the stochastic behavior of chromatin modification systems, as derived in our study, is crucial for determining experimental interventions on molecular players, such as chromatin modifier enzymes, to modulate cell memory. This mechanistic insight is expected to be extremely valuable for applications such as cell fate reprogramming and engineering approaches to cell therapy.
Furthermore, the mathematical results and theoretical tools developed in this paper can be applied beyond the scope of the epigenetic cell memory models analyzed in this research work. In fact, they can be applied to all stochastic models that respect the assumptions considered. Future work will investigate how to generalize these results by removing some of these assumptions, including allowing the Markov chain to have countably many states and $Q(0)$ to have ergodic classes as well as absorbing states. While there is some theory for countably many states, such as in \cite{AAN2004}, the continuous time Markov chains for further applications that we have in mind are not uniformizable and have many transient states for $Q(0)$, and the theory in \cite{AAN2004} needs to be generalized for them.\\

%% file: Appendix_title.tex
\begin{center}
{\LARGE Analysis of singularly perturbed stochastic chemical reaction networks motivated by applications to epigenetic cell memory}\\
\end{center}

\begin{center}
{\large Simone Bruno$^{1,*}$, Felipe A. Campos$^{2,*}$, Yi Fu$^{2}$, Domitilla Del Vecchio$^{1}$, and Ruth J. Williams$^{2}$}\\
\end{center}

\begin{center}
$^1$\textit{\small Department of Mechanical Engineering, Massachusetts Institute of Technology, 77 Massachusetts Avenue, Cambridge, MA 02139. Emails: {\tt\small (sbruno,ddv)@mit.edu}}

$^{2}$\textit{\small Department of Mathematics, University of California, San Diego, 9500 Gilman Drive, La Jolla CA 92093-0112. Email: {\tt\small (fcamposv,yif064,rjwilliams)@ucsd.edu}}

$^{*}$\textit{\small These authors contributed equally: S. Bruno and F. A. Campos}
\end{center}

\section*{Supplementary Information (SI)}

%% file: Appendix_Probability.tex
\subsection{Some results in probability}

    Let $\X$ be a finite set. Recall the notation for matrices introduced in Section \ref{sec:PreliminariesAndNotation}.
    
    \begin{lemma}
    \label{lem:nullityQmatrix}
     Let $X = \{X(t): \: t \geq 0 \}$ be a continuous time Markov chain with state space $\X$ and infinitesimal generator $Q = (Q_{x,y})_{x,y \in \X}$. Then, the number of recurrent classes for $X$ is equal to $\nullity(Q^T)=\nullity(Q)$.
    \end{lemma}
    
    \begin{proof}
    Since $Q$ is a square matrix, the Rank plus Nullity Theorem yields that $\nullity(Q^T)$ $= \nullity(Q)$. Now, consider $\lambda > \max_{x \in \X}|Q_{x,x}|$ and define $P := I + Q/\lambda$, where $I$ is the identity matrix of size $|\X| \times |\X|$. The matrix $P$ is stochastic and such that for every $x \neq y$ in $\X$, $P_{x,y} > 0$ if and only if $Q_{x,y} >0$. As a consequence, the recurrent classes of $X$ are the same as the recurrent classes of $P$. By Theorem IV.2.4 in Isaacson and Madsen \cite{IsaacsonMadsen}, the number of recurrent classes of $P$ is equal to the maximum number of linearly independent left eigenvectors satisfying $\pi P = \pi$. By observing that $\pi P = \pi$ if and only if $\pi Q = 0$, we see that this latter quantity is equal to $\nullity(Q^T)$.
    \end{proof}
   
    The following is Proposition 6.3 in Asmussen \cite{AsmussenQueues2nd}. 
    
    \begin{proposition}
    \label{prop:SpectralRadiusSubstochastic}
    Let $(P_{x,y})_{x,y \in \X}$ be a nonnegative substochastic matrix ($P \one \leq \one$) such that for each $x \in \X$ there are $z_1,\ldots,z_m,y \in \X$ such that $P_{x,z_1}P_{z_1,z_2}\ldots P_{z_m,y} >0$ and $\sum_{z \in \X} P_{y,z} < 1$. Then, $\spr(P) < 1$.
    \end{proposition}
    
    We use Proposition \ref{prop:SpectralRadiusSubstochastic} in order to obtain invertibility for some matrices, as in the next result.
    
    \begin{lemma}
	\label{lem:IMinusPInvertible}
    Let $X = \{X(t) :\: t \geq 0 \}$ be an irreducible continuous time Markov chain with state space $\X$ and with an embedded discrete time Markov chain with transition matrix $P$.  Consider a nonempty set $\B \subseteq \X$ such that $\B \neq \X$ and consider $P^{\Bc}$ to be the matrix obtained by removing the columns and rows of $P$ corresponding to states in $\B$. Then, $I-P^{\Bc}$ is invertible and its inverse is given by the absolutely convergent series $\sum_{k=0}^{\infty} (P^{\Bc})^k$, where $(P^{\Bc})^0 = I$.
    \end{lemma}
    
    \begin{proof}
    Observe that $P^{\Bc}=(P_{x,y})_{x,y \in \Bc}$ is a nonnegative substochastic matrix. Since $X$ is an irreducible continuous time Markov chain, its embedded discrete time Markov chain is also irreducible. Thus, for each $x \in \B^c$, there exist $z_1,\ldots,z_m,y \in \B^c$ and $\tilde{y} \in \B$ such that $P_{x,z_1} P_{z_1,z_2} \ldots P_{z_m,y} P_{y,\tilde{y}} > 0$. Then, $P^{\Bc}_{x,z_1} P^{\Bc}_{z_1,z_2}\ldots P^{\Bc}_{z_m,y} > 0$ and $\sum_{z \in \B^c} P^{\Bc}_{y,z} = \sum_{z \in \B^c} P_{y,z} < 1$ since $P_{y,\tilde{y}} > 0$ and $\sum_{z \in \X} P_{y,z}=1$. By Proposition \ref{prop:SpectralRadiusSubstochastic}, $\spr(P^{\Bc}) < 1$. This fact, together with Theorem 5.6.15 in Horn \& Johnson \cite{HornJohnson} yields the convergence of $\sum_{k=0}^{\infty} (P^{\Bc})^k$. Moreover, $(I-P^{\Bc}) \sum_{k=0}^{\infty} (P^{\Bc})^k = \sum_{k=0}^{\infty} (P^{\Bc})^k (I-P^{\Bc}) = I$, which yields the desired result.
    \end{proof}

    We will use the following continuous time analogue of Proposition \ref{prop:SpectralRadiusSubstochastic}.
    
    \begin{lemma}
    \label{lem:Q_invertible}
    Let $(Q_{x,y})_{x,y \in \X}$ be a matrix such that $Q\one \leq 0$ and such that $Q_{x,x} \leq 0$ for each $x \in \X$ and $Q_{x,y} \geq 0$ for each $x \neq y \in \X$. In addition, suppose that for each $x \in \X$ there are distinct $z_1,\ldots,z_m,y \in \X$ different from $x$ such that $Q_{x,z_1}Q_{z_1,z_2}\ldots Q_{z_m,y} > 0$ and $\sum_{z \in \X} Q_{y,z} < 0$. Then, for every $\upsilon \in \spec(Q)$, the real part of $\upsilon$ is negative. In particular, $Q$ is invertible.
    \end{lemma}
    
    \begin{proof}
    Consider $\lambda > \max_{x \in \X}|Q_{x,x}|$ and define $P := I + Q/\lambda$, where $I$ is the identity matrix of size $|\X| \times |\X|$. The matrix $P$ is nonnegative, substochastic and such that $P_{x,y} = \frac{1}{\lambda}Q_{x,y}$ for every $x \neq y \in \X$. With these elements, we obtain that for each $x \in \X$ there are distinct $z_1,\ldots,z_m,y \in \X$ such that $P_{x,z_1}P_{z_1,z_2}\ldots P_{z_m,y} > 0$ and where $\sum_{z \in \X} P_{y,z} = 1  + \frac{1}{\lambda}\sum_{z \in \X} Q_{y,z}  < 1$. Proposition \ref{prop:SpectralRadiusSubstochastic} yields that $\spr(P) < 1$. By observing that $\upsilon \in \spec(Q)$ implies that $1 + \frac{\upsilon}{\lambda} \in \spec(P)$, we obtain that $1 > |1 + \frac{\upsilon}{\lambda}| \geq |1 + \frac{\Re(\upsilon)}{\lambda}|$ where $\Re(\upsilon)$ is the real part of $\upsilon$. The latter inequality implies that $\Re(\upsilon) < 0$.
    \end{proof}
    
    Consider a nonempty set $\B \subseteq \X$ such that $\B \neq \X$ and a $Q$-matrix written as
    \begin{equation}
        Q= \begin{blockarray}{ccc}
        & \B & \Bc \\
        \begin{block}{c(c|c)}
         \B & Q^{\nBc} & S \\
         \cline{2-3}
          \B^c & R & Q^{\Bc} \\
        \end{block}
        \end{blockarray}
        \; .
    \end{equation} 
    Consider a process $X = \{X(t) :\: t \geq 0 \}$ defined on a measurable space $(\Omega, \F)$ and a collection of probability measures $\{\PP_x :\: x \in \X\}$ on $(\Omega,\F)$ such that for every $x \in \X$, $X$ is a continuous time Markov chain under $\PP_x$ with infinitesimal generator given by $Q$ and such that $\PP_x[X(0)=x]=1$. Consider the stopping time $\tau_{\B} := \inf\{ t \geq 0 :\: X(t) \in \B \}$.
    
    \begin{lemma}
    \label{lem:InvertibilityQ^B}
    The following are equivalent:
        \begin{enumerate}
            \item[(i)]
            For every $x \in \B^c$, there exists a $z \in \B$ such that $x$ leads to $z$ under $Q$, i.e., there are distinct $x_1,...,x_m \in \B^c$, different from $x$, such that $Q_{x,x_1}$, $Q_{x_1,x_2},...,Q_{x_m,z}>0$.
            \item[(ii)]
            $Q^{\Bc}$ is invertible.
            \item[(iii)]
            $\E_x[\tau_{\B}] < \infty$ for every $x \in \B^c$.
        \end{enumerate}
    If any of $(i)-(iii)$ hold, then
        \begin{equation}
        \label{eq:MatrixRepForProbabilities}
            \PP_x\left[X(\tau_{\B})=y\right] =  (-(Q^{\Bc})^{-1}R)_{x,y}
        \end{equation}
    for every $x \in \B^c$ and $y \in \B$. Moreover, if $\B^c$ is a set of transient states for $X$, then $(i)-(iii)$ hold.
    \end{lemma}
    
    Part of the results in Lemma  \ref{lem:InvertibilityQ^B} appear as Lemma 1 in Gaver et al. \cite{Latouche} for the case where $Q^{\B}$ and $S$ are the zero matrix. For completeness, we provide a self-contained proof here.
    
    \begin{proof}
    The implication $(i) \Rightarrow (ii)$ is a consequence of Lemma \ref{lem:Q_invertible} with $\B^c,Q^{\B^c}$ in place of $\X,Q$ there.
    
    In the following, recall that for every $x \in \X$ and function $f:\X \longrightarrow \R$, the process
    \begin{equation}
    \label{eq:MartingaleMarkovChain}
        M^{\B}_f(t) := f(X(t \wedge \tau_{\B}))-f(X(0))-\int_0^{t \wedge \tau_{\B}} \L f(X(s))ds, \qquad t \geq 0,
    \end{equation}
    is a martingale under $\PP_x$ (see Theorem 3.32 in \cite{Liggett}), where $\L f(y) := \sum_{z \in \X} Q_{yz}f(z)$ for $y \in \X$.
    
    For $(ii) \Rightarrow (iii)$, consider the function $f(y) := -[(Q^{\Bc})^{-1}\one]_y\one_{\B^c}(y)$ for $y \in \X$. The reader may verify that $\L f(y) = -1$ for every $y \in \B^c$. Therefore, for an $x \in \B^c$, taking expectations in \eqref{eq:MartingaleMarkovChain} yields $f(x) - \E_x[f(X(t \wedge \tau_{\B}))] = \E_x[t \wedge \tau_{\B}]$ for every $t \geq 0$. Hence, $\E_x[t \wedge \tau_{\B}] \leq 2\sup_{x \in \X}|f(x)|$ for every $t \geq 0$ and we conclude the desired result by letting $t \to \infty$.
    
    For $(iii) \Rightarrow (i)$, we prove that not $(i)$ implies not $(iii)$. Suppose there exists $x \in \B^c$ such that no point of $\B$ is accessible from $x$. Then, $\tau_{\B} = \infty$ $\PP_x$-a.s., so $(iii)$ does not hold.

    For \eqref{eq:MatrixRepForProbabilities}, consider $x \in \B^c$, $y \in \B$ and the martingale $M_f^{\B}$ with $f(z) = (-(Q^{\Bc})^{-1}R)_{z,y} \one_{\B^c}(z)$ $ + \one_{\{y\}}(z)$ for $z \in \X$. The reader may verify that $\L f(x) = 0$ for $x \in \B^c$, which yields that $\E_x[f(X(t \wedge \tau_{\B}))] = f(x)$ for every $t \geq 0$. If any of $(i)-(iii)$ hold, then $\tau_{\B} < \infty$, $\PP_x$-a.s., and on letting $t \to \infty$ and using bounded convergence, we obtain $\E_x[f(X(\tau_{\B}))] = f(x)$, which implies $\PP_x[X(\tau_{\B}) = y] = (-(Q^{\Bc})^{-1}R)_{x,y}$.

    Now, suppose that every $x \in \B^c$ is transient. Then, $\PP_x(\tau_{\B}<\infty)=1$ for each $x \in \B^c$. For a proof by contradiction, suppose that $Q^{\Bc}$ is not invertible, which implies the existence of a nonzero vector $v = (v(x))_{x \in \B^c} \neq 0$ such that $Q^{\Bc}v = 0$. Then, consider the function $f(y) = v(y)\one_{\B^c}(y)$, for which $\L f(y) = 0$ for $y \in \B^c$. Consider an $x \in \B^c$ such that $v(x) \neq 0$, then $M^{\B}_f(t) = f(X(t \wedge \tau_{\B})) - v(x)$ is a bounded martingale. On taking expectations we have $\E_x[f(X(t \wedge \tau_{\B}))] = v(x)$. Since the states in $\B^c$ are transient, $X(\cdot)$ will $\PP_x$-a.s. leave $\B^c$. Then letting $t \to \infty$ and using bounded convergence yields $0 = v(x)$ which is the desired contradiction. Hence $(ii)$ (as well as $(i),(iii)$) must hold.
    \end{proof}
    
    Lemma \ref{lem:InvertibilityQ^B} has a useful consequence in terms of occupations times. In the above context, consider the occupation time of $\B$ by the Markov chain $X$ up to time $t \geq 0$: $\chi_{\B}(t) = \int_{0}^{t} \one_{\B}(X(s))ds$. Denote by $\chi_{\B}(\infty) = \lim_{t \to \infty} \chi_{\B}(t) = \int_0^{\infty} \one_{\B}(X(s))ds$.
    
    \begin{lemma}
    \label{lem:OccupationTimesAreInfinity}
    Suppose that
    \begin{equation}
    \label{eqn:finiteHittingTimeToB}
    \PP_y [\tau_\B < \infty]=1 \text{ for all } y \in \B^c.
    \end{equation}
    Then $\PP_x[\chi_{\B}(\infty)=\infty] = 1$ for every $x \in \X$.
    \end{lemma}

    \begin{remark}
    If any of the conditions (i)-(iii) in Lemma \ref{lem:InvertibilityQ^B} holds, then  \eqref{eqn:finiteHittingTimeToB} holds.
    \end{remark}
    
    \begin{proof}
        Fix $x \in \X$. Let $\sigma_{-1}=0$ and $\sigma_0=\inf\{t \geq \sigma_{-1}: X(t) \in \B \}$. Then, inductively define for $k=0,1,2,\dots$, $\sigma_{2k+1} = \inf\{t \geq \sigma_{2k}: X(t) \in \B^c \}$ and $\sigma_{2k+2} = \inf\{t \geq \sigma_{2k+1}: X(t) \in \B \}$. Using \eqref{eqn:finiteHittingTimeToB} and the strong Markov property of $X$, we have
            \begin{equation}
            \label{eq:SigmaEvenIsFinite}
            \sigma_{2k} < \infty \: \PP_x\text{-a.s. on } \{\sigma_{2k-1} < \infty\}    
            \end{equation}
        for $k=0,1,2,\dots$, and
            \begin{equation*}
            \chi_{\B}(\infty) = \sum_{k=0}^\infty \one_{\{\sigma_{2k} < \infty\}} (\sigma_{2k+1} - \sigma_{2k}) = \lim_{N \rightarrow \infty} \sum_{k=0}^N \one_{\{\sigma_{2k} < \infty\}} (\sigma_{2k+1} - \sigma_{2k}),  
            \end{equation*}
        where terms in the sum indexed by $k:\sigma_{2k} = \infty$ are taken to be zero.
        Now, $\PP_x$-a.s.,
            \begin{equation*}
            \prod_{k=0}^N \exp(-\one_{\{\sigma_{2k} < \infty\}}(\sigma_{2k+1} - \sigma_{2k}))  = \prod_{k=0}^N \one_{\{\sigma_{2k} < \infty\}} \exp (-(\sigma_{2k+1} - \sigma_{2k})),
            \end{equation*}
       where we used \eqref{eq:SigmaEvenIsFinite} and the fact that $\sigma_{-1}=0$, to conclude that the product is zero $\PP_x$-a.s., if $\{\sigma_{2k} = \infty\}$ for any $k \in \{0,1,\ldots,N\}$. Hence,
        \small
        \begin{eqnarray*}
        && \E_x[\exp(-\chi_{\B}(\infty))]\\
        &=& \lim_{N \rightarrow \infty} \E_x\left[\prod_{k=0}^N \one_{\{\sigma_{2k} < \infty\}} \exp (-(\sigma_{2k+1} - \sigma_{2k})) \right]\\
        &=& \lim_{N \rightarrow \infty} \E_x\left[\prod_{k=0}^{N-1} \one_{\{\sigma_{2k} < \infty\}} \exp (-(\sigma_{2k+1} - \sigma_{2k})) \one_{\{\sigma_{2N} < \infty\}} \E_x \left[ \exp (-(\sigma_{2N+1} - \sigma_{2N})) | X(\sigma_{2N})\right] \right].
        \end{eqnarray*}
        \normalsize
        
        On $\{\sigma_{2N} < \infty\}$, we have $X(\sigma_{2N}) \in \B$ and $X(t) \in \B$ for $\sigma_{2N} \leq t < \sigma_{2N+1}$. Hence, for $a > \max_{y \in \B} |Q_{y,y}|$, using the strong Markov property, on $\{\sigma_{2N} < \infty\}$, conditioned on $X(\sigma_{2N})$, $\sigma_{2N+1} - \sigma_{2N}$ stochastically dominates an exponential random variable with parameter $a$ and so 
        \begin{equation*}
            \E_x\left[ \exp (-(\sigma_{2N+1} - \sigma_{2N})) | X(\sigma_{2N})\right] \leq \int_0^\infty e^{-t} a e^{-at}dt = \frac{a}{1+a}.
        \end{equation*}
        Similarly, for $k=N-1,\dots,0$, on $\{\sigma_{2k} < \infty\}$, 
        \begin{equation}
        \label{eq:BoundOnConditionalExpectation}
        \E_x\left[ \exp (-(\sigma_{2k+1} - \sigma_{2k})) | X(\sigma_{2k})\right] \leq \frac{a}{1+a}.   
        \end{equation}
        Then,
        \begin{eqnarray*}
        \E_x[\exp(-\chi_{\B}(\infty))] &\leq& \limsup_{N \rightarrow \infty} \E_x\left[\prod_{k=0}^{N-1} \one_{\{\sigma_{2k} < \infty\}} \exp (-(\sigma_{2k+1} - \sigma_{2k})) \one_{\{\sigma_{2N} < \infty\}} \frac{a}{1+a} \right]\\
        &\leq& \limsup_{N \rightarrow \infty} \E_x\left[\prod_{k=0}^{N-1} \one_{\{\sigma_{2k} < \infty\}} \exp (-(\sigma_{2k+1} - \sigma_{2k})) \frac{a}{1+a} \right].
        \end{eqnarray*}
    Repeatedly conditioning on $\{\sigma_{2k} < \infty\}$, for $k=N-1,\ldots,0$ and using \eqref{eq:BoundOnConditionalExpectation}, we obtain
        \begin{equation*}
        \E_x[\exp(-\chi_{\B}(\infty))] \leq \limsup_{N \rightarrow \infty} \left(\frac{a}{1+a}\right)^{N+1} =0.     
        \end{equation*}
        Hence, $\PP_x[\chi_{\B}(\infty)=\infty] = 1$.
    \end{proof}
    
   \begin{lemma}
    \label{lem:ClassifyTransient}
        Let $X=\{X(t): t \geq 0\}$ be a continuous time Markov chain with state space $\X$ and infinitesimal generator $Q=(Q_{x,y})_{x,y \in \X}$. Suppose there is an absorbing state $y \in \X$. If there are distinct states $z_1,\ldots,z_m \in \X$ different from $x$ and $y$ such that $Q_{x,z_1}Q_{z_1,z_2} \dots Q_{z_m,y}>0$, then $x$ is a transient state for $X$.
    \end{lemma}
    
    \begin{proof}
        Since $Q_{x,z_1}Q_{z_1,z_2} \dots Q_{z_m,y}>0$, we have $\PP_x[X(t_0)=y]>0$ for some $t_0 > 0$ (see Theorem 3.2.1 in \cite{Norris2}). Thus $\PP_x [X(t) \neq x \text{ for all } t\ge t_0] \geq \PP_x[X(t_0)=y] \PP[X(t)=y \text{ for all } t>t_0 | X(t_0)=y] = \PP_x[X(t_0)=y]>0$, which means that $x$ is a transient state.
    \end{proof}

%% file: Appendix_StatDist.tex
\

\

\subsection{Additional results for stationary distributions}
\label{sec:AppendixHigherOrder}

\subsubsection{Higher order terms for linear perturbations}
\label{sec:StatDistHigherOrder}

    Under the assumption of the perturbation being linear (which is the case for our chromatin modification circuit models), we now provide an iterative procedure for computing all of the terms in the series expansion of $\pi(\cdot)$. Additional results for characterizing some of these terms are given in SI - Sections \ref{sec:AdditionalCharactStatDist}, \ref{sec:PartialBalance}.

    
    \begin{theorem}
    \label{thm:HigherOrderTerms_for_stationary distribution}
    Suppose Assumptions \ref{assumption:transient_absorbing}, \ref{assumption:SingleRecurrentClass} and \ref{assumption:linear_case} hold. Then, the following hold for the sequence $\{\pi^{(k)}:\: k \geq 0\}$ in \eqref{eq:ExpansionPiEps}:
        \begin{enumerate}
            \item[(i)]
            $\pi^{(0)}=[\alpha^{(0)},\beta^{(0)}]=[\alpha,0]$ where $\alpha$ is the unique probability vector on $\A$ such that $\alpha Q_{\A}=0$,
            \item[(ii)]
            for every $k \geq 1$, $\pi^{(k)}=[\alpha^{(k)},\beta^{(k)}]$, where
                \begin{equation}
                \label{eq:FormulaBeta^k}
                    \beta^{(k)} = (\alpha^{(k-1)}S_1 + \beta^{(k-1)}T_1)(-T_0)^{-1}
                \end{equation}
        and $\alpha^{(k)}$ is the unique vector such that
            \begin{align}
            \label{eqn:findAlphaK}
            \alpha^{(k)}Q_{\A} &= -\beta^{(k)}[R_1 + T_1(-T_0)^{-1}R_0], \\
             \label{eq:findAlphaK_II}
             \alpha^{(k)}\one &= -\beta^{(k)}\one.  
            \end{align}
      Moreover, if $|\A| \geq 2$, for every $k \geq 1$ we obtain
            \begin{equation}
            \label{eq:FormulaForAlphaK}
            \alpha^{(k)} = \tilde{\alpha}^{(k)} + (- \beta^{(k)}\one - \tilde{\alpha}^{(k)}\one)\alpha,
            \end{equation}
       where $\tilde{\alpha}^{(k)} := - \beta^{(k)}(R_1 +T_1(-T_0)^{-1}R_0) Q_{\A}^\dagger$ for $k \geq 1$ and $Q_{\A}^\dagger$ is a generalized inverse of $Q_{\A}$ \footnote{A generalized inverse $Q_{\A}^\dagger$ of $Q_{\A}$ is such that $Q_{\A} Q_{\A}^\dagger Q_{\A} = Q_{\A}$. The Moore-Penrose inverse is a generalized inverse. The deviation matrix for $\hat{X}_\A$ is $D=(-Q_{\A}+\one \alpha)^{-1}-\one \alpha$, and $-D$ is also a generalized inverse of $Q_{\A}$. Meyer \cite{Meyer19752} suggested that $-D$ is a better generalized inverse to use than the Moore-Penrose inverse since it can be computed efficiently and embeds answers concerning the transitory behavior of the Markov chain. Avrachenkov et al. \cite{AvrachenkovFilarHowlett2}, in the context of discrete time Markov chains, use a suitable deviation matrix when they need a generalized inverse. Here, if we take $Q_{\A}^\dagger=-D$, then the term $\tilde{\alpha}^{(k)}\one$ in \eqref{eq:FormulaForAlphaK} is equal to zero since $D \one = 0$ and then $\alpha^{(k)} = \beta^{(k)}((R_1 +T_1(-T_0)^{-1}R_0)D - \one \alpha)$.}.
        \end{enumerate}
    \end{theorem}

    The proof of Theorem \ref{thm:HigherOrderTerms_for_stationary distribution} is given in SI - Section \ref{proofs}. 
    

\subsubsection{Proofs of Lemmas \ref{lem:alphaL_A0}, \ref{lem:welldefined_QA}, and \ref{lemma3} and Theorem \ref{thm:HigherOrderTerms_for_stationary distribution}}
\label{proofs}
\ 

\noindent \textbf{Proof of Lemma \ref{lem:alphaL_A0}}

    \begin{proof}
    It has already been established before Lemma \ref{lem:alphaL_A0} that $\pi(0)=[\alpha,0]$. By equating to zero the coefficients of the terms $\eps^m$ for $m=0,1, \ldots$ in the series $(\sum_{k=0}^{\infty}\eps^k\pi^{(k)})(\sum_{k=0}^{\infty}\eps^kQ^{(k)})$, we obtain that $\sum_{k=0}^{m} \pi^{(k)}Q^{(m-k)} = 0$ for every $m \geq 0$. In particular, $\pi^{(1)}Q^{(0)} + \pi^{(0)}Q^{(1)} = 0$, which yields,
    \begin{equation*}
    [\alpha^{(1)}, \beta^{(1)}]\left( \begin{array}{c | c }
                0 & 0 \\
               \hline
                R_0 & T_0
    \end{array}    \right) +   [\alpha, 0]\left( \begin{array}{c | c }
                A_1 & S_1 \\
                \hline
                R_1 & T_1
    \end{array}    \right) = 0.    
    \end{equation*}
    From this, we obtain two equations:
    \begin{equation}
    \label{eq:StatDist0order_I}
     \beta^{(1)} R_0 + \alpha A_1 = 0
    \end{equation}
    and
    \begin{equation}
    \label{eq:StatDist0order_II}
    \beta^{(1)} T_0 + \alpha S_1 =0. 
    \end{equation}
    Since $T_0$ is invertible, from \eqref{eq:StatDist0order_II} we obtain \eqref{eq:Beta1InTermsOfAlpha}. We conclude by substituting this formula for $\beta^{(1)}$ in \eqref{eq:StatDist0order_I}.
    \end{proof}
    
\noindent \textbf{Proof of Lemma \ref{lem:welldefined_QA}}
   
    \begin{proof}
   By following the proof of Lemma 2 in  Gaver, Jacobs \& Latouche \cite{Latouche} and using formula \eqref{eq:MatrixRepForProbabilities}, we can prove that the transition rates between $x \neq y \in \A$ for $\hat{X}_{\A}$ are given by $(Q_{\A})_{x,y}$. In essence, the argument is as follows. From the state $x\in \A$, the Markov chain $\tilde{X}$ may move to $y \in \A$ in two ways that lead to a one-step transition for $\hat{X}_{\A}$. First, it could happen that $\tilde{X}$ jumps directly to $y$ at a rate of $(A_1)_{x,y}$. Second, the chain $\tilde{X}$ may go to some state $z \in \T$ at a rate $(S_1)_{x,z}$ and from there, jump between states in $\T$ until getting back to $\A$ at the state $y \in \A$. By \eqref{eq:MatrixRepForProbabilities}, this happens with probability $((-T_0)^{-1}R_0)_{z,y}$. Putting this all together, the rate of transition for $\hat{X}_{\A}$ from $x$ to $y$ will be 
   \begin{equation}
   \label{formulaforlemma43}
     (A_1)_{x,y} + \sum_{z \in \T} (S_1)_{x,z}((-T_0)^{-1}R_0)_{z,y} = (Q_{\A})_{x,y}.
   \end{equation}
    \end{proof}

\noindent \textbf{Proof of Lemma \ref{lemma3}}
   
    \begin{proof}
    Consider $x\neq y \in \A$. Then, there exists a sequence of states $x_0=x,x_1,\ldots,x_m=y$ in $\CC$ such that $\tilde{Q}_{x,x_1}\tilde{Q}_{x_1,x_2}\ldots \tilde{Q}_{x_{m-1},y} > 0$. Roughly speaking, the proof follows by erasing the times that $\tilde{X}$ is outside of $\A$. We now give the details. Consider $i \in \{0,1,\ldots,m-1\}$ with $x_i \in \A$. If $x_{i+1} \in \A$, then, by \eqref{formulaforlemma43}, $(Q_{\A})_{x_i,x_{i+1}} \geq (A_1)_{x_{i},x_{i+1}} = \tilde{Q}_{x_i,x_{i+1}} > 0$. If $x_{i+1} \in \T$, consider the path of states $x_i, x_{i+1},\ldots, x_{k}$ for $0\leq i < k \leq m$ such that $x_i,x_k \in \A$ and $x_{i+1},\ldots, x_{k-1} \in \T$. Since $x_{i+1}$ leads to $x_k$ for $\tilde{X}$, then  $\PP_{x_{i+1}}[\tilde{X}(\tau_{\A}) = x_{k}] >0$ where $\tau_{\A} := \inf\{ t \geq 0:\: \tilde{X}(t) \in \A \}$. By \eqref{eq:MatrixRepForProbabilities}, this means that $((-T_0)^{-1}R_0)_{x_{i+1},x_{k}} > 0$ which yields $(Q_{\A})_{x_i,x_k} \geq (S_1)_{x_i,x_{i+1}}((-T_0)^{-1}R_0)_{x_{i+1},x_k} = \tilde{Q}_{x_i,x_{i+1}}((-T_0)^{-1}R_0)_{x_{i+1},x_k} > 0$. These observations yield a sequence of states $x_0=x,x_{i_1},\ldots,x_{i_j}=y$ in $\A$ such that $(Q_{\A})_{x,x_{i_1}}(Q_{\A})_{x_{i_1},x_{i_2}}\ldots (Q_{\A})_{x_{i_{j-1}},y} > 0$.
    \end{proof}


\noindent \textbf{Proof of Theorem \ref{thm:HigherOrderTerms_for_stationary distribution}}

    \begin{proof}
    Point $(i)$ was established in Theorem \ref{thm:alpha_unique}. For $(ii)$, we equate to zero the coefficients of the terms $\eps^m$ for $m=0,1,2, \ldots$ in the terms of the series, $(\sum_{k=0}^{\infty}\eps^k\pi^{(k)})(Q^{(0)} + \eps Q^{(1)})$
     to obtain that $ \pi^{(0)}Q^{(0)} = 0$ and $\pi^{(k)}Q^{(0)} + \pi^{(k-1)}Q^{(1)} = 0$ for every $k \geq 1$. The latter requires that for all $k \geq 1$,
     \begin{equation*}
         [\alpha^{(k)}, \beta^{(k)}]\left( \begin{array}{c | c }
                0 & 0 \\
                \hline
                R_0 & T_0
    \end{array}    \right) +   [\alpha^{(k-1)}, \beta^{(k-1)}]\left( \begin{array}{c | c }
                A_1 & S_1 \\
                \hline
                R_1 & T_1
    \end{array}    \right) = 0.
     \end{equation*}
    Now, this yields two equations:
        \begin{equation}
        \label{eq:pi_k_1}
           \beta^{(k)} R_0 + \alpha^{(k-1)} A_1 + \beta^{(k-1)}R_1 = 0, 
        \end{equation}
        \begin{equation}
        \label{eq:pi_k_2}
            \beta^{(k)} T_0 +\alpha^{(k-1)} S_1 + \beta^{(k-1)}T_1 = 0.
        \end{equation}
        For $\beta^{(k)}$, we obtain the relation \eqref{eq:FormulaBeta^k} directly from \eqref{eq:pi_k_2} for all $k \geq 1$. For $\alpha^{(k)}$, let's see first that it satisfies \eqref{eqn:findAlphaK}. From \eqref{eq:FormulaBeta^k}, for all $k \geq 1$, $\beta^{(k+1)} = (\alpha^{(k)} S_1 + \beta^{(k)}T_1)(-T_0)^{-1}$ and using this in \eqref{eq:pi_k_1} (with $k$ replaced by $k+1$) we obtain for all $k \geq 1$
        \begin{equation}
        \label{eqS11}
         (\alpha^{(k)} S_1 + \beta^{(k)}T_1)(-T_0)^{-1}R_0 + \alpha^{(k)} A_1 + \beta^{(k)}R_1 = 0.  
        \end{equation}
    By rearranging \eqref{eqS11} and using \eqref{eq:QA}, we obtain \eqref{eqn:findAlphaK} for all $k \geq 1$. On the other hand, since $\inn{\pi^{(k)},\one} = 0$ for every $k \geq 1$, we obtain \eqref{eq:findAlphaK_II}.
    
    For the uniqueness of $\alpha^{(k)}$, for all $k \geq 1$, if $|\A| = 1$, then $\alpha^{(k)}$ has only one entry and it is determined uniquely by \eqref{eq:findAlphaK_II}. If $|\A| \geq 2$, consider another solution $\gamma^{(k)}$ of \eqref{eqn:findAlphaK} and \eqref{eq:findAlphaK_II}, where $(\gamma^{(k)})^T \in \R^{|\A|}$. By Assumption \ref{assumption:SingleRecurrentClass} and Lemma \ref{lem:nullityQmatrix}, $\dim(\ker((Q_{\A})^T))=1$ and therefore, by \eqref{eqn:findAlphaK}, $\alpha^{(k)} - \gamma^{(k)} = c\alpha$ for some $c \in \R$. Using \eqref{eq:findAlphaK_II}, then $0=\alpha^{(k)}\one - \gamma^{(k)}\one = c\alpha\one = c$, and therefore $c=0$, and $\alpha^{(k)} = \gamma^{(k)}$.
    
    For existence of a solution $\alpha^{(k)}$ of \eqref{eqn:findAlphaK}-\eqref{eq:findAlphaK_II}, 
    using the properties
    \begin{eqnarray*}
    &&R_0 \one + T_0 \one =0\quad \hbox{ and }\quad
    R_1 \one + T_1 \one =0 ,
    \end{eqnarray*}
    we have that
    \begin{eqnarray*}
    (R_1  + T_1 (-T_0)^{-1} R_0 )\one 
    &=&R_1 \one + T_1 (-T_0)^{-1} R_0 \one \\
    &=&-T_1 \one + T_1 (-T_0)^{-1} R_0 \one \\
    &=& T_1 (-T_0)^{-1}(T_0 \one + R_0 \one )\\
    &=&0.
    \end{eqnarray*}
    Then, since $\dim(\ker(Q_{\A}))=\dim(\ker((Q_{\A})^T))=1$ and $\one \in \ker (Q_{\A})$, we have 
    \begin{equation*}
        (- \beta^{(k)} (R_1   + T_1 (-T_0)^{-1} R_0 ))^T \in \ker (Q_{\A})^\perp = \text{range} ((Q_{\A})^T),
    \end{equation*}
    and so \eqref{eqn:findAlphaK} has a solution and
    \eqref{eq:findAlphaK_II} will determine the multiple of $\alpha$ to add to any particular solution to obtain the unique solution $\alpha^{(k)}$ of both equations.
    
    Furthermore, 
    if $Q_{\A}^\dagger$ is a generalized inverse of $Q_{\A}$, then
        \begin{equation}
        \tilde{\alpha}^{(k)} := - \beta^{(k)}(R_1 +T_1(-T_0)^{-1}R_0) Q_{\A}^\dagger    
        \end{equation}
    is a solution to \eqref{eqn:findAlphaK} (see \cite{James1978} for an exposition). Similar to the uniqueness argument, $\alpha^{(k)} - \tilde{\alpha}^{(k)} = c\alpha$ for some $c \in \R$. By \eqref{eq:findAlphaK_II},  $c=-\beta^{(k)}\one -\tilde{\alpha}^{(k)}\one $ and we obtain \eqref{eq:FormulaForAlphaK}.
    \end{proof}
    
\subsubsection{Additional characterization of zeroth and first order terms for linear perturbations via restricted processes}        
\label{sec:AdditionalCharactStatDist}        

    In this section, assume that Assumptions \ref{assumption:transient_absorbing} and \ref{assumption:linear_case} hold. We will also sometimes assume Assumptions \ref{assumption:A_is_recurrent} or \ref{assumption:tildeX_irreducible} hold. We will explore additional characterizations of $\alpha$ and $\beta^{(1)}$. Under Assumptions \ref{assumption:transient_absorbing} and \ref{assumption:linear_case}, $A(\eps)$ (defined in \eqref{eq:Q_eps_decomposition}) corresponds to $\eps A_1$ for every $0 \leq \eps < \eps_0$. Since $Q(\eps)$ is irreducible for every $0 < \eps < \eps_0$, from Lemma \ref{lem:InvertibilityQ^B} (with $\A^c$ in place of $\B$ and $Q(\eps)$ in place of $Q$ there), we obtain that $\eps A_1$ is invertible for $0 < \eps < \eps_0$, and therefore $A_1$ is invertible. This will be an important fact for the coming results.
    
    Consider the matrix $\tilde{Q}$ introduced in \eqref{eq:DefTildeQ}. For a continuous time Markov chain $\tilde{X}$ with infinitesimal generator $\tilde{Q}$, denote by $\chi_{\T}(t)$ the occupation time of $\T$ by the Markov chain $\tilde{X}$ up to time $t \geq 0$, with its associated limit $\chi_{\T}(\infty) = \lim_{t \to \infty} \chi_{\T}(t)= \int_0^{\infty}\one_{\T}(X(t))dt$. Since $A_1$ is invertible, by Lemma \ref{lem:InvertibilityQ^B} (with $\B=\T$ and $Q=\tilde{Q}$) and Lemma \ref{lem:OccupationTimesAreInfinity} we have that $\PP[\chi_{\T}(\infty)=\infty]=1$.
    
    Consider the process $\hat{X}_{\T}$ as in \eqref{eq:RestrictiontoAbsorbing}, but with $\A$ replaced by $\T$, which corresponds to observing $\tilde{X}$ only on the time intervals where $\tilde{X}$ is in $\T$. The process $\hat{X}_{\T}$ is a continuous time Markov chain on $\T$. Consider the matrix
    \begin{equation}
         \label{eq:Q_T_linearcase}
         Q_{\T} := T_0 + R_0(-A_1)^{-1}S_1,
        \end{equation}
    which by Lemma \ref{lem:Q_A_Qtilde_are_Qmatrix} is a $Q$-matrix. Similarly to Lemma \ref{lem:welldefined_QA}, we can show that $Q_{\T}$ is the infinitesimal generator of $\hat{X}_{\T}$. Our previous assumptions relate to $\hat{X}_{\T}$ in the following way. 
    
    \begin{lemma}
    \label{lem:SingleRecurrentClassQ_T}
    Suppose Assumptions \ref{assumption:transient_absorbing}, \ref{assumption:A_is_recurrent} and \ref{assumption:linear_case} hold. Then $\hat{X}_{\T}$ has a single recurrent class. Moreover, if Assumption \ref{assumption:tildeX_irreducible} holds, the process $\hat{X}_{\T}$ is irreducible.
    \end{lemma}

    
    \begin{proof}
    Let $\D \subseteq \T$ be a non-empty recurrent class for $\hat{X}_{\T}$ (there must be at least one), and let $\CC \subseteq \X$ be the communicating class under $\tilde{X}$ described in Assumption \ref{assumption:A_is_recurrent}. We will prove that $\D = \CC \setminus \A$, which yields the uniqueness of recurrent classes for $\hat{X}_{\T}$.  If Assumption \ref{assumption:tildeX_irreducible} holds, then $\CC=\X$, which combined with the relation $\D = \CC \setminus \A$, implies that $\D = \X \setminus \A = \T$ and the conclusion follows.

    In order to prove $\D = \CC \setminus \A$, we start by making some observations. First, we prove that there exist states $\tilde{x} \in \D$ and $\tilde{y} \in \A$ such that $\tilde{Q}_{\tilde{x},\tilde{y}} >0$. In fact, if this was not the case, then for every $x \in \D$ and $z \in \A$ we would have $\tilde{Q}_{x,z} = (R_0)_{x,z} = 0$. This yields that for $x \in \D$,
    \begin{equation}
    (Q_{\T})_{x,y} = (T_0)_{x,y} + \sum_{z \in \A}(R_0)_{x,z}[(-A_1)^{-1}S_1]_{z,y} = (T_0)_{x,y},  
    \end{equation}
    for all $y \in \T$. Since $\D$ is a closed class under $Q_{\T}$, $(Q_{\T})_{x,y}=0$ for $y \in \T \setminus \D$ and so $\sum_{y \in \D}(Q_{\T})_{x,y} = \sum_{y \in \T}(Q_{\T})_{x,y}=0$, since $Q_{\T}$ is a $Q$-matrix. Combining this with the previous equation, we obtain that $\sum_{y \in \D} (T_0)_{x,y} = 0$ for all $x \in \D$, which implies that $\D$ is closed under $\tilde{Q}$. This contradicts the fact that $T_0$ is invertible by point $(i)$ in Lemma \ref{lem:InvertibilityQ^B} (with $\B^c = \T$ and $Q=\tilde{Q}$).
    
    Second, we observe that there exist a $\hat{y} \in \A$ and $\hat{x} \in \D$ such that $\tilde{Q}_{\hat{y},\hat{x}} > 0$. In fact, since $A_1$ is invertible, by Lemma \ref{lem:InvertibilityQ^B} (with $\B^c = \A$ and $Q=\tilde{Q}$) there has to be a $\hat{y} \in \A$ and $\hat{x} \in \T$ such that $\tilde{Q}_{\hat{y},\hat{x}} > 0$. To show that $\hat{x} \in \D$, consider that $\tilde{y} \in \A$ communicates with $\hat{y} \in \A$ under $\tilde{Q}$, by Assumption \ref{assumption:A_is_recurrent}, and therefore $\tilde{x}$ leads to $\hat{x}$ under $\tilde{Q}$ and therefore under $Q_{\T}$. Since $\D$ is closed under $Q_{\T}$, $\hat{x} \in \D$. 

    We now prove that $\D \subseteq \CC \setminus \A$. For $x \in \D$, since $T_0$ is invertible there exists a state $y \in \A$ such that $x$ leads to $y$ under $\tilde{Q}$. By Assumption \ref{assumption:A_is_recurrent}, $y$ and $\hat{y}$ are in $\A \subseteq \CC$ and so they communicate under $\tilde{Q}$. It follows that $x$ leads to $\hat{y}$ under $\tilde{Q}$. On the other hand, $\tilde{Q}_{\hat{y},\hat{x}} > 0$ and since $\D$ is a communicating class under $Q_{\T}$, $\hat{x}$ leads to $x$ under $\tilde{Q}$. Thus, $x$ leads to $\hat{y}$ and $\hat{y}$ leads to $x$ under $\tilde{Q}$ and so $x \in \CC$. Thus, $\D \subseteq \CC$ and $\D \subseteq \T=\A^c$, and so $\D \subseteq \CC \setminus \A$. 
    
    To prove that $\CC \setminus \A \subseteq \D$, let $x \in \CC \setminus \A$. Since $\D \subseteq \CC$, then $x$ communicates with the element $\tilde{x}$ in $\D$ under $\tilde{Q}$. This implies that they communicate under $Q_{\T}$ and since $\D$ is a communicating class under $Q_{\T}$, then $x \in \D$. Combining the above we see that $\D=\CC \setminus \A$.
    \end{proof}

    
    When the continuous time Markov chain $\hat{X}_{\T}$ has a single recurrent class $\D$, there is a unique probability vector $\nu$ in $\R^{|\T|}$ such that $\nu Q_{\T} = 0$ and $\nu$ will be the stationary distribution for $\hat{X}_{\T}$ with non-zero entries only for entries corresponding to states in $\D$. We use the vector $\nu$ to characterize $\alpha$ and $\beta^{(1)}$.

    In the following theorem, we use the fact that $A_1$ is invertible. This follows from Lemma S.4 because $A_1=Q^{\A}(1)$, where $Q(1)$ is positive recurrent and so the condition (i) of Lemma S.4 holds with $B_c = \A$.
  
    \begin{theorem}
    \label{thm:alpha_nu}
    Suppose Assumptions \ref{assumption:transient_absorbing}, \ref{assumption:A_is_recurrent} and \ref{assumption:linear_case} hold. Denote by $\nu$ the unique probability vector in $\R^{|\T|}$ such that $\nu Q_{\T} =0$. Then, $\pi(0)=[\alpha,0]$ where $\alpha$ is given by
            \begin{equation}
            \label{eq:AltCharacAlpha}
                \alpha = c\nu R_0(-A_1)^{-1},
            \end{equation}
    and where $c$ is given by $c = (\nu R_0(-A_1)^{-1}\one)^{-1}$. Moreover, $\pi^{(1)}=[\alpha^{(1)},\beta^{(1)}]$ where
        \begin{equation}
        \label{eq:AltCharacBeta}
        \beta^{(1)} = c \nu.
        \end{equation}
    \end{theorem}
    
    
    \begin{proof}
    Following the proof of Theorem \ref{thm:HigherOrderTerms_for_stationary distribution}, equations \eqref{eq:pi_k_1} and \eqref{eq:pi_k_2} yield that
        \begin{equation}
        \label{eq:AltCharac1}
         \beta^{(1)} R_0 + \alpha A_1 = 0, 
        \end{equation}
    and
        \begin{equation}
        \label{eq:AltCharac2}
        \beta^{(1)} T_0 +\alpha S_1 = 0.
        \end{equation}
    From \eqref{eq:AltCharac1} we obtain that  $\alpha = \beta^{(1)}R_0(-A_1)^{-1}$. We substitute this expression in \eqref{eq:AltCharac2} to obtain that $\beta^{(1)}(T_0 + R_0(-A_1)^{-1}S_1) =0$, which is exactly $\beta^{(1)}Q_{\T}=0$. By Lemma \ref{lem:SingleRecurrentClassQ_T} combined with Lemma \ref{lem:nullityQmatrix}, we obtain that $\beta^{(1)} = \tilde{c}\nu$ for some constant $\tilde{c} \in \R$ and therefore $\alpha = \tilde{c}\nu R_0(-A_1)^{-1}$. To show that $\tilde{c}=c$, we observe that since $\alpha\one=1$, then $\tilde{c}(\nu R_0(-A_1)^{-1}\one)=1$ and the desired result follows.
    \end{proof}

    Under the assumptions of Theorem \ref{thm:alpha_nu}, $\beta^{(1)}_x > 0$ for every $x \in \D$, the single recurrent class of $\hat{X}_{\T}$, while $\beta^{(1)}_x=0$ for $x \in \T^{\setminus (1)}$. In fact, 
    %
    %
    using first step analysis, one can show that the entry $(-A_1)^{-1}_{i,j}$ is the expected time that the process $\tilde X$ spends at $j$ when started at $i$, before exiting $\A$. Hence, these entries are non-negative and so does $\nu R_0(-A_1)^{-1}$. This implies that the constant $c$ is positive and the conclusion follows from \eqref{eq:AltCharacBeta} and the properties of $\nu$.

\subsubsection{Additional characterization of zeroth and first order terms via partial balance}
\label{sec:PartialBalance}

 For the last part of this section, we consider an additional characterization for $\beta^{(1)}$ based on the idea of truncated processes and partial balance relations (see Section 9.4 in \cite{FrankKelly}).
    
    Consider a continuous time Markov chain $X = \{X(t):\: t \geq 0 \}$ with infinitesimal generator $Q$ on a finite state space $\X$. Let $\Y$ a non-empty set in $\X$. Define the matrix $\bar{Q} = (\bar{Q}_{x,y})_{x,y \in \Y}$ by $\bar{Q}_{x,y}=Q_{x,y}$ for $x \neq y$ and $\bar{Q}_{x,x}=Q_{x,x} + \sum_{y \notin \Y}Q_{x,y}$. A continuous time Markov chain $\bar{X} = \{\bar{X}(t):\: t \geq 0 \}$ with state space $\Y$ and infinitesimal generator $\bar{Q}$ will be called a \textbf{truncation} of $X$ to $\Y$.

    \begin{assumption}
   \label{assumption:partial_balance}
    For every $0 < \eps < \eps_0$, the truncation of $X^{\eps}$ to $\T$, denoted by $\bar{X}^{\eps}$, is irreducible. In addition, the following partial balance condition holds on $\T$ for every $0 < \eps < \eps_0$:
        \begin{equation}
        \label{eq:partial_balance}
            \pi_x(\eps)\sum_{y \in \A}Q_{x,y}(\eps) = \sum_{y \in \A}\pi_y(\eps)Q_{y,x}(\eps), \qquad  \text{ for every } x \in \T.
        \end{equation}
    \end{assumption}
    
    
    \medskip
    Under Assumption \ref{assumption:partial_balance}, the process $\bar{X}^{\eps}$ has a stationary distribution $\eta(\eps)$ for every $0 < \eps < \eps_0$,
    %
    given by \begin{equation}
            \label{eq:nu_is_conditioned}
                \eta_x(\eps) = \frac{\pi_x(\eps)}{\sum_{y \in \T} \pi_y(\eps)},\;\;\;\;x \in \T
            \end{equation}
  (see Theorem 9.5 in Kelly \cite{FrankKelly}). The following is our main theorem.
  
    
    \begin{theorem}
    \label{thm:alpha_beta_eta}
    Suppose Assumptions \ref{assumption:transient_absorbing}, \ref{assumption:linear_case} and \ref{assumption:partial_balance} hold. Then, the following hold:
        \begin{enumerate}
            \item[(i)]
            the limit $\eta := \lim_{\eps \to 0} \eta(\eps)$ exists and it is a probability vector on $\T$ such that
                \begin{equation}
                \eta\bar{Q}(0)=0,
                \end{equation}
            \item[(ii)]
            the vectors $\alpha$ and $\beta^{(1)}$ can be characterized by
             \begin{equation}
            \beta^{(1)} = c\eta, \qquad \qquad \alpha = c\eta  R_0(-A_1)^{-1},
            \end{equation}
            where $c = (\eta R_0(-A_1)^{-1}\one)^{-1}$, and
            \item[(iii)]
            $\eta Q_{\T}=0$. 
            
            If, in addition, Assumption \ref{assumption:A_is_recurrent} or Assumption \ref{assumption:tildeX_irreducible} holds, then $\eta=\nu$ is the unique stationary distribution for $\hat{X}_{\T}$.
        \end{enumerate}
    \end{theorem}

    \begin{remark}
    Although we know that $\bar{X}^{0}$ is well defined, we do not know a priori whether the process is irreducible or it has a single recurrent class. If the truncation process $\bar{X}^{0}$ has a single recurrent class (or is irreducible), it will have a unique stationary distribution, which we would call $\eta(0)$. But we do not know if such a vector exists. This non existence is what led us to express Theorem \ref{thm:alpha_beta_eta} in terms of the limit $\eta$ which solves $\eta \bar{Q}(0)=0$. If the truncation process $\bar{X}^{0}$ has a single recurrent class, as in the 1D and 2D models, the probability vector $\eta$ is characterized uniquely by solving $\eta\bar{Q}(0)=0$ and $\eta=\eta(0)$.  
    \end{remark}
    
    
        \begin{proof}
        We will first show that $\beta^{(1)}\one > 0$. From \eqref{eq:Beta1InTermsOfAlpha} we know that $\beta^{(1)} = \alpha S_1(-T_0)^{-1}$, which yields $\beta^{(1)}\one = \alpha S_1(-T_0)^{-1}\one$. Since all of the entries in $\alpha$, $S_1$, and $(-T_0)^{-1}$ are nonnegative, it suffices to show $\beta^{(1)}\one \neq 0$. For a proof by contradiction, suppose that $\beta^{(1)}\one=0$. This means that
            \begin{equation}
            \sum_{y \in \T}\sum_{x \in \A} \alpha_x(S_1)_{x,y}((-T_0)^{-1}\one)_y  =0.
            \end{equation}
        All of the entries in the sum are nonnegative, so this means that $\alpha_x(S_1)_{x,y}((-T_0)^{-1}\one)_y=0$ for every $x \in \A$ and $y \in \T$. Now, $((-T_0)^{-1}\one)_y$ is the mean first passage time to $\A$, for the Markov chain that starts at $y$ with infinitesimal generator $Q(0)$, and so $((-T_0)^{-1}\one)_y \geq \frac{1}{|Q(0)_{y,y}|} > 0$. Hence, $\alpha_x(S_1)_{x,y} = 0$ for every $x \in \A$ and $y \in \T$. This yields that $\alpha S_1 =0$ 
        and substituting this in \eqref{eqsatisfied} yields that $\alpha A_1 =0$. Since $A_1$ is invertible, this is a contradiction. 
        
       Since we know that $\beta(0)=\beta^{(0)}=0$, we obtain that
            \begin{align*}
             \eta_x(\eps) &= \frac{\pi_x(\eps)}{\sum_{y \in \T} \pi_y(\eps)} = \frac{\sum_{k=1}^{\infty}\eps^{k} \beta_x^{(k)}}{\sum_{y \in \T} \sum_{k=1}^{\infty}\eps^{k}\beta_y^{(k)}} =  \frac{\sum_{k=1}^{\infty}\eps^{k}\beta_x^{(k)}}{\sum_{k=1}^{\infty}\eps^{k} \sum_{y \in \T} \beta_y^{(k)}}= \frac{\sum_{k=1}^{\infty}\eps^{k-1}\beta_x^{(k)}}{\sum_{k=1}^{\infty}\eps^{k-1} \sum_{y \in \T} \beta_y^{(k)}} \\
             & \to  \frac{\beta_x^{(1)}}{\sum_{y \in \T} \beta_y^{(1)}} = \frac{\beta_x^{(1)}}{\beta^{(1)}\one}.
        \end{align*}
       We then obtain that $\eta$ exists and $\eta_x = \frac{\beta_x^{(1)}}{\beta^{(1)}\one}$ for every $x \in \T$, which is a probability vector on $\T$. Or letting $\eps \rightarrow 0$ in $\eta(\eps)\bar{Q}(\eps) = 0$, we obtain that $\eta \bar{Q} =0$.
       %
        %
 We already know that $\beta^{(1)}=c\eta$. To obtain a value for $c$ that depends only on $\eta$, note that from \eqref{eqsatisfied} and \eqref{eq:Beta1InTermsOfAlpha}, we have $\alpha = \beta^{(1)}R_0(-A_1)^{-1} =c\eta R_0(-A_1)^{-1}$, where $c\eta R_0(-A_1)^{-1}\one =1$ and so $c = (\eta R_0(-A_1)^{-1}\one)^{-1}$.
       
       By following the proof of Theorem \ref{thm:alpha_nu}, we obtain $\beta^{(1)}Q_{\T}=0$ and therefore $\eta Q_{\T}=0$. The other conclusions follow readily.
    \end{proof}
        
    The following criterion offers a practical way to establish \eqref{eq:partial_balance}.
    
        \begin{lemma}
        Let $\A= \{a_1, \ldots, a_n\}$. Suppose there exist distinct states $x_1, \ldots , x_n$ in $\T$ such that for every $0 < \eps < \eps_0$ and for every $k \in \{1,\ldots,n\}$.
            \begin{enumerate}
                \item
                $Q_{a_kx_k}(\eps),Q_{x_ka_k}(\eps) > 0$,
                \item
                $Q_{a_ky}(\eps) = Q_{ya_k}(\eps) = 0$ for every $y \notin \{x_k,a_k\}$.
            \end{enumerate}
        Then \eqref{eq:partial_balance} holds.    
        \end{lemma}
        
            \begin{proof}
            Denote by $\NN = \{x_1,\ldots,x_n\}$. Let $0 < \eps < \eps_0$. For $x \in \T \setminus \NN$, we have that $\pi_x(\eps)\sum_{y \in \A}Q_{x,y}(\eps) =0$ and $\sum_{y \in \A}\pi_y(\eps)Q_{y,x}(\eps) = 0$. Then, equation \eqref{eq:partial_balance} holds for $x \in \T \setminus \NN$.

            For $x_{k} \in \NN$, $\pi_{x_k}(\eps)\sum_{y \in \A}Q_{x_k,y}(\eps)=\pi_{x_{k}}(\eps)Q_{x_{k}a_{k}}(\eps)$. 
            
            On the other hand, $\sum_{y \in \A}\pi_y(\eps)Q_{yx_{k}}(\eps) $ $= \pi_{a_k}(\eps)Q_{a_kx_k}(\eps)$. Since $\pi(\eps)Q(\eps) = 0$, we have
            \begin{align*}
            0=(\pi(\eps)Q(\eps))_{a_k}&=\sum_{x \in \X} \pi_{x}(\eps) Q_{x,a_k}(\eps) \\            &=\pi_{x_k}(\eps) Q_{x_k,a_k}(\eps) + \pi_{a_k}(\eps) Q_{a_k,a_k}(\eps)= \pi_{x_k}(\eps) Q_{x_k,a_k}(\eps) - \pi_{a_k}(\eps) Q_{a_k,x_k}(\eps).
            \end{align*}
            
            Hence, $ \pi_{x_k}(\eps) Q_{x_k,a_k}(\eps) = \pi_{a_k}(\eps) Q_{a_k,x_k}(\eps)$ and \eqref{eq:partial_balance} holds for $x \in \NN$ as well.
            \end{proof}

\subsubsection{Lemma \ref{lem:Q_A_Qtilde_are_Qmatrix}}

    \begin{lemma}
    \label{lem:Q_A_Qtilde_are_Qmatrix}
    Under Assumption \ref{assumption:transient_absorbing}, the matrices $Q_{\A}$ and $\tilde{Q}$ are $Q$-matrices of sizes $|\A| \times |\A|$ and $|\X| \times |\X|$ respectively. If in addition, Assumption  \ref{assumption:linear_case} holds, then $Q_{\T}$ is a $Q$-matrix of size $|\T| \times |\T|$. 
    \end{lemma}
 
    \begin{proof}
    First, observe that $\lim_{\eps \to 0} \frac{1}{\eps}Q_{x,y}(\eps) = (A_1)_{x,y}$ if $x,y \in \A$, while $\lim_{\eps \to 0} \frac{1}{\eps}Q_{x,y}(\eps) = (S_1)_{x,y}$ if $x \in \A$ and $y \in \T$. Then, since $Q(\eps)$ is a $Q$-matrix, $S_1$ has nonnegative entries, $(A_1)_{x,y} \geq 0$ for $x \neq y \in \A$, and 
        \begin{equation}
        \label{eq:EstbalishingQMatrices}
         \sum_{y \in \A} (A_1)_{x,y} + \sum_{y \in \T} (S_1)_{x,y} = 0 \text{ for every } x \in \A.   
        \end{equation}
    For $x \neq y \in \A$, $(Q_{\A})_{x,y} = (A_1)_{x,y} + \sum_{z \in \T} (S_1)_{x,z}((-T_0)^{-1}R_0)_{z,y}$ is nonnegative since $(A_1)_{x,y} \geq 0$ and by \eqref{eq:MatrixRepForProbabilities}, $((-T_0)^{-1}R_0)_{z,y}=0$ for $z \in \T$. For $x \in \A$, $ \sum_{y \in \A} (Q_{\A})_{x,y}$ is equal to
    \begin{align*} 
    & \sum_{y \in \A} (A_1)_{x,y} +  \sum_{y \in \A} \sum_{z \in \T} (S_1)_{x,z}((-T_0)^{-1}R_0)_{z,y} = \sum_{y \in \A} (A_1)_{x,y} +  \sum_{z \in \T} (S_1)_{x,z} \sum_{y \in \A} ((-T_0)^{-1}R_0)_{z,y} \\
    &= \sum_{y \in \A} (A_1)_{x,y} + \sum_{z \in \T} (S_1)_{x,z} = 0,
    \end{align*}
    where we used \eqref{eq:MatrixRepForProbabilities} and \eqref{eq:EstbalishingQMatrices}. Hence $Q_{\A}$ is a $Q$-matrix.
    
    For $\tilde{Q}$, for $x \neq y \in \X$, if $x \in \A$, then $\tilde{Q}_{x,y}$ corresponds to an off diagonal term in $A_1$ or a term in $S_1$, both of which are nonnegative. If $x \in \T$,  then $\tilde{Q}_{x,y}$ corresponds to an off diagonal term in $T_0$ or a term in $R_0$, which are both nonegative since $Q(0)$ is a $Q$-matrix. To check that the row-sums of $\tilde{Q}$ are zero, first consider when $x \in \A$. Then, $\sum_{y \in \X}\tilde{Q}_{x,y} = \sum_{y \in \A} (A_1)_{x,y} + \sum_{y \in \T} (S_1)_{x,y} = 0$ by \eqref{eq:EstbalishingQMatrices}. If $x \in \T$, then $\sum_{y \in \X}\tilde{Q}_{x,y} = \sum_{y \in \A} (R_0)_{x,y} + \sum_{y \in \T} (T_0)_{x,y} = 0$, since this corresponds to summing across a row of $Q(0)$.

    The case of $Q_{\T}$ follows similarly to that for $Q_{\A}$.
    \end{proof}

%% file: Appendix_orderCTMCalgorithm.tex

\bigskip

\subsection{Algorithm to find the order of the pole of the MFPT}
\label{sec:Algstat}
 		
	\small
	
	\noindent \textbf{Input:} $\B \subset \mathcal{X}$, and $k_{xy}$, the order of $Q_{x,y} (\eps)$ for each $(x,y) \in E_0$.
	
	\noindent \textbf{Output:} $p(x)$, the order of the pole of the mean first passage from $x \in \B^c$ to $\B$.
	
	($p$ will also be defined for condensed nodes in the course of the algorithm)
	
	\noindent \textbf{Step 1 (Set up the initial graph $(V,E)$)}
	
	\noindent Construct a directed graph $G=(V,E)$ where $V = \X$ and $E = E_0$.
	
	\noindent Set, for each $u \in V$, $p (u) \leftarrow \min \{ k_{uv}:\: (u,v) \in E \}$.
	
	\noindent Set, for each $(u,v) \in E$, $\mathcal{K}_{uv} \leftarrow k_{uv} - p (u)$.
	
	\noindent \textbf{Step 2 (Condense $\B$ into a single node $a$)}
	
	\noindent Introduce a new node $a$.
	
	\noindent Set, for each $w \in \B^c$ such that $(w,v) \in E$ for some $v \in \B$, $\mathcal{K}_{wa} \leftarrow \min \{ \mathcal{K}_{wv}:\: v \in \B \text{ and } (w,v) \in E\}$.
	
	\noindent Update $V \leftarrow \B^c \cup \{a\}$ and
	
	\noindent $\quad \quad \quad E \leftarrow \{(u,v) \in E:\: u \in \B^c \text{ and } v \in \B^c \} \cup \{(w,a):\: (w,v) \in E \text{ for some } w \in \B^c \text{ and } v \in \B \}$.
	
	\noindent Set, for each $u \in V \setminus \{a\}$, $S(u) \leftarrow  \{u\}$, and $S(a) \leftarrow \B$.
	
	\noindent \textbf{Step 3 (Condense r-connected sets)}
	
	\noindent Repeat the following until $G$ contains no r-connected sets:
	
	Let $C \subset V$ be an r-connected set and $c$ be a new node representing the r-connected set $C$.
	
	Set $p (c) \leftarrow \max_{u \in C} p (u) + \min \{\mathcal{K}_{uv}:\: u \in C, v \notin C \text{ and } (u,v) \in E \}$.
	
	Set, for each $w \in V \setminus C$ such that $(u,w) \in E$ for some $u \in C$,
	
	$\quad \quad \quad \mathcal{K}_{cw} \leftarrow \min \{\mathcal{K}_{uw}:\: u \in C \text{ and } (u,w) \in E \} - \min \{\mathcal{K}_{uv}:\: u \in C, v \notin C \text{ and } (u,v) \in E \}$
	
	Set, for each $w \in V \setminus C$ such that $(w,v) \in E$ for some $v \in C$, 
	
	$\quad \quad \quad \mathcal{K}_{wc} \leftarrow \min \{\mathcal{K}_{wv}:\: v \in C \text{ and } (w,v) \in E \}$.
	
	Update $V \leftarrow (V \setminus C) \cup \{c\}$ and 
	
	$\quad \quad \quad E \leftarrow \{(u,v) \in E:\: u \notin C \text{ and } v \notin C\} \cup \{(c,w):\: (u,w) \in E \text{ for some } u \in C\text{ and } w \notin C\}$
	
	$\quad \quad \quad \quad \quad \quad \quad \cup \{(w,c):\: (w,v) \in E \text{ for some }  w \notin C \text{ and } v \in C\}$.
	
	Set $S(c) \leftarrow \cup_{u \in C} S(u)$.
	
	\noindent \textbf{Step 4 (Compute $p(x)$ where $x \in \B^c$)}
	
	\noindent Repeat the following until $V = \{a\}$:
	
	Let $v^* \in V \setminus \{a\}$ be such that $p(v^*) = \max_{u \in V \setminus \{a\}} p(u)$ and break the tie arbitrarily. 
	
	For each $x \in S(v^*) \subset \B^c$, the order of the pole of the mean first passage time from $x$ to $\B$ is $p(x) \leftarrow p(v^*)$. 
	
	Update, for each $u \in V$ such that $(u,v^*) \in E$, $p(u) \leftarrow \max \{p(u),p(v^*)-\mathcal{K}_{uv^*}\}$.
	
	Update $V \leftarrow V \setminus \{v^*\}$ and $E \leftarrow \{(u,v) \in E: u \neq v^* \text{ and } v \neq v^*\}$.
	
	\normalsize

\subsection{Graphs for the algorithm to find the order of the MFPT}
\label{sec:construction}

Here we elaborate on the graphs that we use in the algorithm and the definitions that we gave in Section \ref{orderMFPT} and Section \ref{sec:Algstat}. While, in the algorithm statement, we used the same notation for the updated graphs as in the original graph, it will be clearer for the justification given in Section \ref{sec:justification} if we specify which copy of the graph we are looking at for each step. Accordingly, we provide a more detailed version of the definitions of these graphs and associated notation in this section.

\noindent \textbf{Step 1: Graph $G$}

For each $\eps \in (0,\eps_{\max})$ and $x \in \X$, the exponential parameter satisfies
\begin{equation*}
    q_x (\eps) = - Q_{x,x}(\eps) = \sum_{y \neq x \in X} Q_{x,y} (\eps) = \sum_{(x,y) \in E_0} Q_{x,y} (\eps) > 0.
\end{equation*} 
Since the order of $Q_{x,y} (\eps)$ is $k_{xy}$ for each $y \in \X$ such that $(x,y) \in E_0$, the order of $q_x (\eps)$ is $p_0(x)=\min \{k_{xy} : (x,y) \in E_0 \}$. For each $(x,y) \in E_0$, the transition probability $P_{x,y} (\eps)$ for the embedded discrete time Markov chain is
\begin{equation*}
P_{x,y} (\eps) = \frac{Q_{x,y} (\eps)}{q_x (\eps)},
\end{equation*} 
the order of which is 
\begin{equation}
\label{eq:Kxy}
\mathcal{K}_{xy} = k_{xy}-\min \{k_{xy} : (x,y) \in E_0 \} = k_{xy} - p_0(x).
\end{equation}
We start with a weighed graph $G=(V,E)$ where $V = \X$ and $E = E_0$. For each $u \in V$, the node weight of $u$ is $p_0(u)$, which is the \textbf{order of the pole of the expected sojourn time} at state $x$ until escape from $x$ for $X^\eps$. For each $(u,v) \in E$, the edge weight of $(u,v)$ is $\mathcal{K}_{uv}$, which is the \textbf{order of the transition probability} from $u$ to $v$.

\noindent \textbf{Step 2: Graph $G^{(0)}$}

If $x \in \B^c$ is such that $(x,y) \in E$ for some $y \in \B$, then the transition probability from $x$ to $\B$ is positive and is given by
\begin{equation*}
P_{x,\B} (\eps) = \sum_{\substack{y \in \B:\\ (x,y) \in E}} P_{x,y} (\eps).
\end{equation*}
For such $x$, the order of $P_{x,\B} (\eps)$ is 
\begin{equation}
\label{eq:Kxa}
\mathcal{K}_{xa} = \min \{\mathcal{K}_{xy} : y \in \B, (x,y) \in E\},
\end{equation}
where $\mathcal{K}_{xy}$ is the order of $P_{x,y} (\eps)$ for each $(x,y) \in E$.

Now, we are ready to specify the graph $G^{(0)} = (V^{(0)}, E^{(0)})$ which serves as the base case for Step 3. We group the nodes in $\B$ into a single node, denoted by $a$, so the set of nodes becomes $V^{(0)} = (V \setminus \B) \cup \{a\} = \B^c \cup \{a\}$. All of the edges starting from or going to a node in $\B$ are then removed. If there was an edge from $x \in \B^c$ to a node in $\B$, then we add back an edge $(x,a)$. We leave out all edges from a node in $\B$ to a node in $\B^c$, since we are interested in the mean first passage time to the set $\B$. The resulting edge set is $E^{(0)} = \{(u,v) \in E: u \in \B^c \text{ and } v \in \B^c \} \cup \{(w,a): (w,v) \in E \text{ for some } w \in \B^c \text{ and } v \in \B \}$. 

Let $S(u) = \{u\}$ for all $u \in V^{(0)} \setminus \{a\}$ and $S(a) = \B$. Note that $\{S(u): u \in V^{(0)}\}$ is a partition of the state space $\mathcal{X}$, denoting the grouping of nodes in $V^{(0)}$. For each $u \in V^{(0)} \setminus \{a\}$ and $x \in S(u) = \{u\}$, we define $p^x_0 (u)$ to be the \textbf{order of the pole of the expected sojourn time} in $S(u)$ before exiting $S(u)$ when starting the process at the state $x$. For each $(u,v) \in E^{(0)}$ and $x \in S(u) = \{u\}$, we define $\mathcal{K}_{uv}^x$ to be the \textbf{order of the probability of a transition} to $S(v)$ upon exiting from $S(u)$ when the process is started at the state $x$. For these terms, $p^x_0 (u) = p_0(u)$ and $\mathcal{K}_{uv}^x = \mathcal{K}_{uv}$, which is the base case for Lemma \ref{thm:MFPTAlgStep3}.

\noindent \textbf{Step 3: Graphs $\{G^{(N)}\}_{N=0}^M$}

In Step 3, we define a sequence of graphs $\{G^{(N)} = (V^{(N)},E^{(N)})\}_{N=0}^M$ recursively, where the exact value of $M \geq 0$ is not pre-determined and is only revealed when an exit condition for the recursion is satisfied. We know this recursion will end after a finite number of iterations because the number of nodes in $V^{(N)}$ is strictly decreasing with $N$. The weight $p_0$ of each node and the weight $\mathcal{K}$ of each edge are also defined iteratively, and each is defined only once. 

We have already defined $G^{(0)}$ in Step 2. Fix $N \in \{ 1, 2, \dots, M+1 \}$, where the value of $M<\infty$ is defined below. At the $N^{th}$ iteration, an edge $(u,v) \in E^{(N-1)}$ is called an \textbf{r-edge} if its edge weight $\mathcal{K}_{uv}$ is $0$; a directed path in $G^{(N-1)}$ is called an \textbf{r-path} if it consists of r-edges only. A set $C \subset V^{(N-1)}$ is called an \textbf{r-connected set} in $G^{(N-1)}$ if $|C|>1$ and there exists an r-path from $u$ to $v$ for any $u \neq v \in C$. Here we use the qualifier ``r'' to indicate that these edges, paths and cycles are ``regular''. If there is no r-connected set in $G^{(N-1)}$, the iteration stops. We set the value of $M$ to the first value of $N-1$ such that $G^{(N-1)}$ does not have any r-connected set. At that time point, the iteration stops and we move to Step 4 where $G^{(M)}$ will be the initial graph for Step 4. Otherwise, $N \in \{1,\dots,M\}$, and we let $C_N$ be an r-connected set in $G^{(N-1)}$, which is condensed to a new node $c_N$ in $G^{(N)}$. 
Then, we define the graph $G^{(N)} = (V^{(N)}, E^{(N)})$, where $V^{(N)} = (V^{(N-1)} \setminus C_N) \cup \{c_N\}$, and $E^{(N)} = \{(u,v) \in E^{(N-1)}: u \notin C_N \text{ and } v \notin C_N \} \cup \{(c_N,w): (u,w) \in E^{(N-1)} \text{ for some } u \in C_N \text{ and } w \notin C_N \} \cup \{(w,c_N): (w,v) \in E^{(N-1)} \text{ for some } w \notin C_N \text{ and } v \in C_N \}$. Let $S(c_N) = \cup_{u \in C_N} S(u)$. Note that $\{S(u): u \in V^{(N)}\}$ is again a partition of the state space $\mathcal{X}$, denoting the grouping of nodes in $V^{(N)}$. 

We define $p^x_0(c_N)$ to be the \textbf{order of the pole of the expected sojourn time} in $S(c_N)$ until the first exit from $S(c_N)$ when the process is started at the state $x \in S(c_N)$. In Lemma \ref{thm:MFPTAlgStep3}, we will show that the value of $p^x_0(c_N)$ is independent of the state $x \in S(c_N)$, and we define $p_0(c_N) = p^x_0(c_N)$. For each $w \in V^{(N-1)} \setminus C_N$ such that $(w,c_N) \in E^{(N)}$, define $\mathcal{K}_{wc_N}^x$ to be the \textbf{order of the probability of a transition} to $S(c_N)$ upon exiting from $S(w)$ when the process is started at the state $x \in S(w)$. 
In Lemma \ref{thm:MFPTAlgStep3}, we will show that the value of $\mathcal{K}_{wc_N}^x$ is independent of the state $x \in S(w)$, and $\mathcal{K}_{w c_N}^x = \mathcal{K}_{w c_N}$ where
\begin{equation}
\label{eq:Kwc}
\mathcal{K}_{w c_N} = \min\{ \mathcal{K}_{wv}: v \in C_N \text{ and } (w,v) \in E^{(N-1)}\}.
\end{equation}
For each $w \in V^{(N-1)} \setminus C_N$ such that $(c_N,w) \in E^{(N)}$, define $\mathcal{K}_{c_N w}^x$ to be the \textbf{order of the probability of a transition} to $S(w)$ upon exiting from $S(c_N)$ when the process is started at the state $x \in S(c_N)$. In Lemma \ref{thm:MFPTAlgStep3}, we will show that the value of $\mathcal{K}_{c_N w}^x$ is independent of the state $x \in S(c_N)$, and $\mathcal{K}_{c_N w}^x = \mathcal{K}_{c_N w}$ where
\begin{eqnarray}
&& \mathcal{K}_{c_N w} = \min\{\mathcal{K}_{uw}: u \in C_N \text{ and } (u,w) \in E^{(N-1)} \} \quad \quad \quad \quad \quad \quad \quad \nonumber\\
&& \quad \quad \quad \quad - \min\{\mathcal{K}_{uv}: u \in C_N, v \notin C_N \text{ and } (u,v) \in E^{(N-1)} \}. \label{eq:Kcw}
\end{eqnarray}
We note that for each $N \in \{1,\dots,M\}$, since there is no edge in $E^{(N-1)}$ that leads from $a$, the node $a$ is never part of any r-connected set, and so $a \in V^{(N)}$ and there is at least one other node in $V^{(N)}$ besides $a$. Also, the irreducibility of $X^\eps$ when $0 < \eps < \eps_0$ implies that there is a path from $x$ to $y$ in $G$ for each $x \in \X \setminus \B$ and $y \in \X$. This implies that if $u_N \neq v_N \in V^{(N)}$ for some $N \in \{0,1,\dots,M\}$ such that $x \in S(u_N)$ and $y \in S(v_N)$, then there is a path from $u_N$ to $v_N$ in $G^{(N)}$. Therefore, for each $N \in \{0,1,\dots,M\}$, there is always an outgoing edge from some $u' \in C_N$ ($u'$ cannot be $a$) to some $v' \in V^{(N-1)} \setminus C_N$ in $G^{(N-1)}$. In addition, as can be seen from the definition of the $\mathcal{K}$'s in \eqref{eq:Kxy}, \eqref{eq:Kxa}, \eqref{eq:Kwc} and \eqref{eq:Kcw}, for each $N=0,1,\dots,M$ and $u'' \in V^{(N)} \setminus \{a\}$, there exists an r-edge $(u'',v'') \in E^{(N)}$ for some $v'' \in V^{(N)}$.
For $G^{(M)}$, $|V^{(M)}|\geq 2$. Furthermorer, if we only look at r-edges (and ignore the other edges), $G^{(M)}$ is an acyclic graph (as it contains no r-connected set). It follows that the node $a$ is the only sink because for each $u \in V^{(M)} \setminus \{a\}$, there is an outgoing r-edge, and thus there is an r-path from $u$ to $a$ for each $u \in V^{(M)}\setminus \{a\}$.

\noindent \textbf{Step 4: Graphs $\{G^{(M,N)}\}_{N=0}^{|V^{(M)}|-1}$}

In Step 4, we define a sequence of graphs $\{G^{(M,N)} = (V^{(M,N)},E^{(M,N)})\}_{N=0}^{|V^{(M)}|-1}$ recursively, where $G^{(M,0)} = G^{(M)}$. In each iteration, the weight of one of the nodes in $V^{(M,N-1)} \subset V^{(M)}$ is finalized and determines the value of $p$ there, and the weights $p_{N-1}$ of other nodes in $V^{(M,N)}$ are updated to $p_{N}$.

Fix $N \in \{ 1, \ldots, |V^{(M)}|-1\}$. At the $N^{th}$ iteration, let $v_N \in V^{(M,N-1)} \setminus \{a\}$ be such that
\begin{equation}
\label{eq:pN-1vN}
p_{N-1} (v_N) = \max_{u \in V^{(M,N-1)} \setminus \{a\}} p_{N-1} (u) = \vcentcolon p(v_N),
\end{equation}
where we break the tie arbitrarily. 

Now, we define the graph $G^{(M,N)} = (V^{(M,N)},E^{(M,N)})$ where 
$V^{(M,N)} = V^{(M,N-1)} \setminus \{v_N \}$ and $E^{(M,N)} = \{(u,v) \in E^{(M,N-1)}: u \neq v_N \text{ and } v \neq v_N \}$. For each $u \in V^{(M,N)}$, let
\begin{equation}
\label{eq:Step4pN}
p_N (u) = 
\begin{cases}
\max \{p_{N-1} (u), p_{N-1} (v_N) - \mathcal{K}_{uv_N} \}, & \text{ for } (u,v_N) \in E^{(M,N-1)}, \\
p_{N-1} (u), & \text{ for } (u,v_N) \notin E^{(M,N-1)}.
\end{cases}
\end{equation}

In Theorem \ref{thm:MFPTAlgStep4}, we will show that for each $x \in S(v_N)$, the order of the pole of the mean first passage time from $x$ to the set $\B$ is $p(x)=p(v_N)$.

\subsection{Justification for the algorithm to find the order of the MFPT}
\label{sec:justification}

Recall that we defined $\{G^{(N)} \}_{N=0}^M = \{(V^{(N)}, E^{(N)})\}_{N=0}^M$ and $G^{(M,0)} = G^{(M)}$ in Section \ref{sec:construction}. Each $G^{(N)}$ defines a partition $\{S(u): u \in V^{(N)}\}$ of $\X$, which will be used in our proofs. Since $X^\eps$ is an irreducible continuous time Markov chain for all $\eps \in (0,\eps_0)$, each $G^{(N)}$ is weakly connected and has the property that any node in $V^{(N)} \setminus \{a\}$ has an out-going edge starting from the node.

In this section, we will provide the justification for the algorithm. Step 1 of the algorithm sets up the original continuous time Markov chain using a skeleton chain. Step 2 of the algorithm serves as the base case for Step 3, and Lemma \ref{thm:MFPTAlgStep3} justifies Steps 2 and 3. Theorem \ref{thm:MFPTAlgStep4} shows that Step 4 works, which gives our main result for the order of the pole of the mean first passage time from each state $x \in \B^c$ to $\B$.

We will start with Sections \ref{sec:bigtheta} and \ref{sec:stoppingtimes}, in which we describe in more detail the Big Theta notation used in this section and define some useful stopping times that will be used in our proof. 

\subsubsection{More on Big Theta notation}
\label{sec:bigtheta}

In Section \ref{orderMFPT}, we have defined orders for analytic functions using Big Theta notation. Here we give a few more definitions and remarks for inequalities involving the Big Theta notation, on how to compare the orders of analytic functions and on arithmetic for orders. These conventions streamline the proofs in the following subsections.

\begin{definition}
    Given $\eps_0 > 0$ and a function $f: (0,\eps_0) \rightarrow \R_{>0}$, we say $f \leq \Theta (\eps ^k)$ if there exist $k \in \Z$ and a strictly positive $M_f \in \R_{> 0}$ such that, for all $0 < \eps < \eps_0$, 
    \begin{equation*}
    f(\eps) \leq M_f \eps ^k.
    \end{equation*}
    We say $f \geq \Theta (\eps ^k)$ if there exist $k \in \Z$ and a strictly positive $m_f \in \R_{> 0}$ such that, for all $0 < \eps < \eps_0$, 
    \begin{equation*}
    f(\eps) \geq m_f \eps^k.
    \end{equation*}
\end{definition}

\begin{remark}
    Let $k,k_1,k_2 \in \Z$ and $k_1 \leq k \leq k_2$. If $f = \Theta (\eps^{k})$, then $f \leq \Theta (\eps^{k_1})$ and $f \geq \Theta (\eps^{k_2})$.
\end{remark}

\begin{remark}
    For functions $f$ and $g$ mapping $(0,\eps_0)$ into $\R_{>0}$, we write $f = g \cdot \Theta (\eps^{k})$ if $\frac{f}{g} = \Theta (\eps^{k})$, $f \leq g \cdot \Theta (\eps^{k})$ if $\frac{f}{g} \leq \Theta (\eps^{k})$, $f \geq g \cdot \Theta (\eps^{k})$ if $\frac{f}{g} \geq \Theta (\eps^{k})$. 
\end{remark}

\begin{lemma}
	\label{thm:OrderArithmetics}
	Let $k_1, k_2 \in \Z$, $\eps_0 > 0$ and $f,g: (0,\eps_0) \rightarrow \R_{>0}$. If $f = \Theta (\eps ^{k_1})$ and $g = \Theta (\eps ^{k_2})$, then 
	\begin{equation*}
 \frac{1}{f} =  \Theta(\eps^{-k_1}), \quad \quad \quad f + g = \Theta (\eps ^{\min\{k_1,k_2\}}), \quad \quad \quad f \cdot g = \Theta (\eps ^{k_1+k_2}),
 \end{equation*}
 \begin{equation*}
 \max \{f, g\} = \Theta (\eps^{\min \{k_1,k_2\}}), \quad \quad \quad \min \{f, g\} = \Theta (\eps^{\max \{k_1,k_2\}}).
  \end{equation*}
\end{lemma}
We leave the proof of Lemma \ref{thm:OrderArithmetics} to the reader.

\subsubsection{Stopping times $\tau_n^{\eps,N}$}
\label{sec:stoppingtimes}

For each graph $G^{(N)} = (V^{(N)}, E^{(N)})$, $N \in \{0,1,\dots,M\}$, recall that $\{S(u): u \in V^{(N)} \}$ is a partition of the state space $\X$. 
We define the series of stopping times $\{\tau_n^{\eps,N}\}_{n=0}^\infty$, which captures times of transitions of $X^\eps$ between sets in the partition $\{S(u): u \in V^{(N)} \}$ of $\X$. Formally, we let $\tau_0^{\eps,N} = 0$, and for $n = 1, 2, \ldots$, we successively define 
 \begin{equation*}
 \tau_{n}^{\eps,N} = \inf \left\{t \geq \tau_{n-1}^{\eps,N} : X^\eps (t) \notin S(v_{n-1}) \right\},
  \end{equation*}
where $v_{n-1}$ is the element in $V^{(N)}$ such that $X^\eps (\tau_{n-1}^{\eps,N}) \in S(v_{n-1})$.

\subsubsection{Justification for Step 3 of the algorithm}
\label{sec:HHThm4.5}

\begin{lemma}
\label{thm:MFPTAlgStep3}
\begin{enumerate}
	\item[(i)] For $N=0$ in Step 3, 
	
	\begin{enumerate}
		\item for each $u \in V^{(0)} \setminus \{a\}$ and $x \in S(u)$, $\E_x [\tau_1^{\eps,0}] = \Theta (\eps^{-p_0 (u)})$.
		\item for each $(u,v) \in E^{(0)}$ and $x \in S(u)$, $\PP_x [X^\eps (\tau_1^{\eps,0}) \in S(v)]	= \Theta (\eps^{\mathcal{K}_{uv}})$.
	\end{enumerate}
		
	\item[(ii)] For $N\in \{1,2,\dots,M\}$ in Step 3, let 
     \begin{equation*}
     k = \min\{\mathcal{K}_{uv}: u \in C_N, v \notin C_N \text{ and } (u,v) \in E^{(N-1)} \}.
      \end{equation*}
    (We note that $k$ depends on $N$ although we will not indicate that in the notation.)
    \begin{enumerate}
        \item For each $x \in S(c_N)$, $\E_x [\tau_1^{\eps,N}] = \Theta (\eps^{-p^x_0 (c_N)})$ where $p^x_0 (c_N) = p_0(c_N)$ and
        \begin{equation*}
        p_0 (c_N) = \max \{ p_0(u) : u \in C_N \} + k,
        \end{equation*}
        \item For each $x \in S(c_N)$ and $w \in V^{(N-1)} \setminus C_N$ such that $(u,w) \in E^{(N-1)}$ for some $u \in C_N$, $\PP_x [X^\eps (\tau_1^{\eps,N}) \in S(w)] = \Theta (\eps^{\mathcal{K}_{c_N w}^x})$ where $\mathcal{K}_{c_N w}^x = \mathcal{K}_{c_N w}$ and
        \begin{equation*}
        \mathcal{K}_{c_N w} = \min\{\mathcal{K}_{uw}: u \in C_N \text{ and } (u,w) \in E^{(N-1)} \} - k,
        \end{equation*}
        \item For each $x \in S(w)$ where $w \in V^{(N-1)} \setminus C_N$ is such that $(w,v) \in E^{(N-1)}$ for some $v \in C_N$, $\PP_x [X^\eps (\tau_1^{\eps,N}) \in S(c_N)]	= \Theta (\eps^{\mathcal{K}_{wc_N}^x})$ where $\mathcal{K}_{wc_N}^x = \mathcal{K}_{wc_N}$ and
        \begin{equation*}
        \mathcal{K}_{wc_N} = \min\{ \mathcal{K}_{wv}: v \in C_N \text{ and } (w,v) \in E^{(N-1)}\}.
        \end{equation*}
	\end{enumerate}
\end{enumerate}
\end{lemma}

\begin{proof}
    Our proof proceeds by induction. The base case ($N=0$) is established in Section \ref{sec:construction}.
    
    For fixed $1 \leq N \leq M$, assume that {\it (i) (a)-(b) and (ii) (a)-(c)} hold with $N$ replaced by $0,1,\dots,N-1$. We abbreviate $\tau_n^{\eps,N-1}$ as $\tau_n^\eps$ for $n = 0, 1, 2, \ldots$. Let \begin{equation*}
    \delta_{out} (C_N) = \{(u,v) \in E^{(N-1)} : u \in C_N \text{ and } v \notin C_N \}
    \end{equation*}
    denote all out-going boundary edges of $C_N$ so that 
    \begin{equation*}
    k = \min\{\mathcal{K}_{uv}: (u,v) \in \delta_{out} (C_N) \}.
    \end{equation*}
		
    First, consider the discrete time process $\{X^\eps (\tau_n^\eps)\}_{n=0}^\infty$, which is not necessarily a Markov process. We will derive a lower bound and an upper bound for 
    \begin{equation}
    \label{eqn:ExpectedVisitInCN}
    \E_x \left[ \sum_{n=0}^{\infty}  \one_{\{X^\eps (\tau_m^\eps) \in S(c_N) \text{ for } 0 \leq m \leq n \}} \right],
    \end{equation}
    which is the expected amount of time that $\{Y_n^\eps = X^\eps (\tau_n^\eps) \}_{n=0}^\infty$ spends in $S(c_N)$ before exiting from there, when started from a fixed state $x \in S(c_N)$.
		 
    For the lower bound, let
    \begin{equation*}
    \rho_1 = \max_{y \in S(c_N)} \PP_y [X^\eps (\tau_1^\eps) \notin S(c_N)],
    \end{equation*}
    the maximum over $y \in S(c_N)$ of the probability that for $X^\eps$ started at $y$, when $X^\eps$ exits $S(u_y)$, where $u_y \in C_N$ such that $y \in S(u_y)$, $X^\eps$ exits outside of $S(c_N)$. By the induction hypothesis, $\PP_y [X^\eps (\tau_1^\eps) \in S(v)] = \Theta (\eps^{\mathcal{K}_{uv}})$ for each $u,v$ such that $(u,v) \in \delta_{out} (C_N)$ and $y \in S(u)$. Thus, using Lemma \ref{thm:OrderArithmetics}, we have
    \begin{equation*}
    \rho_1 =\max_{y \in S(c_N)} \sum_{\substack{(u,v) \in \delta_{out} (C_N):\\ y \in S(u) }} \PP_y [X^\eps (\tau_1^\eps) \in S(v)] = \Theta (\eps^{\min \{ \mathcal{K}_{uv}: (u,v) \in \delta_{out} (C_N)\}}) = \Theta (\eps^{k}).
    \end{equation*}
    For $x \in S(c_N)$, let $\phi_{n}(x) = \PP_x [X^\eps (\tau^\eps_m) \in S(c_N) \text{ for } 0 \leq m \leq n]$ for $n = 0,1,2,\dots$. Then, $\phi_{0} (x) =1$, $\phi_{1} (x) =\PP_x [X^\eps (\tau^\eps_1) \in S(c_N)] \geq 1-\rho_1$, and by the strong Markov property, for $n \geq 2$, 
    \begin{eqnarray*}
        \phi_{n} (x) &=& \sum_{y \in S(c_N)} \PP_x [X^\eps (\tau^\eps_m) \in S(c_N) \text{ for } 0 \leq m \leq n-2; X^\eps (\tau^\eps_{n-1}) =y] \ \PP_y [X^\eps (\tau^\eps_1) \in S(c_N)] \\
        &\geq& \phi_{n-1} (x) (1-\rho_1).
    \end{eqnarray*}
    Hence, $\phi_{n} (x) \geq (1-\rho_1)^n$ for $n=0,1,2,\dots$. Then, for $x \in S(c_N)$,
    \begin{equation}
    \label{eq:rCycleLower}
    \E_x \left[ \sum_{n=0}^{\infty} \one_{\{X^\eps (\tau_m^\eps) \in S(c_N) \text{ for } 0 \leq m \leq n \}} \right] = \sum_{n=0}^{\infty} \phi_{n} (x) \geq \sum_{n=0}^{\infty} (1 - \rho_1)^n = \frac{1}{\rho_1}  = \Theta (\eps^{-k}).
    \end{equation}
    and so \eqref{eqn:ExpectedVisitInCN} is bounded below by $\Theta (\eps^{-k})$.

    For the upper bound, let $w_0 \in C_N$ be such that 
    \begin{equation*}
    \min \{ \mathcal{K}_{w_0 v}: (w_0,v) \in \delta_{out} (C_N)\} = \min \{ \mathcal{K}_{uv}: (u,v) \in \delta_{out} (C_N)\} = k.
    \end{equation*}
    Since the order of the probability $\PP_x [X^\eps (\tau_1^\eps) \notin S(c_N)]$ might equal $k'>k$ for some $w \neq w_0 \in C_N$ and $x \in S(w)$, such a smaller order probability of directly exiting $S(c_N)$ from $S(w)$ makes it seem possible that \eqref{eqn:ExpectedVisitInCN} could be $\Theta(\eps^{-k'})$ for some $k'>k$. Indeed, using a similar approach to the one we used for the lower bound, we can show that \eqref{eqn:ExpectedVisitInCN} is bounded above by $\Theta (\eps^{\max \{ \mathcal{K}_{uv}: (u,v) \in \delta_{out} (C_N)\}}) \geq \Theta(\eps^k)$. However, we would like a more stringent upper bound. To achieve this, we will show below that from $S(w)$, $X^\eps$ can exit $S(c_N)$ at least as quickly by means of a transition from $S(w)$ to $S(w_0)$ via the r-connected set and then from $S(w_0)$ to $V^{(N-1)} \setminus S(c_N)$.
    
    Let $\zeta_0^\eps = 0$, and for $n = 1, 2, \ldots$, we successively define 
    \begin{equation*}
    \eta_{n-1}^{\eps} = \inf \left\{t \geq \zeta_{n-1}^{\eps} : X^\eps (t) \notin S(v) \text{ where } v \in V^{(N-1)} \text{ and } X^\eps(\zeta_{n-1}^{\eps}) \in S(v) \right\},
    \end{equation*}
    \begin{equation*}
    \zeta_n^{\eps} = \inf \left\{t \geq \eta_{n-1}^{\eps} : X^\eps (t) \in S(w_0) \text{ or } X^\eps (t) \notin S(c_N) \right\}.
    \end{equation*}
    Note that $\{\zeta_n^\eps\}_{n=0}^\infty$ and $\{\eta_n^{\eps}\}_{n=0}^\infty$ depend on $N$. 
    Let
    \begin{equation*}
    \rho_2 = \min_{y \in S(w_0)} \PP_y [X^\eps (\zeta_1^\eps) \notin S(w_0)] = \min_{y \in S(w_0)} \PP_y [X^\eps (\zeta_1^\eps) \notin S(c_N)],
    \end{equation*}
    By the induction hypothesis, $\PP_y [X^\eps (\tau_1^\eps) \in S(v)] = \Theta (\eps^{\mathcal{K}_{w_0 v}})$ for each $v$ such that $(w_0,v) \in \delta_{out} (C_N)$ and $y \in S(w_0)$, and so
	\begin{eqnarray}
        \label{eqn:rho2}
        \rho_2 &\geq& \min_{y \in S(w_0)} \PP_y [X^\eps (\tau_1^\eps) \notin S(c_N)] = \Theta (\eps^{\min \{ \mathcal{K}_{w_0v}: (w_0,v) \in \delta_{out} (C_N)\}}) = \Theta (\eps^{k}),
	\end{eqnarray}
    where the inequality holds since starting from any $y \in S(w_0)$, if $X^\eps (\tau^\eps_1) \notin S(c_N)$, $X^\eps$ exits outside of $S(c_N)$ after leaving $S(w_0)$ and so $X^\eps (\zeta^\eps_1) \notin S(w_0)$. 

    For $x \in S(c_N)$, $\PP_x$-a.s., the sum
    \begin{equation}
    \label{eqn:NumVisitToW0}
    \sum_{n=0}^\infty \one_{\{X^\eps (\tau_m^\eps) \in S(c_N) \text{ for } 0 \leq m < n; X^\eps (\tau_n^\eps) \in S(w_0)\}}
    \end{equation}
    counts the number of distinct visits to $S(w_0)$, including the initial start there if $x \in S(w_0)$, before  $X^\eps$ escapes from $S(c_N)$. By the definition of the $\{\zeta_n^\eps\}_{n=0}^\infty$, $\PP_x$-a.s., the sum
    \begin{equation*}
    \one_{\{X^\eps (0) \in S(w_0) \}} + \sum_{n=1}^{\infty} \one_{\{X^\eps (\zeta_m^\eps) \in S(w_0) \text{ for } 1 \leq m \leq n \}}
    \end{equation*}
    counts the same quantity. 
    Thus, for $x \in S(w_0)$, using the strong Markov property and \eqref{eqn:rho2},
    \begin{eqnarray}
        \psi (x)&\vcentcolon =& \E_x \left[\sum_{n=0}^\infty \one_{\{X^\eps (\tau_m^\eps) \in S(c_N) \text{ for } 0 \leq m < n; X^\eps (\tau_n^\eps) \in S(w_0)\}} \right] \nonumber \\
        &=& 1 + \E_x \left[\sum_{n=1}^{\infty} \one_{\{X^\eps (\zeta_m^\eps) \in S(w_0) \text{ for } 1 \leq m \leq n \}} \right] \nonumber \\
        &=& 1 + \E_x \left[ \one_{\{X^\eps (\zeta^\eps_1) \in S(w_0)\}} \E_{X^\eps (\zeta^\eps_1)} \left[ \sum_{n=0}^{\infty} \one_{\{X^\eps (\zeta_m^\eps) \in S(w_0) \text{ for } 0 \leq m \leq n \}} \right] \right] \nonumber \\
        &\leq& 1 + (1 - \rho_2) \max_{y \in S(w_0)} \psi (y). \label{eq:PsixIneq}
    \end{eqnarray}
    Note that $\max_{y \in S(w_0)} \psi (y) < \infty$ because the state space is finite and $X^\eps$ is positive recurrent. Hence, by \eqref{eq:PsixIneq}, $\max_{y \in S(w_0)} \psi (y) \leq \frac{1}{\rho_2}$.
    Then, for $x \in S(c_N) \setminus S(w_0)$, by the strong Markov property,
    \begin{equation*}
    \psi (x) \leq \PP_x [X^\eps (\zeta^\eps_1) \in S(w_0)\}] \max_{y \in S(w_0)} \psi (y) \leq \frac{1}{\rho_2}.
    \end{equation*}
    Thus, for any $x \in S(c_N)$, 
    \begin{equation}
    \label{eq:rCycleUpper0}
        \E_x \left[ \sum_{n=0}^{\infty} \one_{\{X^\eps (\tau_m^\eps) \in S(c_N) \text{ for } 0 \leq m < n; X^\eps (\tau_n^\eps) \in S(w_0)\}} \right] \leq \frac{1}{\rho_2} = \Theta (\eps^{-k}). 
    \end{equation}
 
    Let $w_1, w_2 \in C_N$ be such that $(w_1,w_2) \in E^{(N-1)}$ and $\mathcal{K}_{w_1 w_2} = 0$. By the induction hypothesis, $\PP_y [X^\eps (\tau_1^\eps) \in S(w_2)] = \Theta(1)$ for all $y \in S(w_1)$. Then, for $x \in S(c_N)$,
	\small
	\begin{eqnarray}
		&& \E_x \left[ \sum_{n=0}^{\infty} \one_{\{X^\eps (\tau_m^\eps) \in S(c_N) \text{ for } 0 \leq m < n; X^\eps (\tau_n^\eps) \in S(w_2)\}} \right] \nonumber\\
		&\geq& \sum_{n=0}^{\infty} \PP_x [X^\eps (\tau_m^\eps) \in S(c_N) \text{ for } 0 \leq m < n; X^\eps (\tau_n^\eps) \in S(w_1), X^\eps (\tau_{n+1}) \in S(w_2)]\nonumber\\
		&=& \sum_{n=0}^{\infty} \sum_{y \in S(w_1)} \PP_x [X^\eps (\tau_m^\eps) \in S(c_N) \text{ for } 0 \leq m < n; X^\eps (\tau_{n}^\eps) = y] \cdot \PP_y [X^\eps (\tau_1^\eps) \in S(w_2)] \nonumber\\
		&\geq& \sum_{n=0}^{\infty} \PP_x [X^\eps (\tau_m^\eps) \in S(c_N) \text{ for } 0 \leq m < n; X^\eps (\tau_{n}^\eps) \in S(w_1)] \cdot \min_{y \in S(w_1)} \PP_y [X^\eps (\tau_1^\eps) \in S(w_2)] \nonumber \\
		&=& \E_x \left[ \sum_{n=0}^{\infty} \one_{\{X^\eps (\tau_m^\eps) \in S(c_N) \text{ for } 0 \leq m < n; X^\eps (\tau_n^\eps) \in S(w_1)\}} \right] \cdot \Theta (1) \label{eqn:rCycleIneq}.
 	\end{eqnarray}
 	\normalsize
	where the first equality holds from the strong Markov property of $X^\eps$. Since $C_N$ is an r-connected set, we can start from the node $w_1$ and the order inequality \eqref{eqn:rCycleIneq} can be passed from node to node in $C_N$ and back to the node $w_1$ ($w_0$ is included in the path) so that we will actually have equality in \eqref{eqn:rCycleIneq} and for all $v \in C_N$,
    \begin{eqnarray}
		&& \E_x \left[ \sum_{n=0}^{\infty} \one_{\{X^\eps (\tau_m^\eps) \in S(c_N) \text{ for } 0 \leq m < n; X^\eps (\tau_n^\eps) \in S(v)\}} \right] \nonumber\\
            &=& \E_x \left[ \sum_{n=0}^{\infty} \one_{\{X^\eps (\tau_m^\eps) \in S(c_N) \text{ for } 0 \leq m < n; X^\eps (\tau_n^\eps) \in S(w_0)\}} \right] \cdot \Theta (1).\label{eqn:rCycleEq}
    \end{eqnarray}
    Therefore, combining \eqref{eq:rCycleUpper0} and \eqref{eqn:rCycleEq}, and since there are only finitely many nodes in $C_N$, we can obtain by summing over $v \in C_N$ that \eqref{eqn:ExpectedVisitInCN} is bounded above by $\Theta (\eps^{-k})$. Combining with \eqref{eq:rCycleLower}, we have that, for $x \in S(c_N)$, \eqref{eqn:ExpectedVisitInCN} is $\Theta (\eps^{-k})$.
    Moreover, by \eqref{eqn:rCycleEq}, for each $x \in S(c_N)$ and each $v \in C_N$,
    \begin{equation}\label{eq:SojournInEachNode}
    \E_x \left[ \sum_{n=0}^{\infty} \one_{\{X^\eps (\tau_m^\eps) \in S(c_N) \text{ for } 0 \leq m < n; X^\eps (\tau_n^\eps) \in S(v)\}} \right] = \Theta(\eps^{-k}).
    \end{equation}
    
    To prove {\it (ii) (a)}, fix $x \in S(c_N)$. By the induction hypothesis, $\E_y [\tau_1^\eps] = \Theta (\eps^{-p_0 (u)})$ for each $y \in S(u)$ where $u \in C_N$. Thus, the expected sojourn time in $S(c_N)$ is 
	\small
	\begin{eqnarray}
		&& \E_x [\tau_1^{\eps,N}] = \E_x \left[ \sum_{n=0}^{\infty} \one_{\{X^\eps (\tau_m^\eps) \in S(c_N) \text{ for } 0 \leq m \leq n \}} \cdot (\tau_{n+1}^\eps - \tau_n^\eps) \right] \nonumber\\
		&=& \sum_{n=0}^{\infty} \sum_{u \in C_N} \sum_{y \in S(u)} \E_x [\one_{\{X^\eps (\tau_m^\eps) \in S(c_N) \text{ for } 0 \leq m < n; X^\eps (\tau_n^\eps) = y\}}] \cdot \E_y [\tau_1^\eps - \tau_0^\eps] \nonumber\\
		&=& \sum_{u \in C_N} \E_x \left[ \sum_{n=0}^{\infty} \one_{\{X^\eps (\tau_m^\eps) \in S(c_N) \text{ for } 0 \leq m < n; X^\eps (\tau_n^\eps) \in S(u) \}} \right] \cdot \Theta (\eps^{-p_0 (u)}) \nonumber\\
		&=& \sum_{u \in C_N} \Theta (\eps^{-k-p_0 (u)}) = \Theta (\eps^{-k-\max \{p_0 (u): u \in C_N\} }) = \Theta (\eps^{-p_0 (c_N)}). \nonumber
	\end{eqnarray}
	\normalsize
    where the first equality holds from the strong Markov property of $X^\eps$, we used \eqref{eq:SojournInEachNode} for the third equality, and we used Lemma \ref{thm:OrderArithmetics} for the fourth equality. 
			
	To prove {\it (ii) (b)}, fix $x \in S(c_N)$ and $w \in V^{(N-1)} \setminus C_N$. By the induction hypothesis, $\PP_y [X^\eps (\tau_1^\eps) \in S(w)] = \Theta (\eps^{\mathcal{K}_{uw}})$ for each $y \in S(u)$ where $u \in C_N$. Thus, starting from $x$, the probability of exiting $S(c_N)$ by means of a transition from a state in $S(c_N)$ to a state in $S(w)$ is given by
	\small
	\begin{eqnarray*}
		&& \PP_x [X^\eps (\tau_1^{\eps,N}) \in S(w)] = \sum_{n=0}^{\infty} \PP_x [X^\eps (\tau_m^\eps) \in S(c_N) \text{ for } 0 \leq m \leq n; X^\eps (\tau_{n+1}^\eps) \in S(w)] \\
		&=& \sum_{n=0}^{\infty} \sum_{u \in C_N} \sum_{y \in S(u)} \E_x [\one_{\{X^\eps (\tau_m^\eps) \in S(c_N) \text{ for } 0 \leq m < n; X^\eps (\tau_n^\eps) = y\}}] \cdot \PP_y [X^\eps (\tau_1^\eps) \in S(w)] \\
		&=& \sum_{\substack{u \in C_N:\\ (u,w) \in E^{(N-1)}}} \E_x \left[ \sum_{n=0}^{\infty} \one_{\{X^\eps (\tau_m^\eps) \in S(c_N) \text{ for } 0 \leq m < n; X^\eps (\tau_n^\eps) \in S(u) \}} \right] \cdot \Theta (\eps^{\mathcal{K}_{uw}}) \\
		&=& \sum_{\substack{u \in C_N:\\ (u,w) \in E^{(N-1)}}} \Theta (\eps^{- k + \mathcal{K}_{uw}}) = \Theta (\eps^{- k + \min\{ \mathcal{K}_{uw} : u \in C_N \text{ and } (u,w) \in E^{(N-1)} \} }) = \Theta (\eps^{\mathcal{K}_{c_N w}}),
	\end{eqnarray*}
	\normalsize
    where we used \eqref{eq:SojournInEachNode} for the third equality, and for the second equality, we used the fact that there must be an edge in $E^{(N-1)}$ between $u$ and $w$ if $\PP_y [X^\eps (\tau^\eps_1) \in S(w)] > 0$ for some and hence all $y \in S(u)$.
			
	To prove {\it (ii) (c)}, fix $x \in S(w)$ where $w \in V^{(N-1)} \setminus C_N$. By the induction hypothesis, $\PP_x [X^\eps (\tau_1^\eps) \in S(v)] = \Theta (\eps^{\mathcal{K}_{wv}})$ for each $v \in C_N$. Thus, starting from $x$, the probability of entering $S(c_N)$ by means of a transition from a state in $S(w)$ to a state in $S(c_N)$ is 
	\begin{eqnarray}
		&& \PP_x [X^\eps (\tau_1^{\eps,N}) \in S(c_N)] = \sum_{v \in C_N} \PP_x [X^\eps (\tau_1^\eps) \in S(v)] \nonumber\\
		&=& \sum_{\substack{v \in C_N:\\ (w,v) \in E^{(N-1)}}} \Theta (\eps^{\mathcal{K}_{wv}}) = \Theta (\eps^{\min\{ \mathcal{K}_{wv}: v \in C_N \text{ and } (w,v) \in E^{(N-1)}\} }) = \Theta(\eps^{\mathcal{K}_{wc_N}}). \nonumber
	\end{eqnarray}
 \
\end{proof}

\subsubsection{Justification for Step 4 of the algorithm}
\label{sec:HHLem4.6}

\begin{lemma}
\label{thm:Step4Prep}
    Fix $w \in V^{(M,0)} \setminus \{a\}$. 
    Let $\tau_n^\eps = \tau_n^{\eps,M}$, for $n = 0, 1, 2, \ldots$, as defined in Section \ref{sec:stoppingtimes}.
    Then, starting from $x \in S(w)$, the expected number of distinct visits to $S(w)$, including the initial start there, before $X^\eps$ enters $S(a)$ is
    \begin{equation*}
    \E_x \left[ \sum_{n=0}^{\infty} \one_{\{X^\eps (\tau_m^\eps) \notin S(a) \text{ for } 0 \leq m < n; X^\eps (\tau_n^{\eps}) \in S(w) \}} \right] = \Theta(1).
    \end{equation*}
\end{lemma}

\begin{proof}
Let $\zeta_0^{\eps} = 0$, and for $n = 1, 2, \ldots$, successively define 
\begin{equation*}
\eta_{n-1}^{\eps} = \inf \left\{t \geq \zeta_{n-1}^{\eps} : X^\eps (t) \notin S(v) \text{ where } v \in V^{(M,0)} \text{ and } X^\eps(\zeta_{n-1}^{\eps}) \in S(v)  \right\},
\end{equation*}
\begin{equation*}
\zeta_n^{\eps} = \inf \left\{t \geq \eta_{n-1}^{\eps} : X^\eps (t) \in S(w) \cup S(a) \right\}.
\end{equation*}
Note that for $x \in S(w)$, $\PP_x$-a.s., 
\begin{equation}
\label{eqn:SojournAtW}
\sum_{n=0}^{\infty} \one_{\{X^\eps (\tau_m^\eps) \notin S(a) \text{ for } 0 \leq m < n; X^\eps (\tau_n^{\eps}) \in S(w) \}} 
 = \sum_{n=0}^{\infty} \one_{\{X^\eps (\zeta_m^{\eps}) \in S(w) \text{ for } 0 \leq m \leq n \}},
\end{equation}
since they both count the number of distinct visits to $S(w)$, including the initial start there, before $X^\eps$ enters $S(a)$.

Recall from Section \ref{sec:construction} that for each $u \in V^{(M,0)}\setminus \{a\}$, there is an r-path from $u$ to $a$. Let such an r-path from $w$ to $a$ be $w \rightarrow w_1 \ldots \rightarrow w_d \rightarrow a$ where $w, w_1, \ldots, w_d, a$ are distinct. By definition, an edge $(u,v) \in V^{(M,0)}$ is an r-edge implies that $\PP_z [X^\eps (\tau_1^\eps) \in S(v)] = \Theta (1)$ for all $z \in S(u)$. Thus, for any $y \in S(w)$, using the strong Markov property of $X^\eps$, we have
\small
\begin{eqnarray}
\Theta(1) = 1 &\geq& \PP_y [X^\eps (\zeta_1^\eps) \notin S(w)] = \PP_y [X^\eps (\zeta_1^\eps) \in S(a)] \nonumber\\
&\geq& \PP_y [X^\eps (\tau_1^\eps) \in S(w_1), \ldots, X^\eps (\tau_d^\eps) \in S(w_d), X^\eps (\tau_{d+1}^\eps) \in S(a)] \nonumber\\
&=& \sum_{z \in S(w_d)} \PP_x [X^\eps (\tau_1^\eps) \in S(w_1), \ldots, X^\eps (\tau_d^\eps) \in S(w_d), X^\eps (\tau_d^\eps) = z] \cdot \PP_z [X^\eps (\tau_1^\eps) \in S(a)] \nonumber\\
&=& \PP_y [X^\eps (\tau_1^\eps) \in S(w_1), \ldots, X^\eps (\tau_d^\eps) \in S(w_d)] \cdot \Theta(1) = \ldots = \Theta(1) \cdot \ldots \cdot \Theta(1) = \Theta(1). \nonumber
\end{eqnarray}
\normalsize

Using a similar approach to that used in Section \ref{sec:HHThm4.5}, we can show that for $x \in S(w)$,
\begin{equation*}
\E_x \left[ \sum_{n=0}^{\infty} \one_{\{X^\eps (\zeta_m^{\eps}) \in S(w) \text{ for } 0 \leq m \leq n \}} \right] \geq \frac{1}{\max_{y \in S(w)} \PP_y [X^\eps (\zeta_1^{\eps}) \notin S(w)]} = \Theta(1),
\end{equation*}
and
\begin{equation*}
\E_x \left[ \sum_{n=0}^{\infty} \one_{\{X^\eps (\zeta_m^{\eps}) \in S(w) \text{ for } 0 \leq m \leq n \}} \right] \leq \frac{1}{\min_{y \in S(w)} \PP_y [X^\eps (\zeta_1^{\eps}) \notin S(w)]} = \Theta(1).
\end{equation*}
Combining these inequalities with \eqref{eqn:SojournAtW} yields the desired result.
\end{proof}

\begin{theorem}
\label{thm:MFPTAlgStep4}
    Let $\tau_\B^\eps = \inf \{t \geq 0: X^\eps (t) \in \B \}$ be the first passage time to $\B$ for $X^\eps$. For each $N = 1, \dots, |V^{(M)}|-1$ and $x \in S(v_N)$, we have
    \begin{equation}
    \label{eq:orderMFPT}
    \E_x [\tau_\B^\eps] = \Theta (\eps^{-p(v_N)}).
    \end{equation}
\end{theorem}

\begin{proof}
    It suffices by iteration to prove that for each fixed $1 \leq N \leq |V^{(M)}|-1$, if 
    \begin{equation}
    \label{eq:orderMFPTassumption}
    \E_y [\tau_\B^\eps] = \Theta (\eps^{-p(v_k)}) \text{ for all } y \in S(v_k) \text{ and } 1 \leq k \leq N-1,
    \end{equation}
    then \eqref{eq:orderMFPT} holds for all $x \in S(v_N)$. By convention, \eqref{eq:orderMFPTassumption} holds automatically for $N=1$. 
    
    For the iteration step, fix $1 \leq N \leq |V^{(M)}|-1$, and assume that \eqref{eq:orderMFPTassumption} holds. If $N=1$, let $A=V^{(M,0)} \setminus \{a\}$,  and if $N>1$, let $A = V^{(M,0)} \setminus \{v_1, \ldots, v_{N-1}, a\}$. Recall that for $w \in A$ and $1 \leq k \leq N-1$, we have $w \in V^{(M,k)}$, and so by \eqref{eq:Step4pN},
    \begin{equation}
    \label{eq:pkw}
    p_{k} (w) = 
    \begin{cases}
    \max \{p_{k-1}(w), p_{k-1} (v_k) - \mathcal{K}_{w v_k} \} & \text{ for } (w,v_k) \in E^{(M,k-1)}, \\
    p_{k-1}(w), & \text{ for } (w,v_k) \notin E^{(M,k-1)}.
    \end{cases}
    \end{equation}
    Note that $(w,v_k) \in E^{(M,k-1)}$ if and only if $(w,v_k) \in E^{(M,0)}$. Since $p_{k-1} (v_k) = \vcentcolon p (v_k)$ for $1 \leq k \leq N-1$, by iterating \eqref{eq:pkw}, we can obtain
    \begin{equation}
    \label{eq:pw}
    p_{N-1} (w) = \max \{p_0(w), \max \{ p(v_k) - \mathcal{K}_{w v_k}: 1 \leq k \leq N-1 \text{ and } (w,v_k) \in E^{(M,0)} \} \},
    \end{equation}
    where we make the convention that a maximum over an empty set is $-\infty$. In particular, 
    since $v_N \in A$, we have
    \small
    \begin{equation}
    \label{eq:pvN}
    p(v_N) \vcentcolon= p_{N-1} (v_N) = \max \{p_0(v_N), \max \{ p(v_k) - \mathcal{K}_{v_N v_k}: 1 \leq k \leq N-1 \text{ and } (v_N,v_k) \in E^{(M,0)} \} \}.
    \end{equation}
    \normalsize
    
    Fix $x \in S(v_N)$. We will derive a lower bound and an upper bound for $\E_x [\tau_\B^\eps]$.
    For the lower bound, let $\tau^\eps = \inf \{ t \geq 0: X^\eps (t) \notin S(v_N)\}$.
    Recall that $\PP_x$-a.s., $\tau^\eps = \tau_1^{\eps,M}$ as defined in Section \ref{sec:stoppingtimes}. By Lemma \ref{thm:MFPTAlgStep3}, for each $y \in S(w)$ where $w \in V^{(M,0)}=V^{(M)}$ is such that $(v_N,w) \in E^{(M,0)}=E^{(M)}$, $\PP_x [X^\eps (\tau^\eps) = y] = \Theta (\eps^{\mathcal{K}_{v_N w}})$ and $\E_x [\tau^\eps] = \Theta (\eps^{-p_0 (v_N)})$. By first step analysis,
    \begin{eqnarray}
	\E_x [\tau_\B^\eps] &=& \E_x [\tau^\eps] + \sum_{(v_N,w) \in E^{(M,0)}} \sum_{y \in S(w)} \PP_x [X^\eps (\tau^\eps) = y] \cdot \E_y [\tau_\B^\eps] \nonumber \\
	&\geq& \E_x [\tau^\eps] + \sum_{\substack{1 \leq k \leq N-1: \\(v_N,v_k) \in E^{(M,0)}}} \sum_{y \in S(v_k)} \PP_x [X^\eps (\tau^\eps) = y] \cdot \E_y [\tau_\B^\eps] \nonumber \\
	&=& \Theta (\eps^{-p_0 (v_N)}) + \sum_{\substack{1 \leq k \leq N-1: \\(v_N,v_k) \in E^{(M,0)}}} \Theta (\eps^{ \mathcal{K}_{v_N v_k}}) \cdot \Theta (\eps^{-p (v_k)}) = \Theta (\eps^{-p_{N-1} (v_N)}). \quad \label{eqn:MFPTLowerBound}
    \end{eqnarray}
    where we used \eqref{eq:orderMFPTassumption} in the second last equality, and used Lemma \ref{thm:OrderArithmetics} and \eqref{eq:pvN} for the last equality.

    For the upper bound, let $\eta^\eps = \inf \{t \geq 0: X^\eps (t) \notin \bigcup_{u \in A} S(u) \}$. 
    Let $\tau_n^\eps = \tau_n^{\eps,M}$, for $n = 0, 1, 2, \ldots$, as defined in Section \ref{sec:stoppingtimes}. Then, using Lemma \ref{thm:Step4Prep} and the strong Markov property, for $w \in A \subset V^{(M,0)} \setminus \{a\}$,
    \begin{eqnarray}
    && \E_x \left[ \sum_{n=0}^{\infty} \one_{\{X^\eps (\tau_m^\eps) \in \bigcup_{u \in A} S(u) \text{ for } 0 \leq m < n; X^\eps (\tau_n^{\eps}) \in S(w) \}} \right] \nonumber\\
    &\leq& \PP_x [ X^\eps (\zeta^\eps) \in S(w)] \max_{y \in S(w)} \E_y \left[ \sum_{n=0}^{\infty} \one_{\{X^\eps (\tau_m^\eps) \notin S(a) \text{ for } 0 \leq m < n; X^\eps (\tau_n^{\eps}) \in S(w) \}} \right] \leq \Theta (1).  \label{eq:numVisitTowInA} 
    \end{eqnarray}
    where $\zeta^\eps = \inf \{t \geq 0: X^\eps (t) \in S(w) \}$.

    For $1 \leq k \leq N-1$ such that there exists $w \in A$ where $(w,v_k) \in E^{(M,0)}$,
	\small
	\begin{eqnarray}
		&& \PP_x [X^\eps (\eta^\eps) \in S(v_k) ] = \sum_{n=0}^{\infty} \PP_x [X^\eps (\tau_m^\eps) \in \bigcup_{u \in A} S(u) \text{ for } 0 \leq m \leq n; X^\eps (\tau_{n+1}^{\eps}) \in S(v_k)] \nonumber \\
		&=& \sum_{n=0}^{\infty} \sum_{w \in A} \sum_{y \in S(w)} \PP_x [X^\eps (\tau_m^\eps) \in \bigcup_{u \in A} S(u) \text{ for } 0 \leq m < n; X^\eps (\tau_n^{\eps}) = y] \cdot \PP_y [X^\eps (\tau_1^\eps) \in S(v_k)] \nonumber \\
		&=& \sum_{\substack{w \in A:\\(w,v_k) \in E^{(M,0)}}} \E_x \left[ \sum_{n=0}^{\infty} \one_{\{X^\eps (\tau_m^\eps) \in \bigcup_{u \in A} S(u) \text{ for } 0 \leq m < n; X^\eps (\tau_n^{\eps}) \in S(w) \}} \right] \cdot \Theta (\eps^{\mathcal{K}_{wv_k}}) \nonumber \\
		&\leq& \sum_{\substack{w \in A:\\(w,v_k) \in E^{(M,0)}}} \Theta (1) \cdot \Theta (\eps^{\mathcal{K}_{wv_k}}), \label{eq:PPeta}
	\end{eqnarray}
	\normalsize
    where the second equality holds from strong Markov property of $X^\eps$, the third equality uses Lemma \ref{thm:MFPTAlgStep3}, and we used \eqref{eq:numVisitTowInA} for the last inequality. 
    Using Lemma \ref{thm:MFPTAlgStep3} and \eqref{eq:numVisitTowInA},
	\small
	\begin{eqnarray}
		&& \E_x [\eta^\eps] = \E_x \left[ \sum_{n=0}^{\infty} \one_{\{X^\eps (\tau_m^\eps) \in \bigcup_{u \in A} S(u) \text{ for } 0 \leq m \leq n\}} \cdot (\tau_{n+1}^\eps - \tau_n^\eps) \right] \nonumber \\
		&=& \sum_{n=0}^{\infty} \sum_{w \in A} \sum_{y \in S(w)} \PP_x [X^\eps (\tau_m^\eps) \in \bigcup_{u \in A} S(u) \text{ for } 0 \leq m < n; X^\eps (\tau_n^{\eps}) = y] \cdot \E_y [\tau_1^\eps - \tau_0^\eps] \nonumber\\ 
		&=& \sum_{w \in A} \E_x \left[ \sum_{n=0}^{\infty} \one_{\{X^\eps (\tau_m^\eps) \in \bigcup_{u \in A} S(u) \text{ for } 0 \leq m < n; X^\eps (\tau_n^{\eps}) \in S(w) \}} \right] \cdot \Theta (\eps^{-p_0 (w)}) \nonumber \\
		&\leq& \sum_{w \in A} \Theta (1) \cdot \Theta (\eps^{-p_0 (w)}) \leq \Theta (\eps ^{-p_{N-1} (v_N)}), \label{eq:Eeta}
	\end{eqnarray}
	\normalsize
	where we have used \eqref{eq:pw} and \eqref{eq:pN-1vN} to conclude that $p_0 (w) \leq p_{N-1} (w) \leq p_{N-1} (v_N)$ for all $w \in A$. 
    Therefore, using first step analysis, we have
    \begin{eqnarray}
	\E_x [\tau_\B^\eps]  &=& \E_x [\eta^\eps] + \sum_{1 \leq k \leq N-1} \sum_{y \in S(v_k)} \PP_x [X^\eps (\eta^\eps) = y] \cdot \E_y [\tau_\B^\eps] \nonumber\\
	&=& \E_x [\eta^\eps] + \sum_{1 \leq k \leq N-1} \PP_x [X^\eps (\eta^\eps) \in S(v_k)] \cdot \Theta (\eps^{-p(v_k)}) \nonumber\\
	&\leq& \Theta (\eps ^{-p_{N-1} (v_N)}) + \sum_{\substack{w \in A, \\ 1 \leq k \leq N-1: \\(w,v_k) \in E^{(M,0)}}}  \Theta (\eps^{\mathcal{K}_{wv_k}}) \cdot \Theta (\eps^{-p(v_k)}) \nonumber\\
	&\leq& \Theta (\eps ^{-p_{N-1} (v_N)}) + \Theta (\eps^{- \max\{p_{N-1} (w): w \in A \}} ) = \Theta (\eps^{- p_{N-1} (v_N)}), \quad \quad \quad \label{eqn:MFPTUpperBound}
    \end{eqnarray}
    where we used \eqref{eq:orderMFPTassumption} for the second equality, \eqref{eq:PPeta} and \eqref{eq:Eeta} for the first inequality, and \eqref{eq:pw} and Lemma \ref{thm:OrderArithmetics} for the second inequality.
    
    By \eqref{eqn:MFPTLowerBound} and \eqref{eqn:MFPTUpperBound}, we conclude that $\E_x [\tau_\B^\eps] = \Theta (\eps^{-p_{N-1}(v_N)}) =  \Theta (\eps^{-p(v_N)})$.
\end{proof}

\subsection{Application of the algorithm to the 2D, 3D and 4D models}
\label{sec:alg_3D4D}
\begin{figure}[H]
		\centering \includegraphics[width=\textwidth]{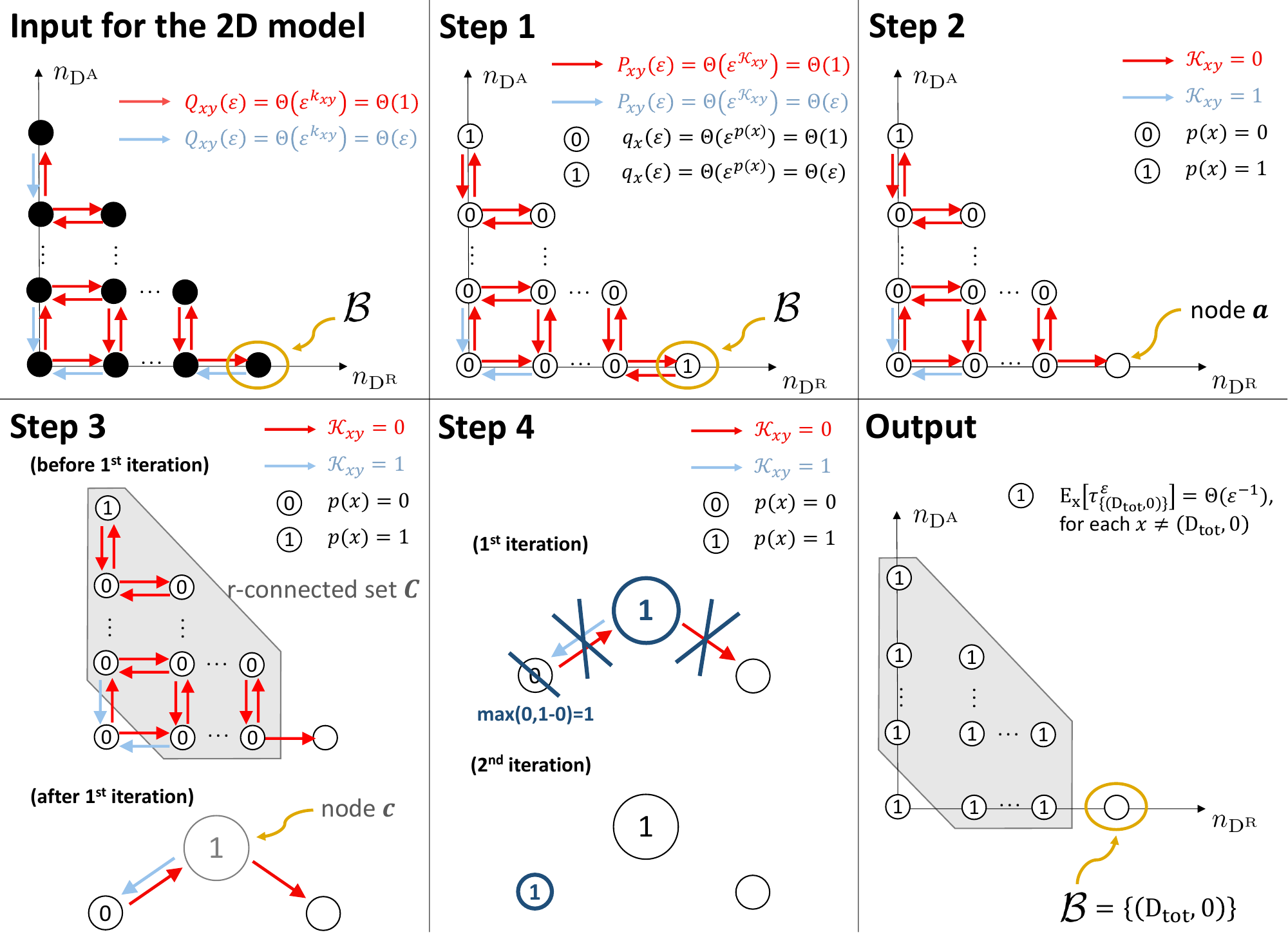}
		\caption{\small { \bf Key steps of the algorithm for the 2D model.} 
		}
		\label{fig:algorithm_2D}
	\end{figure}

\begin{figure}[H]
	\centering \includegraphics[width=\textwidth]{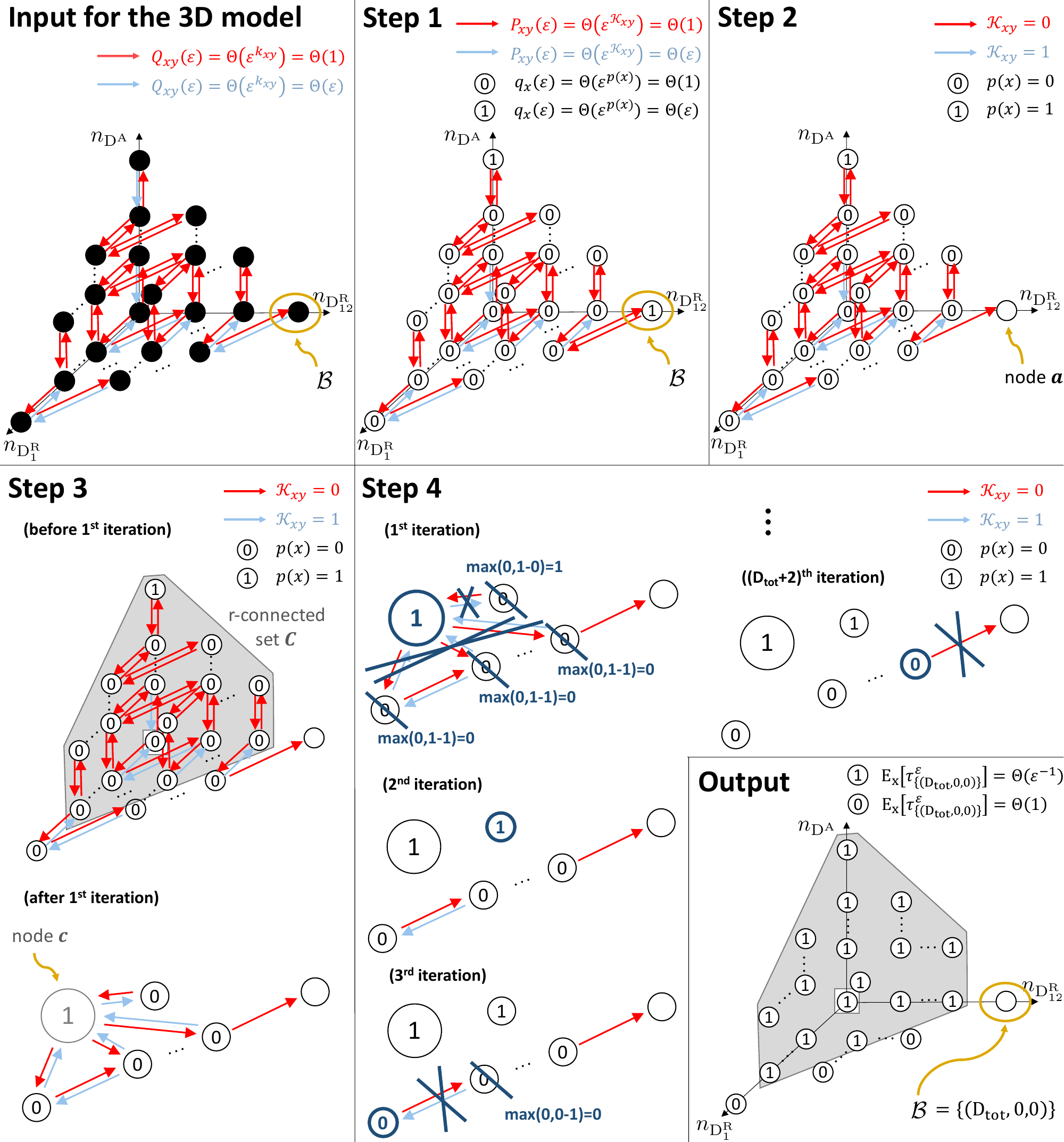}
	\caption{{\bf Key steps of the algorithm for the 3D model (from the fully active state to the fully repressed state).}}
    \label{fig:algorithm_3DAtoR}
\end{figure}

\begin{figure}[H]
	\centering \includegraphics[width=\textwidth]{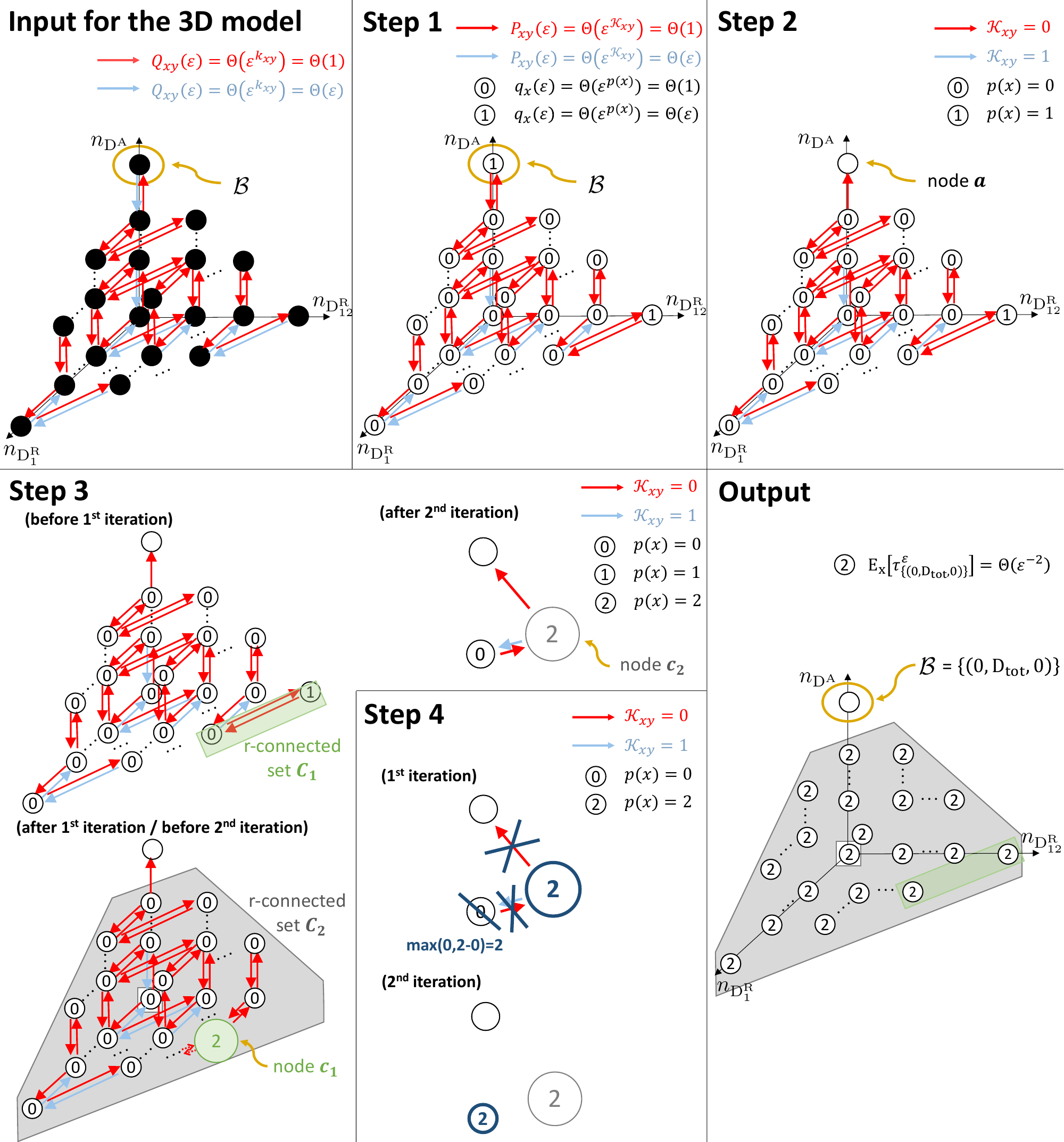}
	\caption{{\bf Key steps of the algorithm for the 3D model (from the fully repressed state to the fully active state).}}
	\label{fig:algorithm_3DRtoA}
\end{figure}

\begin{figure}[H]
	\centering \includegraphics[width=\textwidth]{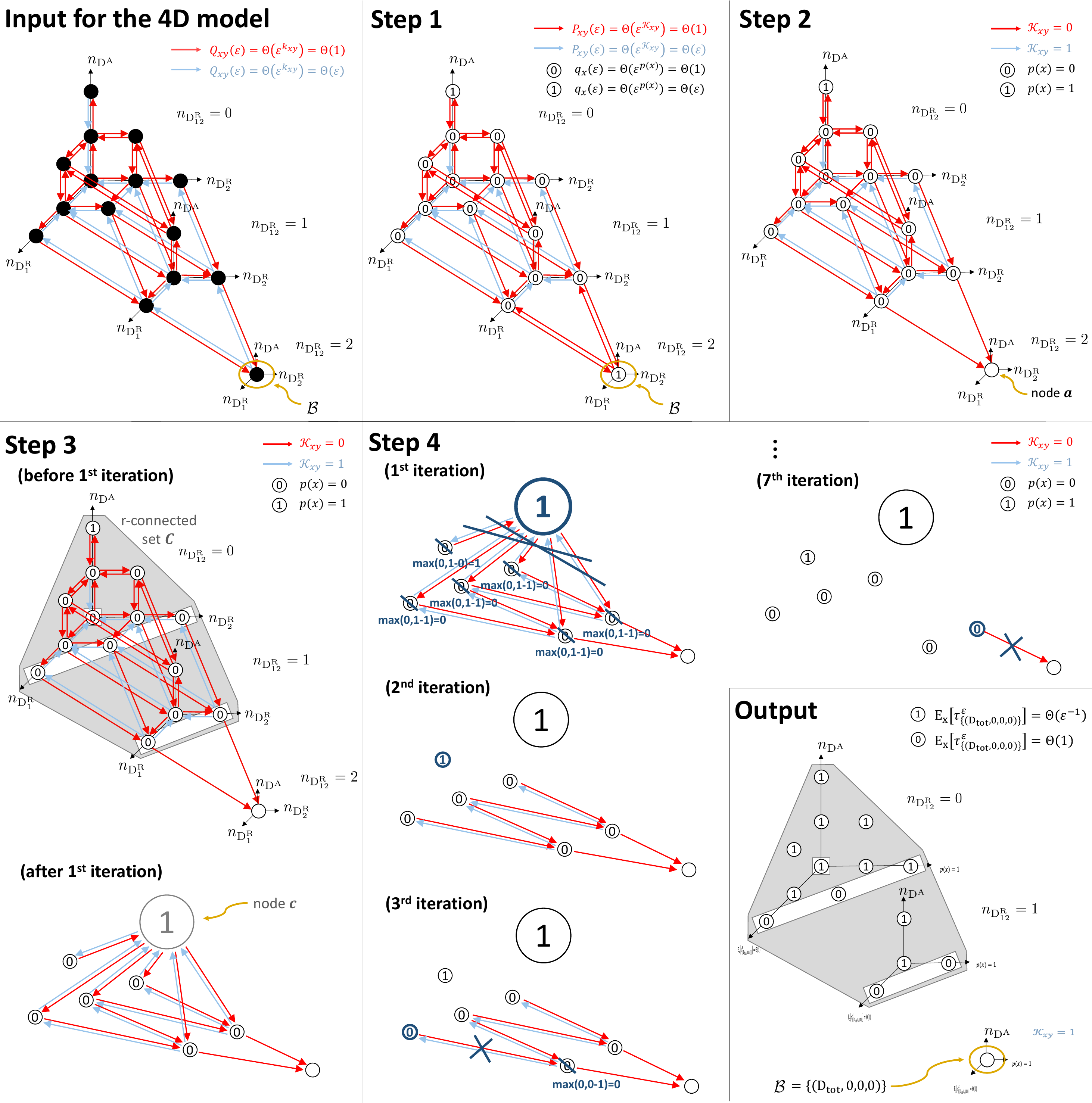}
	\caption{{\bf Key steps of the algorithm for the 4D model (from the fully active state to the fully repressed state).}}
	\label{fig:algorithm_4DAtoR}
\end{figure}

\begin{figure}[H]
	\centering \includegraphics[width=\textwidth]{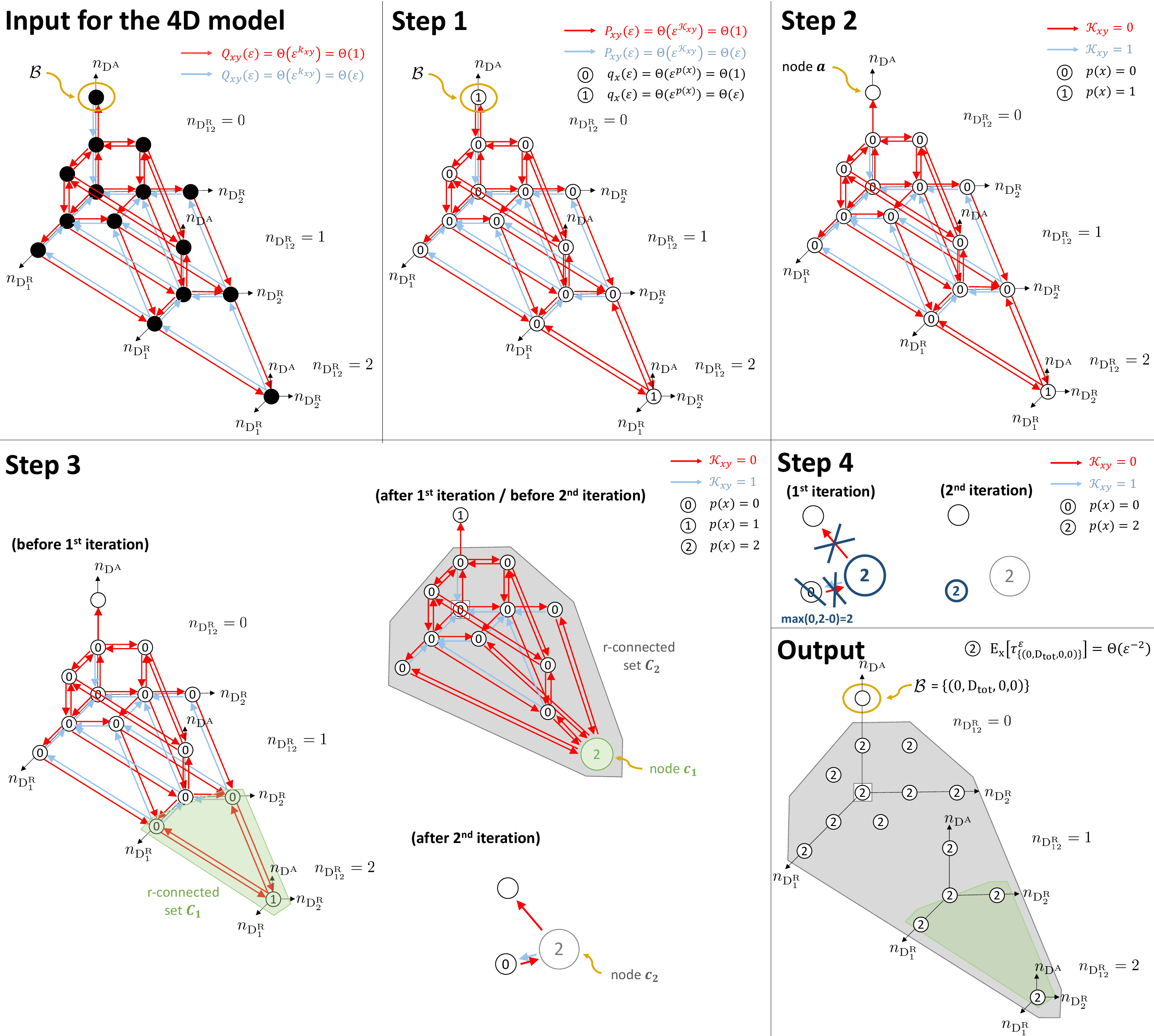}
	\caption{{\bf Key steps of the algorithm for the 4D model (from the fully repressed state to the fully active state).}}
	\label{fig:algorithm_4DRtoA}
\end{figure}

The algorithm is described in Section \ref{orderMFPT}, and it finds the order of the pole of the mean first passage time to $\emptyset \neq \B \subset \X$ from each state in $\B^c$. In this section, we will apply the algorithm to the {\color{red} 2D, } 3D and 4D models and find the order of the poles of the mean first passage times of interest to the fully repressed state and the fully active state (Figure \ref{fig:algorithm_2D} -- \ref{fig:algorithm_4DRtoA}). For each figure, the ``Input'' panel shows the order of each of the non-zero off-diagonal entries in $Q(\eps)$ and the set $\B$ which contains a single state, which is either the fully repressed state or the fully active state. The orders of these non-zero entries in $Q(\eps)$ are represented by colored arrows in the graph (red for order $0$ and blue for order $1$). Step 1 transforms the orders in the infinitesimal generator $Q(\eps)$ into orders for the transition matrix $P(\eps)$ and the exponential parameters $q(\eps)$ to give an equivalent construction for the continuous time Markov chain. The orders of the non-zero entries in $P(\eps)$ are given by $\mathcal{K}$ and represented by colored arrows in the graph, and the number in the circle at a state $x \in \B^c$ is the order of the pole $p(x)$ of $\frac{1}{q_x(\eps)}$ (the mean sojourn time at the state $x$). In Step 2, the set $\B$ contains only one state and is just relabeled as the node $a$. All transitions from $a$ to $\B^c$ are then removed. While the Input, Step 1 and Step 2 are universal across all the figures in this section, we explain the Step 3, Step 4 and Output panels separately for each application below since they are more distinct.

\noindent {\bf 2D model (from the fully active state to the fully repressed state)}: see Figure \ref{fig:algorithm_2D}. The explanation of the panels for Input, Step 1 and Step 2 is given above with $\B = \{(\Dtot,0)^T\}$. 
Step 3 for the 2D model involves only one iteration, where the collection of all nodes except the node $a$ and the origin $0$ (called an r-connected set $C$) is condensed to a single node $c$, and the order of the pole at $c$ is $p(c) = \max_{u \in C} p (u) + \min \{\mathcal{K}_{uv}:\: u \in C, v\notin C \text{ and } (u,v) \in E \} = 1 + 0 = 1$, where $E$ denotes the edge set of the graph in Step 3 before the first iteration. Moreover, $\mathcal{K}_{c0} = \min \{\mathcal{K}_{u0}:\: u \in C \text{ and } (u,0) \in E \} - \min \{\mathcal{K}_{uv}:\: u \in C, v\notin C \text{ and } (u,v) \in E \} = 1-0=1$, $\mathcal{K}_{ca} = \min \{\mathcal{K}_{ua}:\: u \in C \text{ and } (u,a) \in E \} - \min \{\mathcal{K}_{uv}:\: u \in C, v\notin C \text{ and } (u,v) \in E \} = 0-0=0$, and $\mathcal{K}_{0c} = \min \{\mathcal{K}_{0v}:\: v \in C \text{ and } (0,v) \in E \} = 0$. Step 4 involves two iterations. In the first iteration, we fix the node with the largest value of $p$, which is $c$ in our case. At any node other than $a$ that is connected to $c$ (i.e., the origin $0$), the value of $p$ is updated to $p(0) = \max \{ p(0), p(c)-\mathcal{K}_{0c}\} = \max \{0,1-0\} = 1$, and then any edges leading to or from $c$ are removed. In the second iteration, of the remaining nodes, we fix the node with the largest value of $p$, which is the origin. When all of the nodes other than $a$ have been fixed, the order of the pole of the mean first passage time from each state in $\B^c$ to $\B$ is given by the fixed value of the node to which the state belongs.

\noindent {\bf 2D model (from the fully repressed state to the fully active state)} Because of the symmetry in the input graph in Figure \ref{fig:algorithm_2D}, the orders of the poles of the mean first passage times to the fully repressed state can be obtained in the same way as above.

\noindent {\bf 3D model (from the fully active state to the fully repressed state)}: see Figure \ref{fig:algorithm_3DAtoR}.
The explanation of the panels for Input, Step 1 and Step 2 is given above with $\B = \{(\Dtot,0,0)^T\}$. A state represents $(n_\mathrm{D_{12}^R},n_\mathrm{D^A},n_\mathrm{D_{1}^R})^T$. Step 3 involves only one iteration, where the collection of all nodes except for $(0,0,0)^T$, $(0,0,\Dtot)^T$, $(1,0,\Dtot-1)^T$, $(2,0,\Dtot-2)^T$,$\ldots$, $(\Dtot-2,0,2)^T$, and $(\Dtot-1,0,1)^T$ (called an r-connected set $C$) is condensed to a single node $c$. The order of the pole of the sojourn time at $C$ is $p(c) = \max_{u \in C} p (u) + \min \{\mathcal{K}_{uv}: u \in C, v\notin C \text{ and } (u,v) \in E \} = 1 + 0 = 1$, where $E$ denotes the edge set of the graph in Step 3 before the first iteration. Moreover, 
$\mathcal{K}_{c,(0,0,0)^T} = \min \{\mathcal{K}_{u,(0,0,0)^T}: u \in C \text{ and } (u,(0,0,0)^T) \in E \} - \min \{\mathcal{K}_{uv}: u \in C, v\notin C \text{ and } (u,v) \in E \} = 1-0 = 1$, $\mathcal{K}_{(0,0,0)^T,c} = \min \{\mathcal{K}_{(0,0,0)^T,v}: v \in C \text{ and } ((0,0,0)^T,v) \in E \} = 0$, $\mathcal{K}_{c,(0,0,\Dtot)^T} = \min \{\mathcal{K}_{u,(0,0,\Dtot)^T}: u \in C \text{ and } (u,(0,0,\Dtot)^T) \in E \} - \min \{\mathcal{K}_{uv}: u \in C, v\notin C \text{ and } (u,v) \in E \} = 0-0 = 0$, $\mathcal{K}_{(0,0,\Dtot)^T,c} = \min \{\mathcal{K}_{(0,0,\Dtot)^T,v}: v \in C \text{ and } ((0,0,\Dtot)^T,v) \in E \} = 1$, $\ldots$, and $\mathcal{K}_{c,(\Dtot-1,0,1)^T} = \min \{\mathcal{K}_{u,(\Dtot-1,0,1)^T}: u \in C \text{ and } (u,(\Dtot-1,0,1)^T) \in E \} - \min \{\mathcal{K}_{uv}: u \in C, v\notin C \text{ and } (u,v) \in E \} = 0-0 = 0$, $\mathcal{K}_{(\Dtot-1,0,1)^T,c} = \min \{\mathcal{K}_{(\Dtot-1,0,1)^T,v}: v \in C \text{ and } ((\Dtot-1,0,1)^T,v) \in E \} = 1$.
Step 4 involves $(\Dtot+2)$ iterations. In the first iteration, we fix the node with the largest value of $p$, which is $c$ in our case. At any node $u$ other than $a$ 
that is connected to $c$, the value of $p$ is updated according to the formula $p(u)=\max\{p(u),p(c)-\mathcal{K}_{uc}\}$, and then any edges leading to or from $c$ are removed. In the second iteration, the node $(0,0,0)^T$ has the largest value of $p$ among the remaining nodes, and thus is fixed. There is no other nodes connected to $(0,0,0)^T$ at this point, so we move to the next iteration. In the third iteration, the node $(0,0,\Dtot)^T$ is fixed. The node $(1,0,\Dtot-1)^T$ is connected to it, and thus the $p((1,0,\Dtot-1)^T)$ is updated to be $\max\{p((1,0,\Dtot-1)^T),p((0,0,\Dtot)^T)-\mathcal{K}_{(1,0,\Dtot-1)^T,(0,0,\Dtot)^T}\}=0$. Then, any edges leading to or from $(0,0,\Dtot)^T$ are removed. The remaining iterations will be similar to the third one. When all of the nodes other than $a$ have been fixed, the order of the pole of the mean first passage time from each state in $\B^c$ to $\B$ is given by the fixed value of the node to which the state belongs.

\noindent {\bf 3D model (from the fully repressed state to the fully active state)}: see Figure \ref{fig:algorithm_3DRtoA}.
The explanation of the panels for Input, Step 1 and Step 2 is given above with $\B = \{(0,\Dtot,0)^T\}$. Step 3 involves two iterations. In the first iteration, the collection of nodes consisting of $(\Dtot-1,0,1)^T$ and $(\Dtot,0,0)^T$ (called an r-connected set $C_1$) is condensed into a single node $c_1$. The order of the pole of the sojourn time in $C_1$ is $p(c_1) = \max_{u \in C_1} p (u) + \min \{\mathcal{K}_{uv}: u \in C_1, v\notin C_1 \text{ and } (u,v) \in E \} = 1 + 1 = 2$, where $E$ denotes the edge set of the graph in Step 3 before the first iteration. Moreover, $\mathcal{K}_{c_1,(\Dtot-2,0,2)^T} = \min \{\mathcal{K}_{u,(\Dtot-2,0,2)^T}: u \in C_1 \text{ and } (u,(\Dtot-2,0,2)^T) \in E \} - \min \{\mathcal{K}_{uv}: u \in C_1, v\notin C_1 \text{ and } (u,v) \in E \} = 1- 1 = 0$, $\mathcal{K}_{(\Dtot-2,0,2)^T,c_1} = \min \{\mathcal{K}_{(\Dtot-2,0,2)^T,v}: v \in C_1 \text{ and } ((\Dtot-2,0,2)^T,v) \in E \} = 0$, $\mathcal{K}_{c_1,(\Dtot-1,0,0)^T} = \min \{\mathcal{K}_{u,(\Dtot-1,0,0)^T}: u \in C_1 \text{ and } (u,(\Dtot-1,0,0)^T) \in E \} - \min \{\mathcal{K}_{uv}: u \in C_1, v\notin C_1 \text{ and } (u,v) \in E \} = 1- 1 = 0$ and $\mathcal{K}_{(\Dtot-1,0,0)^T,c_1} = \min \{\mathcal{K}_{(\Dtot-1,0,0)^T,v}: v \in C_1 \text{ and } ((\Dtot-1,0,0)^T,v) \in E \} = 0$. In the second iteration of Step 3, the collection of all nodes except for $(0,0,0)^T$ and $a$ (called an r-connected set $C_2$) is condensed to a single node $c_2$. The order of the pole of the sojourn time in $C_2$ is $p(c_2) = \max_{u \in C_2} p (u) + \min \{\mathcal{K}_{uv}: u \in C_2, v\notin C_2 \text{ and } (u,v) \in E \} = 2 + 0 = 2$, where $E$ denotes the edge set of the graph in Step 3 before the second iteration. Moreover, $\mathcal{K}_{c_2,(0,0,0)^T} = \min \{\mathcal{K}_{u,(0,0,0)^T}:\: u \in C_2 \text{ and } (u,(0,0,0)^T) \in E \} - \min \{\mathcal{K}_{uv}:\: u \in C_2, v\notin C_2 \text{ and } (u,v) \in E \} = 1-0=1$, $\mathcal{K}_{c_2,a} = \min \{\mathcal{K}_{ua}:\: u \in C_2 \text{ and } (u,a) \in E \} - \min \{\mathcal{K}_{uv}:\: u \in C_2, v\notin C_2 \text{ and } (u,v) \in E \} = 0-0=0$, and $\mathcal{K}_{(0,0,0)^T,c_2} = \min \{\mathcal{K}_{(0,0,0)^T,v}:\: v \in C_2 \text{ and } ((0,0,0)^T,v) \in E \} = 0$. Step 4 involves two iterations. In the first iteration, we fix the node with the largest value of $p$, which is $c_2$ in our case. At any node other than $a$ that is connected to $c_2$ (i.e., the origin $(0,0,0)^T$), the value of $p$ is updated to $p((0,0,0)^T) = \max \{ p((0,0,0)^T), p(c_2)-\mathcal{K}_{(0,0,0)^T,c_2}\} = \max \{0,2-0\} = 2$, and then any edges leading to or from $c_2$ are removed. In the second iteration, among the remaining nodes, we fix the node with the largest value of $p$, which is the origin. When all of the nodes other than $a$ have been fixed, the order of the pole of the mean first passage time from each state in $\B^c$ to $\B$ is given by the fixed value of the node to which the state belongs. 

\noindent {\bf 4D model (from the fully active state to the fully repressed state)}: see Figure \ref{fig:algorithm_4DAtoR}.
We illustrate how to use the algorithm for the 4D model when $\Dtot=2$; for larger $\Dtot$, the methodology will be the same. A state represents $(n_\mathrm{D_{12}^R},n_\mathrm{D^A},n_\mathrm{D_{1}^R},n_\mathrm{D_{2}^R})^T$. The explanation of the panels for Input, Step 1 and Step 2 is given above with $\B=\{(2,0,0,0)^T\}$. Step 3 involves only one iteration, where the collection of all nodes except for $(0,0,0,0)^T$, $(0,0,2,0)^T$, $(0,0,1,1)^T$, $(0,0,0,2)^T$, $(1,0,1,0)^T$ and $(1,0,0,1)^T$ (called an r-connected set $C$) is condensed to a single node $c$. The order of the pole of the sojourn time in $C$ is $p(c) = \max_{u \in C} p (u) + \min \{\mathcal{K}_{uv}: u \in C, v\notin C \text{ and } (u,v) \in E \} = 1 + 0 = 1$, where $E$ denotes the edge set of the graph in Step 3 before the first iteration. Moreover, the value of $\mathcal{K}$ for edges between $c$ and an original node $w$ that is not in $C$ are defined according to the formula $\mathcal{K}_{c,w} = \min \{\mathcal{K}_{uw}: u \in C \text{ and } (u,w) \in E \} - \min \{\mathcal{K}_{uv}: u \in C, v\notin C \text{ and } (u,v) \in E \}$, $\mathcal{K}_{w,c} = \min \{\mathcal{K}_{w,v}: v \in C \text{ and } (w,v) \in E \}$. Step 4 involves seven iterations. In the first iteration, we fix the node with the largest value of $p$, which is $c$ in our case. At any node $u$ other than $a$ 
that is connected to $c$, the value of $p$ is updated according to the formula $p(u)=\max\{p(u),p(c)-\mathcal{K}_{uc}\}$, and then any edges leading to or from $c$ are removed. In the second iteration of Step 4, the node $(0,0,0,0)^T$ has the largest value of $p$ among the remaining nodes, and then is fixed. There are no other nodes connected to $(0,0,0,0)^T$ at this point, so we move to the next iteration. In the third iteration, the node $(0,0,2,0)^T$ is fixed. The node $(1,0,1,0)^T$ is connected to it, and thus $p((1,0,1,0)^T)$ is updated to be $\max\{p((1,0,1,0)^T),p((0,0,2,0)^T)-\mathcal{K}_{(1,0,1,0)^T,(0,0,2,0)^T}\}=0$. Then, any edges leading to or from $(0,0,2,0)^T$ are removed. The remaining iterations will be similar to the third one. When all of the nodes other than $a$ have been fixed, the order of the pole of the mean first passage time from each state in $\B^c$ to $\B$ is given by the fixed value of the node to which the state belongs. 

\noindent {\bf 4D model (from the fully repressed state to the fully active state)}: see Figure \ref{fig:algorithm_4DRtoA}.
We again illustrate how to use the algorithm for the 4D model when $\Dtot=2$; for larger $\Dtot$, the methodology will be the same. The explanation of the panels for Input, Step 1 and Step 2 is given above with $\B=\{(0,2,0,0)^T\}$. Step 3 involves two iterations. In the first iteration, the collection of the nodes $(1,0,1,0)^T$, $(1,0,0,1)^T$ and $(2,0,0,0)^T$ (called an r-connected set $C_1$) is condensed into a single node $c_1$. The order of the pole of the sojourn time in $C_1$ is $p(c_1) = \max_{u \in C_1} p (u) + \min \{\mathcal{K}_{uv}: u \in C_1, v\notin C_1 \text{ and } (u,v) \in E \} = 1 + 1 = 2$, where $E$ denotes the edge set of the graph in Step 3 before the first iteration. Moreover, the value of $\mathcal{K}$ of edges between $c_1$ and an original node $w$ that is not in $C_1$ are defined if there is an edge between some node $u \in C_1$ and $w$ and according to the formula $\mathcal{K}_{c_1,w} = \min \{\mathcal{K}_{uw}: u \in C_1 \text{ and } (u,w) \in E \} - \min \{\mathcal{K}_{uv}: u \in C_1, v \notin C_1 \text{ and } (u,v) \in E \}$, $\mathcal{K}_{w,c_1} = \min \{\mathcal{K}_{w,v}: v \in C_1 \text{ and } (w,v) \in E \}$. In the second iteration of Step 3, the collection of all nodes except for $(0,0,0,0)^T$ and $a$ (called an r-connected set $C_2$) is condensed to a single node $c_2$. The order of the pole of the sojourn time in $C_2$ is $p(c_2) = \max_{u \in C_2} p (u) + \min \{\mathcal{K}_{uv}: u \in C_2, v\notin C_2 \text{ and } (u,v) \in E \} = 2 + 0 = 2$, where $E$ denotes the edge set of the graph in Step 3 before the second iteration. Moreover, $\mathcal{K}_{c_2,(0,0,0,0)^T} = \min \{\mathcal{K}_{u,(0,0,0,0)^T}:\: u \in C \text{ and } (u,(0,0,0,0)^T) \in E \} - \min \{\mathcal{K}_{uv}:\: u \in C_2, v\notin C_2 \text{ and } (u,v) \in E \} = 1-0=1$, $\mathcal{K}_{c_2,a} = \min \{\mathcal{K}_{ua}:\: u \in C_2 \text{ and } (u,a) \in E \} - \min \{\mathcal{K}_{uv}:\: u \in C_2, v\notin C_2 \text{ and } (u,v) \in E \} = 0-0=0$, and $\mathcal{K}_{(0,0,0,0)^T,c_2} = \min \{\mathcal{K}_{(0,0,0,0)^T,v}:\: v \in C_2 \text{ and } ((0,0,0,0)^T,v) \in E \} = 0$. Step 4 involves two iterations. In the first iteration, we fix the node with the largest value of $p$, which is $c_2$ in our case. At any node other than $a$ that is connected to $c_2$ (i.e., the origin $(0,0,0,0)^T$), the value of $p$ is updated to $p((0,0,0,0)^T) = \max \{ p((0,0,0,0)^T), p(c_2)-\mathcal{K}_{(0,0,0,0)^T,c_2}\} = \max \{0,2-0\} = 2$, and then any edges leading to or from $c_2$ are removed. In the second iteration, of the remaining nodes, we fix the node with the largest value of $p$, which is the origin. When all of the nodes other than $a$ have been fixed, the order of the pole of the mean first passage time from each state in $\B^c$ to $\B$ is given by the fixed value of the node to which the state belongs. 

 \clearpage
 
\clearpage

%% file: Appendix_DeviationMatrix.tex
\subsection{Leading coefficient for the MFPT}

\subsubsection{Proof of Theorem \ref{thm:MFPTcoefficient}}
\label{proofthmleadcoeff}

	\begin{proof}
		Fix $\lambda = \max \{q_x (\eps): x \in \X, 0 \leq \eps < \eps_0\}$. The $\lambda$ here should not be confused with other rates $\lambda$ with subscripts and/or superscripts used elsewhere of this paper. In the following, we use the breve symbol to denote notation associated with discrete time Markov chains defined below. 
		
		Let $Y^\eps = \{Y^\eps(n) : n \in \Z_+\}$ be a discrete time Markov chain with transition matrix $\labelDTMC{P}(\eps) = I + \frac{1}{\lambda} Q(\eps)$ for each $0 \leq \eps < \eps_0$\footnote{
  In general, when $0<\eps<\eps_0$, the discrete time Markov chain $Y^\eps$ is different from the embedded discrete time Markov chain described in Section \ref{sec:pertCTMC}. In particular, the discrete time Markov chain used here can have self loops, whereas the embedded discrete time Markov chain has no self loops.}. Note that $Y^\eps$ is a singularly perturbed discrete time Markov chain under the definition of Avrachenkov et al. \cite{AvrachenkovFilarHowlett2}. Let $\labelDTMC{\Pi}(\eps)$ be the ergodic projection of $Y^\eps$ and $\labelDTMC{H}(\eps)$ be the deviation matrix of $Y^\eps$ (see definitions in SI - Section \ref{sec:deviationmatrix}).
		The ergodic projection of $Y^0$ is $\labelDTMC{\Pi}(0) = \labelDTMC{W} \labelDTMC{M}$, where $I$ is the $|\A| \times |\A|$ identity matrix and
		\begin{equation*}
		    \labelDTMC{W} = 
		\left( \begin{array}{c}
		I \\\hline
		-T_0^{-1} R_0
		\end{array} \right)
		\quad \text{ and } \quad
		\labelDTMC{M} = 
		\left( \begin{array}{c | c }
		I & 0
		\end{array} \right).
		\end{equation*}
		Then,
		\begin{equation*}
            \labelDTMC{M} \Big(\frac{1}{\lambda} Q^{(1)}\Big) \labelDTMC{W} = \frac{1}{\lambda}
		\left( \begin{array}{c | c }
		I & 0
		\end{array} \right)
		\left( \begin{array}{c | c }
		A_1 & S_1 \\\hline
		R_1 & T_1
		\end{array} \right)
		\left( \begin{array}{c}
		I \\\hline
		-T_0^{-1} R_0
		\end{array} \right) = \frac{1}{\lambda} (A_1 +  S_1 (-T_0)^{-1} R_0) = \frac{1}{\lambda} Q_\A.
		\end{equation*}
		Assumptions \ref{assumption:transient_absorbing}, \ref{assumption:SingleRecurrentClass} and Lemma \ref{lem:nullityQmatrix} imply that the null space of this matrix is one dimensional.
		 
		Using the computational algorithm in Section 6.3.1 of \cite{AvrachenkovFilarHowlett2}, the generator\footnote{The transition matrix for a discrete time Markov chain with generator $\mathcal{G}$ is $\mathcal{P}=I+\mathcal{G}$.} for an aggregated discrete time Markov chain is $\labelDTMC{M} \Big(\frac{1}{\lambda} Q^{(1)}\Big) \labelDTMC{W} = \frac{1}{\lambda} Q_\A$, whose null space is one dimensional. Then, by the computational algorithm on page 176-177 of \cite{AvrachenkovFilarHowlett2} the deviation matrix $\labelDTMC{H} (\eps)$ has a Laurent series expansion with order of the pole equal to one:
	\begin{equation*}
        \labelDTMC{H}(\eps) = \frac{1}{\eps} \labelDTMC{H}^{(-1)} + \labelDTMC{H}^{(0)} + \eps \labelDTMC{H}^{(1)} + \ldots.
        \end{equation*}
		Since the aggregated Markov chain has a single recurrent class by Assumption \ref{assumption:SingleRecurrentClass}, the ergodic projection of the aggregated Markov chain is $\one \alpha$, where $\alpha$ is a row vector denoting the unique stationary distribution of the aggregated discrete time Markov chain. The deviation matrix of this aggregated Markov chain is $\labelDTMC{D}=(-\frac{1}{\lambda} Q_\A + \one \alpha)^{-1} - \one \alpha$. By Theorem 6.7 in \cite{AvrachenkovFilarHowlett2}, 
	\begin{equation*}
            \labelDTMC{H}^{(-1)} = \labelDTMC{W}\labelDTMC{D}\labelDTMC{M}= 
		\left( \begin{array}{c}
		I \\\hline
		-T_0^{-1} R_0
		\end{array} \right)
		\labelDTMC{D}
		\left( \begin{array}{c | c }
		I & 0
		\end{array} \right) =
		\left( \begin{array}{c | c }
		\labelDTMC{D} & 0 \\\hline
		T_0^{-1} R_0 \labelDTMC{D} & 0
		\end{array} \right).
        \end{equation*}
		
		For each $0 < \eps < \eps_0$, let $\labelDTMC{h}_{x,y}(\eps)$ be the mean first passage time from $x$ to $y$ in $Y^\eps$. Then, the mean first passage time from $x$ to $y$ in $X^\eps$ is 
		\begin{eqnarray*}
		h_{x,y}(\eps) = \frac{1}{\lambda} \labelDTMC{h}_{x,y}(\eps) &=& \frac{1}{\lambda} \frac{\labelDTMC{H}_{y,y} (\eps) - \labelDTMC{H}_{x,y} (\eps)}{\pi_y (\eps)}\\
		&=& \frac{1}{\lambda} \frac{(\frac{1}{\eps} \labelDTMC{D}_{y,y} + O(1)) - (\frac{1}{\eps} \labelDTMC{D}_{x,y} + O(1))}{\pi_y^{(k_y)} \eps^{k_y}_y +O(\eps^{k_y+1})} \\
		&=& \frac{D_{y,y} - D_{x,y}}{\pi_y^{(k_y)}} \frac{1}{\eps^{k_y+1}} + O \left(\frac{1}{\eps^{k_y}}\right),
		\end{eqnarray*}
            where we used \eqref{eq:deviationmatrix} to show that $\labelDTMC{D}=\lambda D$.
		The above equations use the properties of the deviation matrix given in SI - Section \ref{sec:deviationmatrix}.
		
		When $\hat{X}_\A$ is irreducible, then $\pi^{(0)} = \alpha$ has all strictly positive entries. Again, by SI - Section \ref{sec:deviationmatrix}, the mean first passage time from $x$ to $y$ in $\hat{X}_\A$ is finite and positive, and it is $\frac{1}{\lambda} \frac{\labelDTMC{D}_{y,y} - \labelDTMC{D}_{x,y}}{\alpha_y} = \frac{D_{y,y} - D_{x,y}}{\pi_y^{(0)}}$.
		In this case, the order of the pole of $h_{x,y}(\eps)$ is one and the leading coefficient is the mean first passage time from $x$ to $y$ in $\hat{X}_\A$.
	\end{proof}

\subsubsection{Properties of the deviation matrix for a discrete time Markov chain}
\label{sec:deviationmatrix}

In this section, we will start with a few results stated in Section 6.1 of Avrachenkov et al. \cite{AvrachenkovFilarHowlett2} about discrete time Markov chains with finite state space. These include the definitions and properties of the ergodic projection, the fundamental matrix and the deviation matrix. Then, we show one more fact about the deviation matrix. Lastly, Theorem 4.4.7 of Kemeny and Snell \cite{KemenySnell} gave a formula for mean first passage times for irreducible 
discrete time Markov chains in terms of the fundamental matrix and the stationary distribution, which is also briefly mentioned in \cite{AvrachenkovFilarHowlett2}. We will write this in terms of the deviation matrix and the stationary distribution with a simple modification.

Suppose $Y = \{Y(n) : n \in \Z_+\}$ is a discrete time Markov chain with a finite state space $\Y$. Suppose the state space $\Y$ is partitioned into $m$ ergodic classes (possibly including absorbing states) and a set of transient states, and accordingly, the transition matrix $\labelDTMC{P}$ is
\begin{equation*}
\labelDTMC{P}=\left( \begin{array}{ccc|c}
\labelDTMC{A}_1 & \ldots & 0 & 0 \\
\vdots & \ddots & \vdots & \vdots \\
0 & \ldots & \labelDTMC{A}_m & 0 \\\hline
\labelDTMC{R}_1 & \ldots & \labelDTMC{R}_m & \labelDTMC{T}
\end{array} \right).
\end{equation*}
The ergodic projection of $Y$ is given by the \textit{Cesaro} limit,
\begin{equation*}
\labelDTMC{\Pi}= \lim_{N\rightarrow \infty} \frac{1}{N+1} \sum_{n=0}^N \labelDTMC{P}^n.
\end{equation*}
It follows that $\labelDTMC{\Pi} (I-\labelDTMC{P})=0$ and $\labelDTMC{\Pi} \labelDTMC{\Pi} = \labelDTMC{\Pi}$. The ergodic projection $\labelDTMC{\Pi}$ is the eigenprojection of the transition matrix $\labelDTMC{P}$ corresponding to its maximal eigenvalue $1$. That is, if $\labelDTMC{\pi}_i$ is the unique stationary distribution for the discrete time Markov chain with transition matrix $\labelDTMC{A}_i$ for $1 \leq i \leq m$, then $\labelDTMC{\Pi}=\labelDTMC{W} \labelDTMC{M}$ with
\begin{equation*}
\labelDTMC{W}=
\left( \begin{array}{ccc}
\one & \ldots & 0 \\
\vdots & \ddots & \vdots \\
0 & \ldots & \one \\\hline
(I-\labelDTMC{T})^{-1} \labelDTMC{R}_1 \one & \ldots & (I-\labelDTMC{T})^{-1} \labelDTMC{R}_m \one
\end{array} \right)
\quad \text{ and } \quad
\labelDTMC{M} = 
\left( \begin{array}{ccc|c}
\labelDTMC{\pi}_1 & \ldots & 0 & 0 \\
\vdots & \ddots & \vdots & \vdots \\
0 & \ldots & \labelDTMC{\pi}_m & 0
\end{array} \right),
\end{equation*}
where $\labelDTMC{W}$ and $\labelDTMC{M}$ form bases for the right and left eigenspaces, respectively, which implies that $\labelDTMC{P}\labelDTMC{W}=\labelDTMC{W}$ and $\labelDTMC{M}\labelDTMC{P}=\labelDTMC{M}$.
One can see that $v (I-\labelDTMC{P}+\labelDTMC{\Pi})=0$ implies that $v=0$ and so $(I-\labelDTMC{P}+\labelDTMC{\Pi})$ is invertible.
The fundamental matrix $\labelDTMC{Z}$ and the deviation matrix $\labelDTMC{H}$ of $Y$ are well-defined:
\begin{equation*}
\labelDTMC{Z} = \lim_{N\rightarrow \infty} \frac{1}{N+1} \sum_{n=0}^N \sum_{\ell=0}^n (\labelDTMC{P}-\labelDTMC{\Pi})^\ell = (I-\labelDTMC{P}+\labelDTMC{\Pi})^{-1},
\end{equation*}
\begin{equation*}\labelDTMC{H} 
= \labelDTMC{Z} - \labelDTMC{\Pi} = (I-\labelDTMC{P}+\labelDTMC{\Pi})^{-1} - \labelDTMC{\Pi}.
\end{equation*}
We also have that $\labelDTMC{H} \labelDTMC{\Pi} = (\labelDTMC{Z} -\labelDTMC{\Pi}) \labelDTMC{\Pi}= 0$ since $\labelDTMC{P}  \labelDTMC{\Pi} = \labelDTMC{P}\labelDTMC{W}\labelDTMC{M}=\labelDTMC{W}\labelDTMC{M}=\labelDTMC{\Pi}$, $\labelDTMC{\Pi}= \labelDTMC{Z} (I-\labelDTMC{P}+\labelDTMC{\Pi}) \labelDTMC{\Pi} = \labelDTMC{Z} (\labelDTMC{\Pi}-\labelDTMC{\Pi}+\labelDTMC{\Pi})= \labelDTMC{Z} \labelDTMC{\Pi}$, and $\Pi^2=\Pi$.

Now, we show a property of the deviation matrix that is not in \cite{AvrachenkovFilarHowlett2} and is useful in Section \ref{sec:leadingMFPT}. Suppose $Q$ is an infinitesimal generator for a continuous time Markov chain on $\Y$ and $|Q_{y,y}| \leq \lambda$ for all $y \in \Y$. Then, $\labelDTMC{P}=I+\frac{1}{\lambda} Q$ defines a transition matrix for a discrete time Markov chain. The associated ergodic projection and deviation matrix $\labelDTMC{\Pi}$ and $\labelDTMC{H}$ for $\labelDTMC{P}$ satisfy $\labelDTMC{\Pi} Q = \lambda (\labelDTMC{\Pi} \labelDTMC{P} - \labelDTMC{\Pi}) = 0$ and
\begin{equation*}
\labelDTMC{H} = \left(I- \left(I+\frac{1}{\lambda} Q \right)+\labelDTMC{\Pi} \right)^{-1} - \labelDTMC{\Pi} = \left(-\frac{1}{\lambda} Q+\labelDTMC{\Pi} \right)^{-1} -\labelDTMC{\Pi},
\end{equation*}
and so
\begin{eqnarray}
I &=& (\labelDTMC{H} + \labelDTMC{\Pi}) \left(-\frac{1}{\lambda} Q+\labelDTMC{\Pi} \right) = -\frac{1}{\lambda} \labelDTMC{H} Q + \labelDTMC{H} \labelDTMC{\Pi} - \frac{1}{\lambda} \labelDTMC{\Pi} Q+\labelDTMC{\Pi} \labelDTMC{\Pi} \nonumber\\
&=& -\frac{1}{\lambda} \labelDTMC{H} Q + \labelDTMC{\Pi} = -\frac{1}{\lambda} \labelDTMC{H} Q + \frac{1}{\lambda} \labelDTMC{H} \labelDTMC{\Pi} - \labelDTMC{\Pi} Q + \labelDTMC{\Pi} \labelDTMC{\Pi} = \left(\frac{1}{\lambda} \labelDTMC{H} + \labelDTMC{\Pi} \right)(-Q+\labelDTMC{\Pi}).\nonumber 
\end{eqnarray}
Thus,
\begin{equation}
\label{eq:deviationmatrix}
\frac{1}{\lambda} \labelDTMC{H} = ((-Q+\labelDTMC{\Pi})^{-1}-\labelDTMC{\Pi}),
\end{equation}
where we have used the fact that $-Q+\labelDTMC{\Pi}$ is invertible because $\labelDTMC{\Pi}$ is the eigenprojection of $Q$ corresponding to the eigenvalue $0$.

Lastly, assume that $Y$ is irreducible. Then, $Y$ has a unique stationary distribution $\labelDTMC{\pi}$, which is a row vector, and the ergodic projection of $Y$ is $\one \labelDTMC{\pi}$. By Theorem 4.4.7 in \cite{KemenySnell}, the mean first passage time from $x \in \Y$ to $y \in \Y$ is $\frac{\labelDTMC{Z}_{y,y}-\labelDTMC{Z}_{x,y}}{\labelDTMC{\pi}_y}$. Since $\labelDTMC{Z} = \labelDTMC{H} + \one \labelDTMC{\pi}$, the mean first passage time from $x \in \Y$ to $y \in \Y$ is
\begin{equation*}
\frac{(\labelDTMC{H}_{y,y}+(\one \labelDTMC{\pi})_{y,y})-(\labelDTMC{H}_{x,y}+(\one \labelDTMC{\pi})_{x,y})}{\labelDTMC{\pi}_y}=\frac{(\labelDTMC{H}_{y,y}+\labelDTMC{\pi}_y)-(\labelDTMC{H}_{x,y}+\labelDTMC{\pi}_y)}{\labelDTMC{\pi}_y}=\frac{\labelDTMC{H}_{y,y}-\labelDTMC{H}_{x,y}}{\labelDTMC{\pi}_y}.
\end{equation*}

%% file: Appendix_1DModel.tex
 \subsection{1D Model: additional mathematical details}
 \label{sec:Appendix1Dmodel}

\

\noindent \textbf{Verification of Assumption \ref{assumption:transient_absorbing}.} 
In order to show that Assumption \ref{assumption:transient_absorbing} holds, consider the states $a =0$ and $r=\Dtot$ and the set $\T = \{1, \ldots, \Dtot-1\}$ defined in Section \ref{SS1}. Since $\Dtot \geq 2$, $\T \neq \emptyset$. From \eqref{eq:Qmatrix_eps}, we can see that $Q_{a,a+1}(0)=Q_{a,a-1}(0)=Q_{r,r+1}(0)=Q_{r,r-1}(0)=0$. As a consequence, both $a$ and $r$ are absorbing states under $Q(0)$. To see that the states in $\T$ are transient under $Q(0)$, consider a state $x \in \T$. Since $Q_{z,z+1}(0)= \frac{k^A_E}{V}(\Dtot - z) z >0$ for all $z \in \{1,\dots,\Dtot-1\}$, we have $Q_{x,x+1}(0) \dots Q_{\Dtot-1,\Dtot}(0)>0$. By Lemma \ref{lem:ClassifyTransient} and the fact that $r$ is an absorbing state, we have that $x$ is a transient state for $X^0$.

\noindent \textbf{Verification of Assumptions \ref{assumption:tildeX_irreducible} and \ref{assumption:SingleRecurrentClass}.} 
By Lemma \ref{lemma3}, it suffices to show Assumption \ref{assumption:tildeX_irreducible} holds. From \eqref{eq:Qmatrix_eps}, we can see that $\tilde{Q}_{a,a+1} > 0$. From the analysis made to prove Assumption \ref{assumption:transient_absorbing}, we know that there is a positive probability for $\tilde{X}$ to transition from $x \in \X \setminus \{a,r\}$ to $r$. It follows that any state $x \in \X \setminus \{r\}$ leads to $r$ under $\tilde{X}$. Now, we would like to show that there is a positive probability for transition from $r$ to $x \neq r \in \X$ for the process $\tilde{X}$. This is because $\tilde{Q}_{r,r-1}=b\frac{k_E^A}{V}\Dtot^2 > 0$ and $\tilde{Q}_{z,z-1}=Q_{z,z-1}(0)= \mu \frac{k^A_E}{V}(\Dtot - z)z >0$ for all $z \in \{1,\dots,\Dtot-1\}$. Thus, $r$ leads to any state in $\X \setminus \{r\}$ under $\tilde{X}$. Combining the above, we see that $\tilde{X}$ is irreducible and Assumption \ref{assumption:tildeX_irreducible} holds.

 
\noindent \textbf{Stationary distribution.} 
    Let us consider a one-dimensional finite state continuous time Markov chain in which the state space $\X =\{0,1,\ldots, \upb \}$ and the off-diagonal entries of the infinitesimal generator $Q$ are all zero except for the following positive rates:
    \begin{equation}
    \label{1DQmatrix}
     \begin{aligned}
     Q_{x,x+1} &= \lambda_x  & \text{ if } x \in \{0,\ldots,\upb-1\}, \\
     Q_{x,x-1} &= \gamma_x  & \text{ if } x \in \{1,\ldots,\upb\}.       
    \end{aligned}   
    \end{equation}
	Thus, the continuous time Markov chain is a birth-and-death process, it satisfies detailed balance (see \cite{gardiner1994}) and so the stationary distribution $\pi = (\pi_x)_{x \in \{0,1,...,\upb\}}$ satisfies
    \begin{equation*}
    \label{detbal}
    \pi_x=\frac{\lambda_{x-1}}{\gamma_{x}}\pi_{x-1}, \qquad \text{ for } x \in \{1,\ldots,\upb\}.
    \end{equation*}
    Applying this equality recursively, we can express $\pi_x$, $x\in\{1,\ldots,\upb\}$, as a function of $\pi_0$, obtaining
    \begin{equation}
    \label{detbal7}
    \pi_x= \pi_0\prod_{i=1}^{x} \frac{\lambda_{i-1}}{\gamma_{i}}.
    \end{equation}
    Using the fact that$\sum_{j=0}^{\upb}\pi_j =1$, we obtain 
    \begin{equation}
    \label{detbal8}
    \pi_0=\frac{1}{1+\sum_{j=1}^{\upb}\left(\prod_{i=1}^{j} \frac{\lambda_{i-1}}{\gamma_{i}}\right)}.
    \end{equation}
    Substituting \eqref{detbal8} in \eqref{detbal7}, we obtain
    \begin{equation}
    \label{detbalgeneral}
    \pi_x =\frac{\prod_{i=1}^{x}\frac{\lambda_{i-1}}{\gamma_{i}}}{1+\sum_{j=1}^{\upb}\left(\prod_{i=1}^{j} \frac{\lambda_{i-1}}{\gamma_{i}}\right)}\qquad \text{ for } x\in\{1,\ldots,\upb\}.
    \end{equation}

    Now, consider the one-dimensional continuous time Markov chain introduced in Section \ref{SS1} with state space $\X= \{0,1,\ldots,\mathrm{D_{tot}}\}$ and infinitesimal generator as defined in \eqref{eq:Qmatrix_eps}, which has nonzero off-diagonal entries given, for $\eps >0$, by
	   \begin{equation}
    \label{1Drates}
     \begin{aligned}
    \lambda^\eps_x &:= Q_{x,x+1}(\eps) = \left(\frac{k^A_E}{V}x+\eps \frac{k^A_E}{V} \mathrm{D_{tot}}\right)(\Dtot - x)   & \text{ if } x \in \{0,\ldots,\Dtot-1\}, \\
    \gamma^\eps_x &:= Q_{x,x-1}(\eps) = \mu \left(\frac{k^A_E}{V}(\Dtot - x)+b \eps \frac{k^A_E}{V} \mathrm{D_{tot}}\right)x   & \text{ if } x \in \{1,\ldots,\Dtot\}.
    \end{aligned}   
    \end{equation}
    By substituting the expressions for the rates in \eqref{1Drates} into \eqref{detbal8}-\eqref{detbalgeneral}, and suitably rearranging the terms, we obtain that
    \begin{equation}
         \pi_x(0) = \lim_{\eps \rightarrow 0}\pi_x(\eps) = \begin{cases}
                    \frac{b\mu^{\mathrm{D_{tot}}}}{1 + b\mu^{\mathrm{D_{tot}}}} & \text{ if } x=0 \\
                    0 & \text{ if } x \in \{1,\ldots, \Dtot-1\} \\
                    \frac{1}{1 + b\mu^{\mathrm{D_{tot}}}} & \text{ if } x = \Dtot.\nonumber
        \end{cases}
        \end{equation}

\noindent \textbf{Mean first passage time.} 
	Consider the one-dimensional, finite state, continuous time Markov chain introduced in \eqref{1DQmatrix}. We will determine an analytical expression for the MFPT from $x=\upb$ to $x=0$ and from $x=0$ to $x=\upb$ for this chain. We first focus on the former. For this, we exploit first step analysis (see Equation 3.1 of \cite{Norris2}), proceeding in a similar manner to that for \eqref{MFPTsystem}, to obtain
	\begin{equation}
    \label{1stepanalysisbis}
	\begin{cases}
	h_{0,0}=0, \\ 
	h_{x,0}=\frac{1}{\lambda_x+\gamma_x}+\frac{\lambda_x}{\lambda_x+\gamma_x}h_{x+1,0}+\frac{\gamma_x}{\lambda_x+\gamma_x}h_{x-1,0} & \mbox{if } x\in \{1,\ldots,\upb-1\}, \\
	h_{\upb,0}=\frac{1}{\gamma_{\upb}}+h_{\upb-1,0},\\
	\end{cases}
	\end{equation}
 where for $x,y \in \X$, $h_{x,y}=\E_x [\tau_y]$, $\tau_y = \inf \{t \geq 0:\: X(t)=y\}$, $X$ is the continuous time Markov chain with infinitesimal generator given by \eqref{1DQmatrix}.
Now, defining $\Delta h_{x,x-1}=h_{x,0}-h_{x-1,0}$ for $x\in \{1,...,\upb\}$, we can rewrite (\ref{1stepanalysisbis}) in the following way:
	\begin{equation}\label{1stepanalysistris}
	\begin{cases}
	h_{0,0}=0, \\ 
	\Delta h_{x,x-1} = \frac{1}{\gamma_x}+\frac{\lambda_x}{\gamma_x}\Delta h_{x+1,x} & \mbox{if } x\in  \{1,\ldots,\upb-1\},\\
	\Delta h_{\upb,\upb-1}=\frac{1}{\gamma_{\upb}}.\\
	\end{cases}
	\end{equation}
 From (\ref{1stepanalysistris}), we have an explicit formula for $\Delta h_{\upb,\upb-1}$ and any $\Delta h_{x,x-1}$ can be expressed as a function of $\Delta h_{x+1,x}$. Furthermore, if we sum the $\Delta h_{x,x-1}$ for $x= 1,\ldots,\upb$, we obtain
	\begin{equation}\label{repTTML1D}
	\begin{aligned}
	h_{\upb,0}=h_{\upb,0}-h_{0,0}=\sum_{x=1}^{\upb}\left(\Delta h_{x,x-1}\right)&= \Delta h_{1,0} + \Delta h_{2,1} + ... + \Delta h_{\upb-1,\upb-2} +\Delta h_{\upb,\upb-1}.
	\end{aligned}
	\end{equation}
	Thus, to evaluate the MFPT from $x=\upb$ to $x=0$, we can calculate $\Delta h_{x,x-1}$ for $x=\upb,\upb-1,...,1$ and then sum all of the terms. Defining $r_j=\frac{\lambda_1 \lambda_2...\lambda_j}{\gamma_1 \gamma_2 ... \gamma_j}$, for $j=1,...,\upb$, we obtain
	\begin{equation}
    \label{formulaDto0}
	\begin{aligned}
	h_{\upb,0} &=\frac{1}{\gamma_{\upb}}\left(1+\frac{\lambda_{\upb-1}}{\gamma_{\upb-1}}+\frac{\lambda_{\upb-1}\lambda_{\upb-2}}{\gamma_{\upb-1}\gamma_{\upb-2}}+...+r_{\upb-1}\right)\\
 &+\frac{1}{\gamma_{\upb-1}}\left(1+\frac{\lambda_{\upb-2}}{\gamma_{\upb-2}}+\frac{\lambda_{\upb-2}\lambda_{\upb-3}}{\gamma_{\upb-2}\gamma_{\upb-3}}+...+r_{\upb-2}\right)+...+\frac{1}{\gamma_{1}}\\
 &=\frac{r_{\upb-1}}{\gamma_{\upb}}\left(1+\sum_{i=1}^{\upb-1}\frac{1}{r_i}\right)+\sum_{i=2}^{\upb-1}\left[\frac{r_{i-1}}{\gamma_i}\left(1+\sum_{j=1}^{i-1}\frac{1}{r_j}\right)\right]+\frac{1}{\gamma_1}. 
	\end{aligned}
	\end{equation}
	With a similar procedure, we can obtain the MFPT from $x=0$ to $x=\upb$. More precisely, defining $\tilde r_j=\frac{\gamma_{\upb-1}\gamma_{\upb-2}...\gamma_{\upb-j}}{\lambda_{\upb-1}\lambda_{\upb-2}...\lambda_{\upb-j}}$, we have
	\begin{align}
\notag	h_{0,\upb}&=\frac{1}{\lambda_{0}}\left(1+\frac{\gamma_1}{\lambda_1}+\frac{\gamma_1\gamma_2}{\lambda_1\lambda_2}+...+\tilde r_{\upb-1}\right)\\
\label{formula3} &+\frac{1}{\lambda_{1}}\left(1+\frac{\gamma_2}{\lambda_2}+\frac{\gamma_2\gamma_3}{\lambda_2\lambda_3}+...+\tilde r_{\upb-2}\right)+...+\frac{1}{\lambda_{\upb-1}}\\
\notag	&=\frac{\tilde r_{\upb-1}}{\lambda_{0}}\left(1+\sum_{j=1}^{\upb-1}\frac{1}{\tilde r_i}\right)+\sum_{i=2}^{\upb-1}\left[\frac{\tilde r_{i-1}}{\lambda_{\upb-i}}\left(1+\sum_{j=1}^{i-1}\frac{1}{\tilde r_j}\right)\right]+\frac{1}{\lambda_{\upb-1}}. 
	\end{align}
	A more detailed derivation of the $h_{0,\upb}$ and $h_{\upb,0}$ is given in \cite{BrunoDelVecchio2}.
	
	Let us consider the one-dimensional continuous time Markov chain introduced in Section \ref{SS1}, with state space $\X =  \{0,1,\ldots,\mathrm{D_{tot}}\}$ and infinitesimal transition rates that can be written as in \eqref{1Drates}. Since all of the transition rates are $O(1)$, except for $\lambda^\eps_0$ and $\gamma^\eps_{\mathrm{D_{tot}}}$ which are $O(\eps)$, then both $h_{\mathrm{D_{tot}},0}(\eps)$ and $h_{0,\mathrm{D_{tot}}}(\eps)$ are $O(1/ \eps)$. This means that in the limit as $\eps \rightarrow 0$, $h_{\mathrm{D_{tot}},0}(\eps)$ and $h_{0,\mathrm{D_{tot}}}(\eps)$, which correspond to the time to memory loss of the repressed and active states, respectively, tend to infinity. Substituting parameters in \eqref{formulaDto0} and \eqref{formula3} yields \eqref{formulaDto0INTRO} and \eqref{formula3INTRO}, respectively.
%

%% file: Appendix_2DModel.tex
 \subsection{2D Model: additional mathematical details}
 \label{sec:Appendix2Dmodel}

\ 
 
 \noindent \textbf{Verification of Assumption \ref{assumption:transient_absorbing}.} 
 In order to show that Assumption \ref{assumption:transient_absorbing} holds, consider the states $a =(0,\Dtot)^T$ and $r=(\Dtot,0)^T$ and the set $\T = \{i_1,\ldots,i_m\}$ defined in Section \ref{sec:IllustrativeStationaryDistribution}. From \eqref{rates2D}, we can see that $Q_{a,a+v_j}(0)=Q_{r,r+v_j}(0)=0$ for every $1 \leq j \leq 4$. As a consequence, both $a$ and $r$ are absorbing states under $Q(0)$. To see that the states in $\T$ are transient under $Q(0)$, consider a state $x=(x_1,x_2)^T \in \T$. First, suppose $x_1 \neq 0$. By having the one-step transition along $v_2=(0,-1)^T$ occurring $x_2$ times where $Q_{z,z+v_2}(0) = \frac{k^A_E}{V} z_2 x_1 > 0$ for all $z=(x_1,z_2)^T$ and $1 \leq z_2 \leq x_2$, and having one-step transition along $v_3=(1,0)^T$ occurring $\Dtot-x_1$ times where $Q_{z,z+v_3} (0) = (\Dtot -z_1)\left(k_{W0}^R+k_{W}^R + \frac{k_{M}^R}{V}z_1\right) >0$ for all $z=(z_1,0)^T$ and $x_1 \leq z_1 \leq \Dtot-1$, we have a positive probability of transition from $x$ to $r$ under $Q(0)$. By Lemma \ref{lem:ClassifyTransient} and the fact that $r$ is an absorbing state, we have that $x$ is a transient state for $X^0$. On the other hand, suppose $x_1 = 0$. Since $x = (0,x_2) \in \T$, we have $0 \leq x_2 \leq \Dtot-1$. We can first have a one-step transition along $v_3=(1,0)^T$, where $Q_{x,x+v_3} (0) = (\Dtot -x_2)\left(k_{W0}^R+k_{W}^R\right) > 0$, to reach the state $(1,x_2)^T$ and then take the steps for the $x_1 \neq 0$ case to reach $r$. In this way, there is a positive probability of transition from $x$ to the absorbing state $r$ under $Q(0)$, and thus $x$ is transient by Lemma \ref{lem:ClassifyTransient}.

\noindent \textbf{Verification of Assumptions \ref{assumption:A_is_recurrent} and \ref{assumption:SingleRecurrentClass}.} By Lemma \ref{lemma3}, it suffices to show Assumption \ref{assumption:A_is_recurrent} holds. 
From \eqref{rates2D}, we can see that $\tilde{Q}_{a,a+v_2} > 0$. From the analysis made to prove Assumption \ref{assumption:transient_absorbing}, we know that there is a positive probability to transition from all $x \in \X \setminus \{a,r\}$ to $r$. Now, we would like to show that there is a positive probability to transition from $r$ to $x = (x_1,x_2)^T \in \X \setminus \{(0,0)^T\}$ for the process $\tilde{X}$. We first can have a one-step transition along $v_4=(-1,0)^T$ where $Q_{r,r+v_4}^{(1)} = \mu b \frac{k_{M}^A}{V} \Dtot^2 >0$, then have a one-step transition along $v_1=(0,1)^T$ where $Q_{r+v_4,r+v_4+v_1}^{(0)}= k_{W0}^A+k_{W}^A >0$, then have one-step transitions along $v_4=(-1,0)^T$ occurring $\Dtot-x_1-1$ times where $Q_{z,z+v_4}^{(0)} = \mu \frac{k^A_E}{V} z_1 >0$ for all $z=(z_1,1)^T$ and $x_1+1 \leq z_1 \leq \Dtot-1$. If $x_2 \neq 0$, we finally have one-step transitions along $v_1=(0,1)^T$ occurring $x_2-1$ times where $Q_{z,z+v_1}^{(0)} = (\Dtot -(x_1+z_2))\left(k_{W0}^A+k_{W}^A + \frac{k_{M}^A}{V}z_2\right) >0$ for all $z=(x_1,z_2)^T$ and $1 \leq z_2 \leq x_2-1$; if $x_2 = 0$ and $x_1 \neq 0$, we will make a one-step transition along $v_2=(0,-1)^T$ to $(x_1,0)^T$ where $Q_{z,z+v_2}^{(0)}= \frac{k^A_E}{V} x_1 >0$ with $z=(x_1,1)^T$ and $x_1 \geq 1$. Therefore, we have that there is a positive probability of transition from $r$ to each $x \in \X \setminus \{(0,0)^T\}$. Since $Q_{-v_j,(0,0)^T}^{(0)}=0$ for $j \in \{1,2,3,4\}$ such that $-v_j \in \X$, we conclude that $\CC= \X \setminus \{(0,0)^T\}$ is a closed communicating class under $\tilde{Q}$ and since it contains $\A$, Assumption \ref{assumption:A_is_recurrent} holds. Note that Assumption \ref{assumption:tildeX_irreducible} does not hold.

\noindent \textbf{Stationary distribution.} Here, we derive the expression for $\pi^{(1)}_x$, $x\in \T= \{i_1, \ldots, i_m\}$, for the case $\Dtot=2$. In this case $T(\eps) = T_0+\eps T_1$, with 

\begin{equation*}
T_0=\left( \begin{array}{c c c c}
	-q_3 & k_{W0}^R+k_{W}^R & 0 & 0 \\
	\mu \frac{k^A_E}{V}& -\frac{k^A_E}{V}(1+\mu) & 0 & \frac{k^A_E}{V} \\
	2(k_{W0}^A+k_{W}^A)&  0& -q_5 & 2(k_{W0}^R+k_{W}^R) \\
	0&  k_{W0}^A+k_{W}^A& 0 & -q_6 \end{array}    \right),
 \end{equation*}

\begin{equation*}
T_1=\left( \begin{array}{c c c c}
	-2 \frac{k_{M}^A}{V} & 0 & 2 \frac{k_{M}^A}{V} & 0 \\
	2\frac{k_{M}^A}{V} \mu b  & -2 \frac{k_{M}^A}{V}(1+\mu b)& 0 &2 \frac{k_{M}^A}{V} \\
	0  & 0 & 0 &0 \\
	0  & 0 & 2\frac{k_{M}^A}{V} \mu b &-2\frac{k_{M}^A}{V}    \end{array}    \right),
 \end{equation*}

\noindent in which $q_3=(k_{W0}^A+k_{W}^A + \frac{k_{M}^A}{V}+k_{W0}^R+k_{W}^R )$, $q_5=2(k_{W0}^A+k_{W}^A+k_{W0}^R+k_{W}^R)$ and $q_6=(k_{W0}^A+k_{W}^A+k_{W0}^R+k_{W}^R + \frac{k_{M}^R}{V})$. Then, by \eqref{eq:Beta1InTermsOfAlpha}, $\beta^{(1)} =\pi^{(1)}_x$, $x\in \T$, is given by $\beta^{(1)} = \alpha S_1(-T_0)^{-1}$, where $\alpha = (\pi^{(0)}_{(0,2)},\pi^{(0)}_{(2,0)})$ was derived in Section \ref{sec:IllustrativeStationaryDistribution} (Eq. \eqref{pi2DmodelRESULTS}). After some calculations, $\pi^{(1)}_x$ can be written for $\Dtot=2$ as
    \begin{align}
\notag    \pi^{(1)}_{i_1}&=\frac{4b \eps^2\left(\frac{k^A_M}{V}\right)^2\bar k_{W}^A\mu^2\left(\bar k^R_W(\bar k^R_W+\frac{k^R_M}{V})+(\bar k^A_W+\frac{k^A_M}{V})\left((1+\mu)(\bar k^R_W+\frac{k^R_M}{V})+\mu\frac{k^A_M}{V}\right)\right)}{d_1d_2},\\
\notag \pi^{(1)}_{i_2}&=\frac{4b \eps^2\left(\frac{k^A_M}{V}\right)^2\bar k_{W}^A\bar k^R_W\mu\left((\bar k^R_W+\frac{k^R_M}{V})(\bar k^R_W+ \bar k^A_W+k^A_M)+\mu(\bar k^A_W+\frac{k^A_M}{V})(\bar k^R_W+ \bar k^A_W+k^R_M)\right)}{\frac{k^A_E}{V}d_1d_2},\\
\notag \pi^{(1)}_{i_3}&=0,\\
\notag \pi^{(1)}_{i_4}&=\frac{4b \eps^2\left(\frac{k^A_M}{V}\right)^2\bar k_{W}^R\mu\left(
    (\bar k^R_W+\frac{k^R_M}{V})\left((1+\mu)(\bar k^A_W+\frac{k^A_M}{V})+\mu\bar k^R_W\right)+\mu\bar k^A_W(\bar k^A_W+\frac{k^A_M}{V})\right)}{d_1d_2},
    \end{align}
in which $\bar k^A_W=k_{W0}^A+k_{W}^A$, $\bar k^R_W=k_{W0}^R+k_{W}^R$, $d_1=\eps\frac{k^A_M}{V}(\bar k^R_W(\bar k^R_W+\frac{k^R_M}{V})+b\bar k^A_W\mu^2(\bar k^A_W+\frac{k^A_M}{V}))$, and $d_2=((\bar k^A_W+\frac{k^A_M}{V})((1+\mu)(\bar k^R_W+\frac{k^R_M}{V})+\mu \bar k^A_W)+\bar k^R_W(\bar k^R_W+\frac{k^R_M}{V}))$ and in which $i_1=(0,1)^T,i_2=(1,1)^T,i_3=(0,0)^T$, and $i_4=(0,2)^T$.

%% file: Appendix_3DModel.tex
\subsection{3D Model: additional mathematical details}
 \label{sec:Appendix3Dmodel}

\noindent \textbf{Verification of Assumption \ref{assumption:transient_absorbing}.} 
In order to show that Assumption \ref{assumption:transient_absorbing} holds, consider the states $a =(0,\Dtot,0)^T$ and $r=(\Dtot,0,0)^T$ and the set $\T = \{i_1,\ldots,i_m\}$ defined in Section \ref{SD3D}. From \eqref{rates3D}, we can see that $Q_{a,a+v_j}(0)=Q_{r,r+v_j}(0)=0$ for every $1 \leq j \leq 6$. As a consequence, both $a$ and $r$ are absorbing states under $Q(0)$. To see that the states in $\T$ are transient under $Q(0)$, consider a state $x=(x_1,x_2,x_3)^T \in \T$. First, suppose $x_1+x_3 \neq 0$. By having the one-step transitions along $v_2=(0,-1,0)^T$ occurring $x_2$ times where $Q_{z,z+v_2}(0) = \frac{k^A_E}{V}(x_3+2x_1) z_2 > 0$ for all $z=(x_1,z_2,x_3)^T$ and $1 \leq z_2 \leq x_2$, then having one-step transitions along $v_3=(0,0,1)^T$ occurring $\Dtot-x_1-x_3$ times where $Q_{z,z+v_3} (0) = (\Dtot -(x_1+z_3))\left(k^1_{W0}+k^1_{W} + \frac{k^{'}_{M}}{V}x_1\right)>0$ for all $z=(x_1,0,z_3)^T$ and $x_3 \leq z_3 \leq \Dtot-x_1-1$ and finally having one-step transitions along $v_5=(1,0,-1)^T$ occurring $\Dtot-x_1$ times where $Q_{z,z+v_5}(0) = (\Dtot-z_1) \left(k^2_{W0}+ \frac{k_{M}}{V}z_1 + \frac{\bar k_{M}}{V}\frac{\Dtot+z_1-1}{2}\right)>0$ for all $z=(z_1,0,\Dtot-z_1)^T$ and $x_1 \leq z_1 \leq \Dtot-1$, we have a positive probability of transition from $x$ to $r$ under $Q(0)$. By Lemma \ref{lem:ClassifyTransient} and the fact that $r$ is an absorbing state, we have that $x$ is a transient state for $X^0$. On the other hand, suppose $x_1+x_3 = 0$. Since $x = (0,x_2,0)^T \in \T$, we have $0 \leq x_2 \leq \Dtot-1$. We can first have a one-step transition along $v_3=(0,0,1)^T$, where $Q_{x,x+v_3}(0) = (\Dtot -x_2) \left(k^1_{W0}+k^1_{W}\right) >0$, to reach the state $(0,x_2,1)^T$ and then take the steps in the $x_1+x_3 \neq 0$ case. In this way, there is a positive probability of transition from $x$ to the absorbing state $r$, and thus $x$ is transient by Lemma \ref{lem:ClassifyTransient}.

\noindent \textbf{Verification of Assumption \ref{assumption:SingleRecurrentClass}.} 
To show that Assumption \ref{assumption:SingleRecurrentClass} holds, consider the continuous time Markov chain $\tilde{X}$ with infinitesimal generator $\tilde{Q}$ as described in \eqref{eq:DefTildeQ} and shown in Fig. \ref{fig:3Dmodel}(d). We will first see that $\{i_m,r\}$ forms a closed class under $\tilde{Q}$. For this, we see that $Q_{r,r+v_j}(\eps)$ vanishes for every $1 \leq j \leq 5$ and $\eps \geq 0$, while $Q_{r,r+v_6}(\eps)= \eps\mu b\frac{k_{M}^A}{V}\Dtot^2$. Therefore, the only transition from $r$ under $\tilde{Q}$ is given by $\tilde{Q}_{r,r+v_6}= \mu b\frac{k_{M}^A}{V}\Dtot^2 > 0$, where $r+v_6 = i_m$. From \eqref{rates3D}, we can see that $Q_{i_m,i_m+v_j}(0)=0$ for every $j \in \{1,2,3,4,6\}$ and $Q_{i_m,i_m+v_5}(0)= k^2_{W0}+ \frac{k_{M}}{V}(\Dtot -1) + \frac{\bar k_{M}}{V}(\Dtot -1) > 0$. Since $i_m+v_5=r$ we see that $\tilde{Q}_{i_m,r} > 0$. Therefore, $\{i_m,r\}$ forms a closed class under $\tilde{Q}$. The fact that $\hat{X}_{\A}$, shown in Fig. \ref{fig:3Dmodel}(e) consists of erasing the times from $\tilde{X}$ in which the process is in $\T$, together with Lemma \ref{lem:welldefined_QA}, yields that $r$ is an absorbing state under $Q_{\A}$. From \eqref{rates3D}, we can see that $\tilde{Q}_{a,i_1} > 0$. From the analysis made to prove Assumption \ref{assumption:transient_absorbing} we obtain that $i_1$ leads to $r$ under $\tilde{Q}$ which is part of a closed class. By interpreting $\hat{X}_{\A}$ again as a time-change of $\tilde{X}$, by Lemma \ref{lem:welldefined_QA} we obtain that $a$ is transient under $Q_{\A}$. As a consequence, $Q_{\A}$ has a single recurrent class consisting of the state $r$, and so Assumption \ref{assumption:SingleRecurrentClass} holds, and furthermore, $\alpha=[\alpha_a,\alpha_r]$ with $\alpha_a=0$ and $\alpha_r=1$. In addition, the previous arguments show that neither Assumption \ref{assumption:A_is_recurrent} nor \ref{assumption:tildeX_irreducible} holds for this model. 
    
\noindent \textbf{Stationary distribution.} 
Here, we derive an expression for $\pi^{(1)}_x$, $x \in \T= \{i_1, \ldots, i_m\}$.
Matrices $A_1$, $S_1$, and $R_0$ can be written as 
\begin{equation}\nonumber
     A_1 = \left( \begin{array}{cc}
        - s_1 & 0  \\
            0 & -s_2
    \end{array}    \right),\;
         S_1 = \left( \begin{array}{ccccc}
                s_1 & 0 & \ldots & \ldots & 0 \\
                0 & \ldots & \ldots & 0 & s_2
    \end{array}    \right),\;
         R_0 = \left( \begin{array}{cc}
                r_1 & 0 \\
                0 & 0                \\
                \vdots & \vdots       \\
                0 & 0                \\
                0 & r_2
    \end{array}    \right),
    \end{equation}    
    with $s_1=\frac{k_{M}^A}{V}\mathrm{D^2_{tot}}$, $s_2=\mu b \frac{k_{M}^A}{V}\mathrm{D^2_{tot}}$, $r_1=(k_{W0}^A+k_{W}^A + \frac{k_{M}^A}{V}(\Dtot-1))$, and $r_2=(k^2_{W0} + (\frac{k_{M}}{V} + \frac{\bar k_{M}}{V})(\Dtot - 1))$. From \eqref{eq:Beta1InTermsOfAlpha}, $\beta^{(1)}=[\pi^{(1)}_{i_1},...,\pi^{(1)}_{i_m}]= \alpha S_1(-T_0)^{-1}$, and so, given that the last row of $T_0$ is made of all zeros except for the last element, that is $(T_0)_{i_m,i_m}$ and given that $(T_0)_{i_m,i_m}=k^2_{W0}+ (\frac{k_{M}}{V} + \frac{\bar k_{M}}{V})(\Dtot - 1)$,  $\beta^{(1)}=[0,...,0,\pi^{(1)}_{i_m}]$, with 
    \begin{equation*}
    \pi^{(1)}_{i_m}=\frac{\mu b\frac{k_{M}^A}{V}\Dtot^2}{k^2_{W0}+ (\frac{k_{M}}{V} + \frac{\bar k_{M}}{V})(\Dtot - 1)}.
    \end{equation*}
    Now, $\alpha^{(1)}=[\pi^{(1)}_{a},\pi^{(1)}_{r}]$ is the unique vector such that
%
            \begin{align}
            \alpha^{(1)}Q_{\A} = -\beta^{(1)}[R_1 + T_1(-T_0)^{-1}R_0],\;\;\;\;\; \alpha^{(1)}\one = -\beta^{(1)}\one.  
            \end{align}
    
    For an illustration, suppose $\Dtot=2$. Then $A_0 =0$, $S_0=0$ and matrices $A_1\in \R^{2 \times 2}$ and $S_1\in \R^{2 \times 8}$ are given by 
\begin{equation}\nonumber
     A_1 = \left( \begin{array}{cc}
        -4\frac{k^A_M}{V} & 0  \\
            0 & -4\frac{k^A_M}{V}\mu b 
    \end{array}    \right),\;
         S_1 = \left( \begin{array}{ccccc}
                4\frac{k^A_M}{V} & 0 & \ldots & \ldots & 0 \\
                0 & \ldots & \ldots & 0 & 4\frac{k^A_M}{V}\mu b 
    \end{array}    \right).
    \end{equation}
Furthermore, $R_1=0$ and $R_0\in \R^{8 \times 2}$ can be written as
\begin{equation}\nonumber
    R_0 = \left( \begin{array}{cc}
                \bar k_{W}^A + \frac{k_{M}^A}{V} & 0 \\
                0 & 0                \\
                \vdots & \vdots        \\
                0 & 0                \\
                0 &  k^2_{W0}+ \frac{k_{M}}{V} + \frac{\bar k_{M}}{V}
\end{array}    \right).
\end{equation}
Finally, matrices $T_0$ and $T_1$ can be written as 

\begin{equation*}
	\resizebox{0.98\textwidth}{!}{
		\begin{minipage}{\linewidth}
			\begin{flalign*}
				T_0 &= \left( \begin{array}{cccccccc}
					-\bar k^A_W  - \frac{k_{M}^A}{V} -\bar k^1_W & 0 & 0 & 0 & \bar k^1_W & 0 & 0 & 0\\
					0 & -\frac{k^A_E}{V}(2+\mu) & 0 & 2\frac{k^A_E}{V}\ & \mu \frac{k^A_E}{V} & 0 & 0 & 0\\
					2\bar k^A_W & 0 & -2(\bar k^A_W  +\bar k^1_W) & 0 & 0 & 2 \bar k^1_W & 0 & 0\\
					0 & \bar k^A_W & 0 & -(\bar k^A_W+\bar k^1_W+\frac{k^{'}_{M}}{V}) & 0 & 0 & 0 & \bar k^1_W+\frac{k^{'}_{M}}{V}\\
					\mu'\frac{k^A_E}{V} & k^2_{W0} & 0 & 0 & -(\frac{k^A_E}{V}(1+\mu')+k^2_{W0}) & \frac{k^A_E}{V} & 0 & 0\\
					0 & 0 & 0 & k^2_{W0} & \bar k^A_W & -(k^2_{W0} + \bar k^A_W + \bar k^1_W) & \bar k^1_W & 0\\
					0 & 0 & 0 & 0 & 0 & 0 & -2(k^2_{W0}+\frac{\bar k_{M}}{V}) & 2(k^2_{W0}+\frac{\bar k_{M}}{2V})\\
					0 & 0 & 0 & 0 & 0 & 0 & 0 & -2(k^2_{W0}+\frac{k_{M}}{V}+\frac{\bar k_{M}}{V})\\
				\end{array} \right),\\
				T_1 &= \left( \begin{array}{cccccccc}
					-2\frac{k_{M}^A}{V} & 0 & 2\frac{k_{M}^A}{V} & 0 & 0 & 0 & 0 & 0\\
					0 & -2\frac{k_{M}^A}{V}(1+b\mu) & 0 & 2\frac{k_{M}^A}{V}\ & 2b\mu\frac{k_{M}^A}{V} & 0 & 0 & 0\\
					0 & 0 & 0 & 0 & 0 & 0 & 0 & 0\\
					0 & 0 & 0 & -2b\mu\frac{k_{M}^A}{V}  & 0 & 2b\mu\frac{k_{M}^A}{V} & 0 & 0\\
					2\mu'\beta \frac{k_{M}^A}{V} & 0 & 0 & 0 & -2\frac{k_{M}^A}{V}(1+\beta\mu') & 2\frac{k_{M}^A}{V} & 0 & 0\\
					0 & 0 & 2\mu' \beta \frac{k_{M}^A}{V} & 0 & 0 & -2\mu' \beta \frac{k_{M}^A}{V} & 0 & 0\\
					0 & 0 & 0 & 0 & 0 & 4\mu' \beta \frac{k_{M}^A}{V} & -4\mu' \beta \frac{k_{M}^A}{V} & 0\\
					0 & 0 & 0 & 2\mu' \beta \frac{k_{M}^A}{V} & 0 & 0 & 2\mu b \frac{k_{M}^A}{V} & -2(\mu' \beta \frac{k_{M}^A}{V}+\mu b \frac{k_{M}^A}{V})\\
				\end{array} \right),
			\end{flalign*}
		\end{minipage}
	}
\end{equation*}

\noindent in which $\bar k^A_W=k_{W0}^A+k_{W}^A$, $\bar k^R_W=k_{W0}^R+k_{W}^R$. Now, by applying Theorem \ref{thm:HigherOrderTerms_for_stationary distribution}, we first obtain that $\pi(0)=\pi^{(0)}=[\alpha,0] = \left[\alpha_a,\alpha_r, 0 \ldots ,0 \right]$ where $\alpha$ is the unique probability vector such that $\alpha Q_{\A} =0$. In this case,
    	\begin{equation}\label{QA3Dsi}
    	    Q_\A = \frac{K_1+\mu K_2}{K_3+\mu K_4 + \mu' K_5 + \mu \mu' K_6} \left( \begin{array}{c c}
	    -1 & 1 \\
	    0 & 0 \\
    \end{array} \right),
    \end{equation}
    with 
      \begin{align}
\notag    K_1&=8\bar k^1_W\frac{k_{M}^A}{V} (\bar k^1_W+\frac{k^{'}_{M}}{V})((\bar k^1_W +k^2_{W0})(\frac{k^A_E}{V}+k^2_{W0})+\bar k^A_W k^2_{W0}),\\
\notag     K_2&=4\bar k^1_W\frac{k_{M}^A}{V}\frac{k^A_E}{V}((\bar k^1_W+\frac{k^{'}_{M}}{V})(\bar k^1_W +k^2_{W0})+\bar k^1_W \bar k^A_W),\\
\notag     K_3&=\frac{k^A_E}{V}(2(\bar k^1_W+\bar k^A_W + k^2_{W0}+\frac{k^{'}_{M}}{V}+\frac{k_{M}^A}{V}))(\bar k^1_W)^2\\
 \notag      &+2\bar k^1_W(\bar k^1_W(\bar k^1_W+k^2_{W0}(k^2_{W0}+\frac{k^{'}_{M}}{V}))+\frac{k^{'}_{M}}{V}(k^2_{W0})^2)\\
 \notag      &+2\frac{k_{M}^A}{V}k^2_{W0}(\frac{k^{'}_{M}}{V}(k^2_{W0}+\bar k^A_W)+\bar k^1_W+(2\bar k^1_W+\bar k^A_W+k^2_{W0}))\\
\notag       &+2\frac{k_{M}^A}{V}\frac{k^{'}_{M}}{V}k^2_{W0}(\bar k^1_W+\bar k^A_W)+2\bar k^1_W k^2_{W0}(\frac{k^A_E}{V}\frac{k^{'}_{M}}{V}+\bar k^A_W\frac{k_{M}^A}{V}+2\bar k^A_W\frac{k^{'}_{M}}{V}),\\
\notag     K_4&=\frac{k^A_E}{V}((\bar k^1_W +\frac{k^{'}_{M}}{V})(\bar k^A_W+\frac{k_{M}^A}{V})(\bar k^1_W +k^2_{W0})+(\bar k^1_W +k^2_{W0})\bar k^1_W\frac{k^{'}_{M}}{V})\\
\label{eq:3Dconstants}     &+\frac{k^A_E}{V}\bar k^1_W(\bar k^1_W(\bar k^A_W+k^2_{W0}+1)+\bar k^A_W(\bar k^A_W+\frac{k_{M}^A}{V})),\\
\notag     K_5&=2\frac{k^A_E}{V}(\bar k^1_W +k^2_{W0}+\bar k^A_W)(\bar k^A_W+\frac{k_{M}^A}{V})(\bar k^1_W +\frac{k^{'}_{M}}{V}),\\
\notag     K_6&=\frac{k^A_E}{V}(\bar k^1_W +k^2_{W0}+\bar k^A_W)(\bar k^A_W+\frac{k_{M}^A}{V})(\bar k^A_W+\bar k^1_W +\frac{k^{'}_{M}}{V}),
    \end{align} 
and then $\alpha_a=0$ and $\alpha_r=1$. Let us now derive an expression for $\beta^{(1)}$. Starting from the transient states $\T= \{i_1, \ldots, i_8\}$, we obtain that $\beta^{(1)}=[\pi^{(1)}_{i_1},...,\pi^{(1)}_{i_8}]$ can be determine by $\beta^{(1)} = \alpha S_1(-T_0)^{-1}$, obtaining $\beta^{(1)}=[0,...,0,\pi^{(1)}_{i_8}]$, with \begin{equation*}
\pi^{(1)}_{i_8}=\frac{4\mu b\frac{k_{M}^A}{V}}{k^2_{W0}+ (\frac{k_{M}}{V} + \frac{\bar k_{M}}{V})}.
\end{equation*}
 Finally, $\alpha^{(1)}=[\pi^{(1)}_{a},\pi^{(1)}_{r}]$ is the unique vector such that $\alpha^{(1)}Q_{\A} = -\beta^{(1)}[R_1 + T_1(-T_0)^{-1}R_0]$ and $ \alpha^{(1)}\one = -\beta^{(1)}\one$. After some calculations, we obtain
 \begin{equation}
        \pi^{(1)}_{a}=\frac{\mu^2\mu'^2K_7}{K_8(K_9+K_{10}\mu)},\;\;\;
        \pi^{(1)}_{r}=-\pi^{(1)}_{a}-\pi^{(1)}_{i_8}=-\frac{\mu^2\mu'^2K_7+\mu K_{11}(K_{9}+K_{10}\mu)}{K_{8}(K_{9}+K_{10}\mu)},
        \end{equation}
    with 
   \begin{equation}\nonumber
       \begin{aligned}
       \bar k^1_W&=k^1_{W0}+k^1_{W},\;\;\;\;\;\bar k^A_W=k^A_{W0}+k^A_{W},\\
           K_7&=2\beta b \frac{k_{M}^A}{V} \bar k^A_W (\bar k^A_W+\frac{k_{M}^A}{V})(\bar k^A_W+\bar k^1_W +k^2_{W0}),\\
           K_8&=\bar k^1_W(k^2_{W0}+\frac{k_{M}}{V}+\frac{\bar k_{M}}{V}),\;\;\;\;\;K_{10} = \frac{k^A_E}{V}\left((\bar k^1_W+k^2_{W0})(\bar k^1_W+\frac{k^{'}_{M}}{V})+\bar k^A_W \bar k^1_W\right)\\
           K_9&=2(\bar k^1_{W}+\frac{k^{'}_{M}}{V})(\frac{k^A_E}{V}(\bar k^1_W+k^2_{W0})+ k^2_{W0}(k^2_{W0}+\bar k^1_W+\bar k^A_W)),\;\;\;\;\;
           K_{11}=4b\bar k^1_W\frac{k_{M}^A}{V}.
       \end{aligned}
   \end{equation}

    \noindent \textbf{Time to memory loss.} 
    As a reminder, we define the time to memory loss of the active state as $h_{a,r} (\eps)$ and the time to memory loss of repressed state as $h_{r,a} (\eps)$. Let us start by deriving the order and the leading coefficient of $h_{a,r} (\eps)$ and $h_{r,a} (\eps)$. By \eqref{eqn:orderSD}, we know the order of the stationary distribution at $a$ and $r$ are $k_a=-\min\{1-2,0\}=1$ and $k_r=-\min\{1-1,0\}=0$, respectively. This is consistent with the results in Section \ref{SD3D}. Moreover, the leading coefficient in the stationary distribution for the fully repressed and fully active states are $\pi^{(0)}_r=1$ and $\pi^{(1)}_a>0$, respectively. Now, $\hat{X}_\A$ has the infinitesimal generator in the form of 
    \begin{equation*}
    Q_\A=
	\left( \begin{array}{c c}
	    -(Q_\A)_{a,r} & (Q_\A)_{a,r} \\
	    0 & 0 \\
    \end{array} \right),
    \end{equation*}
    and $\hat{X}_\A$ has a unique stationary distribution $\alpha=[0,1]$.
    By Theorem \ref{thm:MFPTcoefficient},
    \begin{equation*}
    D=\left( \begin{array}{c c}
	    \frac{1}{(Q_\A)_{a,r}} & -\frac{1}{(Q_\A)_{a,r}} \\
	    0 & 0 \\
    \end{array} \right),
    \end{equation*}
    \begin{equation*}
    h_{a,r}(\eps)= \frac{D_{r,r} - D_{a,r}}{\pi_r^{(k_r)}} \frac{1}{\eps^{k_r+1}} + O\left(\frac{1}{\eps^{k_r}}\right) = \frac{D_{r,r} - D_{a,r}}{\alpha_r} \frac{1}{\eps} + O(1) = \frac{1}{(Q_\A)_{a,r}}\frac{1}{\eps} + O(1),
    \end{equation*}
    and

    \begingroup
    \small
    \begin{equation*}
    h_{r,a}(\eps)= \frac{D_{a,a} - D_{r,a}}{\pi_a^{(k_a)}} \frac{1}{\eps^{k_a+1}} + O\left(\frac{1}{\eps^{k_a}}\right) = \frac{D_{a,a} - D_{r,a}}{\pi_a^{(1)}} \frac{1}{\eps^2} + O \left(\frac{1}{\eps}\right)= \frac{1}{(Q_\A)_{a,r} \cdot \pi_a^{(1)}}\frac{1}{\eps^2} + O\left(\frac{1}{\eps}\right).
    \end{equation*}
    \endgroup
    As an example, when $\Dtot=2$, $Q_\A$ and $\pi_a^{(1)}$ are shown in \eqref{QA3D} and \eqref{pi3D}, and we obtain that 
    \begin{equation*}
    h_{a,r}(\eps)=\frac{K_3+\mu K_4 + \mu' K_5 + \mu \mu' K_6}{K_1+\mu K_2} \frac{1}{\eps} + O(1),
    \end{equation*}
    and
    \begin{equation*}h_{r,a}(\eps)=\frac{K_3+\mu K_4 + \mu' K_5 + \mu \mu' K_6}{K_1+\mu K_2} \frac{K_8(K_9+K_{10}\mu)}{\mu^2\mu'^2K_7} \frac{1}{\eps^2} + O\left(\frac{1}{\eps}\right),
    \end{equation*}
    where $K_i$, $i=1,...,11$, are non-negative functions independent of $\mu$ and $\mu'$ as defined in \eqref{eq:3Dconstants}.

Now, let us verify that both conditions $(i)$ and $(ii)$ of Theorem \ref{thm:Compartheorem} hold. To this end, let us first write the directions of the six possible transitions of the continuous time Markov chain $X^{\eps}(t)$, which are $v_1=(0,1,0)^T$, $v_2=(0,-1,0)^T$, $v_3=(0,0,1)^T$, $v_4=(0,0,-1)^T$, $v_5=(1,0,-1)^T$, $v_6=(-1,0,1)^T$, with the associated infinitesimal transition rates that can be written as	       
%
$\rate_1(x)=f_A(x)$, $\rate_2(x)=g^A_\eps(x)$, $\rate_3(x)=f_{R1}(x)$, $\rate_4(x)=g^{\eps}_{R1}(x)$, $\rate_5(x)=f_{R12}(x)$, $\rate_6(x)=g^{\eps}_{R12}(x)$.
%
Define the matrix
    \begin{equation}\label{matrixA3D}
       A= \begin{bmatrix}
        1 & 0 & 0\\
        0 & -1 & 0\\
        1 & 0 & 1
        \end{bmatrix}\nonumber
    \end{equation}
    and, for $x\in \X$, $(K_A+x)\cap \X=\{ w \in \X:\: x \preccurlyeq_A w \}$. Let us also introduce infinitesimal transition rates $\Breve{\rate}_{i}(x)$, $i=1,2,...,6$, defined as for $\rate_{i}(x)$, $i=1,2,...,6$, with all the parameters having the same values except that $\mu'$ is replaced by $\Breve{\mu}'$, with $\mu' \ge \Breve{\mu}'$. Given that $Av_1= (0,-1,0)^T$, $Av_2= (0,1,0)^T$, $Av_3=(0,0,1)^T$, $Av_4=(0,0,-1)^T$, $Av_5=(1,0,0)^T$, $Av_6=(-1,0,0)^T$, condition $(i)$ of Theorem \ref{thm:Compartheorem} holds.

   %

   To verify condition $(ii)$ of Theorem \ref{thm:Compartheorem}, consider $x\in \X$ and $y \in \partial_1(K_A+x) \cap \X=\{ w \in \X:\: x_1=w_1, x_2\ge w_2, x_1 + x_3\le w_1 + w_3\}=\{ w \in \R^3:\: x_1=w_1, x_2\ge w_2, x_3\le w_3\}$. Given that $\inn{A_{1\bullet},v_5}=1$ and $\inn{A_{1\bullet},v_6}=-1$, we need to verify that $\rate_{5}(x)\le\Breve{\rate}_{5}(y)$ and $\rate_{6}(x)\ge\Breve{\rate}_{6}(y)$. Since 
   $x_1=y_1, x_2\ge y_2, x_3\le y_3$, then 
   $\rate_{5}(x)=x_3\left(k^2_{W0}+ \frac{k_{M}}{V}x_1 + \frac{\bar k_{M}}{V}(x_1+\frac{x_3-1}{2})\right)\le y_3\left(k^2_{W0}+ \frac{k_{M}}{V}x_1 + \frac{\bar k_{M}}{V}(y_1+\frac{y_3-1}{2})\right)= \Breve{\rate}_{5}(y)$ and $\rate_{6}(x)=x_1\mu\left(\eps \frac{k_{M}^A}{V}\Dtot b + x_2\frac{k^A_E}{V}\right)$ 
   
   \noindent $\ge y_1\mu\left(\eps \frac{k_{M}^A}{V}\Dtot b+ y_2\frac{k^A_E}{V}\right)= \Breve{\rate}_{6}(y)$.
    Let us now consider $x\in \X$ and $y \in \partial_2(K_A+x) \cap \X =\{ w \in \X:\: x_1\le w_1, x_2=w_2,  x_1 + x_3\le w_1 + w_3\}$. Given that $\inn{A_{2\bullet},v_1}=-1$ and $\inn{A_{2\bullet},v_2}=1$, we need to verify that $\rate_1(x)\ge\Breve{\rate}_1(y)$ and $\rate_2(x)\le\Breve{\rate}_2(y)$. Since
    $x_1\le y_1, x_2=y_2,  x_1 + x_3\le y_1 + y_3$, then
    $\rate_{1}(x)=(\Dtot -(x_1+x_2+x_3))\left(k_{W0}^A+k_{W}^A + \frac{k_{M}^A}{V}x_2\right)\ge (\Dtot -(y_1+y_2+y_3))\left(k_{W0}^A+k_{W}^A + \frac{k_{M}^A}{V}y_2\right)=\Breve{\rate}_{1}(y)$ and $\rate_{2}(x)=x_2\left(\eps \frac{k_{M}^A}{V}\Dtot + \frac{k^A_E}{V}(x_3+2x_1)\right) \le y_2\left(\eps \frac{k_{M}^A}{V}\Dtot + \frac{k^A_E}{V}(y_3+2y_1)\right)=\Breve{\rate}_{2}(y)$.
    Finally, consider $x\in \X$ and $y \in \partial_3(K_A+x) \cap \X=\{ w \in \X \:|\: x_1\le w_1, x_2\ge w_2,  x_1 + x_3= w_1 + w_3\}=\{ w \in \X \:|\: x_1\le w_1, x_2\ge w_2,  x_3 \ge w_3\}$. Given that $\inn{A_{3\bullet},v_3}=1$ and $\inn{A_{3\bullet},v_4}=-1$, we need to check that $\rate_3(x)\le\Breve{\rate}_3(y)$ and $\rate_4(x)\ge\Breve{\rate}_4(y)$. Since $x_1\le y_1, x_2\ge y_2,  x_3 \ge y_3$, then $\rate_3(x)=(\Dtot -(x_1+x_2+x_3))\left(k^1_{W0}+k^1_{W} + \frac{k^{'}_{M}}{V}x_1\right)\le(\Dtot -(y_1+y_2+y_3))\left(k^1_{W0}+k^1_{W} + \frac{k^{'}_{M}}{V}y_1\right)=\Breve{\rate}_3(y)$ and $\rate_4(x)=x_3\mu'\left(\eps \frac{k_{M}^A}{V}\Dtot \beta + x_2\frac{k^A_E}{V}\right)\ge y_3\Breve{\mu}'\left(\eps \frac{k_{M}^A}{V}\Dtot \beta + y_2\frac{k^A_E}{V}\right)=\Breve{\rate}_4(y)$.

We can then conclude that all of the conditions of Theorem \ref{thm:Compartheorem} hold.

%% file: Appendix_4DModel.tex
\subsection{4D Model: additional mathematical details}
 \label{sec:Appendix4Dmodel}

\ 

\noindent \textbf{Verification of Assumption \ref{assumption:transient_absorbing}.} 
In order to show that Assumption \ref{assumption:transient_absorbing} holds, consider the states $a =(0,\Dtot,0,0)^T$ and $r=(\Dtot,0,0,0)^T$ and the set $\T = \{i_1,\ldots,i_m\}$ defined in Section \ref{SD4D}. From \eqref{rates4D}, we can see that $Q_{a,a+v_j}(0)=Q_{r,r+v_j}(0)=0$ for every $1 \leq j \leq 10$. As a consequence, both $a$ and $r$ are absorbing states under $Q(0)$. To see that the states in $\T$ are transient under $Q(0)$, consider a state $x=(x_1,x_2,x_3,x_4)^T \in \T$. First, suppose $x_1+x_3+x_4 \neq 0$. By having the one-step transitions along $v_2=(0,-1,0,0)^T$ occurring $x_2$ times where $Q_{z,z+v_2}(0) = \frac{k^A_E}{V}(x_3+x_4+2x_1)z_2 > 0$ for all $z=(x_1,z_2,x_3,x_4)^T$ and $1 \leq z_2 \leq x_2$, then having one-step transitions along $v_3=(0,0,1,0)^T$ occurring $\Dtot-x_1-x_3-x_4$ times where $Q_{z,z+v_3} (0) = (\Dtot -(x_1+z_3+x_4))\left(k^1_{W0}+k^1_{W} + \frac{k^{'}_{M}}{V}(x_1+x_4)\right) >0$ for all $z=(x_1,0,z_3,x_4)^T$ and $x_3 \leq z_3 \leq \Dtot-x_1-x_4-1$, then having one-step transitions along $v_9=(1,0,0,-1)^T$ occurring $x_4$ times where $Q_{z,z+v_9} (0) = (x_1+x_4-z_1) \left(k^1_{W0} + \frac{k^{'}_{M}}{V}\frac{x_1+x_4+z_1-1}{2}\right) > 0$ for all $z=(z_1,0,\Dtot-x_1-x_4,x_1+x_4-z_1)^T$ and $x_1 \leq z_1 \leq x_1+x_4-1$, and finally having one-step transitions along $v_7=(1,0,-1,0)^T$ occurring $\Dtot-x_1-x_4$ times where $Q_{z,z+v_7} (0) = (\Dtot-z_1) \left(k^2_{W0}+ \frac{k_{M}}{V}z_1  +  \frac{\bar k_{M}}{V}\frac{\Dtot+z_1-1}{2}\right) >0$ for all $z=(z_1,0,\Dtot-z_1,0)^T$ and $x_1+x_4 \leq z_1 \leq \Dtot-1$, we have a positive probability of transition from $x$ to $r$ under $Q(0)$. By Lemma \ref{lem:ClassifyTransient} and the fact that $r$ is an absorbing state, we have that $x$ is a transient state for $X^0$.  On the other hand, suppose $x_1+x_3+x_4 = 0$. Since $x = (0,x_2,0,0)^T \in \T$, we have $0 \leq x_2 \leq \Dtot-1$. We can first have a one-step transition along $v_3=(0,0,1,0)^T$, where $Q_{x,x+v_3} (0) = (\Dtot -x_2) (k^1_{W0}+k^1_{W}) > 0$, to reach the state $(0,x_2,1,0)^T$ and then take the steps in the $x_1+x_3+x_4 \neq 0$ case. In this way, there is a positive probability of transition from $x$ to the absorbing state $r$, and thus $x$ is transient by Lemma \ref{lem:ClassifyTransient}.

\noindent \textbf{Verification of Assumption \ref{assumption:SingleRecurrentClass}.} 
To show that Assumption \ref{assumption:SingleRecurrentClass} holds, consider the continuous time Markov chain $\tilde{X}$ with infinitesimal generator $\tilde{Q}$ as described in \eqref{eq:DefTildeQ} and shown in Fig. \ref{fig:4Dmodel1}(d). We will first see that $\{i_{m-1},i_m,r\}$, with $i_{m-1}=r+v_{10}$ and $i_m=r+v_8$, forms a closed class under $\tilde{Q}$. For this, we see that $Q_{r,r+v_j}(\eps)$ vanishes for every $j=\{1,2,3,4,5,6,7,9\}$ and $\eps \geq 0$, while $Q_{r,r+v_8}(\eps)= \eps\mu b\frac{k_{M}^A}{V}\Dtot^2$ and $Q_{r,r+v_{10}}(\eps)= \eps\mu' \beta\frac{k_{M}^A}{V}\Dtot^2$. Therefore, the only transitions from $r$ under $\tilde{Q}$ are to $i_{m-1}$ with rate $\tilde{Q}_{r,i_{m-1}}= \mu' \beta\frac{k_{M}^A}{V}\Dtot^2 > 0$ and to $i_m$ with rate $\tilde{Q}_{r,i_m}= \mu b\frac{k_{M}^A}{V}\Dtot^2 > 0$. From \eqref{rates4D}, we can see that $Q_{i_{m-1},i_{m-1}+v_j}(0)=0$ for $j \in \{1,2,3,4,5,6,7,8,10\}$, $Q_{i_m,i_m+v_j}(0)=0$ for $j \in \{1,2,3,4,5,6,8,9,10\}$, $Q_{i_{m-1},r}(0)=Q_{i_{m-1},i_{m-1}+v_9}(0)= k^1_{W0}+ \frac{k'_{M}}{V}(\Dtot -1) > 0$ and $Q_{i_m,r}(0)=Q_{i_m,i_m+v_7}(0)= k^2_{W0}+ \frac{k_{M}}{V}(\Dtot -1) + \frac{\bar k_{M}}{V}(\Dtot -1) > 0$. Therefore, $\{i_{m-1},i_m,r\}$ forms a closed class under $\tilde{Q}$. The fact that $\hat{X}_{\A}$ consists of erasing the times from $\tilde{X}$ in which the process is in $\T$, together with Lemma \ref{lem:welldefined_QA}, yields that $r$ is an absorbing state under $Q_{\A}$. From \eqref{rates4D}, we can see that $\tilde{Q}_{a,i_1} > 0$ where $i_1=(0,\Dtot-1,0,0)^T$. From the analysis made to prove Assumption \ref{assumption:transient_absorbing}, we obtain that $i_1$ leads to $r$ under $\tilde{Q}$, which is part of the closed class $\{i_{m-1},i_m,r\}$. By interpreting $\hat{X}_{\A}$ again as a time-change of $\tilde{X}$, by Lemma \ref{lem:welldefined_QA} we obtain that $a$ is transient under $Q_{\A}$. As a consequence, $Q_{\A}$ has a single recurrent class consisting of the state $r$. Thus, Assumption \ref{assumption:SingleRecurrentClass} holds, and furthermore, $\alpha=[\alpha_a,\alpha_r]$ with $\alpha_a=0$ and $\alpha_r=1$. In addition, the previous arguments show that neither Assumption \ref{assumption:A_is_recurrent} nor \ref{assumption:tildeX_irreducible} holds for this model. 

\noindent \textbf{Stationary distribution.} Here, we derive an expression for $\pi^{(1)}_x$, $x \in \T= \{i_1, \ldots, i_m\}$.
Matrices $A_1$, $S_1$, and $R_0$ can be written as
\begin{equation}\nonumber
     A_1 = \left( \begin{array}{cc}
        - s_1 & 0  \\
            0 & -(s_2+s_3)
    \end{array}    \right),\;
         S_1 = \left( \begin{array}{cccccc}
                s_1 & 0 & \ldots & \ldots & 0 & 0 \\
                0 & \ldots & \ldots & 0 & s_2 & s_3
    \end{array}    \right),\;
         R_0 = \left( \begin{array}{cc}
                r_1 & 0 \\
                0 & 0                \\
                \vdots & \vdots       \\
                0 & 0                \\
                0 & r_2                \\
                0 & r_3
    \end{array}    \right),
    \end{equation}  
with $s_1=\frac{k_{M}^A}{V}\mathrm{D^2_{tot}}$, $s_2=
\mu'\beta\eps\frac{k_{M}^A}{V}\mathrm{D^2_{tot}}$, $s_3=
\mu b\eps\frac{k_{M}^A}{V}\mathrm{D^2_{tot}}$, 
$r_1=(k_{W0}^A+k_{W}^A + \frac{k_{M}^A}{V}(\Dtot-1))$,
$r_2=(k^1_{W0} + \frac{k^{'}_{M}}{V}(\Dtot - 1))$ and $r_3=(k^2_{W0} + (\frac{k_{M}}{V} + \frac{\bar k_{M}}{V})(\Dtot - 1))$. Now, we determine $\beta^{(1)}=[\pi^{(1)}_{i_1},...,\pi^{(1)}_{i_m}] = \alpha S_1(-T_0)^{-1}$. Given that the only two elements different from zero in the last two rows of $T_0$ are $(T_0)_{i_{m-1},i_{m-1}}=(k^1_{W0} + \frac{k^{'}_{M}}{V}(\Dtot-1)$ and $(T_0)_{i_m,i_m}=(k^2_{W0}+  (\frac{k_{M}}{V}+\frac{\bar k_{M}}{V}(\Dtot-1))$, we obtain $\beta^{(1)}=[0,...,0,\pi^{(1)}_{i_{m-1}},\pi^{(1)}_{i_m}]$, with 
\begin{equation}\nonumber
\pi^{(1)}_{i_{m-1}}=\frac{\mu' \beta\frac{k_{M}^A}{V}\Dtot^2}{k^1_{W0}+ \frac{k'_{M}}{V}(\Dtot - 1)},\;\;\;\;\;\; \pi^{(1)}_{i_m}=\frac{\mu b\frac{k_{M}^A}{V}\Dtot^2}{k^2_{W0}+ (\frac{k_{M}}{V} + \frac{\bar k_{M}}{V})(\Dtot - 1)}.
\end{equation}
    Now, $\alpha^{(1)}=[\pi^{(1)}_{a},\pi^{(1)}_{r}]$ is the unique vector such that $\alpha^{(1)}Q_{\A} = -\beta^{(1)}[R_1 + T_1(-T_0)^{-1}R_0]$, $\alpha^{(1)}\one = -\beta^{(1)}\one$. 

    As an example, suppose $\Dtot=2$, $\beta=b$, $k^1_{W}=k^2_{W}=k^A_{W}=0$, $k^1_{W0}=k^2_{W0}=k^A_{W0}= k_{W0}$ and $k'_M=\bar k_M = k^A_M= k_M $. Then, we have that $A_0 =0$, $S_0=0$ and matrices $A_1\in \R^{2 \times 2}$ and $S_1\in \R^{2 \times 13}$ are equal to 
\begin{equation}\nonumber
     A_1 = \left( \begin{array}{cc}
        -4\frac{k_M}{V} & 0  \\
            0 & -4\frac{k_M}{V}(\mu b +\mu' \beta)
    \end{array}    \right),\;
         S_1 = \left( \begin{array}{cccccc}
                4\frac{k_M}{V} & 0 & \ldots & \ldots & 0 & 0 \\
                0 & \ldots & \ldots & 0 & 4\frac{k_M}{V}\mu' b & 4\frac{k_M}{V}\mu b 
    \end{array}    \right).
    \end{equation}
Furthermore, $R_1=0$ and $R_0\in \R^{13 \times 2}$ can be written as
\begin{equation}\nonumber
    R_0 = \left( \begin{array}{cc}
                k_{W0} + \frac{k_{M}}{V} & 0 \\
                0 & 0                \\
                \vdots & \vdots        \\
                0 & 0                \\
                0 & k_{W0}+ \frac{k_{M}}{V}\\
                0 &  k_{W0}+ 2\frac{k_{M}}{V}
\end{array}    \right).
\end{equation}
Finally, matrices $T_0$ and $T_1$ can be written as 
\begin{equation}\nonumber
    T_0 = \left( \begin{array}{cc}
        T_0^1 & T_0^2\\
        T_0^3 & T_0^4\\
\end{array}    \right),
    T_1 = \left( \begin{array}{cc}
        T_1^1 & T_1^2\\
        T_1^3 & T_1^4\\
\end{array}    \right),
\end{equation}
\noindent with $T_1^2= 0^{7 \times 6}$ and

\begin{equation*}
	\resizebox{0.99\textwidth}{!}{
		\begin{minipage}{\linewidth}
			\begin{flalign*}
				T_0^1 &= \left( \begin{array}{ccccccc}
					-3k_{W0}-\frac{k_{M}}{V} & k_{W0} & 0 & 0 & 0 & k_{W0} & 0\\
					\mu\frac{k^A_E}{V} & -(1+\mu)\frac{k^A_E}{V}-k_{W0} & 0 & \frac{k^A_E}{V} & 0 & 0 & 0\\
					2k_{W0} & 0 & -6k_{W0} & 2k_{W0} & 0 & 0 & 2k_{W0}\\
					0 & k_{W0} & 0 & -(4k_{W0}+2\frac{k_{M}}{V}) & k_{W0}+\frac{k_{M}}{V} & 0 & 0\\
					0 & 0 & 0 & 0 & -(k_{W0}+\frac{k_{M}}{V}) & 0 & 0\\
					\mu' \frac{k^A_E}{V} & 0 & 0 & 0 & 0 & -(1+\mu')\frac{k^A_E}{V}-k_{W0} & \frac{k^A_E}{V}\\
					0 & 0 & 0 & 0 & 0 & k_{W0} & -(4k_{W0}+\frac{k_{M}}{V})\\
				\end{array} \right),\\
				T_0^2 &= \left( \begin{array}{cccccc}
					0 & 0 & 0 & 0 & 0 & 0\\
					0 & 0 & k_{W0} & 0 & 0 & 0\\
					0 & 0 & 0 & 0 & 0 & 0\\
					k_{W0}+\frac{k_{M}}{V} & 0 & 0 & k_{W0} & 0 & 0\\
					0 & 0 & 0 & 0 & k_{W0}+\frac{k_{M}}{V} & 0\\
					0 & 0 & k_{W0} & 0 & 0 & 0\\
					k_{W0}+\frac{k_{M}}{V} & k_{W0} & 0 & k_{W0} & 0 & 0\\
				\end{array} \right),\;\;
				T_1^3 = \left( \begin{array}{ccccccc} 
					0 & 0 & 0 & 2\mu' b\frac{k_{M}}{V} & 0 & 0 & 2\mu b\frac{k_{M}}{V}\\
					0 & 0 & 0 & 0 & 0 & 0 & 2\mu ' b \frac{k_{M}}{V}\\   
					0 & 2\mu ' b \frac{k_{M}}{V} & 0 & 0 & 0 & 2 \mu b \frac{k_{M}}{V} & 0\\
					0 & 0 & 0 & 2 \mu' b\frac{k_{M}}{V} & 0 & 0 & 2 \mu b\frac{k_{M}}{V}\\
					0 & 0 & 0 & 0 & 2 \mu' b\frac{k_{M}}{V} & 0 & 0\\
					0 & 0 & 0 & 0 & 0 & 0 & 0\\
				\end{array} \right),\\
				T_1^1 &= \left( \begin{array}{ccccccc}
					-2\frac{k_{M}}{V}& 0 & 2\frac{k_{M}}{V} & 0 & 0 & 0 & 0\\
					2\mu b\frac{k_{M}}{V} & -2(\mu b+1)\frac{k_{M}}{V} & 0 & 2\frac{k_{M}}{V} & 0 & 0 & 0\\
					0 & 0 & 0 & 0 & 0 & 0 & 0\\
					0 & 0 & 2\mu b \frac{k_{M}}{V} & -2\mu b \frac{k_{M}}{V} & 0 & 0 & 0\\
					2\mu' b \frac{k_{M}}{V} & 0 & 0 & 0 & 0 & -2(\mu' b+1)\frac{k_{M}}{V} & 2\frac{k_{M}}{V}\\
					0 & 0 & 2 \mu' b\frac{k_{M}}{V}  & 0 & 0 & 0 & -2 \mu' b\frac{k_{M}}{V} \\
				\end{array} \right),\;\;
				T_0^3 = \left( \begin{array}{ccccccc} 
					0 & 0 & 0 & 0 & 0 & 0 & 0\\
					0 & 0 & 0 & 0 & 0 & 0 & 0\\
					0 & \mu ' \frac{k^A_E}{V} & 0 & 0 & 0 & \mu \frac{k^A_E}{V} & 0\\         
					0 & 0 & 0 & 0 & 0 & 0 & 0\\ 
					0 & 0 & 0 & 0 & 0 & 0 & 0\\
					0 & 0 & 0 & 0 & 0 & 0 & 0\\
				\end{array} \right),\\
				T_0^4 &= \left( \begin{array}{cccccc}
					-(2k_{W0}+\frac{k_{M}}{V}) & 0 & 0 & 0 & k_{W0}+\frac{k_{M}}{V} & k_{W0}\\
					0 & -(k_{W0}+\frac{k_{M}}{V}) & 0 & 0 & 0 &  k_{W0}+\frac{k_{M}}{V}\\
					0 & 0 & -(\mu + \mu' + 2)\frac{k^A_E}{V} & 2\frac{k^A_E}{V} & 0 & 0\\
					0 & 0 & k_{W0} & -3(k_{W0}+\frac{k_{M}}{V}) & k_{W0}+2\frac{k_{M}}{V} & k_{W0}+\frac{k_{M}}{V}\\
					0 & 0 & 0 & 0 & -(k_{W0}+\frac{k_{M}}{V}) & 0\\
					0 & 0 & 0 & 0 & 0 & -(k_{W0}+2\frac{k_{M}}{V})\\
				\end{array} \right),\\
				T_1^4 &= \left( \begin{array}{cccccc}
					-2(\mu+\mu')b\frac{k_{M}}{V} & 0 & 0 & 0 & 0 & 0\\
					0 & -2\mu ' b \frac{k_{M}}{V} & 0 & 0 & 0 & 0\\
					0 & 0 & -2(\mu b + \mu' b + 1)\frac{k_{M}}{V} & 2\frac{k_{M}}{V} & 0 & 0\\
					0 & 0 & 0 & -2(\mu + \mu')b\frac{k_{M}}{V}  & 0 & 0\\
					2 \mu b\frac{k_{M}}{V} & 0 & 0 & 2 \mu b\frac{k_{M}}{V} & -2(2\mu + \mu')b\frac{k_{M}}{V} & 0\\
					2 \mu' b\frac{k_{M}}{V} & 2 \mu b\frac{k_{M}}{V} & 0 & 2 \mu' b\frac{k_{M}}{V} & 0 & -2(\mu + 2 \mu')b\frac{k_{M}}{V}\\
				\end{array} \right),
			\end{flalign*}
		\end{minipage}
	}
\end{equation*}

    Now, by applying Theorem \ref{thm:HigherOrderTerms_for_stationary distribution}, we first obtain that $\pi(0)=\pi^{(0)}=[\alpha,0] = \left[\alpha_a,\alpha_r, 0 \ldots ,0 \right]$ where $\alpha$ can be obtained by solving $\alpha Q_{\A}= 0$. In this case, we obtain that
    	\begin{equation}\nonumber
    	    Q_\A = \frac{K_1((\mu')^2K_2+(\mu)^2K_3+\mu\mu'K_4+\mu'K_5+\mu K_6+K_{7})}{\mu'\mu(\mu'+\mu) K_8+(\mu')^2K_{9}+(\mu)^2K_{10}+\mu\mu'K_{11}+\mu'K_{12}+\mu K_{13}+K_{14}} \left( \begin{array}{c c}
	    -1 & 1 \\
	    0 & 0 \\
    \end{array} \right),
    \end{equation}
    with 

 \begingroup
 \scriptsize
\begin{align}
\notag K_1&=4k_{W0}\frac{k_M}{V},\\
\notag K_2&=\left(\frac{k^A_E}{V}\right)^2\left(6\left(\frac{k_M}{V}\right)^3+39k_{W0}\left(\frac{k_M}{V}\right)^2+68\frac{k_M}{V}(k_{W0})^2+32(k_{W0})^3\right),\\
\notag K_3&=\left(\frac{k^A_E}{V}\right)^2\left(6\left(\frac{k_M}{V}\right)^3+36k_{W0}\left(\frac{k_M}{V}\right)^2+64\frac{k_M}{V}(k_{W0})^2+32(k_{W0})^3\right),\\
\notag K_4&=\left(\frac{k^A_E}{V}\right)^2\left(12\left(\frac{k_M}{V}\right)^3+75k_{W0}\left(\frac{k_M}{V}\right)^2+132\frac{k_M}{V}(k_{W0})^2+64(k_{W0})^3\right),\\
\notag K_5&=\left(\frac{k^A_E}{V}\right)^2\left(24\left(\frac{k_M}{V}\right)^3+140k_{W0}\left(\frac{k_M}{V}\right)^2+222\frac{k_M}{V}(k_{W0})^2+96(k_{W0})^3\right)\\
\notag &+\frac{k^A_E}{V}k_{W0}\left(24\left(\frac{k_M}{V}\right)^3+158k_{W0}\left(\frac{k_M}{V}\right)^2+280\frac{k_M}{V}(k_{W0})^2+128(k_{W0})^3\right),\\
\notag K_6&=\left(\frac{k^A_E}{V}\right)^2\left(24\left(\frac{k_M}{V}\right)^3+134k_{W0}\left(\frac{k_M}{V}\right)^2+216\frac{k_M}{V}(k_{W0})^2+96(k_{W0})^3\right)\\
\notag &+\frac{k^A_E}{V}k_{W0}\left(24\left(\frac{k_M}{V}\right)^3+152k_{W0}\left(\frac{k_M}{V}\right)^2+272\frac{k_M}{V}(k_{W0})^2+128(k_{W0})^3\right),\\
\notag K_{7}&=\left(\frac{k^A_E}{V}\right)^2\left(24\left(\frac{k_M}{V}\right)^3+124k_{W0}\left(\frac{k_M}{V}\right)^2+180\frac{k_M}{V}(k_{W0})^2+72(k_{W0})^3\right)\\
\notag &+\frac{k^A_E}{V}k_{W0}\left(48\left(\frac{k_M}{V}\right)^3+284k_{W0}\left(\frac{k_M}{V}\right)^2+456\frac{k_M}{V}(k_{W0})^2+192(k_{W0})^3\right)\\
\notag &+(k_{W0})^2\left(24\left(\frac{k_M}{V}\right)^3+160k_{W0}\left(\frac{k_M}{V}\right)^2+288\frac{k_M}{V}(k_{W0})^2+128(k_{W0})^3\right),\\
\notag K_{8}&=\left(\frac{k^A_E}{V}\right)^2\left(
6\left(\frac{k_M}{V}\right)^4+
48k_{W0}\left(\frac{k_M}{V}\right)^3+
126(k_{W0})^2\left(\frac{k_M}{V}\right)^2+
132\frac{k_M}{V}(k_{W0})^3+
72(k_{W0})^4\right),\\
\notag K_{9}&=\left(\frac{k^A_E}{V}\right)^2\left(
6\left(\frac{k_M}{V}\right)^4+
51k_{W0}\left(\frac{k_M}{V}\right)^3+
146(k_{W0})^2\left(\frac{k_M}{V}\right)^2+
168\frac{k_M}{V}(k_{W0})^3+
64(k_{W0})^4\right),\\
\notag K_{10}&=\left(\frac{k^A_E}{V}\right)^2\left(
6\left(\frac{k_M}{V}\right)^4+
48k_{W0}\left(\frac{k_M}{V}\right)^3+
136(k_{W0})^2\left(\frac{k_M}{V}\right)^2+
160\frac{k_M}{V}(k_{W0})^3+
64(k_{W0})^4\right),\\
\notag K_{11}&=\left(\frac{k^A_E}{V}\right)^2\left(
24\left(\frac{k_M}{V}\right)^4+
191k_{W0}\left(\frac{k_M}{V}\right)^3+
509(k_{W0})^2\left(\frac{k_M}{V}\right)^2+
547\frac{k_M}{V}(k_{W0})^3+
200(k_{W0})^4\right),\\
\notag &+\frac{k^A_E}{V}k_{W0}\left(
12\left(\frac{k_M}{V}\right)^4+
96k_{W0}\left(\frac{k_M}{V}\right)^3+
252(k_{W0})^2\left(\frac{k_M}{V}\right)^2+
264\frac{k_M}{V}(k_{W0})^3+
96(k_{W0})^4\right),\\
\notag K_{12}&=\left(\frac{k^A_E}{V}\right)^2\left(
19\left(\frac{k_M}{V}\right)^4+
149k_{W0}\left(\frac{k_M}{V}\right)^3+
416(k_{W0})^2\left(\frac{k_M}{V}\right)^2+
463\frac{k_M}{V}(k_{W0})^3+
168(k_{W0})^4\right),\\
\notag &+\frac{k^A_E}{V}k_{W0}\left(
18\left(\frac{k_M}{V}\right)^4+
161k_{W0}\left(\frac{k_M}{V}\right)^3+
489(k_{W0})^2\left(\frac{k_M}{V}\right)^2+
588\frac{k_M}{V}(k_{W0})^3+
224(k_{W0})^4\right),\\
\notag K_{13}&=\left(\frac{k^A_E}{V}\right)^2\left(
18\left(\frac{k_M}{V}\right)^4+
143k_{W0}\left(\frac{k_M}{V}\right)^3+
399(k_{W0})^2\left(\frac{k_M}{V}\right)^2+
452\frac{k_M}{V}(k_{W0})^3+
168(k_{W0})^4\right),\\
\notag &+\frac{k^A_E}{V}k_{W0}\left(
18\left(\frac{k_M}{V}\right)^4+
158k_{W0}\left(\frac{k_M}{V}\right)^3+
476(k_{W0})^2\left(\frac{k_M}{V}\right)^2+
576\frac{k_M}{V}(k_{W0})^3+
224(k_{W0})^4\right),\\
\notag K_{14}&=\left(\frac{k^A_E}{V}\right)^2\left(
12\left(\frac{k_M}{V}\right)^4+
98k_{W0}\left(\frac{k_M}{V}\right)^3+
276(k_{W0})^2\left(\frac{k_M}{V}\right)^2+
306\frac{k_M}{V}(k_{W0})^3+
108(k_{W0})^4\right),\\
\notag &+\frac{k^A_E}{V}k_{W0}\left(
24\left(\frac{k_M}{V}\right)^4+
214k_{W0}\left(\frac{k_M}{V}\right)^3+
654(k_{W0})^2\left(\frac{k_M}{V}\right)^2+
780\frac{k_M}{V}(k_{W0})^3+
288(k_{W0})^4\right),\\
\notag &+(k_{W0})^2\left(
12\left(\frac{k_M}{V}\right)^4+
116k_{W0}\left(\frac{k_M}{V}\right)^3+
384(k_{W0})^2\left(\frac{k_M}{V}\right)^2+
496\frac{k_M}{V}(k_{W0})^3+
192(k_{W0})^4\right).
\end{align}
       \endgroup
       \normalsize

    Let us now derive an expression for $\pi^{(1)}$. Starting with the transient states $\T= \{i_1, \ldots, i_{15}\}$, we obtain that $\beta^{(1)}=[\pi^{(1)}_{i_1},...,\pi^{(1)}_{i_{15}}] = \alpha S_1(-T_0)^{-1}$, and so $\beta^{(1)}=[0,...,0,\pi^{(1)}_{i_{14}},\pi^{(1)}_{i_{15}}]$, with 
\begin{equation}\nonumber
\pi^{(1)}_{i_{12}}=\frac{4\mu' \beta\frac{k_{M}^A}{V}}{k^1_{W0}+ \frac{k'_{M}}{V}},\;\;\;\;\;\; \pi^{(1)}_{i_{13}}=\frac{4\mu b\frac{k_{M}^A}{V}}{k^2_{W0}+ (\frac{k_{M}}{V} + \frac{\bar k_{M}}{V})}.
\end{equation}
 Finally, $\alpha^{(1)}=[\pi^{(1)}_{a},\pi^{(1)}_{r}]$ is the unique vector such that $\alpha^{(1)}Q_{\A} = -\beta^{(1)}[R_1 + T_1(-T_0)^{-1}R_0]$ and $ \alpha^{(1)}\one = -\beta^{(1)}\one$. After some calculations, we obtain
       \begin{equation}
    \begin{aligned}\nonumber
        \pi^{(1)}_{a}&=\frac{(\mu\mu')^2K_{15}((\mu+\mu')K_{16}+K_{17})}{K_{20}((\mu')^2K_2+(\mu)^2K_3+(\mu'+\mu)K_4+\mu'K_5+\mu K_6+K_{7})},\\
        \pi^{(1)}_{r}&=-\pi^{(1)}_{a}-\pi^{(1)}_{i_{14}}-\pi^{(1)}_{i_{15}}\\
        &=-\frac{(\mu\mu')^2K_{15}((\mu+\mu')K_{16}+K_{17})}{K_{20}((\mu')^2K_2+(\mu)^2K_3+(\mu'+\mu)K_4+\mu'K_5+\mu K_6+K_{7})}-\mu' K_{18}-\mu K_{19},
     \end{aligned}   
        \end{equation}
       with

       \begingroup
       \small
       \begin{align}
\notag K_{15} &= 2  \frac{k^A_E}{V}\frac{k_M}{V}\left(3\frac{k_M}{V}+2k_{W0}\right)b^2,\;\;\;K_{16} = \frac{k^A_E}{V}\left(2\left(\frac{k_M}{V}\right)^2+16(k_{W0})^2+12\frac{k_M}{V}k_{W0}\right),\\
\notag K_{17} &= \frac{k^A_E}{V}\left(21\frac{k_M}{V}k_{W0}+4\left(\frac{k_M}{V}\right)^2+24(k_{W0})^2\right)+k_{W0}\left(24\frac{k_M}{V}k_{W0}+4\left(\frac{k_M}{V}\right)^2+32(k_{W0})^2\right),\\
\notag K_{18} &=\frac{4 b\frac{k_{M}}{V}}{k_{W}+ \frac{k_{M}}{V}},\;\;K_{19}=\frac{4 b\frac{k_{M}}{V}}{k_{W}+ 2\frac{k_{M}}{V}},\notag K_{20} =\left(2\frac{k_M}{V}+k_{W0}\right).
       \end{align}
       \endgroup
       \normalsize


\noindent \textbf{Time to memory loss.} 
As a reminder, we define the time to memory loss of the active state as $h_{a,r}(\eps)$ and the time to memory loss of the repressed state as $h_{r,a}(\eps)$. Let us start by deriving the order and the leading coefficients of $h_{a,r}(\eps)$ and $h_{r,a}(\eps)$. By \eqref{eqn:orderSD}, the order of the stationary distribution at $a$ and $r$ are $k_a=-\min\{1-2,0\}=1$ and $k_r=-\min\{1-1,0\}=0$, respectively. This is consistent with the results obtained in Section \ref{SD4D}. As obtained for the 3D model, here we obtain $\pi^{(0)}_r=1$ and $\pi^{(1)}_a>0$, and thus
\begin{equation*}
h_{a,r}(\eps)= \frac{1}{(Q_\A)_{a,r}} \frac{1}{\eps} + O(1),\;\;\text{and}\;\;h_{r,a}(\eps)= \frac{1}{(Q_\A)_{a,r} \cdot \pi_a^{(1)}} \frac{1}{\eps^2} + O\left(\frac{1}{\eps}\right).
\end{equation*}

Now, in order to exploit Theorems S.2 and 3.4 from \cite{Monotonicitypaper} and determine how $\mu'$ affects $h_{a,r}(\eps)$ and $h_{r,a}(\eps)$, we introduce a small approximation in the transition rates of $X^{\eps}$, namely, $\frac{x_3-1}{2}\approx x_3 $ and $ \frac{x_4-1}{2}\approx x_4$ in $f_{R121}(x)$ and $f_{R122}(x)$, respectively. This approximation can be justified by introducing the reasonable assumption that each nucleosome characterized by a repressive modification ($\mathrm{D^R_1}$ and $\mathrm{D^R_2}$) has the ability to catalyze the establishment of the opposite repressive mark on itself. Now, let us verify that both conditions $(i)$ and $(ii)$ of Theorem S.2 in \cite{Monotonicitypaper2} hold. These conditions can be written as follows:
	\begin{enumerate}
		\item[(i)]
		For each $1 \leq j \leq n$, the vector $A v_j$ has entries in $\{-1,0,1\}$ only.
		\item[(ii)]
		For each $x \in \X$, $1 \leq i \leq m$ and $y \in \partial_i(K_A+x) \cap \X$ we have that for each $1 \le k \le s$,
		\begin{equation*}
		\sum_{j \in G^{k,-}_i} \Breve{\rate}_j (y) \leq \sum_{j \in G^{k,-}_i} \rate_j(x), \quad \text{ where } G^{k,-}_i = \{j \in G^k \:|\: \inn{A_{i\bullet},v_j} = -1 \},
		\end{equation*}
		and
		\begin{equation*}
		\sum_{j \in G^{k,+}_i} \Breve{\rate}_j(y) \geq \sum_{j \in G^{k,+}_i} \rate_j(x), \quad \text{ where } G^{k,+}_i = \{j \in G^k \:|\: \inn{A_{i\bullet},v_j} = 1 \}.
		\end{equation*}
	\end{enumerate}
To verify that these conditions hold, let us first note the ten possible transitions vectors for the continuous time Markov chain $X^{\eps}(t)$:     
 %
 $v_1=-v_2=(1,0,-1,0)^T$, $v_3=-v_4=(1,0,0,-1)^T$, $v_5=-v_6=(0,1,0,0)^T$, $v_7=-v_8=(0,0,1,0)^T$, $v_9= -v_{10}=(0,0,01)^T$, with the associated infinitesimal transition rates
%
$\rate_1(x)=f_{R121}(x)$, $\rate_2(x)=g^{\eps}_{R121}(x)$, $\rate_3(x)=f_{R122}(x)$, $\rate_4(x)=g^{\eps}_{R122}(x)$, $\rate_5(x)=f_A(x)$, $\rate_6(x)=g^{\eps}_A(x)$, $\rate_7(x)=f_{R1}(x)$, $\rate_8(x)=g^{\eps}_{R1}(x)$, $\rate_9(x)=f_{R2}(x)$, $\rate_{10}(x)=g^{\eps}_{R2}(x)$.
Let
    \begin{equation}\label{matrixA4D}
       A= \begin{bmatrix}
        0 & -1 & 0 & 0\\
        1 & 0 & 1 & 0\\
        1 & 0 & 0 & 1\\
        1 & 0 & 1 & 1
        \end{bmatrix}.\nonumber
    \end{equation}
     Then, $(K_A+x)\cap \X=\{ w \in \X:\: x \preccurlyeq_A w \}$. Consider infinitesimal transition rates $\Breve{\rate}_{i}(x)$, $i=1,2,...,10$, defined as for $\rate_{i}(x)$, $i=1,2,...,10$, with all the parameters having the same values except that $\mu'$ is replaced by $\Breve{\mu}'$, with $\mu' \ge \Breve{\mu}'$. Now, condition $(i)$ of Theorem S.2 in \cite{Monotonicitypaper2} holds since
     $Av_1=-Av_2=(0,0,1,0)^T$, $Av_3=-Av_4=(0,1,0,0)^T$, $Av_5=-Av_6=(-1,0,0,0)^T$, $Av_7=-Av_8=(0,1,0,1)^T$ and $Av_9= -Av_{10}=(0,0,1,1)^T$. Assumption S.1 in \cite{Monotonicitypaper2} holds with $G^1 = \{9,1\}$, $G^2 = \{10,2\}$, $G^3 = \{7,3\}$, $G^4 = \{8,4\}$, $G^5 = \{5\}$, $G^6 = \{6\}$ and $\sigma(1)=9$, $\sigma(2)=1$, $\sigma(3)=10$, $\sigma(4)=2$, $\sigma(5)=7$, $\sigma(6)=3$, $\sigma(7)=8$, $\sigma(8)=4$, $\sigma(9)=5$, $\sigma(10)=6$.
     To verify that also condition $(ii)$ of Theorem S.2 in \cite{Monotonicitypaper2} holds, let us start with considering $x\in \X$ and $y \in 
     \partial_1(K_A+x) \cap \X=\{ w \in \X:\: x_2 =  w_2, x_1 + x_3\le w_1 + w_3,x_1 + x_4\le w_1 + w_4, x_1 + x_3 + x_4\le w_1 + w_3 + w_4\}$. Given that $\inn{A_{1\bullet},v_5}=-1$ and $\inn{A_{1\bullet},v_6}=1$, we must verify that $\rate_5(x)\ge\Breve{\rate}_5(y)$ and $\rate_6(x)\le\Breve{\rate}_6(y)$. Since 
   $x_2 =  y_2, x_1 + x_3\le y_1 + y_3,x_1 + x_4\le y_1 + y_4, x_1 + x_3 + x_4\le y_1 + y_3 + y_4$, then $\rate_{5}(x)=(\Dtot -(x_1+x_2+x_3+x_4))\left(k_{W0}^A+k_{W}^A + \frac{k_{M}^A}{V}x_2\right)\ge (\Dtot -(y_1+y_2+y_3+y_4))\left(k_{W0}^A+k_{W}^A + \frac{k_{M}^A}{V}y_2\right)= \Breve{\rate}_{5}(y)$ and 
   $\rate_{6}(x)= x_2\left(\eps \frac{k_{M}^A}{V}\Dtot + \frac{k^A_E}{V}(x_3+x_4+2x_1)\right) \le
   y_2\left(\eps \frac{k_{M}^A}{V}\Dtot + \frac{k^A_E}{V}(y_3+y_4+2y_1)\right)= \Breve{\rate}_{6}(y)$.
    Let us now consider $x\in \X$ and $y \in \partial_2(K_A+x) \cap \X=
    \{ w \in \X:\: x_2 \ge  w_2, x_1 + x_3= w_1 + w_3,x_1 + x_4\le w_1 + w_4, x_1 + x_3 + x_4\le w_1 + w_3 + w_4\}$. Given that $\inn{A_{2\bullet},v_3}=\inn{A_{2\bullet},v_7}=1$ and $\inn{A_{2\bullet},v_4}=\inn{A_{2\bullet},v_8}=-1$, we need to verify that $\rate_3(x)+\rate_7(x)\le\Breve{\rate}_3(y)+\Breve{\rate}_7(y)$ and $\rate_4(x)+\rate_{8}(x)\ge\Breve{\rate}_4(y)+\Breve{\rate}_{8}(y)$ hold. Since
    $x_2 \ge  y_2, x_1 + x_3= y_1 + y_3,x_1 + x_4\le y_1 + y_4, x_1 + x_3 + x_4\le y_1 + y_3 + y_4$, then 
    $\rate_3(x)+\rate_7(x)=(\Dtot -(x_1+x_2+x_3))\left(k^1_{W0} + \frac{k^{'}_{M}}{V}(x_1+x_4)\right)\le (\Dtot -(y_1+y_2+y_3))\left(k^1_{W0} + \frac{k^{'}_{M}}{V}(y_1+y_4)\right)=\Breve{\rate}_3(y)+\Breve{\rate}_7(y)$ and $\rate_4(x)+\rate_{8}(x)=(x_3+x_1)\mu'\left(\eps \frac{k_{M}^A}{V}\Dtot \beta + x_2\frac{k^A_E}{V}\right)\ge (y_3+y_1)\Breve{\mu}'\left(\eps \frac{k_{M}^A}{V}\Dtot \beta + y_2\frac{k^A_E}{V}\right)=\Breve{\rate}_4(y)+\Breve{\rate}_{8}(y)$.
    Let us now consider $x\in \X$ and $y \in \partial_3(K_A+x) \cap \X=
    \{ w \in \X:\: x_2 \ge  w_2, x_1 + x_3\le w_1 + w_3,x_1 + x_4= w_1 + w_4, x_1 + x_3 + x_4\le w_1 + w_3 + w_4\}$. Given that $\inn{A_{3\bullet},v_1}=\inn{A_{3\bullet},v_9}=1$ and $\inn{A_{3\bullet},v_2}=\inn{A_{3\bullet},v_{10}}=1$, we need to verify that $\rate_1(x)+\rate_9(x)\le\Breve{\rate}_1(y)+\Breve{\rate}_9(y)$ and $\rate_2(x)+\rate_{10}(x)\ge\Breve{\rate}_2(y)+\Breve{\rate}_{10}(y)$ hold. Since
    $x_2 \ge  y_2, x_1 + x_3\le y_1 + y_3,x_1 + x_4= y_1 + y_4, x_1 + x_3 + x_4\le y_1 + y_3 + y_4$, then 
    $\rate_1(x)+\rate_9(x)=(\Dtot -(x_1+x_2+x_4))\left(k^2_{W0} +  \frac{k_{M}}{V}(x_1+x_4)  +  \frac{\bar k_{M}}{V}(x_1+x_3)\right)\le (\Dtot -(y_1+y_2+y_4))\left(k^2_{W0} +  \frac{k_{M}}{V}(y_1+y_4)  +  \frac{\bar k_{M}}{V}(y_1+y_3)\right)=\Breve{\rate}_1(y)+\Breve{\rate}_9(y)$ and $\rate_2(x)+\rate_{10}(x)=(x_1+x_4)\mu\left(\eps \frac{k_{M}^A}{V}\Dtot b + x_2\frac{k^A_E}{V}\right)\ge (y_1+y_4)\mu\left(\eps \frac{k_{M}^A}{V}\Dtot b + y_2\frac{k^A_E}{V}\right)=\Breve{\rate}_2(y)+\Breve{\rate}_{10}(y)$.
    Finally, let us consider $x\in \X$ and $y \in \partial_3(K_A+x) \cap \X=
    \{ w \in \X:\: x_2 \ge  w_2, x_1 + x_3\le w_1 + w_3,x_1 + x_4\le w_1 + w_4, x_1 + x_3 + x_4= w_1 + w_3 + w_4\}$. Given that $\inn{A_{4\bullet},v_7}=\inn{A_{4\bullet},v_9}=1$ and $\inn{A_{4\bullet},v_8}=\inn{A_{4\bullet},v_{10}}=-1$, we need to verify that
    $\rate_7(x)\le\Breve{\rate}_7(y)$, $\rate_9(x)\le\Breve{\rate}_9(y)$,
    $\rate_8(x)\ge\Breve{\rate}_{8}(y)$ and $\rate_{10}(x)\ge\Breve{\rate}_{10}(y)$ hold. Since $x_2 \ge  y_2, x_1 + x_3\le y_1 + y_3,x_1 + x_4\le y_1 + y_4, x_1 + x_3 + x_4= y_1 + y_3 + y_4$, that also imply $x_1\le  y_1, x_3\ge y_3,x_4\ge y_4$, then  $\rate_7(x)=(\Dtot -(x_1+x_2+x_3+x_4))\left(k^1_{W0} + \frac{k^{'}_{M}}{V}(x_1+x_4)\right)\le (\Dtot -(y_1+y_2+y_3+y_4))\left(k^1_{W0} + \frac{k^{'}_{M}}{V}(y_1+y_4)\right)=\Breve{\rate}_7(y)$, $\rate_9(x)=(\Dtot -(x_1+x_2+x_3+x_4))\left(k^2_{W0}+  \frac{k_{M}}{V}(x_1+x_4)  +  \frac{\bar k_{M}}{V}(x_1+x_3)\right)\le (\Dtot -(y_1+y_2+y_3+y_4))\left(k^2_{W0}+  \frac{k_{M}}{V}(y_1+y_4)  +  \frac{\bar k_{M}}{V}(y_1+y_3)\right)=\Breve{\rate}_9(y)$, $\rate_8(x)=x_3\mu'\left(\eps \frac{k_{M}^A}{V}\Dtot \beta + x_2\frac{k^A_E}{V}\right)\ge y_3\Breve{\mu}'\left(\eps \frac{k_{M}^A}{V}\Dtot \beta + y_2\frac{k^A_E}{V}\right)=\Breve{\rate}_{8}(y)$ and $\rate_{10}(x)=x_4\mu\left(\eps \frac{k_{M}^A}{V}\Dtot b + x_2\frac{k^A_E}{V}\right)$
    
    \noindent $\ge y_4\mu\left(\eps \frac{k_{M}^A}{V}\Dtot b + y_2\frac{k^A_E}{V}\right)=\Breve{\rate}_{10}(y)$.
Then, condition $(ii)$ of Theorem S.2 in \cite{Monotonicitypaper2} also holds.

We can then conclude that all of the conditions of Theorem S.2 in \cite{Monotonicitypaper2} hold and so do the conclusions of Theorem 3.4 in \cite{Monotonicitypaper2}, as per the remarks in SI - Section S.3 in \cite{Monotonicitypaper2}.